\documentclass[12pt]{article}
\usepackage{amsfonts}
\usepackage{epsfig}
\title{Outer Billiards on the Penrose Kite: Compactification and Renormalization}
\author{Richard Evan Schwartz \thanks{\hskip 5 pt Supported by 
N.S.F. Research Grant DMS-0072607}}

\newtheorem{theorem}{Theorem}[section]

\newtheorem{lemma}[theorem]{Lemma}

\newtheorem{corollary}[theorem]{Corollary}
\newtheorem{conjecture}[theorem]{Conjecture}

\def\startproof{{\bf {\medskip}{\noindent}Proof: }}

\def\endproof{$\spadesuit$  \newline}

\def\A{\mbox{\boldmath{$A$}}}% 
\def\Q{\mbox{\boldmath{$Q$}}}% 
\def\R{\mbox{\boldmath{$R$}}}% 
\def\T{\mbox{\boldmath{$T$}}}% 
\def\Z{\mbox{\boldmath{$Z$}}}% 

\begin{document}
\maketitle
\begin{abstract}
We give a fairly complete analysis of outer
billiards on the Penrose kite. Our
analysis reveals that this $2$ dimensional dynamical
system has a $3$-dimensional compactification,
a certain polyhedron exchange map defined on the $3$-torus,
and that this $3$-dimensional system admits a
renormalization scheme.  The two features, the
compactification and the renormalization scheme,
allow us to make sharp statements concerning
the distribution, large- and fine-scale geometry,
and hidden algebraic symmetry, of the orbits. One concrete result
is that the union of the unbounded orbits has
Hausdorff dimension $1$.
We establish many of the results with
computer-aided proofs that involve only integer
arithmetic. 
\end{abstract}

\section{Introduction}

\subsection{Background}
\label{backg}

Outer billiards is a dynamical system defined relative to
a convex shape in the plane.
B.H. Neumann [{\bf N\/}] introduced outer billiards in
the late 1950s, and J. Moser [{\bf M1\/}]
popularized the system as a toy model for celestial mechanics.
See [{\bf T1\/}], [{\bf T2\/}], and [{\bf DT1\/}] for 
expositions of outer billiards and many references.

\begin{center}
\resizebox{!}{1.3in}{\includegraphics{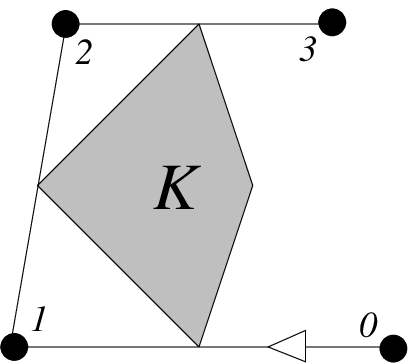}}
\newline
{\bf Figure 1.1:\/} outer billiards relative to $K$.
\end{center}

To define an outer billiards system,
one starts with a bounded convex set $K \subset \R^2$
and considers a point $x_0 \in \R^2-K$.  
One defines $x_1$ to be the point such that
the segment $\overline{x_0x_1}$ is tangent to $K$ at its
midpoint and $K$ lies to the right of the ray
$\overrightarrow{x_0x_1}$. 
The iteration $x_0 \to x_1 \to x_2...$
is called the {\it forwards outer billiards orbit\/} of $x_0$.
It is defined for almost every point of $\R^2-K$.
The backwards orbit is defined similarly.

Outer billiards is an affinely natural system,
in the sense that an affine map $T: P \to Q$
carrying the convex shape $P$ to the convex
shape $Q$ also carries the outer billiards orbits
relative to $P$ to the outer billiards orbit
relative to $Q$.   The reason is that affine
maps carry line segments to line segments
and respect the property of bisection.

One of the central questions about outer billiards
is the {\it Moser-Neumann question\/}, which
asks if an outer billiards system can have unbounded orbits.
Here is an abbreviated list of work on this problem.
\begin{itemize}
\item  J. Moser [{\bf M2\/}] sketches a proof, inspired by
KAM theory, that outer
billiards on $K$ has all bounded orbits provided
that $\partial K$ is at least $C^6$ smooth and positively curved.
R. Douady [{\bf D\/}] gives a complete proof.
\item In Vivaldi-Shaidenko [{\bf VS\/}], Kolodziej [{\bf Ko\/}], and
Gutkin-Simanyi [{\bf GS\/}], it is proved (each with different methods)
that outer billiards on a {\it quasirational polygon\/} has
all orbits bounded. \index{quasirational polygon}  This class of polygons
includes polygons with rational vertices and regular polygons.
In the rational case, all orbits are periodic.  
\item In [{\bf T2\/}], Tabachniikov shows the existence of
aperiodic orbits in the regular pentagon case, and works out
a renormalization scheme to explain their structure. 
 \item D. Genin [{\bf G\/}] shows that all orbits are bounded
for the outer billiards systems associated to trapezoids.
\index{trapezoids} See \S \ref{quad}.
\item In [{\bf S1\/}], we settled Moser-Neumann question by showing that
outer billiards has some unbounded orbits when defined relative to
$K(\phi^{-3})$.   Here $\phi$ is the golden ratio and
$K(A)$ denotes the kite with vertices
\begin{equation}
\label{parameter}
\label{quad}
(-1,0); \hskip 20 pt (0,1) \hskip 20 pt (0,-1) \hskip 20 pt (A,0); \hskip 30 pt
\end{equation}
Figure 1.1 shows an example. 
\item In [{\bf S2\/}], we showed that outer billiards
has unbounded orbits relative to $K(A)$, when $A$ is any
irrational number in $(0,1)$.
\item Dolgopyat and Fayad [{\bf DF\/}] showed
that outer billiards has unbounded orbits relative to the half-disk
and other ``caps'' made from slicing a disk nearly in half.
\end{itemize}
The shapes in [{\bf DF\/}] and [{\bf S2\/}] are the only
ones known to produce unbounded orbits, though certainly
it now seems that unboundedness is a common phenomenon.
Our monograph [{\bf S2\/}] gives more details about the
history of the problem.

The set $\R \times \Z_{\rm odd\/}$ is invariant under
the dynamics on $K(A)$. We call the orbits that lie in
this set {\it special\/}.   In [{\bf S2\/}] we gave quite
a lot of information about special orbits on kites.   For instance,
we gave a formula for the Hausdorff dimension of the set
$U_1(A)$ of unbounded special orbits, in terms of something
akin to the continued fraction expansion of $A$.
As a special case,
\begin{equation}
\label{dim1}
\dim(U_1(\phi^{-3}))=\frac{\log(2)}{\log(\phi^3)}.
\end{equation}
We also showed that every unbounded special orbit
is self-accumulating.  This is to say that every point of
an unbounded special orbit $O$ is an accumulation point
of  $O$.

In the $300$ page [{\bf S2\/}] we only considered the special orbits,
for the sake of ``brevity''. There is quite a bit more
to say about the general orbits, and our purpose
here is to say some of it, at least for the Penrose kite.
The first phenomenon is that
outer billiards on the Penrose kite, an
unbounded $2$ dimensional system, has a
$3$ dimensional compactification.
We saw similar things in [{\bf S1\/}] and
[{\bf S2\/}]. 

 The second phenomenon is that
this higher dimensional compactification
has a renormalization scheme.  The
renormalization scheme is the main
new feature of this paper. Its existence
allows us to get some precise results
about the dynamics.  We think that a
similar scheme exists in great generality,
but we currently don't have any techniques
for investigating it in general.

\subsection{The Distribution of Unbounded Orbits}
\label{basic}

In all that follows, it goes without saying that
our results concern only outer billiards on the Penrose kite.

\begin{theorem}
\label{penrose1}
Every orbit is either periodic or unbounded in both directions,
and the union of unbounded orbits has Hausdorff dimension $1$.
\end{theorem}

It is convenient to state our remaining results in terms
of the square of the outer billiards map, which we
call $\psi$.   The map $\psi$ leaves invariant the
set
\begin{equation}
\R^2_y=\bigcup_{n \in \Z} \R \times (y+2n).
\end{equation}
Each $\R^2_y$ is a discrete countable family of horizontal lines.
We always take $y$ as a point of the circle $\R/2\Z$.
In [{\bf S1\/}] and [{\bf S2\/}] we studied the
orbits in $\R^2_1$.

Let $\Z[\phi]$ denote the ring of elements
$m+n\phi$ where $m,n \in \Z$.  We will sometimes
use the notation
\begin{equation}
\left[\matrix{m\cr n}\right]=m+n\phi.
\end{equation}
We define an equivalence relation on points of
$\R/2\Z$.  We say that $a \sim b$ if
\begin{equation}
\label{mainequiv}
b= \pm \phi^{3k} a + 2m+2n\phi; \hskip 30 pt
k,m,n \in \Z.
\end{equation}
We will explain this equivalence relation in a
more natural way in \S \ref{renormXX}

An even length increasing sequence
$a_1<...<a_{2n}$ canonically defines a
Cantor set $C$, as follows.  Let $T_k$ be
the similarity carrying $[a_1,a_{2n}]$
to $[a_{2k-1},a_{2k}]$.  Then $C$ is the
limit set of the semigroup generated by
$T_1,...,T_n$.  For instance, the
sequence $0<1/3<2/3<1$ defines the usual
middle-third Cantor set.
We let $C^{\#}$ denote the
set obtained from $C$ by
removing the endpoints of all the
complementary regions. 

\begin{theorem}
\label{penrose2}
The set $\R^2_y$ contains unbounded orbits if and only
if $y \sim c$ and $c \in C^{\#}$, where $C$ is the
Cantor set defined by the sequence
$$\left[\matrix{0\cr 0}\right]<
\left[\matrix{2\cr -1}\right]<
\left[\matrix{4\cr -2}\right]<
\left[\matrix{6\cr -3}\right]<
\left[\matrix{-2\cr 2}\right]<
\left[\matrix{0\cr 1}\right].
$$
The set of such $y$ is a dense set of
Hausdorff dimension $\log(3)/\log(\phi^3)$.
\end{theorem}

\subsection{The Distribution of Periodic Orbits}
\label{tiling}

We define the {\it winding number\/} of a periodic orbit $p$
to be half the number of times $\{\psi^n(p)\}$ intersects the
strip 
\begin{equation}
\Sigma=\R \times [-2,2]
\end{equation}
This definition makes sense geometrically. It turns out that,
at least far from the origin, the $\psi$-orbits generally
wind around the origin, nearly following a large octagon,
and return to $\Sigma$ after each half-revolution.
See \S \ref{sob} and \S \ref{returnmap}.   The following
result says that orbits of high winding number are
extremely pervasive.  In particular, it says that
$\R^2_y$ has periodic orbits of arbitrarily
high winding number provided that $y \not \in 2\Z[\phi]$.

\begin{theorem}
\label{winding}
For any integer $N$, there is a
finite subset $W_N \subset 2\Z[\phi]$ such that
$\R^2_y$ has periodic orbits of winding
number greater than $N$ if $y \not \in W_N$.
\end{theorem}

As with any polygonal outer billiards system, every periodic point
is contained in a maximal convex polygon consisting of points which all have
the same period and combinatorial behavior.   We call these maximal polygons
{\it periodic tiles\/}.  For convenience, we include $K$ itself as a
periodic tile.  We call the union of the periodic tiles the
{\it dynamical tiling\/}, even though it is only a tiling of a 
subset of the plane, and we denote it $\cal D$.
Theorem \ref{penrose1} says
that $\cal D$ fills up everything
but $1$-dimensional set.

A clean,
finitary description of $\cal D$ is a beyond our reach.  $\cal D$ is
quite complicated. 
In particular, it contains infinitely many different shapes.
However, we will get a near-complete understanding of the portion of $\cal D$
contained in the triangle $T$ with vertices
\begin{equation}
\label{basictriangle}
(0,1); \hskip 30 pt (0,\phi^{-3}); \hskip 30 pt (\phi^{-1},\phi).
\end{equation}
$T$ is bounded by lines extending $3$ of the sides of $K$,
as in Figure 1.2.  We call $T$ the {\it fundamental triangle\/}.

Figure 1.2 below shows a tiling $\cal T$ of a full measure subset of $T$
by an infinite union of kites and octagons.
the kites are all similar to each other and the octagons
are all similar to each other.  The similarity factors
all have the form $\phi^{3k}$, where $\phi$ is the
golden ratio and $k$ is an integer.
We will describe $\cal T$ precisely in \S \ref{tiling1}.
The fractal set
which is the complement of the polygonal
tiles has Hausdorff dimension $1$. The set
of points in this fractal set, having well
defined orbits, also has Hausdorff dimension $1$.

\begin{center}
\resizebox{!}{7.3in}{\includegraphics{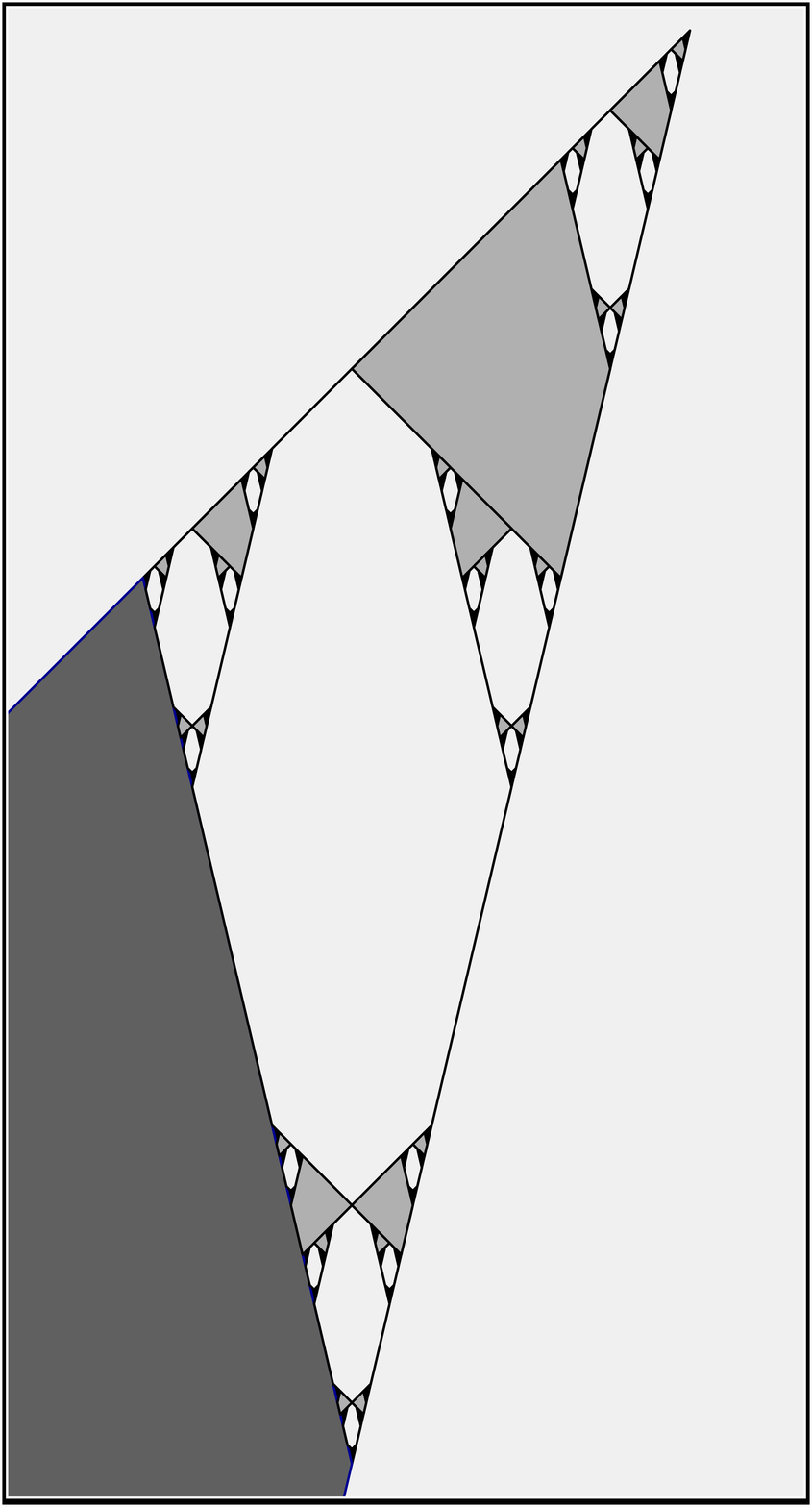}}
\newline
{\bf Figure 1.2:\/} The triangle $T$ and its tiling $\cal T$.
\end{center}

\begin{theorem}
\label{penrose9}
Suppose that $p \in T$ has a well defined
orbit.  Then $p$ has a periodic orbit if and only
if $p$ lies in the interior of a tile of $\cal T$.
Each tile of $\cal T$ is partitioned into
finitely many periodic tiles.  For any $N$,
only finitely many of the tiles contain
periodic orbits having winding number less than $N$.

\end{theorem}

\noindent
{\bf Remarks:\/} \newline
(i)
Theorem \ref{penrose9} is not quite as sharp as
we would like. Experimental evidence suggests that
each tile of $\cal T$ is itself a periodic tile.
Indeed, we plotted Figure 1.2 using our outer billiards
program, which finds the periodic tiles.
\newline
(ii) The set $y$ such that the horizontal line of height
$y$ intersects $T-{\cal T\/}$ in an infinite set is
precisely the Cantor set from Theorem \ref{penrose2}.
This is where the Cantor set comes from.
\newline
(iii) 
The fractal set of aperiodic points for outer billiards
on the regular pentagon has a self-similar structure akin
to the one in Figure 1.2.  In that case, the aperiodic 
orbits are all bounded.
See Tabachnikov [{\bf T2\/}] for details. 

\subsection{Renormalization}
\label{renormXX}

Now we reconsider the equivalence relation
defined in connection with Theorem
\ref{penrose2}.
To say that $a \sim b$ means that $b=T(a)$, where
$T$ is an affine map defined over $\Z[\phi]$
whose multiplying coefficient
is a unit in $\Z[\phi]$ and whose action is
congruent to the identity mod $2\Z[\phi]$.
Such maps form a group, and act on
$\R/2\Z$ with dense orbits.  The equivalence
classes we defined in connection
with Theorem \ref{penrose2} 
are precisely the orbits of this
group action.

We think of $G_2$
as related to $\Gamma_2$, the {\it level 2 congruence
subgroup\/} of $PSL_2(\Z)$, the modular group.
In [{\bf S2\/}] we discovered that, with respect
to the special orbits on arbitrary kites,
there is a kind of hidden $\Gamma_2$-symmetry,
For instance, the dimension
of the set of unbounded orbits on
$\R^2_1$, as a function of the kite parameter, is
a $\Gamma_2$-invariant function.  We were not able
to see the kind of renormalization structure
that we establish here, but the way we think
of things is that the union of all the dynamical
systems defined by outer billiards on kites
is a kind of plane-bundle over the parameter interval.
We think there is a group $\widehat \Gamma_2$
acting (in way that meaningfully relates to
the dynamics) on this bundle in such a way that $\Gamma_2$ 
gives the action on the base space and our
group $G_2$ here is the restriction of
$\widehat \Gamma_2$ to a particular fiber.

Having indulged in some speculation, we now return
to concrete results.

Our renormalization works best for orbits
which we call {\it generic\/}.
Say that a point $p=(x,y) \in \R^2$ is {\it generic\/} if
it does not satisfy any equation of the form
\begin{equation}
ax+by+c=0; \hskip 30 pt a,b,c \in \Z[\phi]; \hskip 30 pt a \not = 0.
\end{equation}
That is, $(x,y)$ is generic if it does not lie on
a non-horizontal line that is defined over $\Z[\phi]$.
The outer billiards map is
entirely defined on the set of generic points, and
preserves this set.  So, it makes sense to speak
of a generic orbit.   Theorem \ref{non-generic} below
explains the sense in which we do not miss much
by ignorning the non-generic orbits.  In \S \ref{tricky}
we discuss the issues surrounding the renormalization
of non-generic orbits.

We define
two kinds of equivalence relations between
orbits.

\begin{itemize}
\item Let $\langle O\rangle$ denote the graph of an orbit
$O=\{p_n\}$, namely the subset  $\{(n,p_n)\} \subset \R^3$.
We call two orbits $O_1$ and $O_2$
{\it coarsely equivalent\/} if there is a
$K$ bi-lipschitz map $h: \R^3 \to \R^3$ such
that $h(\langle O_1 \rangle)$ and $\langle O_2 \rangle$ lie
in $K$-tubular neighborhoods of each other.
We call $K$ the {\it coarse equivalence constant\/}.

\item We call the orbits $O_1$ and $O_2$ are
 {\it locally equivalent\/} if, for each point
$p_1 \in O_1$, there is a point
$p_2 \in O_2$, open disks $\Delta_1$ and $\Delta_2$,
and a similarity $S: \Delta_1 \to \Delta_2$ such that
$S(p_1)=S(p_2)$ and $S$ conjugates the first
return map $\psi|\Delta_1$ to one of the two return maps
$\psi^{\pm 1}|\Delta_2$, at least on generic points.
In particular,
$S(O_1 \cap \Delta_1)=O_2 \cap \Delta_2$.  Also,
$S$ carries ${\cal D\/} \cap \Delta_1$ to
${\cal D\/} \cap \Delta_2$ modulo the operation
of subdividing each tile into finitely many
smaller polygons.
\end{itemize}

\begin{theorem}
\label{ULE}
Suppose that $y_1 \sim y_2$ are two
parameters in the same $G_2$ orbit.  Then
there is a bijection between a
certain subset of the orbits in
$\R^2_{y_1}$ and a certain subset of
the orbits in $\R^2_{y_2}$.  These
subsets contain all generic unbounded
orbits and also all generic periodic
orbits having sufficiently high winding number.
Corresponding orbits
are both locally and coarsely equivalent,
and the coarse equivalence constant only
depends on $(y_1,y_2)$ and not on
the individual orbit.
\end{theorem}

The last statement in Theorem \ref{ULE} is
important for the periodic orbits.  Every
two periodic orbits are coarsely equivalent,
so we need some kind of uniformity in order
to make a meaningful statement.

Here are two applications of
Theorem \ref{ULE}.

\begin{theorem}
\label{self}
Every generic unbounded orbit is self-accumulating in
at least one direction.
\end{theorem}
Theorem \ref{self} is probably true for the non-generic
unbounded orbits as well, but our techniques fall a bit
short of this result.

\begin{theorem}
\label{penrose10}
Any generic unbounded orbit is locally and coarsely equivalent to
a generic unbounded orbit that intersects the fundamental
triangle $T$.  In particular, in a small neighborhood of
any generic point $p$ with an unbounded orbit, the dynamical
tiling in a neighborhood $p$ is isometric to a neighborhood
of $\cal T$, modulo the addition or removal of a countable
set of lines.
\end{theorem}

The ambiguity concerning the countable set of lines
comes from the set containing the non-generic orbit.
The main thing that is missing in Theorem \ref{penrose10}
is a description
of the dynamical tiling in the neighborhood of
points that do not have well-defined orbits.
Figure 1.3 shows what the dynamical tiling looks
like in a certain region whose lowest vertex
is $(3,0)$, a point which turns out to be a
fixed point of renormalization in a sense that
we will make precise in our Fixed Point Theorem from \S 5.

\begin{center}
\resizebox{!}{4in}{\includegraphics{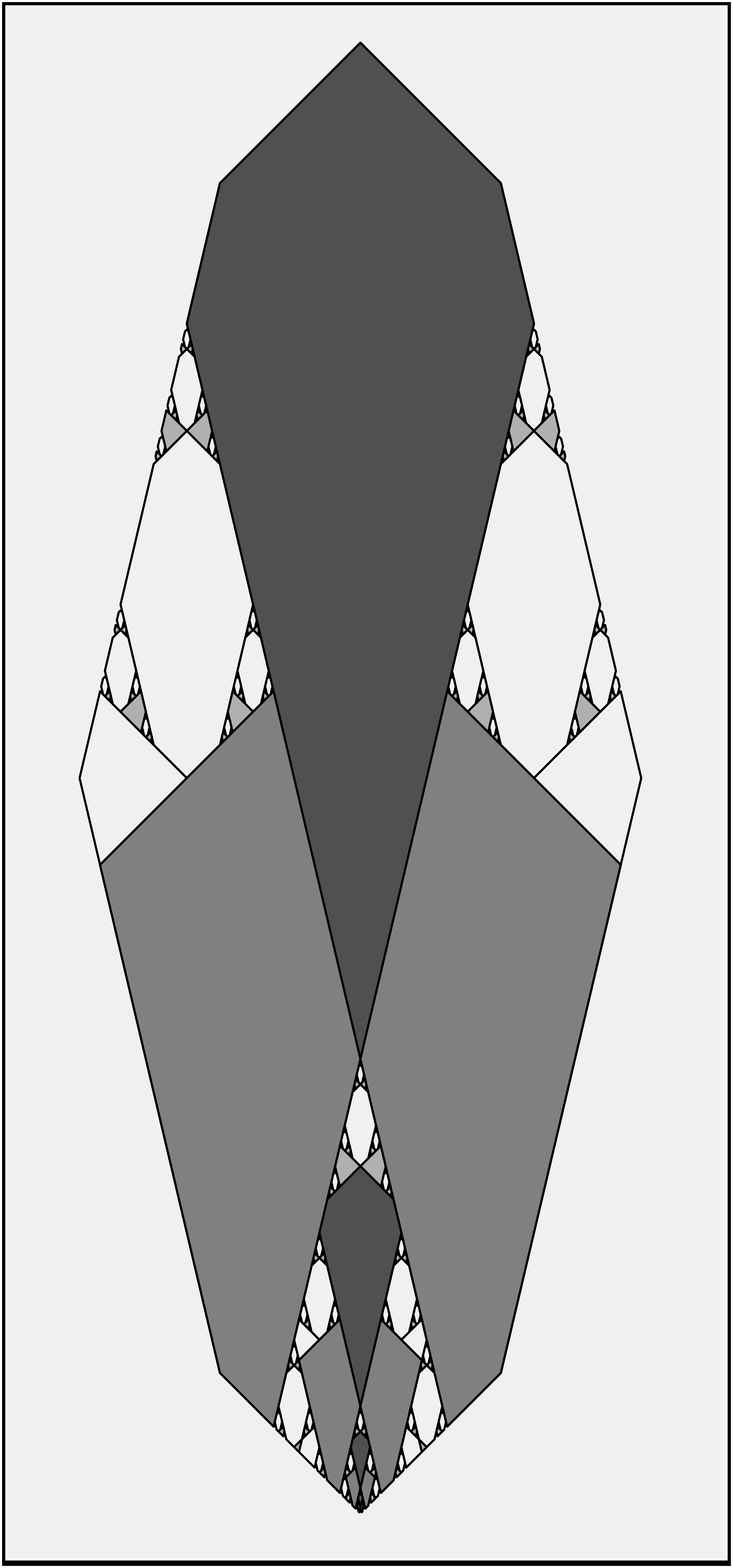}}
\newline
{\bf Figure 1.3\/}: The dynamical tiling near
$(3,0)$.
\end{center}

\subsection{Fine Points}

In this section we present a few more result that seem
a bit more specialized than the ones above.

For all $r>0$, let $Y_r$ denote
those $y \in Y$ such that $\R^2_y$ has an unbounded
orbit that intersects the disk of radius $r$ about
the kite vertex $(\phi^{-3},0)$.  

\begin{theorem}
\label{outside}
For all $r>0$ the set $Y_r$ is a
nowhere dense set having Hausdorff dimension
$\log(3)/\log(\phi^3)$.
\end{theorem}
Comparing Theorem \ref{outside} and
Theorem \ref{penrose2}, we can say qualitatively
that most unbounded orbits
stay far from the origin.  In particular, for any
compact subset, we can find an unbounded orbit
that avoids this set.
This is much different than what happens
for the special orbits. 
For instance, our result [{\bf S2\/}, Erratic Orbit Theorem]
says in particular that every unbounded orbit in $\R^2_1$ 
returns infinitely often to every neighborhood of
the kite vertex $(0,1)$.

Our next result quantifies the sense in which there
are fewer non-generic unbounded orbits than there
are generic unbounded orbits.

\begin{theorem}
\label{non-generic}
The union of non-generic unbounded orbits has
Hausdorff dimension at most $\log(3)/\log(\phi^3)$, and
$\R^2_y$ contains non-generic unbounded orbits only if
it contains generic unbounded orbits.
\end{theorem}

Here is one more result about periodic orbits.
\begin{theorem}
\label{stable}
Suppose that $y=m+n\phi$ with $m$ and $n$
odd integers.  Then
there is some open neighborhood $V$ of the
line $\R \times \{y\}$ that contains no
unbounded orbits.
\end{theorem}

Theorem \ref{stable} says that these special
lines are contained in the interior of the
closure of the union of the periodic tiles.
One might say that these lines are paved over
with periodic tiles.

Here is a nice special case of
Theorem \ref{penrose2}.

\begin{theorem}
\label{ring}
Suppose that $y=m+n\phi$ with $m$ and $n$
integers.  Then $\R^2_y$ has unbounded orbits if and
only if $m$ is odd and $n$ is even.
\end{theorem}

We wonder which
$y \in \Q[\phi]$ are such that
$\R^2_y$ has unbounded orbits.
Theorem \ref{penrose2} reduces
this to an arithmetical question about
the Cantor set $C$.

\subsection{Discussion}

Outer billiards is a non-compact $2$-dimensional
system.  We exhibit a $3$ dimensional
compactification of (a certain first return map of)
the outer billiards system.  The compactification
turns out to be a polyhedron exchange map 
defined on a $3$-torus $\widehat \Sigma$.  This result is
very similar to the Master Picture Theorem in
[{\bf S2\/}].  We call our result the
Compactification Theorem.

Once we find the compactification, we will observe that
it has a renormalizable structure that is responsible
for most of the results presented above.
Specifically, we will find sets $\widehat A, \widehat B \subset \widehat \Sigma$
together with a $3$ to $1$ map
$\widehat R: \widehat A \to \widehat B$ such that $\widehat R$ conjugates the
first return map on $\widehat A$ to the first return map on $\widehat B$.
Our main result here is Theorem \ref{master2}, the
Renormalization Theorem. 

It is worth comparing the renormalization scheme here to the
one [{\bf T2\/}] for regular pentagons.  In that case,
one picks a certain bounded subset $B$ of the plane
and observes that the first return map to $B$ is
renormalizable. What this means, in part, is that there is a subset
$A \subset B$ and a similarity $R: A \to B$ which
conjugates $f_A$ to $f_B$. Here $f_A$ is the first
return map to $A$ and $f_B$ is the first return
map to $B$.  It also means that the subsets $A$ and
$B$ are large enough to capture all the dynamical
properties of the whole map.

A direct renormalization like this would be very difficult 
to establish in the presence of unbounded orbits,
because the first return times to any sufficiently
large compact set would be unbounded.  Also,
thanks to Theorem \ref{outside}, there really
is no compact set that ``captures'' all the
relevant dynamics.
What we do for the penrose kite is show the existence
of the renormalization scheme on a higher dimensional
compactification.  Once we make the compactification,
the renormalization is at least vaguely similar
to what happens for the regular pentagon.

One difference between the renormalization
scheme here and the one for the regular
pentagon is that the one here involves a
$3$-to-$1$ covering map rather than a
similarity.  Another difference is that
the regular pentagon case can be analyzed
by hand, just following the orbits of several
convex polygons.  Here we need to keep
track of about a million polyhedra just to
see that the scheme works. This is what seems to make a
computer-aided proof necessary.

The moral of the story is that if one
wants to find remormalization schemes for
polygonal outer billiards, one should first
compactify.  In some sense, this is a lesson
I learned from John Smillie.  When hearing
about my earlier work on outer billiards,
Smillie guessed that probably there was a
renormalization scheme behind it.

\subsection{Overview of the Paper}

In \S 2 we study a certain $3$-to-$1$ covering
map $R: \R/2\Z \to \R/2\Z$ which governs the
structure of our renormalization map $\widehat R$.
This map is closely related to the Cantor
set in Theorem \ref{penrose2}.
The work is in \S 2 is a microcosm for the
rest of the paper:  We give some theoretical
arguments to reduce the main result to a
finite calculation which is too big to do by
hand, and then we explain a rigorous computer
calculation that finishes the proof.

In \S 3 we study the tiling $\cal T$ shown in Figure 1.2
and prove a number of results about it.

In \S 4 we present some background information
about outer billiards and polyhedron exchange maps.
Most of this information also appears
in [{\bf S1\/}] and [{\bf S2\/}].	

In \S 5 we state the Compactification Theorem,
the Renormalization Theorem, as well as
several useful auxilliary results.  One of
these auxilliary results, the
Fundamental Orbit Theorem, explains the
structure of Figure 1.2.  Another
result, the Fixed Point Theorem, explains
the dilation symmetry in Figure 1.3.
The remaining results, the Near
Reduction Theorem and the Far Reduction Theorem,
explain the sense
in which renormalization ``brings orbits
closer to the origin''.

In \S 6 we put together
the material from \S 2-5 to deduce all the
theorems mentioned in the introduction.
We warn the reader that the order we prove
these results is rather different from
the way we have stated them. For instance,
it takes almost every other result in
order to prove Theorem \ref{penrose1}.

In \S 7 we explain the main computational
ideas we use in our proofs.
These algorithms perform fairly 
standard tasks -- e.g. detecting
whether a point is contained in the interior
of a polyhedron.  Later chapters will
refer back to the methods explained
in \S 6.  Given that the computational
algorithms perform fairly standard tasks,
the reader won't lose much understanding
of the overall proof if they just skim
the material in \S 7.  On the other hand,
we think that \S 7 might be very useful
for someone who would like to reproduce our
results of to prove similar results in a
related stting.

In \S 8 we prove a technical result, known
as the Pinwheel Lemma, which is helpful in 
proving the Compactification Theorem.
Versions of the Pinwheel Lemma also
appear in [{\bf S1\/}],
[{\bf S2\/}], and [{\bf S3\/}], and related
results appear in other works on the subject
by other authors.

In \S 9 we prove the Compactification Theorem.
Our proof is very similar to the proof of
[{\bf S1\/}, Arithmetic Graph Lemma].

In \S 10 we prove the Remormalization Theorem
by reducing it to an explicit calculation involving
the polyhedron exchange map.

In \S 11 we prove the Fundamental Orbit Theorem
and the Fixed Point Theorem.

In \S 12 we prove the Near Reduction Theorem.

In \S 13 we prove the Far Reduction Theorem.

In \S 14 we include coordinates for all 
the polyhedra involved in the polyhedron
exchange map, as well as coordinates
the sets $\widehat A$ and $\widehat B$.

\subsection{Computational Issues}

The general strategy of our paper is to
reduce all the results to statements
about finite partitions of various sets by polygons
and polyhedra.  These partitions sometimes
involve a huge number of individual pieces, on the order of a million,
and so it is necessary for us to use the
computer to deal with them effectively.

One source of potential error in a computer
aided proof is floating-point (or roundoff)
error.  To avoid any possibility 
of floating-point errors, we perform our calculations
using exact arithmetic in the number ring $\Z[\phi]$.
The special nature of our constructions allows
us to do this.  With exact arithmetic calculations,
the one potential hazard is overflow error.
We avoid overflow error by checking the sizes of
the integers involved after every arithmetic
operation.  
We describe the main features of these calculations in
\S \ref{arithmetic}. 

We illustrate our structural results in detail in \S 4,
but there is probably no way for the reader to appreciate the
details of the objects involved in without seeing
explicit (and interactive) computer plots.
We made an extremely detailed and
extensive java applet that lets the reader see everything
in the paper.  This applet is available on my website.
For the reader who would like to do his/her own experiments,
we include enough information in the appendix so that
in principle one could reproduce the calculations.

We would like to comment on some of
the figures in the paper.   To illustrate certain definitions
which make sense for any kite $K(A)$, 
we will draw $K(1/4)$
in place of $K(\phi^{-3})$, because $K(1/4)$ is much
easier to draw by hand.   We draw these pictures mainly
to give the reader a picture of what is going on, and for
these purposes a picture of $K(1/4)$ tells the whole story.
Note that $1/4-\phi^{-3}=.0139...$, so the pictures
we draw for $K(1/4)$ are geometrically quite close
to the ones for $K(\phi^{-3})$.
On the other hand, the computer pictures we draw will
show $K(\phi^{-3})$.

\subsection{Further Results}

There are some other things I've noticed
about outer billiards on the Penrose kite. The
other things have a different character
from the results here.
They have to do with the patterns
one sees in the so-called arithmetic graph associated
to the dynamics.  In \S \ref{freeze1}
I briefly discuss one of the results, something I call
the {\it freezing phenomenon\/}. I would like to
have presented some of these results, this paper
is already long enough and also I have not worked
out proofs for these other results.

\subsection{Acknowledgements}

I would especially like to thank John Smillie, who has for
several years repeatedly
told me to ``get renormalization into the picture'' of
outer billiards.  I would also like to thank
Gordon Hughes, Rick Kenyon, Curt McMullen, and
Sergei Tabachnokov for helpful
conversations about matters related to this work.

\newpage
\section{The Circle Renormalization Map}

In this chapter, we define and then study the
renormalization map $R$ that we mentioned in the
introduction.  We also define and study the
tiling $\cal T$ we mentioned in the introduction
and plotted in Figure 1.2.  The map $R$ and
the tiling $\cal R$ are closely related.

\subsection{Basic Definition}
\label{renorm}

Let 
\begin{equation}
\T=\R/2\Z
\end{equation}
Here we give a precise definition of the
renormalization map $R: \T \to \T$ discussed
in the introduction.

We decompose $[0,2]$ into $5$ intervals.
\begin{itemize}
\item $I_1=[0,\phi^{-2}]=[0,2-\phi]$.
\item $I_2=[\phi^{-2},2\phi^{-2}]=[2-\phi,1-\phi^{-3}]$.
\item $I_3=[2\phi^{-2},2-2\phi^{-2}]=[1-\phi^{-3},1+\phi^{-3}]$.
\item $I_4=[2-2\phi^{-2},2-\phi^{-2}]=[1+\phi^{-3},\phi]$.
\item $I_5=[2-\phi^{-2},2]=[\phi,2]$.
\end{itemize}
We define $R$ as follows.
\begin{itemize}
\item If $y \in I_1$ then $R(y)=\phi^3 y$.  Note that $R(I_1)=I_1 \cup ... \cup I_4$.
\item If $y \in I_2$ then $R(y)=y+\phi - \phi^{-2}$.  Note that $R(I_2)=I_5$.
\item If $y \in I_3$ then $R(y)=\phi^3 y-\phi^3+1$.  Note that $R(I_3)=I_1 \cup ... \cup I_5$.
\item If $y \in I_4$ then $R(y)=y-\phi+\phi^{-2}$.  Note that $R(I_4)=I_1$.
\item If $t \in I_5$ then $R(y)=\phi^3 y -2 \phi^3+2$.  Note that $R(I_5)=I_2 \cup ... \cup I_5$.
\end{itemize}
The map $R$ pieces together correctly on the endpoints of these
intervals, and induces a
degree $3$ covering map $R: \T \to \T$.

Note that $R$ preserves the ring $\Z[\phi]$,
because $\phi^k \in \Z[\phi]$ for all $k \in \Z$.
Note also that the second iterate $R^2$ is strictly expanding.

\subsection{The First Descent Lemma}
\label{descent}

As we mentioned above, $R$ preserves the set
\begin{equation}
\A=\Z[\phi] \cap [0,2]
\end{equation}
The goal of this section is to prove the following
number-theoretic result.

\begin{lemma}[Descent I]
The action of $R$ on $\A$ is the identity
mod $2\Z[\phi]$, and
for any $y \in \A$, there is some $n$ such
that $R^n(y)$ is one of:
$$0; 
\hskip 33 pt 1; 
\hskip 33 pt 
\phi, \hskip 8 pt 2-\phi;
\hskip 33 pt -1+\phi, \hskip 8 pt
-3+3\phi, \hskip 8 pt 3-\phi, \hskip 8pt 5-3\phi.$$
Each collection of points, grouped according to the
parity of their coefficients, forms a periodic cycle for $R$.
\end{lemma}

One can see directly from the formulas that $R$ is
the identity mod $2\Z[\phi]$.  We concentrate on the
proof of the second statement.

Given $y =m+n \phi \in \A$ we define
\begin{equation}
N(y)=\max(|m|,|n|).
\end{equation}

\begin{lemma}
If $y \in \A \cap (I_2 \cup I_4)$ then
$N(R(y)) \leq N(y)+2$.
\end{lemma}

\startproof
If $y=m+n\phi \in I_2$ then
$R(y)=(m-2)+(n+2)\phi.$
In this case the
result is obvious.  If $y \in I_4$, the proof is very similar.
\endproof

\begin{lemma}
\label{descent0}
If $y \in \A \cap (I_1 \cup I_3 \cup I_5)$ then
$N(R(y))<N(y)/2+8$.
\end{lemma}

\startproof
Suppose first that $y \in I_1$.  
Let $y=m+n\phi$.  We compute
$$R(y)=(1+2\phi)(m+n\phi)=(m+2n)+(2m+3n)\phi.$$
Now we observe the following two inequalities.
$$|m+2n|=|(m+\phi n) - (\phi-2)n| \leq 2+\phi^{-2}|n|<2+N(y)/2$$
$$|2m+3n|=|(2m+2 \phi n) - (2\phi-3)n| \leq 4+\phi^{-3}|n|<4+N(y)/2.$$
The last inequalities come from the fact that $y \in [0,2]$.
When $y$ lies in $I_3$ or $I_5$ the argument is the same, except that we have
must either add $2$ or $4$ to our estimate to account for the translational
part of $R$.
\endproof

\begin{corollary}
\label{descent1}
For any $y \in [0,2]$ there is some $n$ such that $N(R^n(y))<20$.
\end{corollary}

\startproof
Suppose that $y \in \A$ has the property that $N(y)>20$.
Note that $R(y)$ and $R^2(y)$ cannot both lie 
in $I_2 \cup I_4$.  Hence, we may combine our last two results to establish
the bound $N(R^2(y))<y/2+10<N(y).$
So $R^2(y)$ has smaller integer norm than does $y$.
\endproof

The above work reduces the proof of the Descent Lemma I
to a finite calculation of what happens to those
$y$ for which $N(y)<20$.  We make the finite
calculation and see that the Descent Lemma I holds
for these values as well.  Indeed, in the next
section, we will describe a much more extensive
calculation.

We mention a variant which we will need in one place.
Let $S$ be the involution $S(y)=2-y$.

\begin{lemma}
\label{descent3}
Let $y \in \Z[\phi]$ and let
$y'=R_1 \circ ... \circ R_N(y)$, where
each $R_k$ is either $R$ or $R \circ S$.
Then $y'$ is one of the values listed in
the Descent Lemma I provided that $N$
is sufficiently large.
\end{lemma}

\startproof
The map $R$ commutes with $S$ and all the
cycles listed in the Descent Lemma I are
invariant under $S$.  (The first cycle
is really $[0]=[2]$.)  The result follows
immediately from these two facts.
\endproof

\subsection{The Second Descent Lemma}
\label{sdl}

Let $G_2$ be the affine group defined in
\S \ref{renormXX}.  This group has two components,
so to speak. We let $G_2^+$ denote the index
$2$ subgroup consisting of maps of the form
\begin{equation}
T(x)=\phi^{3k} x + b; \hskip 30 pt k \in \Z;
\hskip 30 pt b \in 2\Z[\phi].
\end{equation}
This formula differs from the one in
Equation \ref{mainequiv} only in that we do not
allow a minus sign in front of $\phi^{3k}$.

\begin{lemma}[Descent II]
Two elements $y_1,y_2 \in \R/2\Z$ lie in the same
$G_2^+$ orbit iff there are positive
integers $n_1$ and $n_2$ such that
$R^{n_1}(y_1)=R^{n_2}(y_2)$.
\end{lemma}

\startproof
All the maps defining $R$ lie in $G_2^+$.  
So, if  there are positive
integers $n_1$ and $n_2$ such that
$R^{n_1}(y_1)=R^{n_2}(y_2)$ then
$y_1$ and $y_2$ lie in the same
$G_2^+$ orbit.  The converse is
the interesting direction.

For the converse, suppose that $y_1$ and $y_2$
lie in the same $G_2^+$ orbit.  Considering
the action of $R$, we can replace $y_2$ by
some image $y_2'=R^k(y_2)$ such that
$y_1-y_2' \in 2\Z[\phi]$.  So, without loss
of generlity, we can consider the case where
we already know that $y_1-y_2 \in 2\Z[\phi]$.
The Descent Iemma I proves this result whenever 
our points lie in $\Z[\phi]$.  So, it suffices to
consider the case when neither point belongs
to $2\Z[\phi]$.

We find it more convenient to work with a new
map that is closely related to $R$.
We define $\rho=R$ on $I_1 \cup I_3 \cup I_5$ and
$\rho=R^2$ on $I_2 \cup I_4$.  Unlike $R$,
which pieces together continuously across the
endpoints of the intervals in the partition,
the map $\rho$ is not defined on the endpoints.
However, the points we are considering, and
their orbits, never hit these endpoints.
It suffices to prove this result for $\rho$
in place of $R$.

To describe the map $\rho$, we let
$f[m,n]$ denote the map
\begin{equation}
x \to \phi^3 x + m + n \phi.
\end{equation}
Then
\begin{itemize}
\item On $I_1$, we have $\rho=f[0,0]$.
\item On $I_2$, we have $\rho=f[2,-2]$.
\item On $I_3$, we have $\rho=f[0,-2]$.
\item On $I_4$, we have $\rho=f[-2,-2]$.
\item On $I_5$, we have $\rho=f[0,-4]$.
\end{itemize}
This is a short calculation which (after many tries)
we did correctly.  We omit the details.

We write $\langle y_1,y_2 \rangle = \max(|m|,|n|)$, where
$y_1-y_2 = m+n\phi$.  The same argument as in
Lemma \ref{descent0} and Corollary \ref{descent1} shows
(with tons of room to spare) that
\begin{equation}
\label{descent2}
\langle y_1,y_2 \rangle >100 
\hskip 30 pt
\Longrightarrow
\hskip 30 pt
\langle y'_1,y'_2 \rangle<\langle y_1,y_2 \rangle-1.
\end{equation}
Here we have set $y'_k=\rho(y_k)$.
Equation \ref{descent2} reduces this lemma to
the case when $\langle y_1,y_2 \rangle \leq 100$.
In the next section, we will explain our computer-assisted proof that
The Descent Lemma holds for such choices.
\endproof

\subsection{A Dynamical Computation}
\label{micro}

Notice that our statement that the Descent Lemma II holds
for pairs $(y_1,y_2)$ with $\langle y_1,y_2 \rangle \leq 100$
is not {\it obviously\/} a finite calculation, because
it involves infinitely many values.  However, we will explain
how to reduce the problem to a finite calculation, which we
then make. This situation is typical of the results in this
paper.  The challenge is reduce seemingly infinite statements
to finite computations.

Let $Q=[0,2]^2$.  The set of pairs $(y_1,y_2)$ of interest
to us lie on a finite number of line segments of slope $1$
that are contained in $Q$.  By switching the order of
the two points if necessary, we can assume that $y_1<y_2$.
We let $\rho$ act on $Q$ by having $\rho$ act
separately on each coordinate.

The square $Q$ is partitioned into $25$ subsquares
\begin{equation}
Q_{ij}=I_i \times I_j
\end{equation}
on which $\rho$ is entirely defined and a similarity.
Our point of view is that $\rho$ acts separately
on each $Q_{ij}$, and the action on the various
boundaries depends on which square we include
the boundary in.  Anyway, we don't care about
what happens on the boundaries:  As we said above,
the points we consider, and their orbits, never
hit the boundaries.

Say that a {\it diagonal\/} is a segment of slope
$1$ contained in $Q$.  We call the diagonal
{\it small\/} if it is contained in one of the
$25$ subsquares. The endpoints of a small diagonal
might lie in the boundary of the subsquare, but
this is fine with us.  Given a small diagonal
$I$, we can define $\rho(I)$ using the action
of $\rho|_{Q_{ij}}$.  Note that $\rho(I)$ is
another diagonal, but not necessarily a small one.

Each diagonal has a canonical decomposition into
small diagonals:  We just take the intersections
with the $25$ sub-squares.  This we have a kind
of dynamical system defined on lists of
diagonals:  Given a list of diagonals, we
first subdivide each member of the list into
small diagonals.  Next, we let $\rho$ act
on all the small diagonals.  And so on.
As one final nicety, we switch the two
coordinates of the endpoints, if necessary,
so that all our diagonals lie about the
line $y_2-y_1=0$. (We do this simply for
computational convenience.)

Say that a {\it good seed\/} is a diagonal of the form
\begin{equation}
\Delta(m,n)=Q \cap 
\{y_2-y_1 = 2m+2n \phi\},
\end{equation}
where
\begin{equation}
2m+2n\phi \in [0,2]; \hskip 30 pt
\max(|2m|,|2n|) \leq 100.
\end{equation}
We run the dynamical system starting with any good seed
and we find that, after finitely many steps, the only
remaining intervals lie in
\begin{equation}
\Delta(0,0) \cup \Delta(-2,2) \cup \Delta(4,-2).
\end{equation}
Points $(y_1,y_2) \in \Delta(0,0)$ obviously satisfy
$y_1=y_2$.  When $(y_1,y_2) \in \Delta(-2,2)$, it means
that $y_2-y_1 =\phi-\phi^{-2}$.  In this case, we either
have $y_1 \in I_2$ and $y_2=R(y_1)$ or $y_2 \in I_4$ and
$y_1=R(y_2)$.  Compare the definition of the map $R$.

It only remains to deal with those points in
$\Delta(4,-2)$.  For this purpose, we just have
analyze the dynamics more carefully.
Define
\begin{equation}
\Delta'(4,-2) = \Delta(4,-2) \cap Q_{24}.
\end{equation}
When we perform the dynamics on
$\Delta(4,-2)$ we find that the following occurs.
\begin{itemize}
\item $\Delta(4,-2)$ breaks up into $5$ small diagonals,
one of which is $\Delta'(4,-2)$.
\item $\rho$ maps each of the $4$ other small diagonals
into $\Delta(0,0) \cup \Delta(-2,2)$.
\item $\rho$ maps $\Delta'(4,-2)$ back into $\Delta(4,-2)$.
\end{itemize}
This analysis shows that the Renormalization Lemma
can only fail for a pair of points
$(y_1,y_2)$ such that
$\rho^n(y_1) \in I_2$ for all $n$
(and also $\rho^n(y_2) \in I_4$ for all $n$.)
But the fixed point of $\rho_n|_{I_2}$ is
an endpoint of $I_2$ and $\rho$ is an expanding
map.  Since $y_1$ is not this endpoint, we
see that $\rho^n(y_1)$ eventually escapes
$I_2$, and we are done.
\newline
\newline
{\bf Remarks:\/}
\newline
(i)
We perform the calculations with exact arithmetic, as
explained in \S \ref{arithmetic}. \newline
(ii) In \S \ref{overflow} we explain how we eliminate any
possibility of overflow error in our calculations. 
Even without specific guards against overflow error (which
we do have) for all our seeds the dynamical system reaches the
$3$ end-states above very quickly and all integers
remain pretty small.  \newline

\subsection{The Cantor Set}
\label{cset}

In this section, we give some information about
the Cantor set $C$ from Theorem \ref{penrose2}.
First of all, the main property of $C$ is that
both $C$ and $C^{\#}$ are 
forward $R$-invariant.  That is,
$R(C)=C$ and $R(C^{\#})=C^{\#}$.
Indeed, this is how we discovered
the map $R$.    

\begin{lemma}
\label{mod1}
Let $y \in \Z[\phi]$.  Then
$y \in C^{\#}$ only if
$y \equiv 1$ mod $2\Z[\phi]$
and $y \in C-C^{\#}$ only if
$y=m+n\phi$ with $m$ even.
\end{lemma}

\startproof
Of the $8$ values listed in
the conclusion of the Descent Lemma I, we
see that only $1$ lies in
$C^{\#}$.  Given any $y \in \Z[\phi]$
which intersects $C^{\#}$, we simply
note that $R^n(y) \in C^{\#}$ for all
$n$. But there is some $n$ such that
$R^n(y)$ is one of the $8$ values listed
in the Descent Lemma I.  But this means
that $R^n(y)=1$.  Hence $y \sim 1$ mod $G_2^+$.
But this means that $y \equiv 1$ mod $2\Z[\phi]$.

For the second statement, we just have to rule
out the case that $y \equiv -1+\phi$
mod $2\Z[\phi]$.  But we check easily that
$-1+\phi \not \in C$.  But $C$ is forward
$R$-invariant and $R^n(y)=-1+\phi$ for
some $n$. This situation is impossible.
\endproof

To each point $y \in \R/2\Z$ we assign a
{\it renormalization sequence\/} in the digits
$\{1,2,3,4,5\}$.  The sequence is such that
the $n$th term is $k$ if and only if
$R^n(y) \in I_k$.  Not every point has a
unique renormalization sequence.  A
point $y$ has a non-unique renormalization
sequence if and only if $R^n(y) \in \partial I_k$
for some $n$ and some $k$.  Examining the
endpoints of our intervals and also the conclusion
of the Descent Lemma, we see that this happens
if and only if $y=m+n \phi$, where $m$ and $n$
are integers and $n$ is even.  In particular,
all points of $C^{\#}$ have unique
renormalization sequences.

\begin{lemma}
\label{no5}
A point lies in $C^{\#}$ if and only if
its renormalization sequence is unique and
has no $5$'s in it.
\end{lemma}

\startproof  
A point $y$ has a non-unique renormalization
sequence if and only if $R^n(y) \in \partial I_k$
for some $n$ and some $k$.  Examining the
endpoints of our intervals and also the conclusion
of the Descent Lemma, we see that this happens
if and only if $y=m+n \phi$, where $m$ and $n$
are integers and $n$ is even.  In particular,
all points of $C^{\#}$ have unique
renormalization sequences.

Note that the renormalization sequence of a
point in $y \in C^{\#}$ cannot start with
$5$, because $I_5=[\phi,2]$ only shares its
bottom endpoint with $C$, and this endpoint
lies in $C-C^{\#}$.  Since $C^{\#}$ is forward
invariant, we see that the renormalization
sequence cannot have any $5$'s in it at all.

Now we know that $C^{\#}$ only contains points
that have unique renormalization sequences
with no $5$'s in them.  Conversely, suppose
$y$ has a unique renormalization without $5$'s.
Suppose that $y \not \in C$. Since $y \not \in I_5$,
we can say that $y$ lies on one of the 
bounded components of $\R-C$. These components
all have diameter $\phi^{-3k+1}$ for $k=1,2,3...$
 The largest
component has size $\phi^{-2}$, and
is precisely the interior of $I_2$.
But $y \not \in I_2$ because then
$R(y) \in I_5$.  But, iteration
of $R$ increases the sizes of all
gaps except the largest one.
Hence $R^n(y)$ lies in the largest
gap for some $n$. This is a contradiction.

The endpoints of the gaps do not have
unique renormalization sequences.  So, the
same argument rules out the possibility
that $y \in C-C^{\#}$. We conclude that
$y \in C^{\#}$.
\endproof

\newpage
\section{The Fundamental Tiling}

\subsection{Definition of the Tiling}
\label{tiling1}

Let $\cal T$ be the tiling of the fundamental
triangle $T$ shown in Figure 1.2.  In this
section we define $\cal T$ precisely.
As we mentioned in the introduction, 
$T$ is bounded by three
sides of the Penrose kite $K$. 
The reader can most easily understand
the definitions we make by referring back
to Figure 1.2.

Let $f: \R^2 \to \R^2$ be the map that
fixes the top vertex of $T$ and shrinks
distances by a factor of $\phi^3$.  We
define $K_0=f(K)$, where $K$ is the Penrose kite.
$K_0$ is the largest kite in the tiling $\cal T$.
The largest octagon $J_0$ in $\cal T$ has
$4$-fold dihedral symmetry, and $3$ vertices located at
\begin{equation}
(13-8\phi,4-2\phi) \hskip 30 pt
(5-3\phi,6-3\phi) \hskip 30 pt (-3+2\phi,-2+2\phi).
\end{equation}
This is enough information to characterize $J_0$ uniquely.
A calculation shows that 
\begin{equation}
\label{subdivide}
T-K_0-J_0=(T_{11} \cup T_{12}) \cup T_{31} \cup T_{32} \cup T_{41}
\end{equation}
Here $T_{ij}$ is a similar copy of $T$, with scaling factor $\phi^{-3}$.
Our notation is such that $T_{jk} \subset I_j$, the interval used in the
definition of the map $R$. Just to pin thing down exactly, we say that
the centroid of $T_{i1}$ lies to the left of the centroid of $T_{i2}$.
The bottom triangles $T_{11}$ and $T_{12}$,
mirror images of each other, are not disjoint.
  However,  the similarities carrying $T$ to
$T_{11}$ and $T_{12}$ both map $J_0$ to the same smaller octagon.
Therefore, we can compatibly subdivide each of our smaller triangles into a
kite, an octagon, and a union of $5$ smaller triangles.  Continuing
this process forever, we get the tiling $\cal T$.

Let $S$ denote the union of points in
$T$ that do not belong to the interiors of any of the
tiles of $\cal T$.  Let $S^{\#}$ denote those points
of $S$ that do not belong to the boundarty of and
tile of $\cal T$.  Then $S=S^{\#}$ is contained in 
a countable union of lines defined over $\Z[\phi]$.
These are the lines extending the sides of 
the boundaries of the tiles in $\cal T$.
We call $S^{\#}$ the {\it fundamental fractal\/}.
We have the basic relation
\begin{equation}
\label{vertical}
C^{\#} = \pi_2(S^{\#}); \hskip 30 pt
C = \pi_2(S).
\end{equation}
Here $\pi_2$ is projection onto the second coordinate,
and $C$ is the Cantor set from Theorem \ref{penrose2}.

\subsection{Hausdorff Dimension}

The purpose of this section is to prove the following result.
Here $\dim(S)$ refers to the Hausdorff dimension of $S$.
See [{\bf F\/}] for details about Hausdorff dimension.

\begin{lemma}[Dimension]
\label{DIM}
Let $S$ be the fundamental fractal and
let $L$ be any countable collection of non-horizontal lines.
Then $\dim(S-L)=1$.  In particular, 
$\dim(S^{\#})=1$.
\end{lemma}

\startproof
After we remove $J_0$ and $K_0$ from $T$ we
are left with the union in Equation \ref{subdivide}.
We have the equation
\begin{equation}
T_{11} \cup T_{12} - R_{12}^{-1}(K_0)= T_{11} \cup T_{12}',
\end{equation}
where $T_{12}'$ is similar to
$T$, with similarity factor $\phi^{-6}$.
Figure 3.1 shows this operation.

\begin{center}
\resizebox{!}{1.5in}{\includegraphics{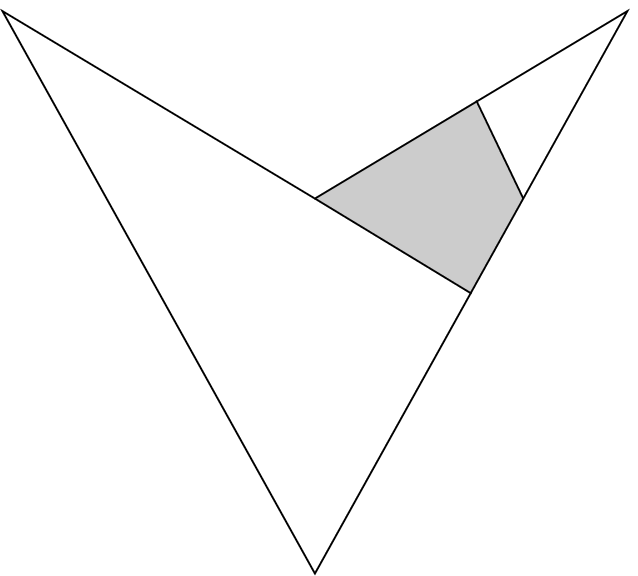}}
\newline
{\bf Figure 3.1:\/} $(T_{11} \cup T_{12})-R_{12}^{-1}(K_0)= T_{11} \cup T_{12}'$.
\end{center}

So, removing two open kites and an open octagon from
$T$ leaves the $5$ (temporarily renamed) disjoint triangles $T_1,...,T_5$.  We have
\begin{equation}
\sum_{k=1}^5 {\rm diam\/}(T_k)=\phi^{-3}(4+\phi^{-3}) {\rm diam\/}(T)={\rm diam\/}(T)
\end{equation}
The set $S$ is the self-similar fractal which is the limit set of the
semigroup generated by the similarities which carry $T$ to each of
the smaller triangles.
From here it is an exercise to show that
$\dim(S)=1$.

Note that $C$, the Cantor set from
Theorem \ref{penrose2} is the set $y \in [0,2]$ such that
the line of height $y$ intersects $S$.
Each line in $L$ intersects $S$ in a set
of dimension $\dim(C)<1$.  Since
$\dim(S)=1$ and $\dim(S \cap L)<1$, we
have $\dim(S-L)=1$.
\endproof

\subsection{The Horizontal Intersections}
\label{horizontal}

Let $C^{\#}$ be as in the previous chapter.  Every
$y \in C^{\#}$ is such that 
\begin{equation}
L_y=\R \times \{y\}
\end{equation}
intersects $S^{\#}$ nontrivially.
The goal of this chapter is to prove the following result.

\begin{lemma}[Horizontal]
For any $y \in C^{\#}$, the set $L_y \cap S$ is
a Cantor set. Hence $L_y \cap S^{\#}$ is obtained
from a Cantor set by removing at most countably
many points.
\end{lemma}

We begin by explaining how the set $S$ may be
constructed recursively, as the nested intersection
of finite unions of (overlapping) triangles.

Recall that $T-J_0-K_0=\bigcup T_{ij}$, a union of
$5$ triangles.  Let $\rho_{ij}$ be the similarity
such that $\rho_{ij}(T_{ij})=T$.   Define
$S_0=T$ and (inductively) $S_n$ such that
\begin{equation}
\rho_{ij}(S_n \cap T_{ij})=S_{n-1}.
\end{equation}
Then $S_n$ is a finite union of triangles and
\begin{equation}
S=\bigcap S_n.
\end{equation}
Our proof involves an analysis of how
these triangles sit with respect to the
horizontal lines.
\newline
\newline
{\bf Remarks:\/} \newline
(i)  Note the
similarity between the maps $\rho_{ij}$ and
the map $\rho$ considered in the previous chapter.
precisely, $\rho_{ij}$ acts on the horizontals
having heights in $I_i$ exactly as $\rho$ acts
on $I_i$.  We will pursue this analogy further
in the next section. \newline
(ii)
We have to be careful with our definition,
because $T_{11}$ and $T_{12}$ overlap.  Here
is the justification for what we do. Assume by
induction that $S_{n-1}$ is well-defined and
that $S_{n-1} \cap (T_{11} \cap T_{12})$
has bilateral symmetry and that $\rho_{11}$ and $\rho_{12}$
have the same action on this intersection.
Then $S_n$ is well-defined and inherits all
these same properties by symmetry.
\newline

Let $\{k_n\}$ be the renormalization sequence
associated to $y$.  By Lemma \ref{no5} this
sequence and has no $5$'s.
Since $S_n$ is a finite union of triangles,
the intersection

\begin{equation}
\Lambda_n=L_y \cap S_n
\end{equation} 
is a finite union of disjoint intervals.

\begin{lemma}
\label{claim0}
Let $J$ be a maximal interval of $\Lambda_n$.
Then $J \cap \Lambda_{n+1}$ is nonempty and
contains $2$ maximal intervals
in case $k_{n}=3$ or 
$(k_n,k_{n+1})=(1,4)$.
\end{lemma}

\startproof
Let $y_0=y$ and
$y_n=R^n(y)$ and $L(n)=L_{y_n}$.
Our sequence $\{k_n\}$ starts with $k_0$.   We have
$y_n \in I_{k_n}$.

Recall that $S_0=T$.
It follows from induction and fact that
$\rho_{ij}(T_{ij})=T$ 
that there is a sequence of maps
$\rho_0,...,\rho_{n-1}$ such that
\begin{equation}
\label{boot2}
\rho_{n-1} \circ ... \circ \rho_0(J \cap \Lambda_{n+m})=L(n) \cap S_m; \hskip 30 pt
m=0,1,...
\end{equation}
Here $\rho_j=\rho_{k_j,m_j}$, where $m_j \in \{1,2\}$.
The sequence $\{\rho_i\}$ is not necessarily unique,
because $T_{11}$ and $T_{12}$ overlap. We don't mind this.

Since $k_n \in \{1,3,4\}$, the set on the right hand side of
Equation \ref{boot2} is a
nontrivial union of intervals.   Hence $J \cap \Lambda_{n+1}$ has at
least one interval.  When $k_n=3$ the set
$L(n) \cap S_1$ contains $2$ intervals.  Hence
$J \cap \Lambda_{n+1}$ contains $2$ intervals in this case.
When $k_{n}=1$ and $k_{n+1}=4$ the line $L(n)$ lies
above the top vertex of $T_{11} \cap T_{12}$.  Hence
$L(n) \cap S_1$ again consists of $2$ intervals.  Hence,
so does $J \cap \Lambda_{n+1}$.
\endproof

\begin{corollary}
\label{claim1}
Let $J$ be a maximal interval of $\Lambda_n$.  Then
$J \cap \Lambda_m$ contains at least $2$ disjoint intervals for $m$ sufficiently large.
\end{corollary}

\startproof
The sequence associated to $y$ cannot terminate in an infinite
string of $1$'s.  Otherwise, there is another sequence associated
to $y$ which terminates in an infinite string of $5$'s.  So,
the associated sequence either has an infinite number of
$3$'s or an infinite number of $(1,4)$'s. The corollary
now follows immediately from Lemma \ref{claim0}.
\endproof

It follows from Corollary \ref{claim1} that
the nested intersection $\bigcap \Lambda_n$ is a
Cantor set.   This completes the proof of
the Horizontal Lemma.

\subsection{Notation}
\label{notation}

We have already mentioned the similarities 
$\rho_{ij}$.  These maps have the property that
$\rho_{ij}(T_{ij})=T$.  These maps would seem to
suit us perfectly well, but it turns out that there
is a slightly more elaborate collection of maps that
are better adapted to the structure of outer billiards.
We introduce these maps here.  These maps have the
advantage that they are all homotheties.

Define the reflection in the vertical line $x=1$:
\begin{equation}
\eta(x,y)=(2-x,y).
\end{equation}
We first change our notation a bit.  We let
$T^+=T$ and $T_{ij}^+=T$, etc.
That is, we attach the $(+)$ superscript to
all the objects associated to the fundamental
triangle.   Next, we define
$X^-=\eta(X^+)$ for any object $X$.

The two triangles
$T^+$ and $T^-$ are mirror images of each other.
We will see eventually that the dynamical tiling
intersects $T^-$ exactly in the tiling ${\cal T\/}^-$.
It turns out that the dynamics on $T^+ \cup T^-$ 
works out more nicely than the dynamics on just $T^+$,
even though ultimately all our results are phrased 
just in terms of $T^+$.

We define the following ten maps.
\begin{itemize}
\item $R_{11}^+=\eta_{11}: T_{11}^+ \to T^+$.
\item $R_{12}^+=\eta \circ \rho_{11}: T_{11}^+\to T^-$.
\item $R_{31}^+=\eta \circ \rho_{31}: T_{31}^+ \to T^-$.
\item $R_{32}^+=\eta_{32}: T_{32}^+ \to T^+$.
\item $R_{41}^+$ is the isometry carrying $T_{41}^+$ to $T_{12}^+$. 
\item $R_{ij}^-=\eta \circ R_{ij}^+ \circ \eta$.
\end{itemize}

The maps $R_{ij}^-$ acts similarly to the map
$R_{ij}^+$.  For instance $R_{31}^-(T_{31}^-)=T^+$.

Recall that the renormalization map $R$ equals
the map $R_i$ on the interval $I_i$.
The maps $R_{ij}^{\pm}$ are all similarities
preserve the horizontal foliation and act on the
horizontal lines intersecting $T_{ij}^{\pm}$ as $R_i$ acts
on $I_i$.  In particular, the maps $R_{41}^{\pm}$ are
better adapted to $R_4$ than
the map $\rho_{41}$, which is really adapted
to the map $R_1 \circ R_4$.

The maps $R_{41}^{\pm}$ are orientation
preserving isometries.   The remaining $8$ maps are
orientation reversing similarities, with expansion
constaint $\phi^3$.

\subsection{The Renormalization Set}
\label{rset}

Given subsets $A,B \subset T^+ \cup T^-$, we write
$A \to B$ if $A \subset T_{ij}^{\pm}$ and 
$B=R_{ij}^{\pm}(A)$.
We use this definition in particular for points.
Suppose that $(p_1,q_1)$ is a pair of
points, both at the same height.  We write
$(p_1,q_1) \to (p_2,q_2)$ if $p_1 \to p_2$
and $q_1 \to q_2$.
We write $p \sim q$ if 
\begin{equation}
(p,q) \leadsto ... \leadsto(p',q'); \hskip 30 pt p'=q'.
\end{equation}
Define
\begin{equation}
\Upsilon(p)=\{q|\ q \sim p\}.
\end{equation}
Note that $\Upsilon(p)$ consists of points that are all on the same
horizontal level as $p$.
We call $\Upsilon(p)$ the {\it renormalization set\/} of $p$.
\newline
\newline
\noindent
{\bf Remark:\/}  
In \S \ref{analy}, we will
see that $p$ and $q$ lie in the same outer billiards orbit
provided that $p$ and $q$ both have unbounded
orbits and $p \sim q$. This is a step in our
proof that generic unbounded orbits are self-accumulating.
\newline

Let $C^{\#}$ be the set studied in the previous chapter.
The goal of this section is to prove the following result.

\begin{lemma}[Density]
Let $\Lambda^{\pm} = L \cap S^{\pm}$, where $L$ is a horizontal
line whose height lies in $C^{\#}$.
For any point $p \in \Lambda^+ \cup \Lambda^-$, the set $\Upsilon(p)$ is
dense in $\Lambda^+ \cup \Lambda^-$.
\end{lemma}

the Density Lemma.
We call $T^+$ and $T^-$ the {\it distinguished
triangles\/} of depth $0$. 
Recall that $S^{\pm}$ is contained in the nested intersection
of sets $S_n^{\pm}$. 
We say that a {\it distinguished triangle\/} is a
maximal triangle of $S_n^{\pm}$.  There
are $2 \times 5^n$ distintuished triangles.
We call $n$ the {\it depth\/} of the distinguished
triangle. 

We say that a distinguished triangle is {\it related\/} to a
renormalization set if it intersects the
horizontal line containing the renormalization set.
Let $P(n)$ be the property that every renormalization
set intersects each related distinguished triangle 
depth $n$.  It suffices to prove that the statement $P(n)$ is
true for all $n$.

If $\tau$ is any distinguished triangle of depth $n$,
then $\tau \to \tau' \to \tau''$ where
one of $\tau'$ or $\tau''$ has depth $n-1$.
From this, and from the definitions,
we see that $P(n-1)$ implies $P(n)$.
It just remains to establish $P(0)$.

Say that points $p,q \in T^+ \cup T^-$ are
{\it distantly placed\/} if one of the points
lies in $T^+$ and the other lies in $T^-$.
Say that $y \in C^{\#}$ is {\it good\/} if $P(0)$ holds
for all renormalizations sets of height $y$.

\begin{lemma}
Suppose that $y_1 \in R^{-1}(y_2)$ and $y_2$ is good.  
Then $y_1$ is good. 
\end{lemma}

\startproof
Let $p_1$ be some point having height $y_1$.
We have $p_1 \to p_2$ for some $p_2$ having
height $y_2$.   Since $y_2$ is good,
there is some $q_2$
such that $p_2 \sim q_2$ and $p_2,q_2$
are distantly placed.
Without loss of generality, assume that $p_2 \in T^+$.
There is a depth $1$ distinguished triangle
$\tau^+$ such that $\tau^+ \to T^+$
and $p_1 \in \tau^+$.  Let
$\tau^-=\rho(\tau^+)$.  Then
$\tau^- \to T^-$ and points in
$\tau^-$ are distantly placed from
points in $\tau^+$.   In particular,
we can find $q_1 \in \tau^-$
such that $q_1 \to q_2$.  But
then $p_1$ and $q_1$ are distantly placed
and $p_1 \sim q_1$.  Since $p_1$ was
chosen arbitrarily, $y_1$ is good.
\endproof

\begin{lemma}
\label{good1}
Suppose that $y \in C^{\#} \cap I_3$.  Then $y$ is good.
\end{lemma}

\startproof
Any horizontal line having a height in $I_3$ intersects
the $4$ disjoint triangles $T_{3j}^{\pm}$ for $j=1,2$.
Let $\rho$ be the union of the $4$ special maps.
Then 
$$\rho(T_{31}^+)=\rho(T_{32}^-)=T^+; \hskip 30 pt
\rho(T_{32}^+)=\rho(T_{31}^-)=T^-.$$
It follows from this equation that every
renormalization set of height $y$ either intersects both
of $(T_{31}^+,T_{32}^-)$ or both
of $(T_{32}^+,T_{31}^-)$.
\endproof

\begin{lemma}
\label{good2}
Suppose that $y \in C^{\#} \cap I_1 \cap R^{-1}(I_4)$.
Then $y$ is good.
\end{lemma}

\startproof
Any horizontal line having a height in $I_3$ intersects
the $4$ triangles $T_{1j}^{\pm}$ for $j=1,2$,
in disjoint intervals.  (This is true even though the
triangles themselves are not disjoint.)  The rest
of the proof is as in Lemma \ref{good1}.
\endproof

Let $y \in C^{\#}$ be arbitrary.  
As in the proof of Corollary \ref{claim1}, there
must be some $n$ such that $R^n(y)$ satisfies
either Lemma \ref{good1} or Lemma \ref{good2}.
But then $R^n(y)$ is good.  But then $y$ is
good as well.  Hence, $P(0)$ holds.  But then
$P(n)$ holds for all $n$.
This completes the proof of the Density Lemma.

\newpage

\section{Preliminaries}

\subsection{Polytope Exchange Maps}
\label{pem}

\noindent
{\bf Definition:\/}
For us, a {\it polytope exchange map\/} is a quadruple
$(M,X_1,X_2,\Psi)$, where $M$ is a flat manifold (possibly with boundary),
$X_1$ and $X_2$ are locally finite partitions of $M$ into convex polytopes,
and $f: M \to M$ is a piecewise isometric bijection which carries $X_1$ to $X_2$.
We mean that $f$ is a translation when restricted to each polytope
$P$ of $X_1$ and $f(P)$ is a polytope of $X_2$.  Technically, $f$ is not
defined on the boundaries of the polytopes.
When $M$ is a compact manifold, we require that the partitions be finite.

As a special case, suppose that
 $M=\R^n/\Lambda$, where $\Lambda \subset \Z^n$
is a discrete group of translations.  We say that a polytope in
$M$ is {\it golden\/} if any lift to $\R^n$ has all vertices
with coordinates in the ring $\Z[\phi]$.  This definition is
independent of lift. 
We call an associated polytope exchange map {\it golden\/} if all the
polytopes are golden, and if all the translations are defined
by vectors in $(\Z[\phi])^n$.
In this paper we will consider two golden polytope exchange maps.
\begin{itemize}
\item A $2$ dimensional non-compact polygon exchange $\Psi$, whose
domain is the infinite strip $\Sigma=\R \times [-2,2]$.
\item A $3$ timensional compact polyhedron exchange $\widehat \Psi$, whose
domain is the torus $\widehat \Sigma=(\R/2\Z)^3$.  
\end{itemize}

\noindent
{\bf Fibered Polyhedron Exchange Maps:\/}
Let $\cal H$ be the foliation of
$\widehat \Sigma$ by horizontal $2$-tori -- those obtained by
holding the third coordinate constant.
We say that a leaf of $\cal H$ is {\it golden\/} if
its height lies in $\Z[\phi]$. We call a golden
polyhedron exchange map on $\widehat \Sigma$ {\it fibered\/} if
it preserves the leaves of $\cal H$ and if
the restriction to each golden leaf is a
golden polygon exchange map. The second condition
does not follow automatically from the first;
there is an auxilliary condition we need on
the edges of the polyhedron partition, as we
now explain.

Let $e$ be a non-horizontal edge of some polyhedron in the
partition.  Let $(x_1,y_1,z_1)$ and $(x_2,y_2,z_2)$
be the endpoints of $e$. The $4$ quantities
\begin{equation}
\label{fibered}
\frac{x_2-x_1}{z_2-z_1}, \hskip 15 pt
\frac{y_2-y_1}{z_2-z_1}, \hskip 15 pt
\frac{x_2z_1-x_1z_2}{z_2-z_1}, \hskip 15 pt
\frac{y_2z_1-y_1z_2}{z_2-z_1}
\end{equation}
all belong on $\Z[\phi]$ if and only if every
nontrivial intersection of $e$ and a golden leaf
has coordinates in $\Z[\phi]$.  We
omit the easy proof of this result.
\newline
\newline
{\bf Periodic Tiles:\/}
Let $(M,X_1,X_2,f)$ be a polytope exchange map, not necessarily
a golden one.  Given a periodic point
$p \in M$, there is a maximal convex polytope $D_p$ consisting
of points $q$ such that $p$ and $q$ have the same period and
the same itinerary.  By {\it itinerary\/} we mean the sequence
$\{n_k\}$, defined so that the $(n_k)$th polytope $P(n_k)$ of $X_1$
contains $f^k(p)$.   We have
\begin{equation}
D_p=\bigcap_k f^{-k}\big(P(n_k)\big)
\end{equation}
Being the finite intersection of convex polytopes, $D_p$ is also a convex
polytope.  We call $D_p$ a {\it periodic tile\/}.   We call the
union of all the periodic tiles the {\it dynamical tiling\/}.
When all points of $M$ are periodic, the dynamical tiling is another
partition of $M$. 

Our terminology is slightly misleading.  In general, the
dynamical tiling only fills a subset of $M$.  This subset
need not even be dense in $M$.  A nice example is
furnished by outer billiards on a generic trapezoid.
See [{\bf G\/}], where the issue is phrased somewhat
differently.  In our examples, thanks to the
results in this paper, the dynamical tilings are dense.
\newline
\newline
{\bf Homogeneity:\/}
As one might expect, the dynamical tiling has some nice
homogeneity properties.
We say that two points $x_1,x_2 \in M$ are {\it tile isometric\/} 
if there are open disks $\Delta_1$ and $\Delta_2$ such that
\begin{itemize}
\item $x_j \in \Delta_j$ for $j=1,2$.
\item ${\cal D\/} \cap \Delta_1$ is isometric
 to ${\cal D\/} \cap \Delta_2$.
\end{itemize}
Note that $x_j$ need not be the center of $\Delta_j$ and the isometry
from from ${\cal D\/} \cap \Delta_1$ to ${\cal D\/} \cap \Delta_2$
need not carry $x_1$ to $x_2$.
We call a subset $S \subset M$ {\it tile homogeneous\/} if every pair
of points in $M$ are tile isometric.

As a related definition, we say that $x_1,x_2 \in M$ are
{\it tile similar\/} if the above properties hold, with the
word {\it similar\/} replacing the word {\it isometric\/}.

\begin{lemma}
Every orbit is tile homogeneous.
\end{lemma}

\startproof
Let $p_0$ and $p_n=f^n(p_0)$ be two points in the orbit.
Let $g=f^n$.  The map $g$ is an isometry in a neighborhood of $p_0$
and $g$ conjugates $f$ to itself.  It follows from this fact that
$p_0$ and $p_n$ are tile isometric.
\endproof

\subsection{The Return Map}
\label{sob}
\label{returnmap}

 The square outer billiards map
$\psi$ is a local translation.
For points far from the origin, the dynamics of $\psi$ is quite simple.
forward iterates of $\psi$ generally circulate counterclockwise around
the kite $K$.   Near the origin, the dynamics of $\psi$ is slightly more
complicated.  Our paper [{\bf S3\/}] considers these dynamics
carefully, for fairly general convex polygons.

\begin{center}
\resizebox{!}{4.5in}{\includegraphics{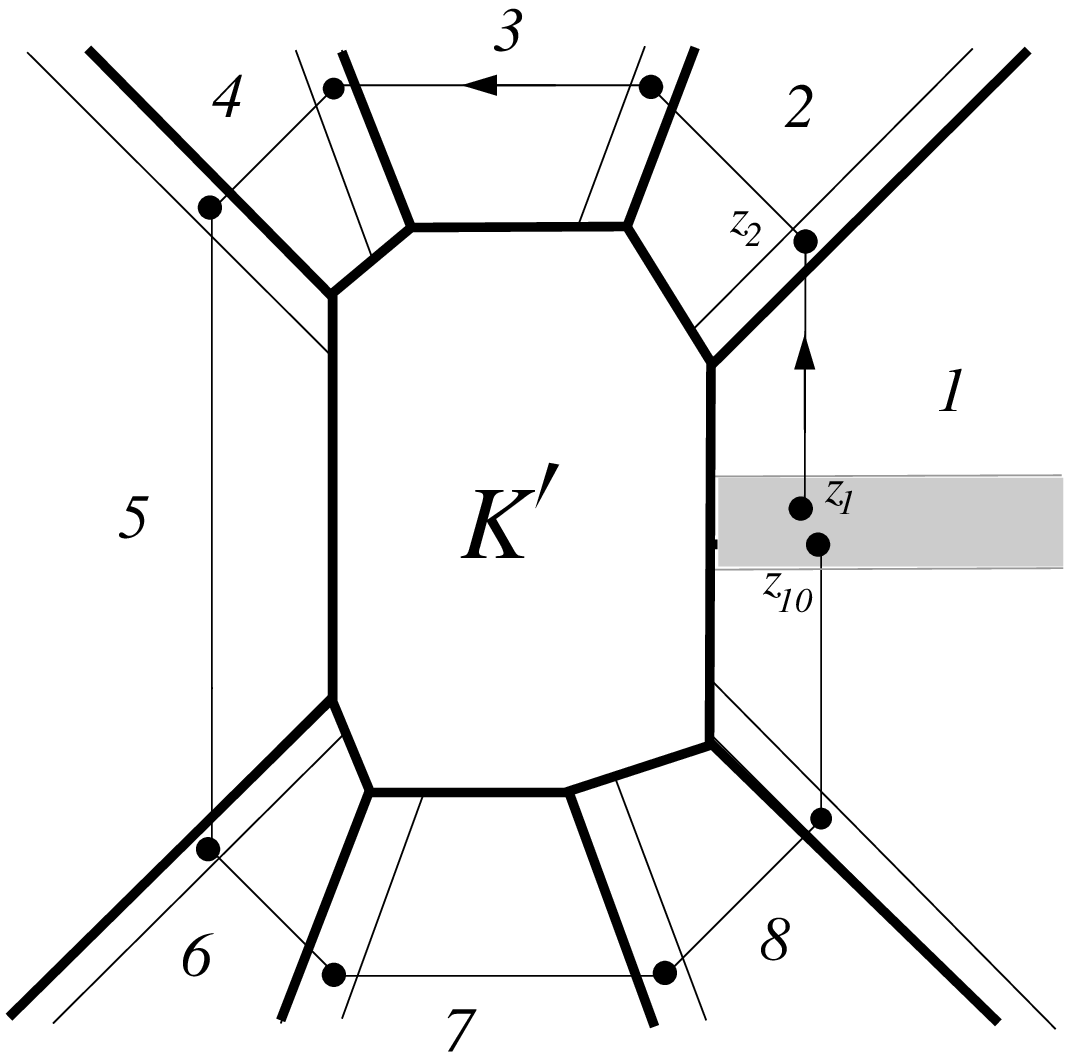}}
\newline
{\bf Figure 4.1:\/} The second return map far from the origin
\end{center}

The set $K'$ shown in Figure 4.1 is a large compact set that
contains $K$ well in its interior.
Given that the $\psi$-orbits generally circulate around the kiet $K$, at
least far from $K$,
it makes sense to consider the return map to each half
of a suitable horizontal strip.  Define
\begin{equation}
\Sigma=\R \times [-2,2].
\end{equation}
half this strip is shaded in Figure 4.1.
What makes this strip canonical is that, far from
$K$ and near the $x$-axis, consecutive iterates
of $\psi$ differ by the vector $(0,\pm 4)$.  So,
$\Sigma$ has just the right width.
\newline
\newline
{\bf Remark:\/}
Sometimes we will want to leave off the bottom
boundary of $\Sigma$ and sometimes we won't.
It turns out that the first return map is the
identity on the boundary of $\Sigma$, so we
can essentially just ignore these points.
\newline

The definition of the return map to ``each half''
of $\Sigma$ presents some problems
for points near $K$.
There are ``bad points'' of $\Sigma$ that are
too close to $K$.    These bad points are in the 
regions labelled $B$ in Figure 4.2.

\begin{center}
\resizebox{!}{2.5in}{\includegraphics{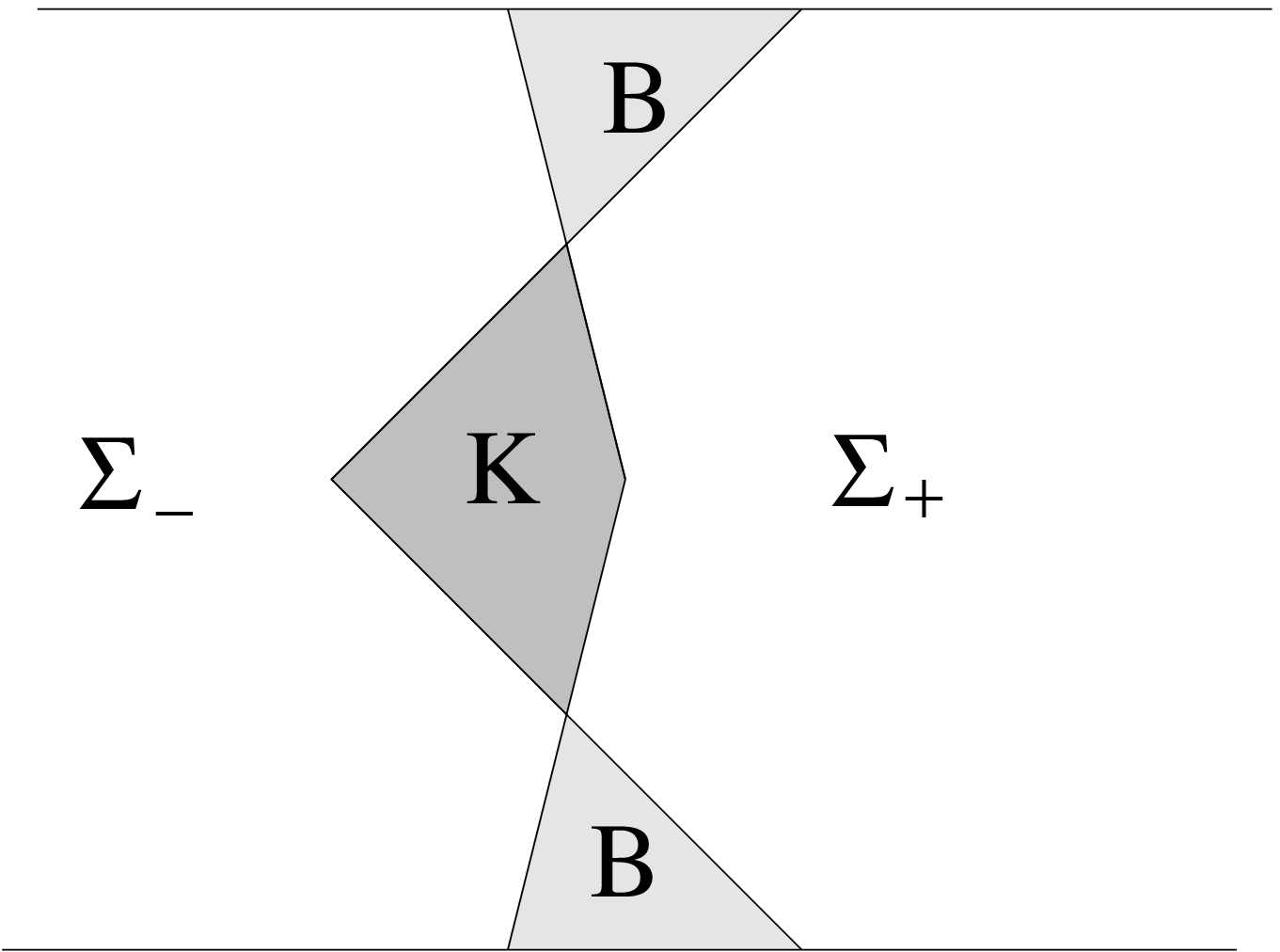}}
\newline
{\bf Figure 4.2:\/} A decomposition of $\Sigma$ into regions.
\end{center}

The problem with the bad points is
that one cannot really say which side of $K$ they are on.
We mention this problem in order to justify a modified version
of the return map, which we now define.
Let $\Sigma_-$ and $\Sigma_+$ be the two components of
$\Sigma-K-B$.   Let
$\Psi_{\pm}: \Sigma_{\pm} \to \Sigma_{\pm}$ to be the first return map
of $\psi$ to $\Sigma_{\pm}$.  Just to be clear,
$\Psi_+$ is the first return to $\Sigma_+$ and $\Psi_-$ is the
first return to $\Sigma_-$. We define $\Psi=\Psi_+ \cup \Psi_-$ to
be the ``union map'', defined on
$\Sigma_+ \cup \Sigma_-$.  Finally, to get a map on all
of $\Sigma$, we define
$\Psi$ to be the identity on $K \cup B$.  Note that this definition
does not correspond with the action of outer billiards on $B$,
but we will deal specially with the points in $B$ whenever necessary.

We call $\Psi: \Sigma \to \Sigma$ the {\it return map\/}.

\begin{lemma}[Return]
\label{return}
The following is true.
\begin{enumerate}
\item Every $\psi$ orbit intersects $\Sigma_+ \cup \Sigma_-$.
\item $\Psi$ is well defined on all points of $\Sigma$ that have a well-defined orbit.
\item There is some $C$ such that $|\Psi(p)-p|<C$ for all $p \in \Sigma$
with a well defined orbit.
\end{enumerate}
\end{lemma}

\startproof
First of all, we proved the same result
in [{\bf S2\/}, \S 2.3], in the context of special
orbits (on $\R^2_1$) for arbitrary kites.  The proof there
works in this setting with only minor changes.

Here we explain a different proof.  The return Lemma is
obvious for points of the form $(x,y)$ with $|x|>20$.
These points just circulate around the kite, nearly
following a giant octagon before coming back to the
strip.  For points in the strip $\Sigma$ having
$|x| \leq 20$, the calculation we make in connection
with the Pinwheel Lemma in \S \ref{pinwheelmap}
in particular establishes the Return Lemma.
The calculation we make there
simply involves covering $[-20,20] \times [-2,2]$
with $572$ convex polygons such that the first
return map exists and is well defined on the interior
of each tile.  We also prove that the remaining
points, the ones in the boundaries of our tiles,
do not have well defined orbits.
\endproof

The Return Lemma allows us to work with $\Psi$ rather than $\psi$.
We state the following result in terms of the action of
$\Psi^+$ on $\Sigma^+$.  The same result holds with
$(-)$ in place of $(+)$.

\begin{lemma}
\label{qi1}
There is a canonical bijection between the unbounded $\psi$ orbits
and the unbounded $\Psi^+$ orbits.  The bijection is such that a
given $\Psi$ orbit corresponds to the $\psi$ orbit that contains it.
Two $\Psi^+$ orbits are coarse equivalent if and only if the
corresponding $\psi$ orbits are coarse equivalent. 
\end{lemma}

\startproof
We associate to each $\Psi^+$ orbit the unique
$\psi$ orbit that contains it.  This injective association
is also surjective, by Statement 1 of the Return Lemma.
Given our description of $\Psi$ for points far from the
origin, it is clear that we can reconstruct the coarse equivalence
class of a $\psi$-orbit from the coarse equivalence
class (defined the same way) for the corresponding
$\Psi^+$ orbit. Up to a uniformly bounded error,
we obtain the graph of the $\psi$ orbit
from the graph of the $\Psi^+$ orbit by attaching a
centrally symmetric octagon of a suitable radius to
each point of $\Gamma$ that is sufficiently far from
the origin.  
\endproof

\subsection{An Unboundedness Criterion}
\label{unbounded criterion}

Here we establish a useful criterion for
the unboundedness of the $\Psi$ orbits.
Let $\Z_+$ denote the set of positive
integers and let $\Z_-$ denote the set of
negative integers.  We say that a subset
$S \subset \Z$ is {\it uniformly dense\/} in
$\Z_+$ if there is some $N$ such that every
point of $\Z_+$ is within $N$ units of $\Z_+$.
We make the same definition relative to
$\Z_-$.

\begin{lemma}
\label{UNB}
Let $O$ be an infinite $\Psi$ orbit.  Suppose that
there is an open horizontal line segment $S$ such that
$O \cap S$ is a nonempty and nowhere dense subset
of $S$.  Then $O$ is unbounded in both directions.
\end{lemma}

\startproof
We will assume that $O$ is bounded in the forward
direction and derive a contradiction.
Let $A=\phi^{-3}$.  Given the locations of the vertices of
the Penrose kite, we have the following formula.
\begin{equation}
\label{pregraph}
\Psi(p)-p=(2n_1A+2n_2,2n_3); \hskip 30 pt n_1+n_2+n_3 \equiv 0 \hskip 5 pt {\rm mod\/} 
\hskip 5 pt 2.
\end{equation}
When $p$ is far from the origin, the orbit of $p$ stays within a uniformly thin
tubular neighborhood of a centrally symmetric octagon, as we mentioned
in connection with Statement 3 of the Return Lemma.  The sides of
this octagon, which depends on $p$, are always integer multiples
of the vectors listed in \S \ref{sob}.    Moreover, opposite sides have
the same length.   For this reason, the vectors entering into the
sum that defines $\Psi(p)-p$ nearly cancel in pairs, and we
find that there is a uniform bound to $\max(|n_1|,|n_2|)$
in Equaton \ref{pregraph}.

Let $(x_0,y_0)$ be some point of $O$.
We can find integers $(a_n,b_n)$ such that
\begin{equation}
\Psi^n(x_0,y_0)=(x_0,y_0)+2a_n A + 2b_n.
\end{equation}
The differences $|a_{n+1}-a_n|$ and $|b_{n+1}-b_n|$ are
uniformly bounded.

The sequence $\Omega=\{a_n A+b_n\}$ is both infinite
and bounded.  Hence, $\Omega_1=\{a_n\}$ has infinitely
many values.  Our uniform bound on
$|a_{n+1}-a_n|$ now implies that
$\Omega_1$ is uniformly dense
in at least one of $\Z_+$ or $\Z_-$.
Combining this fact with the fact that
$\Omega$ is bounded, we see that
there is some $N$ such that the 
union 
\begin{equation}
\Omega^*=\bigcup_{i,j<N} \Omega+(i,j)
\end{equation}
either contains every integer combination of
the form $aA+b \in (0,1)$ with $a>0$,
or every such integer combination with $a<0$.
In either case, $\Omega^*$ is
dense in $(0,1)$.

But $\Omega^*$ is a finite union of translates
of $\Omega$.  We have shown that a finite union
of translates of $\Omega$ is dense in $(0,1)$.
But the set $x_0+2\Omega$ is a subset of $O$.
Hence, a finite union of translates of $O$ is
dense in some line segment.   Since
$O$ is contained in the union of two lines, and
the finite union of nowhere dense linear subsets is
again nowhere dense,
$O$ intersects some line segment
in a set that is not nowhere dense.
Since the map $\Psi$ is a piecewise translation, 
$O \cap S$ is not nowhere dense in $S$.
This contradiction finishes our proof.
\endproof

\subsection{The Arithmetic Graph}
\label{agdefined}

The arithmetic graph gives us a good way to visualize the
orbits of the first return map $\Psi: \Sigma \to \Sigma$.
Unlike our earlier papers
[{\bf S1\/}] and [{\bf S2\/}], the arithmetic graph
does not play an important role in our proofs. However,
we find it very useful to illustrate some concepts
with the arithmetic graph.  

For $y \in (0,2)$ we define
\begin{equation}
\label{hor2}
\Sigma_y=(\R \times \{y\}) \cup (\R \times \{y -2\}) \subset \Sigma
\end{equation}
The set $\Sigma_y$ is a union of $2$ horizontal lines.
(The case $y=0,2$ leads to a trivial picture, and we ignore it.)
In light of Equation \ref{pregraph}, it makes sense to try to understand
our orbits in terms of the triples of integers $(n_1,n_2,n_3)$
rather than in terms of the values $(2n_1A+2n_2,n_3)$. 
The parity of $n_1+n_2$ determines $n_3$.   Thus, the pair of
integers $(n_1,n_2)$ really determines the behavior of
$\Psi$ on $p$.  

Fixing some $y \in (0,2)$, we define a map
$M_y: \Z^2 \to \Sigma_y$ as follows.
We first choose some {\it offset value\/} $\xi \in \R$.
This offset value selects which orbits we focus on.
Given a point $(m,n) \in \Z^2$ we define
\begin{equation}
M_{y,\xi}(m,n)=\big(2mA+2n+\xi,y+\tau(m,n)\big).
\end{equation}
Here $\tau(m,n)=0$ if $m+n$ is odd and $\tau(m,n)=-1$ if $m+n$ is even.

It is a consequence of Equation \ref{pregraph} that
$\Psi$ preserves the image $M_{y,\xi}(\Z^2)$.  That is, given
$(m_0,n_0) \in \Z^2$, there is another $(m_1,n_1) \in \Z^2$ such that
$\Psi \circ M(m_0,n_0) = M(m_1,n_1)$.   
We write $(m_0,n_0) \to (m_1,n_1)$ in this case. Here we have set
$M=M_{y,\xi}$.

Given $(y,\xi)$, we form the arithmetic graph
$\Gamma(y,\xi)$ as follows:  We include $(m,n)$ as
a vertex of $\Gamma$ if and only if $M(m,n)$ has
a well defined orbit.  We connect the vertex
$(m_0,n_0)$ to the vertex $(m_1,n_1)$ if and only if
$(m_0,n_0) \to (m_1,n_1)$.
In this way, we produce 
a $2$-valent directed graph in the plane, whose vertices
lie in $\Z^2$. 
For kites, the graph $\Gamma$ is always embedded.
We proved this for the special orbits (i.e. those on
$\Sigma_1$) relative to any kite parameter in [{\bf S2\/}].
The proof for the general orbit is similar, though we have
not written down the details.

Figure 4.3 shows a portion of the unique unbounded
$\gamma$ component of the graph $\Gamma(1,\phi^{-2})$.
(The straight line segment at the bottom is just for reference.)
The straight line is $M^{-1}(0)$.
In [{\bf S1\/}] we studied this component and
showed that indeed $\gamma$ is unbounded. 

\begin{center}
\resizebox{!}{2.2in}{\includegraphics{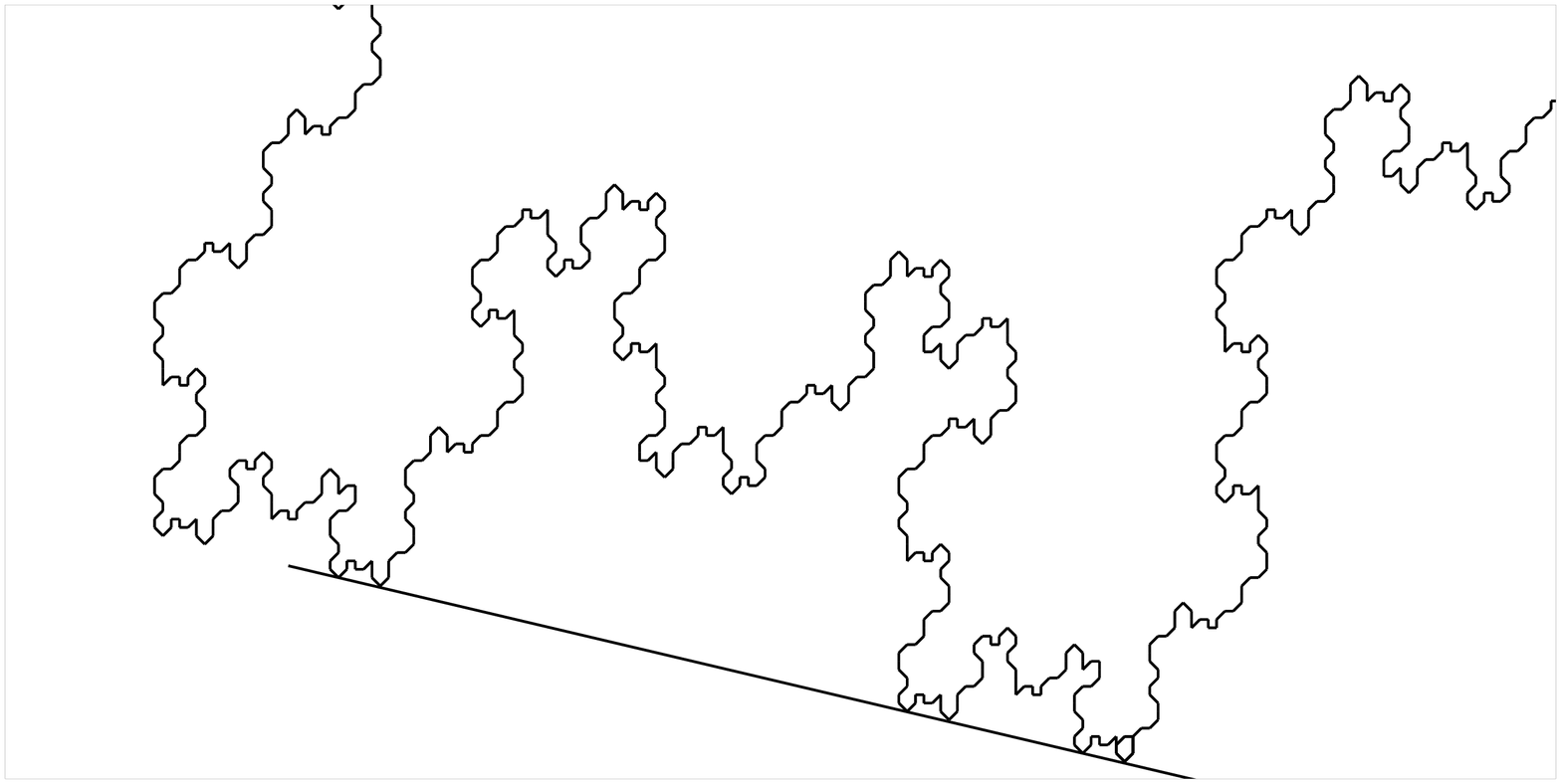}}
\newline
{\bf Figure 4.3:\/} The unbounded component of
$\Gamma(1,\phi^{-2})$.
\end{center}

\begin{center}
\resizebox{!}{2.2in}{\includegraphics{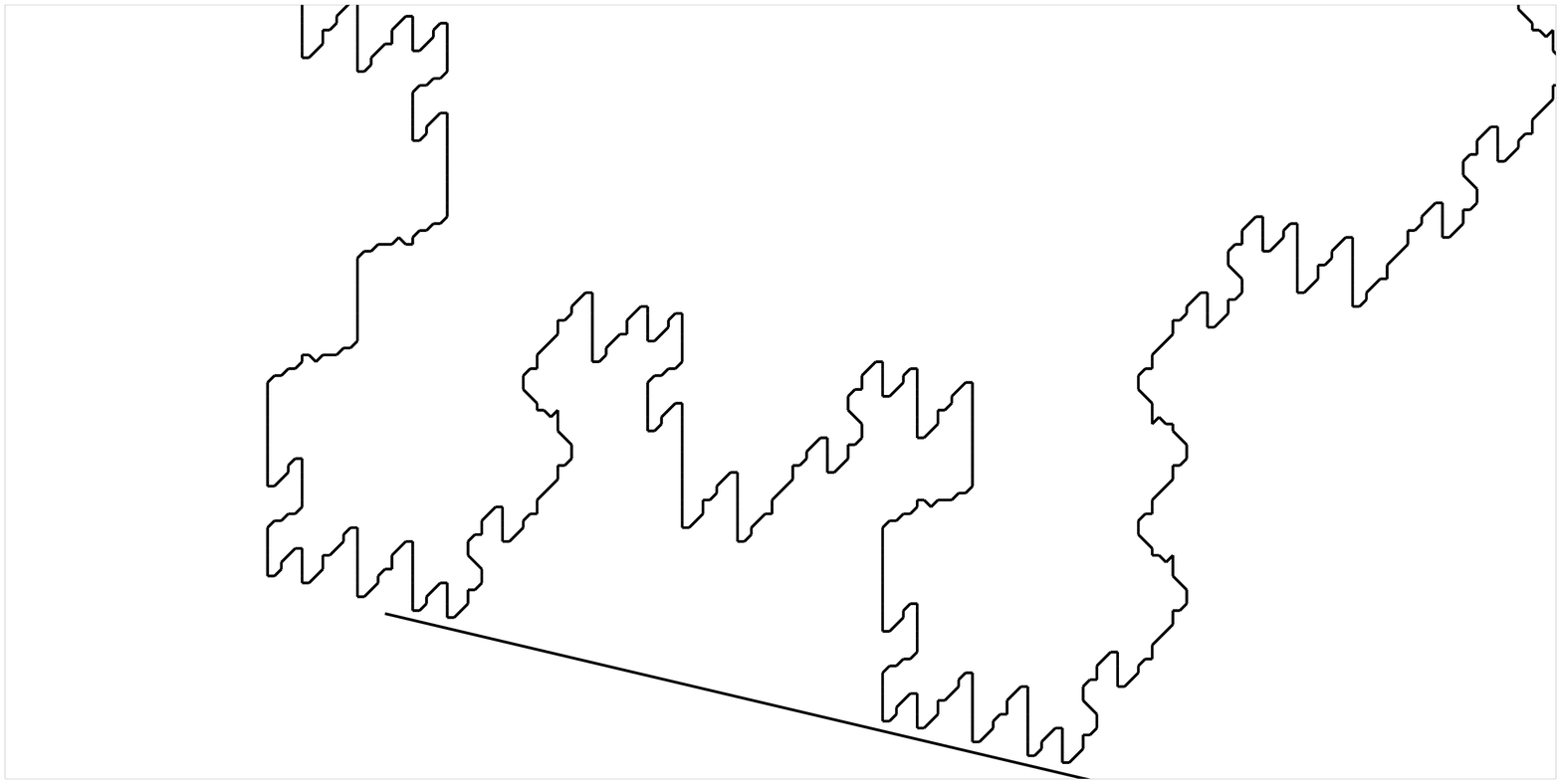}}
\newline
{\bf Figure 4.4:\/} The unbounded component of
$\Gamma(\phi^{-3},\phi^{-5})$.
\end{center}

Figure 4.4 shows the unbounded component of
$\Gamma(\phi^{-3},\phi^{-5})$.  This component
is quasi-isometric to the one in the previous picture but
looks different locally.  This is an illustration of
Theorem \ref{ULE} in action.

\subsection{The Freezing Phenomenon}
\label{freeze1}

The arithmetic graph illustrates an interesting
phenomemon.
Figure 4.5 shows the portion of a graph that
corresponds to the parameter $y=\phi^{-5}$.
Notice the long range linear order.  This
phenomenon becomes more and more extreme
as $y \to 0$:  Longer and longer portions
of the arithmetic graph contain these
nearly linear portions.  The components of
the arithmetic graph follow along the lines,
switching onto a new line at every intersection
so as to avoid collisions.  We call this the
{\it freezing phenomenon\/} because the
orbits seem to freeze into a characeristic shape.

\begin{center}
\resizebox{!}{4.2in}{\includegraphics{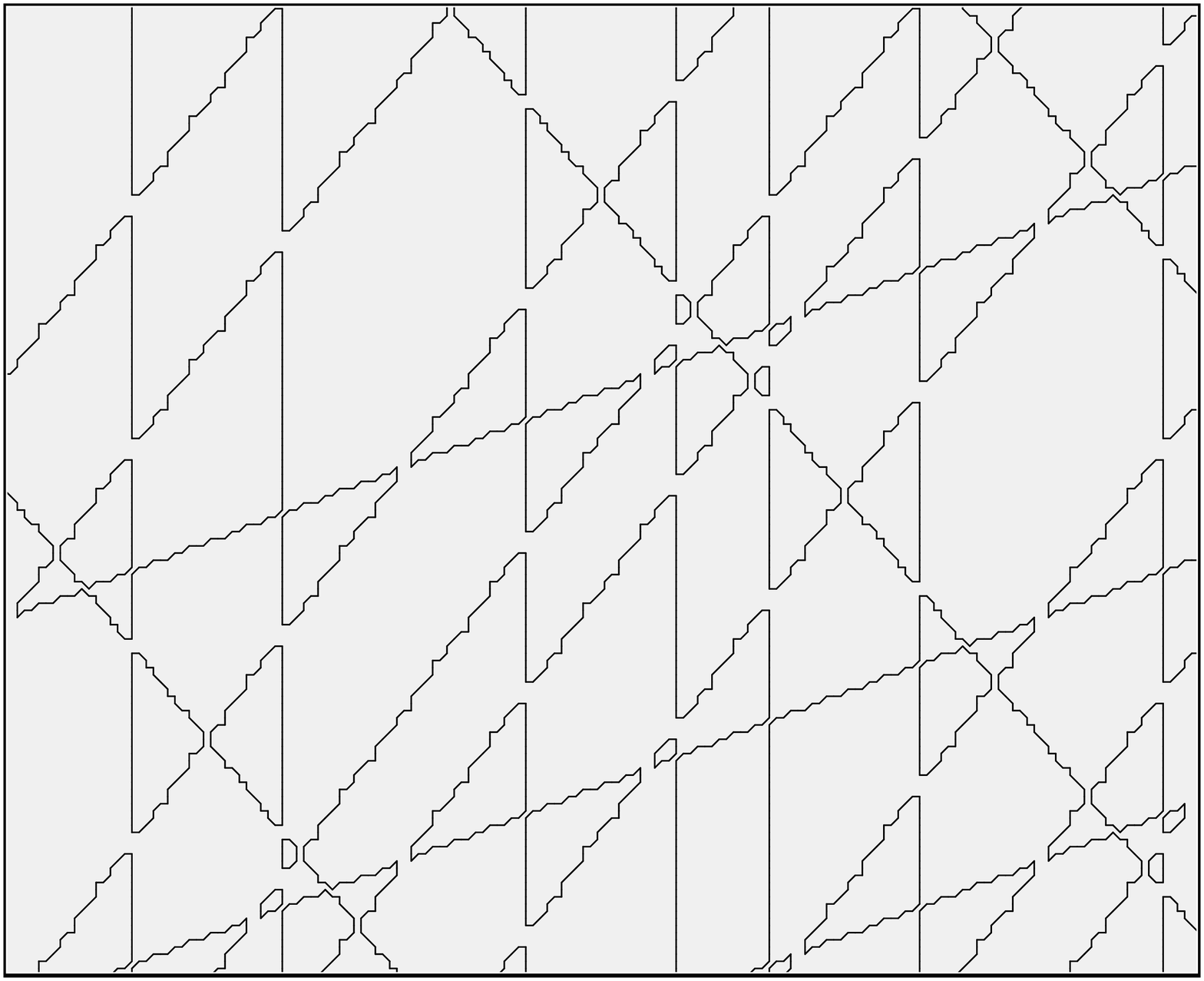}}
\newline
{\bf Figure 4.5:\/} Part of $\Gamma(\phi^{-5},-2\phi^{-4})$.
\end{center}

The freezing phenomenon seems almost to be in contradiction
with Theorem \ref{ULE}, as we now explain.
Let $R$ be the renormalization map.  The map $R^2$ is
expanding, so the full preimage of any point $y \in \R/2\Z$
is dense.  Hence, Theorem \ref{ULE} has the following
corollary.
\begin{corollary}
\label{quasi}
Let $O$ be any generic unbounded orbit.  Then there is a dense set
$y \in (0,2)$ such that $U_y$ contains an orbit that
is both coarsely and locally similar to $O$.  In
particular, this is true for a sequence $\{y_n\}$
converging to $0$.
\end{corollary}

What makes Corollary \ref{quasi} seem at odds
with the freezing phenomenon is that the 
arithmetic graph is definitely changing shape
as $y \to 0$, but somehow the global shape
is always coarsely equivalent to some fixed
arithmetic graph.  The escape from the contradiction
is that the bi-lipschitz constant implicit in the
definition of coarse equivalence tends to $\infty$
as $y$ tends to $0$.

Given that both the freezing phenomenon and
Corollary \ref{quasi} are true, the
renormalization discussed in our theorems
must be expressible directly in terms of the
multigrid system of lines.  This is indeed
the case.  On the level of the multigrid, the
renormalization is reminiscent both of
the renormalization one sees for Sturmian
sequences and the renormalization one sees
for the Penrose tilings, especially when it
is expressed in terms of De Bruijn's pentagrids.
See [{\bf DeB\/}].
So, there really is an underlying connection
between outer billiards on the Penrose kite
and the Penrose tiling. 

The freezing phenomenon works for all kites,
and it is part of a larger phenomenon,
though I don't see the renormalization
scheme for a general kite.  
In [{\bf S3\/}] I wrote some
informal notes describing the connection
betIen the arithmetic graph and the multigrids
of lines, but I did not include any proofs.
I don't currently know any.

\newpage

\section{Structural Results}
\label{structural}

\subsection{Compactification}

The map $\Psi: \Sigma \to \Sigma$ turns out to be an infinite golden polygon
exchange map.   This quasi-periodicity is the driving idea behind our
next result.  Define the flat $3$-torus
\begin{equation}
\widehat \Sigma=\T^3; \hskip 30 pt \T=\R/2\Z.
\end{equation}

\begin{theorem}[Compactification]
\label{master1}
There is a fibered golden polyhedron exchange map
$\widehat \Psi: \widehat \Sigma \to \widehat \Sigma$, and
an injective embedding
$\Theta: \Sigma \to \widehat \Sigma$, given by
the equation
\begin{equation}
\label{thetadef}
\Theta(x,y)=\bigg(1,\frac{1}{2},0\bigg)+\bigg(\frac{x}{\phi},\frac{x-y}{2},y\bigg)
\end{equation}
 which is a semi-conjugacy
between $\Psi$ and $\widehat \Psi$. 
\end{theorem}

Some terms require explanation. 
First, the coordinates of $\Theta$ are interpreted as
living in $\T$.  The translational part of $\Theta$ is
somewhat arbitrary.  We found this choice convenient
for the purpose of drawing pictures.
  To say that $\Theta$ is a
{\it semi-conjugacy\/} is to say that
$\widehat \Psi \circ \Theta = \Theta \circ \Psi$ wherever all
maps are defined.   Since $\Theta$ is injective, we will
often identify $\Sigma$ with the subset
$\Theta(\Sigma) \subset \widehat \Sigma$.  In this
way, we think of $\widehat \Sigma$ as a compactification
of $\Sigma$.   With this interpretation, the semi-conjugacy
condition just says that $\widehat \Psi$ extends $\Psi$.

The set $\Theta(\Sigma)$ is contained in a countable dense union of
parallel planes that are transverse to the horizontal planes.
We define
\begin{equation}
\Sigma_y=\Big(\R \times \{y\}\Big) \cup \Big(\R \times \{y-2\}\Big).
\end{equation}
for $y \in (0,2)$.  Given that $\psi$ preserves
the set $\R \times (y+2\Z)$, the map
$\Psi$ preserves $\Sigma_y$ for each $y$.
The image $\Theta(\Sigma_y)$ is densely contained
in a single horizontal plane, and $\Theta$ is a semi-conjugacy
between the map $\Theta:\Sigma_y \to \Sigma_y$ and a
exchange map on the corresponding horizontal plane.

The semi-conjugacy between 
$\Psi: \Sigma_1 \to \Sigma_1$ and
$\widehat \Psi: \widehat \Sigma_1 \to \widehat \Sigma_1$
is equivalent to the Arithmetic Graph Lemma of 
[{\bf S1\/}].  Here $\widehat \Sigma_1$ is the
horizontal plane of height $1$. Similarly,
Theorem \ref{master1} is closely related to 
the Master Picture Theorem in [{\bf S2\/}].
See \S 7 for a discussion.

\subsection{Structure of the Compactification}
\label{partitionpix}

In terms of raw data, we describe the compactification
precisely in \S \ref{appendix}.  Here we
highlight the important features.
We begin by showing pictures.
Each picture is a different
horizontal slice of $\widehat \Sigma$.  Let
$R$ be the parallelogram with vertices
\begin{equation}
\label{heights}
(-A,0); \hskip 15 pt
(A,2); \hskip 15 pt
(-A+2,0); \hskip 15 pt
(A+2,0); \hskip 30 pt A=\phi^{-3}=\sqrt 5-2.
\end{equation}
The solid $X=R \times [0,2]$ is a fundamental domain
for the action of $2\Z^3$ on $\R^3$ and conveniently
we can choose lifts of our polyhedra so that they
give a partition of $X$.  Figure 5.1 shows the
intersection of $X$ with the plane $z=0$.  The
grey polygons represent the intersections of
the polyhedra in our partition with this plane.

\begin{center}
\resizebox{!}{3in}{\includegraphics{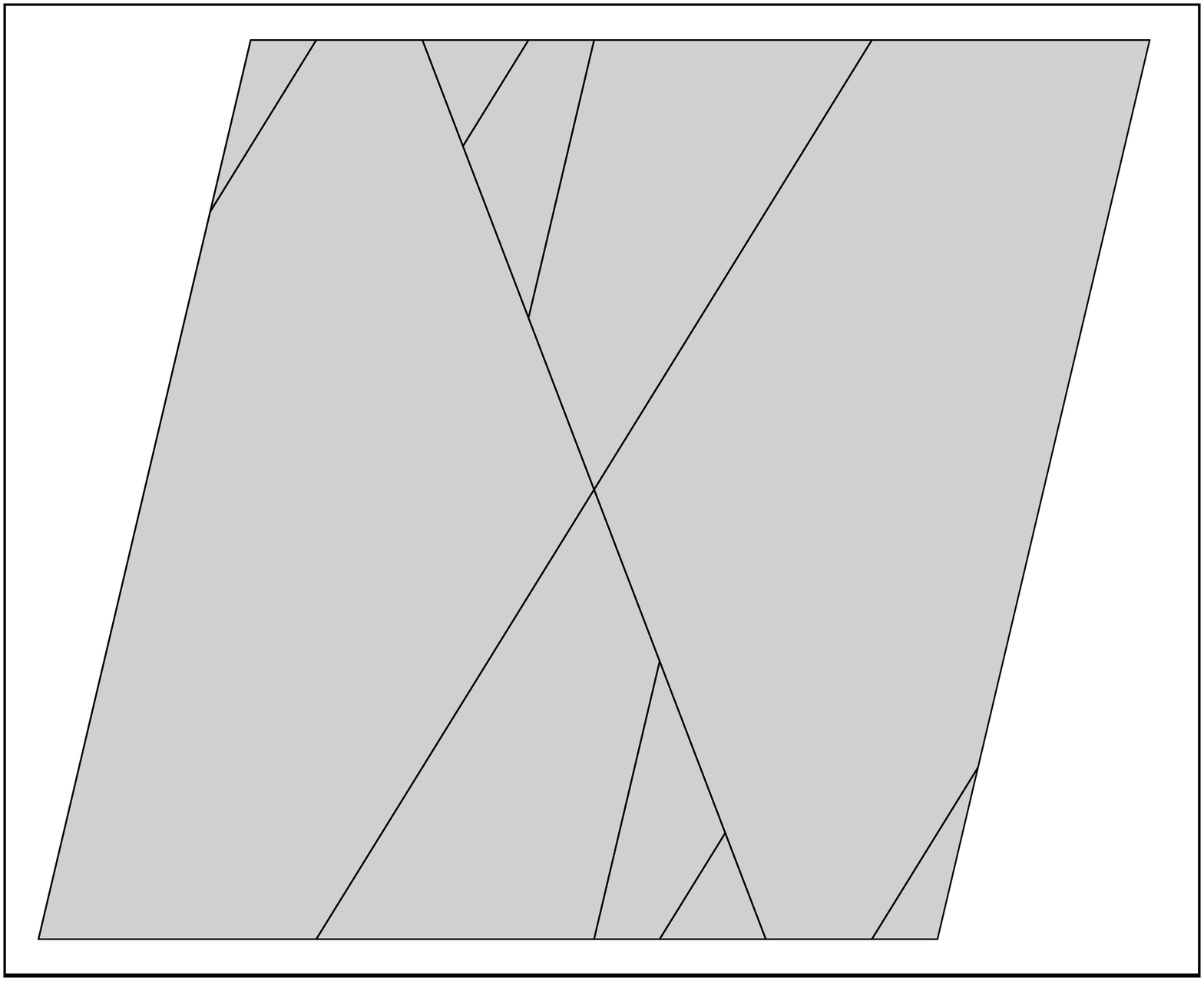}}
\newline
{\bf Figure 5.1\/}: The slice at $0$.
\end{center}

Below we show the slices at the heights
\begin{equation}
0;\hskip 10 pt \phi^{-3}; \hskip 10 pt \phi^{-2} \hskip 10 pt \phi^{-1}
\hskip 10 pt 2\phi^{-2}; \hskip 10 pt 1.
\end{equation}
The slices at heights $t$ and $2-t$ are isometric to each other,
via the isometry which rotates $180$ degrees about the midpoint
of the parallelogram. Indeed, the map
$(x,y,z) \to (2-x,2-y,2-z)$ is an automorphism of the partition.
 The slices values we choose are precisely
the values in $[0,1]$ which contain some polyhedron vertices.
Here are the remaining slices

\begin{center}
\resizebox{!}{3.5in}{\includegraphics{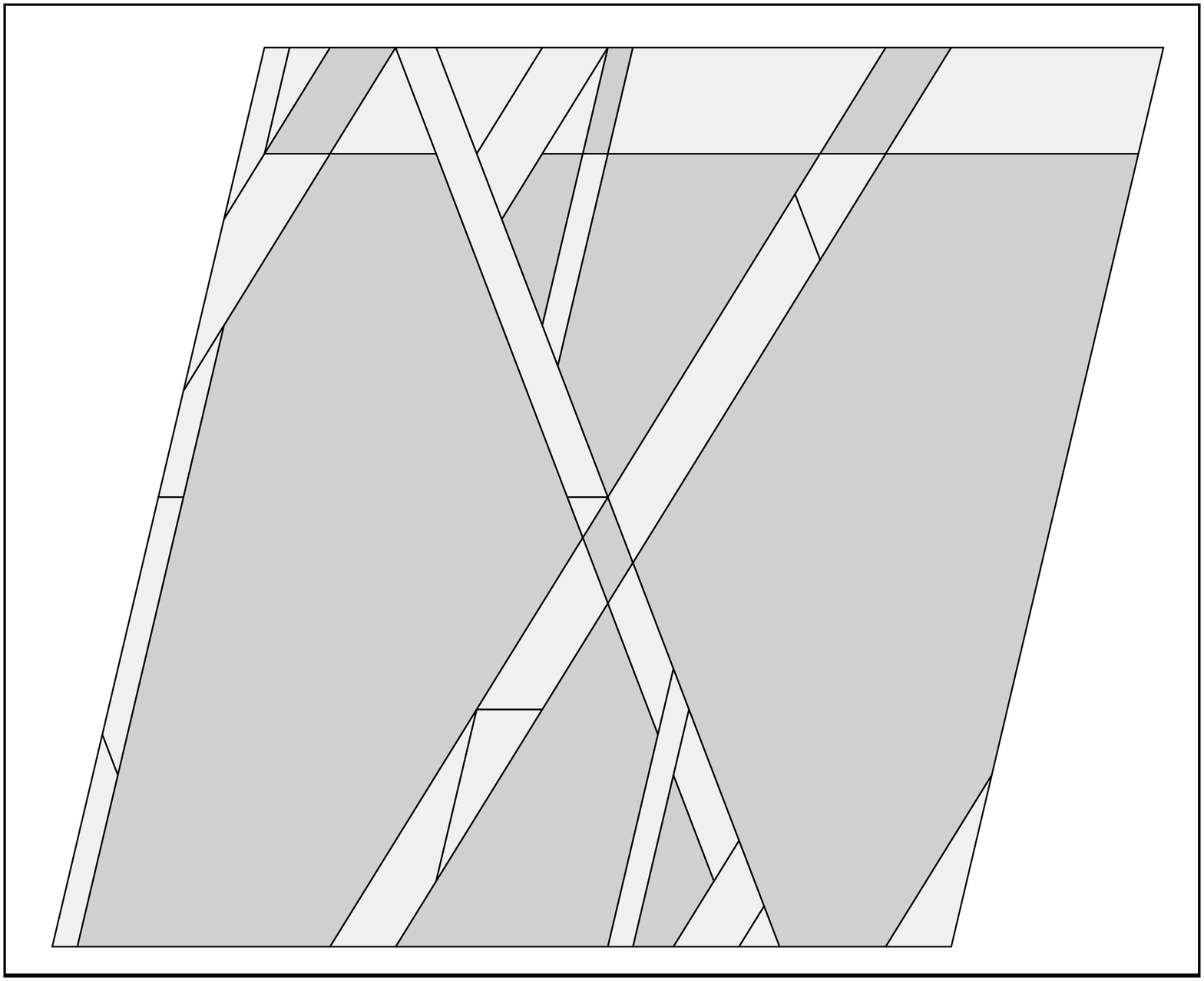}}
\newline
{\bf Figure 5.2\/}: The slice at $\phi^{-3}$.
\end{center}

\begin{center}
\resizebox{!}{3.5in}{\includegraphics{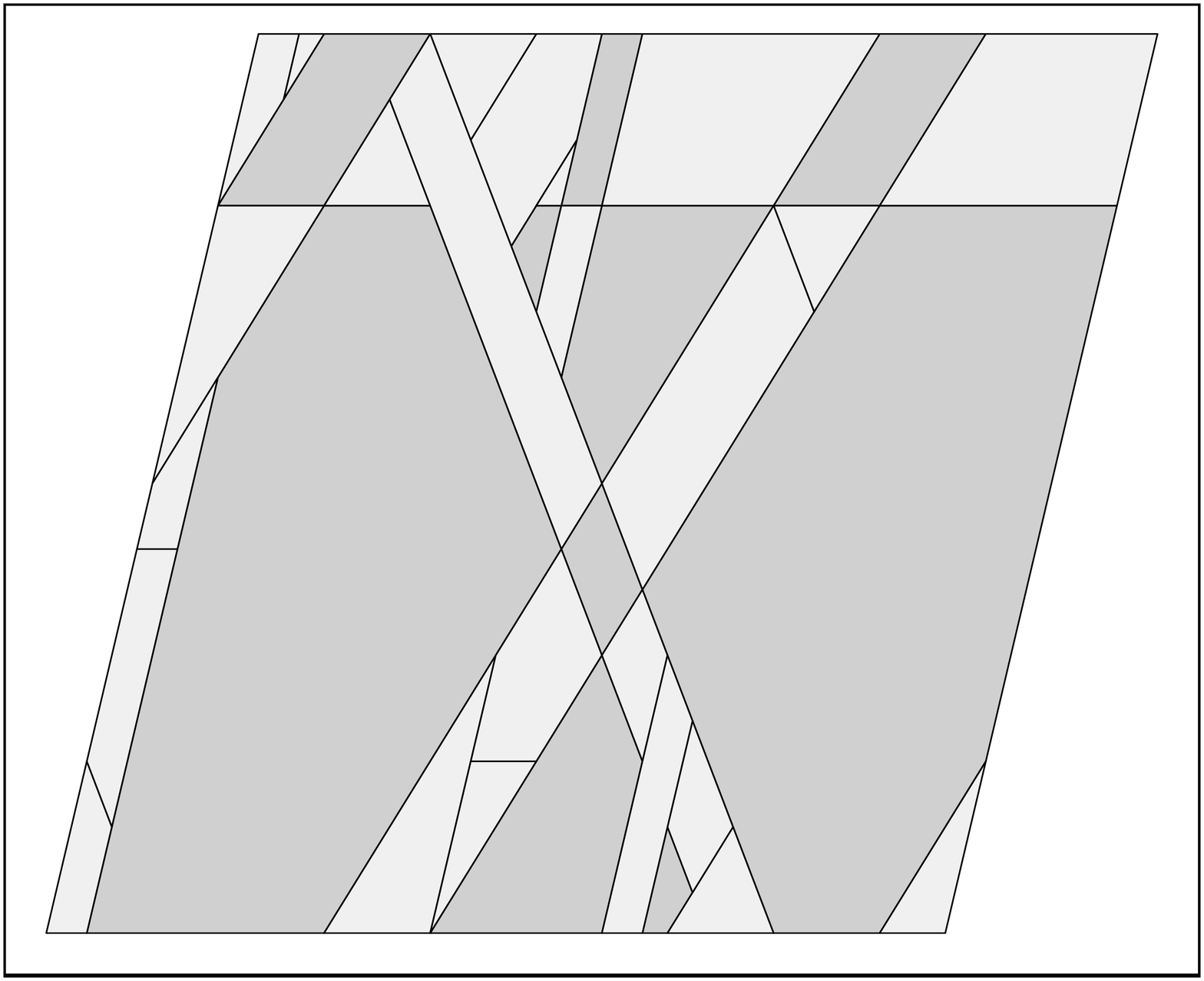}}
\newline
{\bf Figure 5.3\/}: The slice at $\phi^{-2}$.
\end{center}

\begin{center}
\resizebox{!}{3.5in}{\includegraphics{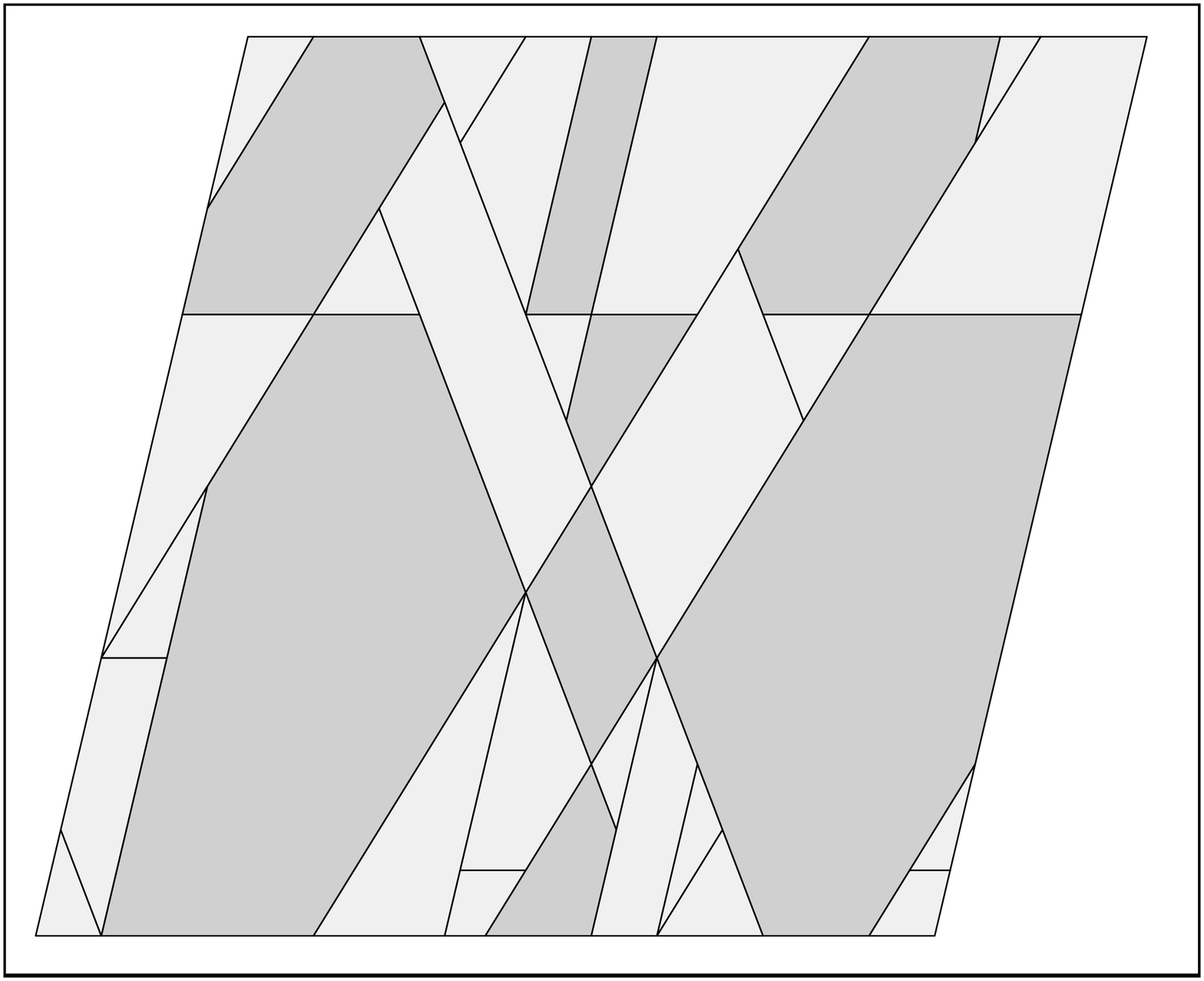}}
\newline
{\bf Figure 5.4\/}: The slice at $\phi^{-1}$.
\end{center}

\begin{center}
\resizebox{!}{3.5in}{\includegraphics{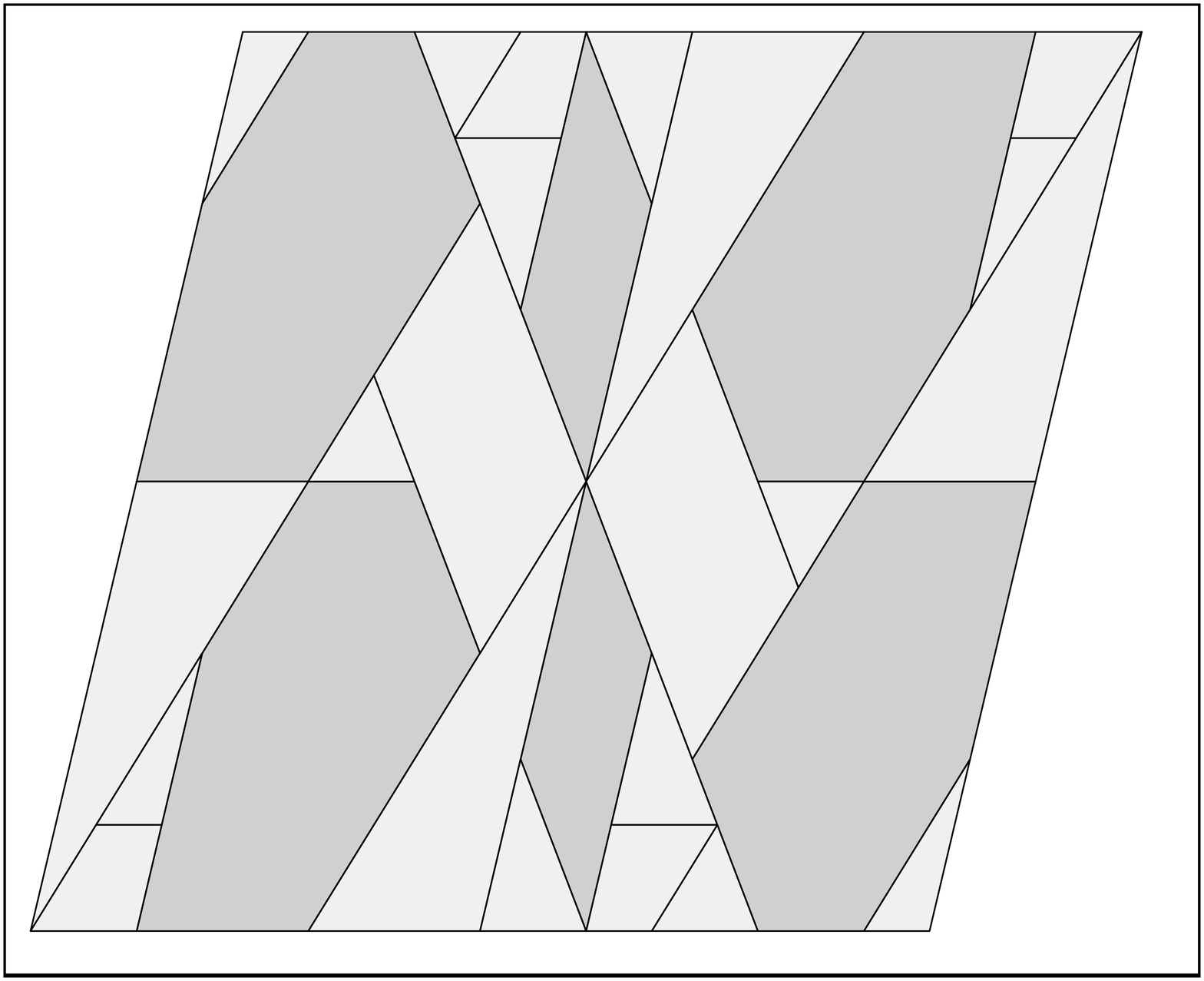}}
\newline
{\bf Figure 5.5\/}: The slice at $2\phi^{-2}$.
\end{center}

\begin{center}
\resizebox{!}{3.3in}{\includegraphics{Pix/slice5.ps}}
\newline
{\bf Figure 5.6\/}: The slice at $1$.
\end{center}

The shading in our pictures has the following explanation.
The map $\widehat \Psi$ is the identity
on a polygon iff the polygon is colored grey.
In particular, $\widehat \Psi$ is the identity
on the slice $z=[0]$.

Now we discuss the action of $\widehat \Psi$ on the
partition. Say that a 
{\it special affine involution\/} of $\widetilde \Sigma$
is an order $2$ affine map that preserves each horizontal slice
and is an isometric rotation in each horizontal slice. 
The map
\begin{equation}
\tau(x,y,z) \to (-x,-3z-y,z)
\end{equation}
is one such map. For each polyhedron $P$ in the partition, there is a
special affine involution $I_P$ such that $I_P(P)$ is also in the partition.
We have
\begin{equation}
\label{PE}
\widehat \Psi|_P=\tau \circ \iota_P.
\end{equation}

\begin{lemma}
$\widehat \Psi$ is a fibered golden
polyhedron exchange map, as in \S \ref{pem}.
\end{lemma}

\startproof
A direct calculation shows that every edge of
every polyhedron in the partition satisfies
Equation \ref{fibered}.  The equations
for the vertices of these polyhedra are
listed in the appendix.
\endproof

There is one other feature of the pictures we would like to
mention.  The lines in Figure 5.1 are present in all the pictures.
These ``persistent lines''
divide $\widehat \Sigma$ into $10$ convex prisms -- i.e.
polyhedra of the form $P \times [0,2]$, where $P$ is a convex polygon
in the plane.  Hence, the $64$ polyhedra in the partition can
divided into $10$ groups, each of which partitions one of the prisms.
Figure 5.7 illustrates this for the slice at height $1$.
Our java program allows the user to highlight the
clusters one at a time.

In particular, the $64$ polyhedra together
partition the fundamental domain $\widehat F$, and not just
the torus.  In other words, the outer edges in all the figures
are actually part of the polyhedra, and not just artifacts of
the way we have drawn the picture.  A moment's reflection reveals
that this issue arises whenever one wants to depict a partition
of a torus by polygons.

\begin{center}
\resizebox{!}{2.7in}{\includegraphics{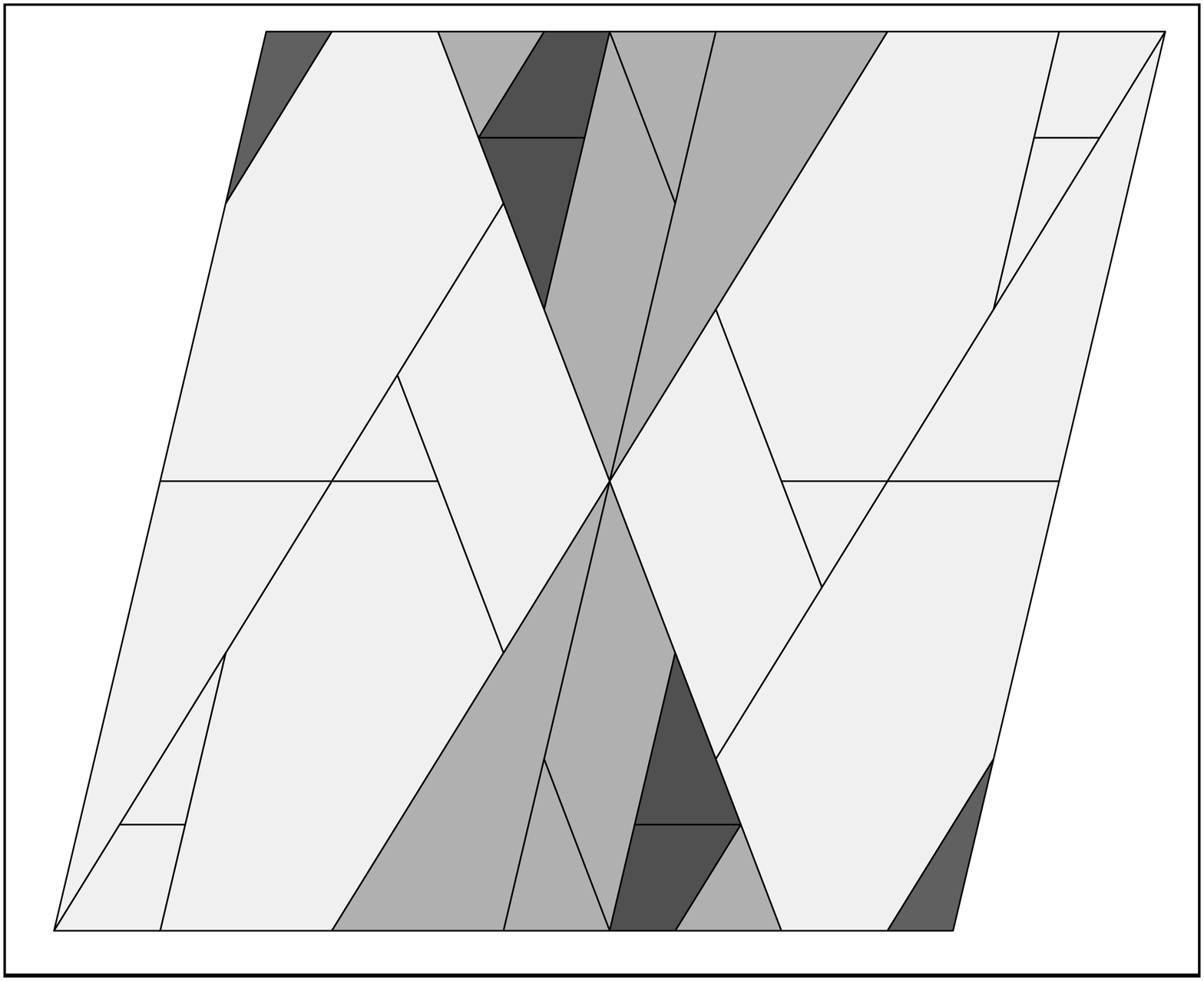}}
\newline
{\bf Figure 5.7\/}: The slice at $1$, divided into prisms.
\end{center}

There is a second way to group the $64$ polyhedra into
$10$ (non-vertical) prisms.  If $P$ is one of the original
$10$ prisms, then $\tau(P)$ is another prism which is
a finite union of some of the pieces.  Here
$\tau$ is the affine involution defined above.
Thus, the partition of $\widehat \Sigma$ into
$64$ convex polyhedra is compatible with two interlocking
partitions of $\widehat \Sigma$ into $10$ prisms each.

The existence of the two families of prisms suggests
that our polyhedron exchange map is actually the
square of a piecewise affine map that is defined
in terms of the two prism partitions.  This is
indeed the case, but we did not find this characterization
useful.

\subsection{The Renormalization Theorem}

We call an open set
$\widehat A \subset \widehat \Sigma$ {\it atomic\/} if
$\widehat A$ has a finite partition into 
golden polyhedra $\widehat P \cup ... \cup \widehat P_n$
such that the first return map $\Psi: \widehat A \to \widehat A$
is entirely defined, and a translation, when restricted
to the interior of each $\widehat P_k$.   
We call the polyhedra $\widehat P_k$ the {\it atoms\/}
of $\widehat A$.  Our definition does not uniquely
define {\it the\/} atoms, but in practice our atoms
will be the maximal ones.

Suppose that $\widehat A$ and $\widehat B$ are two
atomic sets.  We say that a
map $f: \widehat A \to \widehat B$ is an
{\it atomic bijection\/} if $f$
bijectively maps the atoms of $\widehat A$ to the
atoms of $\widehat B$, and the restriction
of $f$ to each atom is a homothety.

We say that a {\it layer\/} of $\widehat A$ is
the set of all points in $\widehat A$ whose
third coordinates belong to some interval.
We call $f: \widehat A \to \widehat B$ an
{\it atomic cover\/} if
$\widehat A$ has a partition into
layers $\widehat A_1,...,A_{dk}$ and
$\widehat B$ has a paritition into
slabs $\widehat B_1,...,\widehat B_k$
such that $f$ is surjective and $d$-to-$1$ and
$f: \widehat A_i \to \widehat B_j$ is
an atomic bijection for each $i$. Here
the index $j$ depends on $i$.  Each index $j$
corresponds to $d$ indices $i$.  Assuming
that the map $f$ is given, we call the
abovementioned layers {\it the\/} layers of 
$\widehat A$ and $\widehat B$.

An orbit (or orbit portion) is {\it generic\/} if it
does not intersect any non-horizontal plane that is
defined over $\Z[\phi]$.

\begin{theorem}[Renormalization]
\label{master2}
There is a pair of atomic sets
$\widehat A, \widehat B \subset \Sigma$
and a $3$-to-$1$ atomic covering map
$\widehat R: \widehat A \to \widehat B$ with
the following properties.
\begin{enumerate}
\item $\widehat R$ conjugates
$\widehat \Psi|\widehat A$ to $\widehat \Psi|\widehat B$
or to $\widehat \Psi^{-1}|\widehat B$, according as
$\widehat R$ acts as a
translation or a dilation. 
\item $\widehat R$ acts on the horizontal
planes exactly as the map $R$ acts on $\R/2\Z$.
\item $\widehat R$ maps
$\Theta(\Sigma) \cap \widehat A$ to
$\Theta(\Sigma) \cap \widehat B$.
\item Any generic orbit portion of length $812$ intersects
$\widehat A$ and any generic orbit portion of
length $109$ intersects $\widehat B$. 
\item For any generic $p \in \widehat A$, the 
orbit of $\widehat R(p) \in \widehat B$ returns to
$\widehat B$ in fewer steps than the orbit of $p$
returns to $\widehat A$.
\item Any generic orbit which intersects $\widehat A$ also intersects $\widehat B$.
\end{enumerate}
\end{theorem}

\noindent
{\bf Remark:\/} The constant $109$ is optimal.  The constant
$812=703+109$ is an artifact of our proof. The smaller constant $703$
would be optimal. See \S \ref{PROOFA} for a
discussion of these matters.

\subsection{Structure of the Renormalization}
\label{renormstruct}

In terms of raw data, we
describe $\widehat A$, $\widehat B$, and
$\widehat R$ precisely in \S \ref{B}
and \S \ref{AR}. In this section we give the
reader a feel for these objects.

\subsubsection{The Set $\widehat B$}

We describe $\widehat B$ precisely in \S \ref{B}.

The partition
\begin{equation}
\label{Bpart}
0< \phi^{-2}<2  \phi^{-2}<1<2-\phi^{-2}<2-2  \phi^{-2}<2.
\end{equation}
divides $\R/2\Z$ into $6$ intervals,
$J_1,J_2,J_{31},J_{32},J_4,J_5$.
This is the same partition that defines the intervals
$I_1,...,I_5$, except that the value $1$ has been inserted,
so as to split $I_3$ into two symmetric halves.
The layers of $\widehat B$ (relative to $\widehat R$)
are the ones corresponding to these intervals.

\begin{center}
\resizebox{!}{1.7in}{\includegraphics{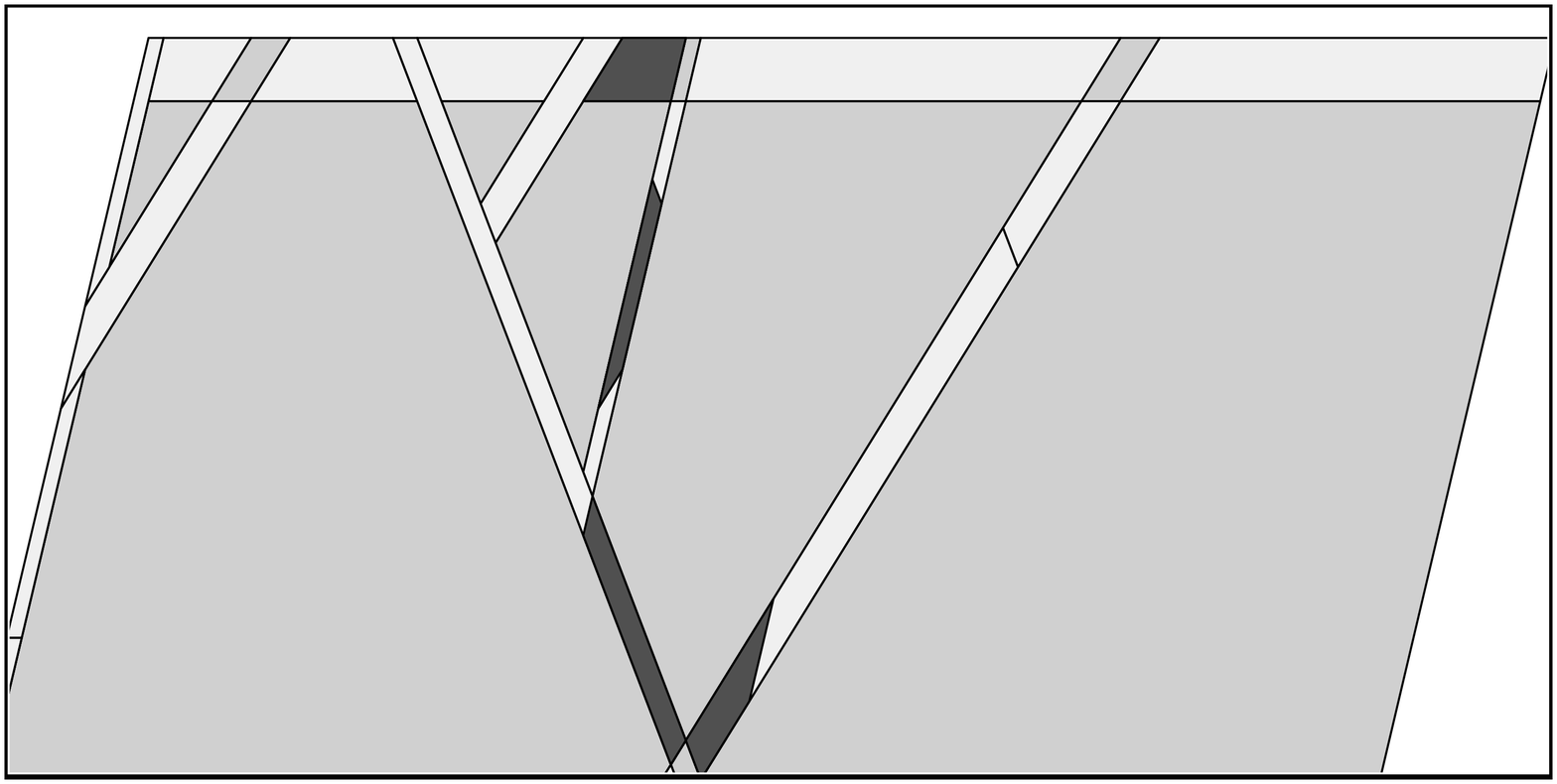}}
\newline
{\bf Figure 5.8\/}: The slice of $\widehat B$ at
$\phi^{-6}$.
\end{center}

Figure 5.8 shows a slice of $\widehat B$ at the
parameter $z=\phi^{-6}$.  The $4$ dark polygons
comprose the slice.   We are showing
the top half of $\widehat \Sigma_{z}$.
  This slice is part of
the first layer $\widehat B_0$. Each layer
$\widehat B_i$ decomposes into $4$ convex polyhedra,
which we call {\it branches\/}. Each branch
further decomposes into between
$9$ and $48$ atoms on which
the first return map is well defined.

The fact that each layer of
$\widehat B$ decomposes into $4$ branches is
practically forced by the structure $\widehat \Psi$,
as we now explain.  For $z$ near $0$, the restriction of
$\widehat \Psi$ to the slice $\widehat \Sigma_z$ is very nearly
the identity. It fails to be the identity only in a thin
neighborhood of $10$ line segments.   There is a set of $4$ slopes
such that each line segment has one of these $4$ slopes.
In this way, the $10$ ``active strips'' of
$\widehat \Sigma_z$ are nearly partitioned into
$4$ groups. (We say {\it nearly\/} because these
strips intersect.)

As $z \to 0$ there exist arbitrarily long orbits
which remain within a single group. If we
want to choose $\widehat B$ so that every sufficiently
long orbit intersects $\widehat B$, we must include
a polyhedron that nontrivially intersects each of the $4$ groups.
The choice of these polyhedra is not uniquely
determined by this requirement, but we have made
a reasonable and efficient choice.  Once we make
our choices for $z$ near zero, the picture
is determined for the entire interval
$[0,\phi^{-2}]$ by a kind of continuation principle:
The polyhedra naturally open up and
follow along as the active strips get wider.

\begin{center}
\resizebox{!}{3in}{\includegraphics{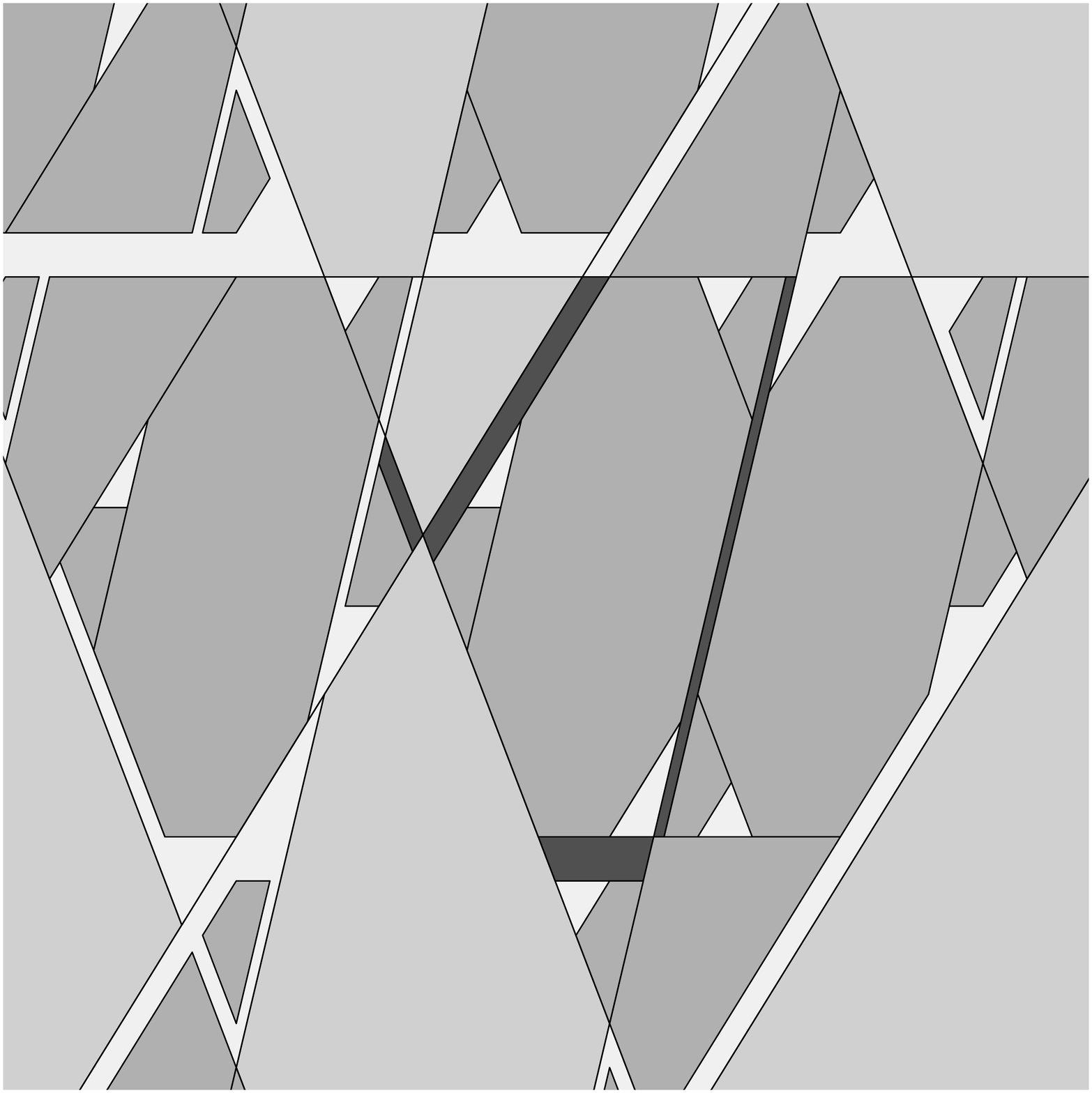}}
\newline
{\bf Figure 5.9\/}: The slice of $\widehat B$ (dark grey) at
$-17+11\phi$.
\end{center}

We say that the {\it freezing phenomenon\/}
is the tendency of the long orbits
to align along thin strips.  This phenomenon
is another manifestation of the one discussed
in \S \ref{freeze1}.
What happens in the first interval happens in each of
the $6$ intervals. The parameters
$2\phi^{-2}$ and $2-2\phi^{-2}$ also exhibit
the freezing phenomenon.
Figure 5.9 shows the slice 
of $\widehat B$ at the parameter
$-17+11\phi$.
This parameter is fairly near $2\phi^{-2}$.
The light grey polygons are some
periodic tiles.
These periodic tiles nearly fill up the
slice, but they leave some thin cracks. These
thin cracks line up along the $4$ basic directions,
and the long orbits accumulate in the
cracks.  The set $\widehat B$ fits inside the
cracks, with one polyhedron per direction.
We will list the vertices of the polyhedra comprising
$\widehat B$ in the appendix.  The reader can see
much better pictures of $\widehat B$ (and $\widehat A$)
using our applet.

\subsubsection{The set $\widehat A$}

The partition associated to $\widehat A$ consists
of the $18$ intervals we get by pulling
back the $J$ partition under the action of $R$.
Precisely:
\begin{itemize}
\item Within $I_1$, the partition is $R^{-1}(J_1),...,R^{-1}(J_4)$.
\item Within $I_2$, the partition is $R^{-1}(J_5)$.
\item Within $I_3$, the partition is $R^{-1}(J_1),...,R^{-1}(J_5)$.
\item Within $I_4$, the partition is $R^{-1}(J_1)$.
\item Within $I_5$, the partition is $R^{-1}(J_2),...,R^{-1}(J_5)$.
\end{itemize}
The $18$ layers of $\widehat A$ correspond to these $18$ intervals.
Again, each natural piece is decomposed into $4$
The map $\widehat R$ carries each branch
of each each layer of $\widehat A$ 
homothetically to the corresponding
branch of the corresponding layer of 
$\widehat B$.  As with the
$\widehat B$-branches, the
$\widehat A$-branches are further decomposed
into atoms.

\subsubsection{The Map $\widehat R$}

Here we highlight the general features.
Let $\widehat R_{ijk}$ denote the homothety that
expands distances by $\phi^3$ and fixes the point $(i,j,k)$.
\begin{enumerate}
\item When resricted to $\widehat A(I_1)$, the map
is $\widehat R_{110}$ or
$\widehat R_{120}$.
\item When restricted to $\widehat A(I_2)$, the map is translation
by $(0,0,2\phi^{-1})$.
\item When resricted to $\widehat A(I_3)$, the map
is one of $R_{1k1}$ for $k \in \{0,1,2\}$.
\item  When restricted to $\widehat A(I_4)$, the map is translation
by $(0,0,-2\phi^{-1})$.
\item When restricted to $A(I_5)$, the map is
one of $\widehat R_{102}$ or $\widehat R_{112}$.
\end{enumerate}

Now we will show the Renormalization Theorem
in action. Figure 5.10 shows a closeup
of the set $\widehat A$ at the slice
$5-3\phi$ together with part of an orbit.

\begin{center}
\resizebox{!}{1.8in}{\includegraphics{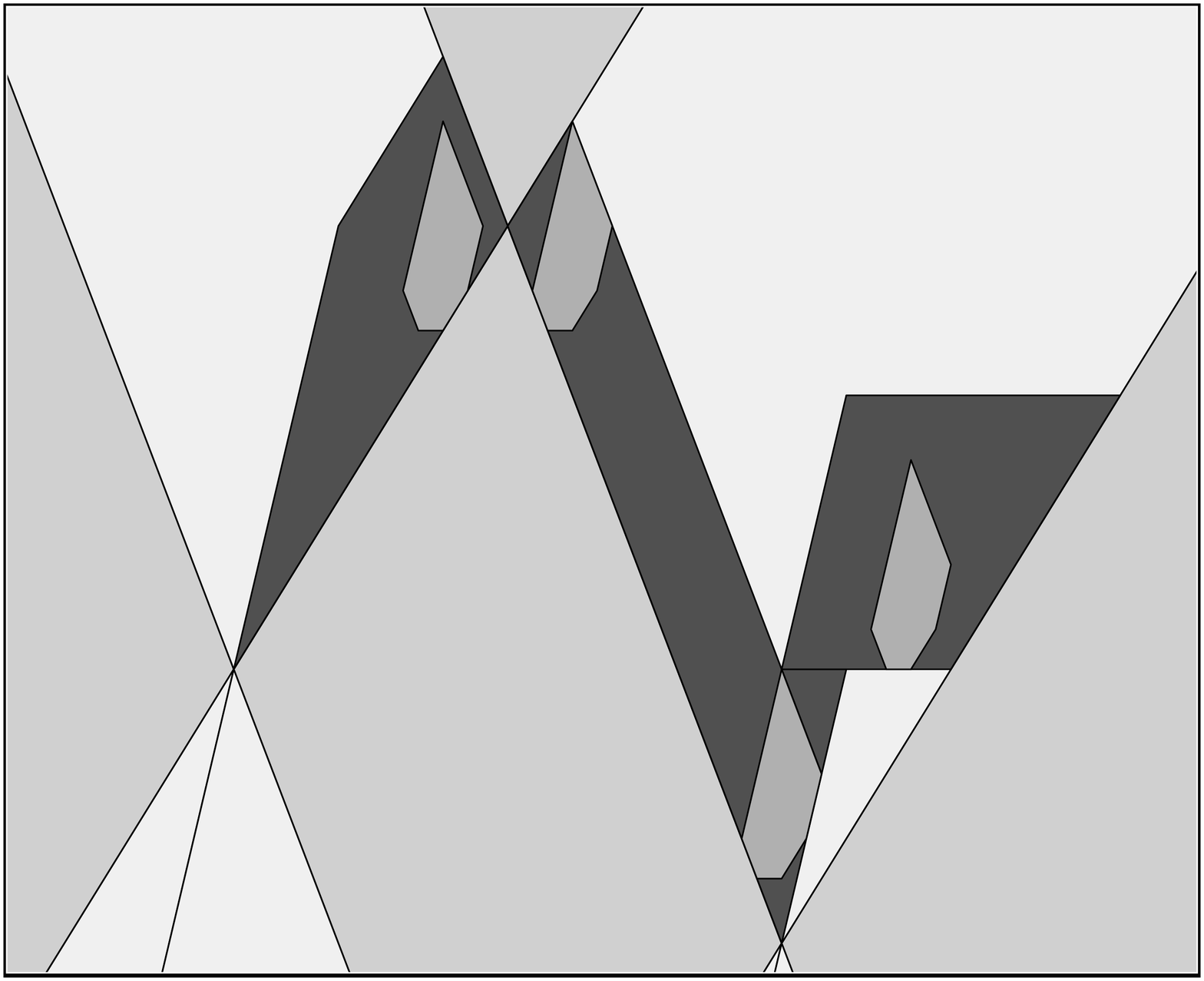}}
\newline
{\bf Figure 5.10\/}: Closeup of an orbit.
\end{center}

Figure 5.11 zooms out, to reveal much more
of the orbit in the same slice.  The set $\widehat A$
is at the very bottom.

\begin{center}
\resizebox{!}{4.5in}{\includegraphics{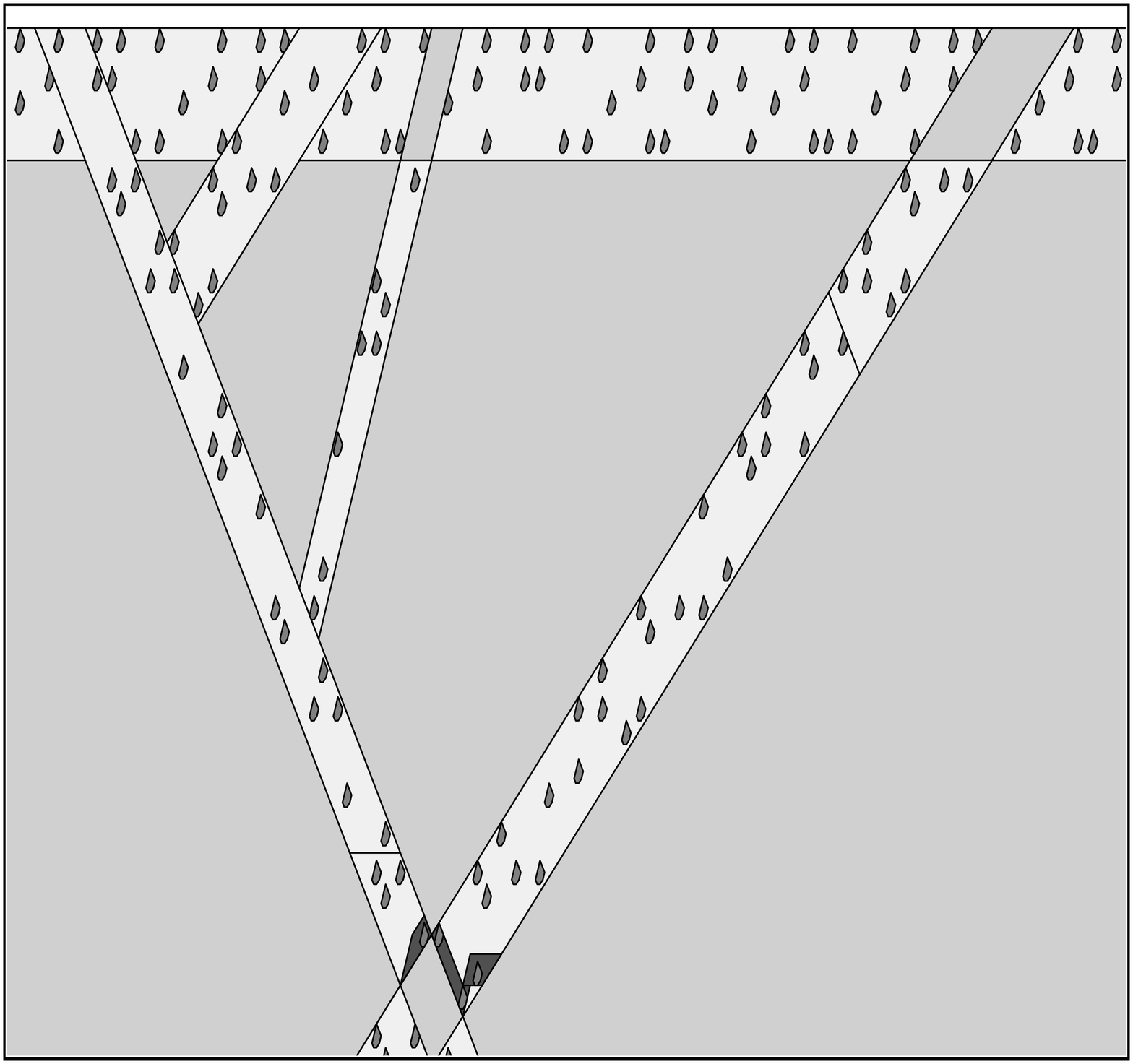}}
\newline
{\bf Figure 5.11\/}: The slice of $\widehat A$ 
at $5-3\phi$, together with an orbit.
\end{center}

Now for the magic trick.  Figure 5.12 shows the
slice of $\widehat B$ at the parameter
$-1+\phi=R(5-3\phi)$, together with an orbit.
The orbit in Figure 5.12 looks different overall
from the orbit in Figure 5.11, but the new
orbit intersects $\widehat B$ in the same way
that the old orbit intersects $\widehat A$.
In other words, were we to take a closeup of
just the renormalization set, we couldn't
tell whether it was the $\widehat A$-set
inside the first slice or the $\widehat B$-set
in the second slice.
The reader can see many more pictures like this
using our applet. 

\begin{center}
\resizebox{!}{3.6in}{\includegraphics{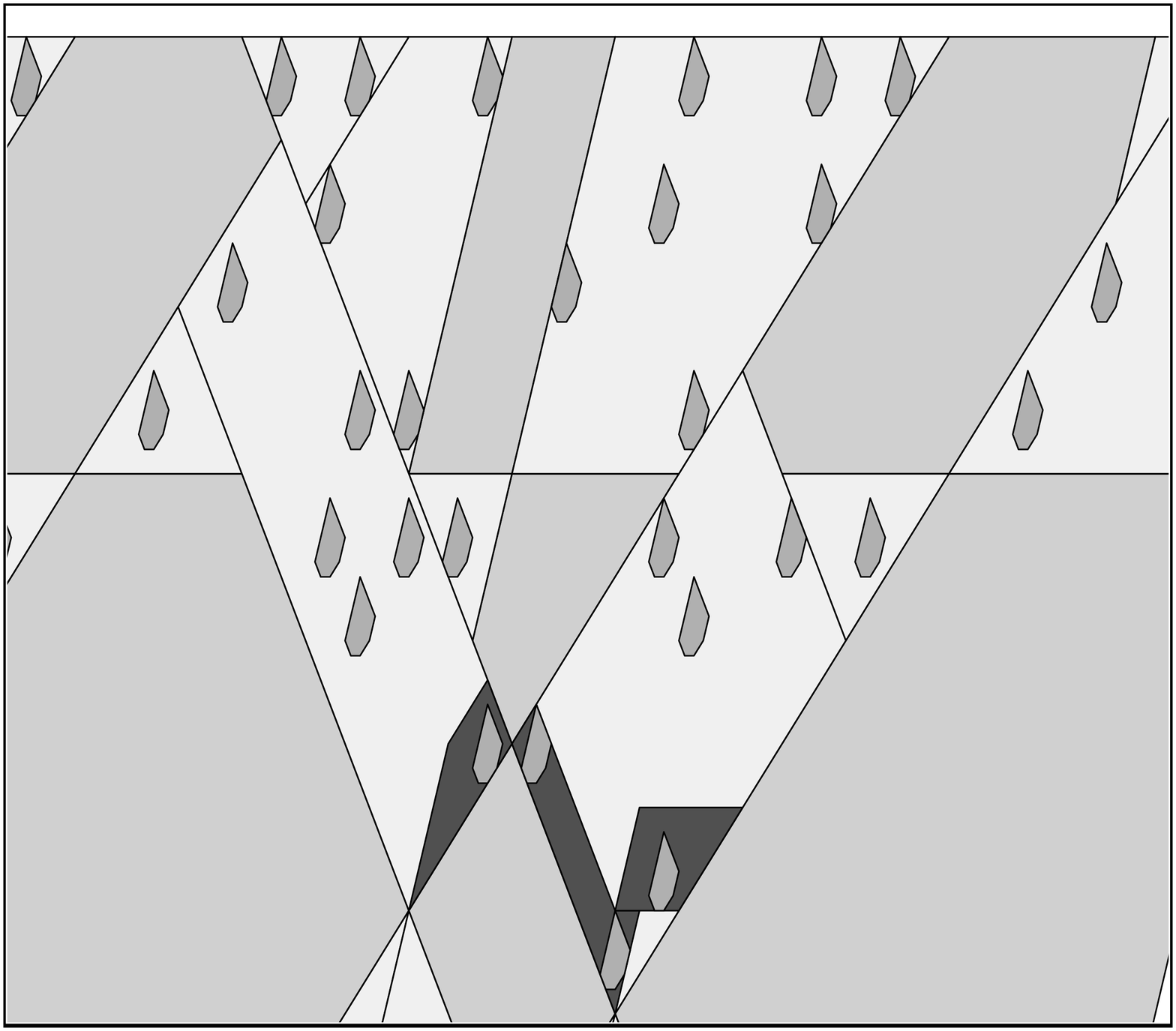}}
\newline
{\bf Figure 5.12\/}: The slice of $\widehat B$ 
at $-1+\phi$, together with an orbit.
\end{center}

The Renormalization Theorem and the freezing phenomenon
together combine to give a nice qualitative
description of the dynamical tiling associated
to the exchange map. Consider a parameter $z_0$
such that all orbits are periodic.  (Any
parameter $z_0$ such that $R^n(z_0)=0$ will
have this property.) At $z_0$, the slice
is tiled by periodic tiles.  As the parameter moves
away from $z_0$, some cracks open up, and
these cracks are filled with smaller
periodic tiles associated to longer periodic
orbits.  If the new parameter is chosen
just right, the slice is again tiled by
periodic tiles: the original large ones
and the much smaller ones that live in 
the cracks.  Moving the parameter a bit,
some new cracks open up and even smaller
tiles fill in these cracks. And so on.

\subsection{Renormalization of Orbits}

\noindent
{\bf Genericity:\/}
We call a point of $\widehat \Sigma$
{\it generic\/} if it does not lie in any non-horizontal
plane that is defined over $\Z[\phi]$.   The map
$\Theta$ carries the set of generic points in
$\Sigma$ into the set of generic points of
$\widehat \Sigma$.  Moreover, the translation
vectors defining $\widehat \Psi$ all lie
in $\Z[\phi]^3$.  Hence $\widehat \Psi$
preserves the set of generic orbits.
All the atoms of $\widehat A$ and $\widehat B$
have their faces in the planes we have
excluded. Hence, any generic point in
$\widehat A$ or $\widehat B$ lies in the
interior of an atom.  See \S \ref{tricky}
below for a further discussion.
\newline
\newline
{\bf Basic Definition:\/}
The Renormalization Theorem allows us to define a
{\it renormalization operation\/} on generic orbits.  Let
$O_1$ be a generic orbit that intersects $\widehat A$.
We know that $O_1$ does not intersect the boundaries of any
of the atoms.  
We start with $p_1 \in O_1 \cap \widehat A$.
We let $p_2=\widehat R(p_1)$ and we let $O_2$ be the orbit of
$p_2$.  This definition is independent of choice of $p_1$ thanks
to Item 1 of the Renormalization Theorem. For instance, if
we choose the first return point
$p_1'=\Psi^k(p_1)$ for a suitable power of $k$,
then $\widehat R(p_1')=p_2'=\widehat \Psi^{\pm k}(p_2)$.
We write 
\begin{equation}
O_1 \leadsto O_2
\end{equation}
when $O_1$ and $O_2$ are related as above.
We call $O_2$ the {\it renormalization\/} of $O_1$.
Given the nature of the map $\widehat R$, the
orbit $O_2$ is also generic.
\newline
\newline
{\bf Distinguished Orbits:\/}
There is one situation where our construction above
is ambiguous.
We call an orbit
$\widehat A$-distinguished
if it lies in a horizontal plane that contains the
top or bottom of one of the $\widehat A$ layers.
These orbits are somewhat of a nuisance.
When $O_1$ is a distinguished orbit, we
define the operation
$O_1 \leadsto O_2$ by including $O_1$ in
one layer of $\widehat A$ or the
other.  We will not take
the trouble to prove the $O_2$ does not
depend on the choice of layer
because in the cases of interest to us,
the uniqueness either doesn't matter or
comes as a byproduct of our proof.
\newline
\newline
{\bf Renormalization of Cores:\/}
We now mention a more precise kind of renormalization.
Say that a $\widehat A$-{\it core\/} is an orbit portion 
$p_1,...,p_n$ such that $p_1$ and $p_{n+1}$
lie in $\widehat A$ but $p_2,...,p_n$ do not.
We define $\widehat B$-cores in the same way.
If $\alpha$ is an $\widehat A$-core and $\beta$ is
a $B$-core, we write $\alpha \leadsto \beta$
if $\widehat R$ maps the endpoints
of $\alpha$ to the endpoints of $\beta$.
In this case, we have $O_1 \leadsto O_2$,
where $O_1$ is the orbit containing the
points of $\alpha$ and $O_2$ is the
orbit containing the points of $\beta$.
\newline
\newline
{\bf Restriction to the Strip:\/}
So far we have been talking about the picture in
the compactification $\widehat \Sigma$, but we
can transfer everything over to $\Sigma$.
By the Compactification Theorem, we can consider
$\Sigma$ as a subset of $\widehat \Sigma$.
With this interpretation, the action of
$\Psi$ is just the restriction of $\widehat \Psi$.
By Item 4 of the Renormalization Theorem, the
operation $O_1 \leadsto O_2$ preserves the set
of generic infinite $\Psi$-orbits.
\newline
\newline
{\bf Coarse Equivalence:\/}
The arithmetic graph illustrates the nature of
the renormalization map on cores.  We say that
a {\it strand\/} is the arithmetic graph
of a core, translated so that one endpoint
is the origin.  The renormalization operation
on cores gives a map from the set of all
$A$-strands to the set of all $B$-strands.
Figure 5.13 plots an $A$-strand in black
and the corresponding $B$-strand in grey.
The two strands start at the origin, which
is the endpoint at right, and they have
the other endpoint in common as well.
This example is what we will
call type-$1$.  It corresponds to 
a layer of $\widehat A$ on which $\widehat R$
is an isometry.

\begin{center}
\resizebox{!}{2.6in}{\includegraphics{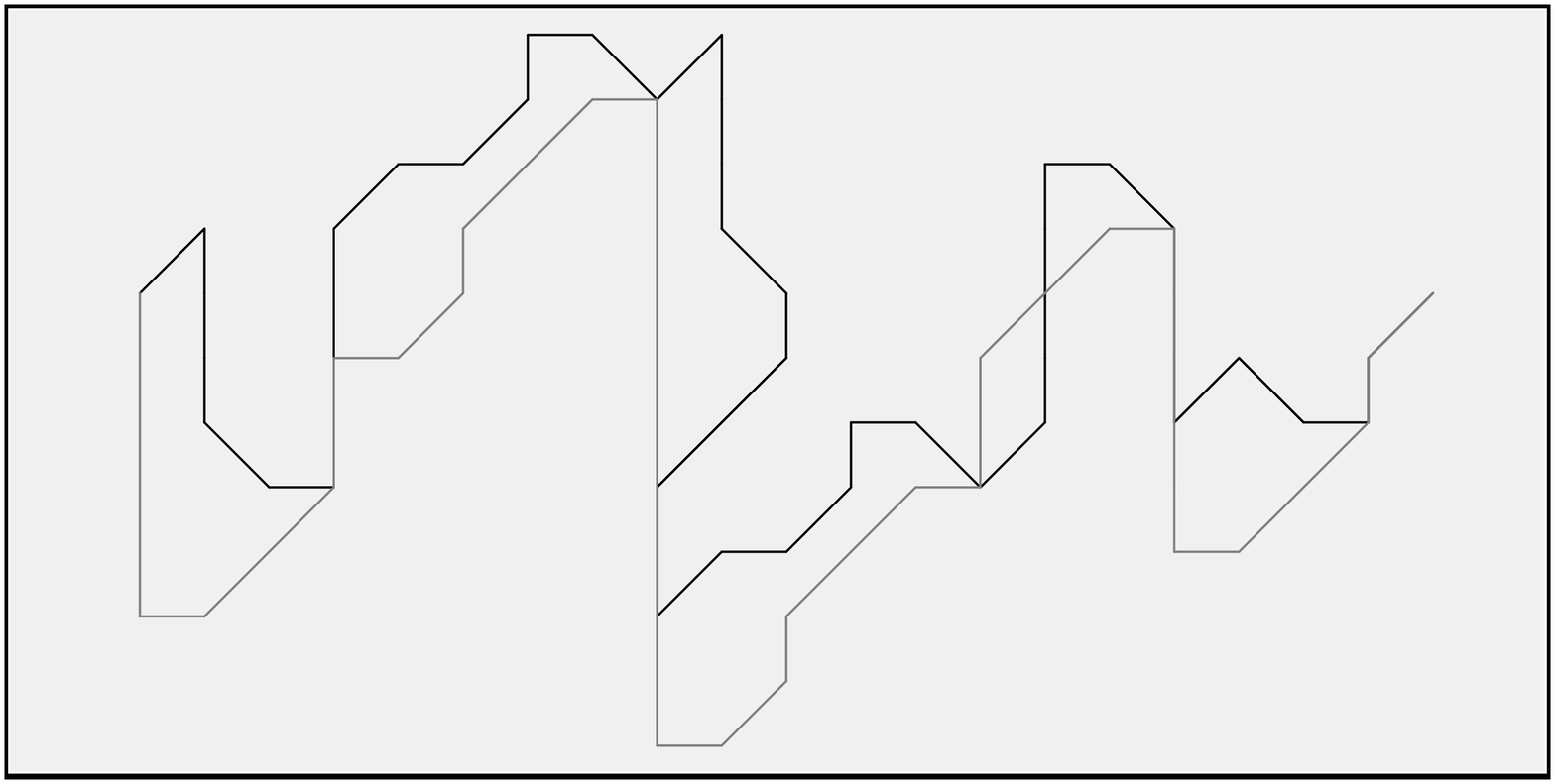}}
\newline
{\bf Figure 5.13:\/} Corresponding $A$ and $B$ strands.
\end{center}

Figure 5.14 shows another example.  This example
corresponds to a layer of $\widehat A$ where
$\widehat R$ is a dilation.  In this example 
the origin is at bottom right, and we have
shrunk the $A$-strand by a factor of
$\phi^{-3}$. Notice the remarkable agreement.
This remarkable agreement is the subject
of the Far Reduction Theorem below, and also
the basis for Theorem \ref{ULE}.

\begin{center}
\resizebox{!}{4.3in}{\includegraphics{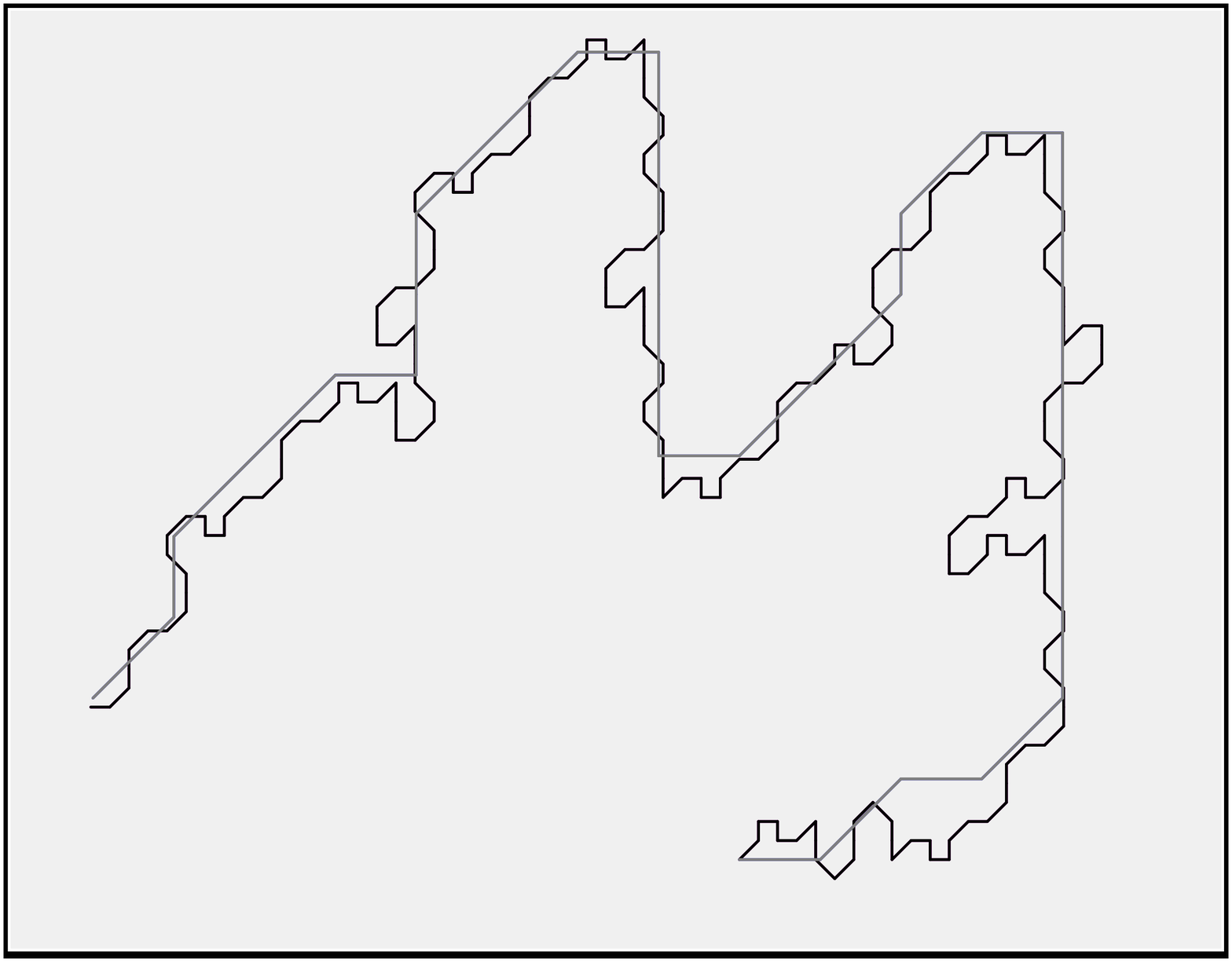}}
\newline
{\bf Figure 5.14:\/} Corresponding $A$ and $B$ strands.
\end{center}

There is a total of $2034=3 \times 678$ $A$-strands and a total of
$678$ $B$-strands.  Using our program, the tireless reader can
see pictures like the ones above for every pair of corresponding
strands.  

\subsection{The Fundamental Orbit Theorem}
\label{fund1}

To state the Fundamental Orbit Theorem,
we use the notation established in \S \ref{notation}.
Let $T$ be the fundamental triangle, the subject
of Theorem \ref{penrose9}.
Recall that $\cal T$ is the tiling of $T$ discussed in \S 3.
Recall that $J_0$ and $K_0$ are the largest octagon and kite in
$\cal T$.  Say that these tiles have {\it depth zero\/}.
Say that a tile of $\cal T$ has {\it depth $n$\/} if it
is similar to $J_0$ or $K_0$ by a factor of
$\phi^{-3n}$.    Let ${\cal T\/}(n)$ denote the finite
union of tiles having depth $n$ or less.

Recall that $\psi$ is the outer billiards map.
It is easy to check that $\psi^{-1}$ is entirely
defined on the interior of $T^+ \cup T^-$.

\begin{theorem}[Fundamental Orbit]
Let $p_1 \in T_{ij}^{\pm}-{\cal T\/}(2)$
be a generic point and let
$$p_2 = \psi^{-1}(R_{ij}^{\pm}(p_1))$$
Let $\alpha_1$ be an $A$-core that contains
$p_1$.  Then $\alpha_1 \leadsto \alpha_2$,
where $\alpha_2$ is the $B$-core containing
$p_2$.  In particular $O_1 \leadsto O_2$,
where $O_k$ is the orbit containing $p_k$.
\end{theorem}

The exponents $(\pm)$ are written loosely in the Fundamental Orbit
Theorem.   The choice of $(+)$ or $(-)$ in every cases is taken so
that the relevant maps have the proper domain and range.

The fundamental orbit says that, on the
level of orbits, the maps
$R_{ij}^{\pm}$ and the renormalization map
have the same action.  The Fundamental Orbit
Theorem is the result behind Theorem \ref{penrose9}.

\subsection{The Reduction Theorems}

Here we present the two results which
explain the qualitative action of
the action of $\widehat R$ on orbits.

We call two orbits $O_1$ and $O_2$ 
{\it associates\/} if one of $4$ things
holds.
\begin{itemize}
\item $O_1=O_2$.
\item $O_1=\overline{O}_2$, the complex conjugate orbit.
\item $O_1=\psi'(O_2)$. Here $\psi'$ is the outer billiards map.
\item $O_1=\psi'(\overline O_2)$.
\end{itemize}
Associate orbits are clearly locally and coarsely similar.
equivalent.

We write $O_1 \to O_2$ if $O_1' \leadsto O_2'$, where
$O_k'$ is an associate of $O_k$.  The two relations
$(\leadsto)$ and $(\to)$ are practically the same.
We use the latter because it allows us to prove our
results with significantly less computation.

Define
\begin{equation}
\Sigma_{24}=[-24,24] \times [-2,2]
\end{equation}

\begin{theorem}[Near Reduction]
Let $O_1$ be any generic 
infinite orbit that intersects
$\Sigma_{24}$.  Then there is some $m$ such
that $O_1 \to ... \to O_m$, and
$O_m$ intersects $T$, the fundamental triangle.
More precisely, let $y \in (0,2)$ be 
any point such that $\{R^n(y)\}$ does
not contain $0$.  Then there is some $N$ with the
following property.  If $O_1$ contains a point within
$1/N$ of the segment $[-24,24] \times \{y\}$, then
$O_1 \to ... \to O_m$ and $O_m$ intersects $T^+$ and
$m<N$.
\end{theorem}

Given an $A$-core $\alpha$, we define
\begin{equation}
\label{minx}
|\alpha|_x=\min_{(x,y) \in \alpha} |x|.
\end{equation}
That is, $|\alpha_x|$ measures how close the $x$-coordinates
of $\alpha$ come to $0$.
We make the same definition for $\beta$.

We say that $\alpha$ has type-$1$ if the renormalization
map $\widehat R$ is a piecewise translation on
the layer containing $\alpha$.  We say that
$\alpha$ has type-$2$ if $\widehat R$ is a piecewise
dilation by $\phi^3$ on the layer containing $\alpha$.

\begin{theorem}[Far Reduction]
\label{renorm2}
Let $\alpha$ and $\beta$
be $\widehat A$ and $\widehat B$ strands
respectively.  Suppose that
$\alpha \leadsto \beta$.  Then
the following is true.
\begin{enumerate}
\item If $\alpha$ has type $1$ then
$|\alpha_x|-C<|\beta|_x<|\alpha|_x+12$.
\item If $\alpha$ has type $2$ then
$\phi^{-3}|\alpha_x|-C<|\beta|_x<\phi^{-3}|\alpha|_x+15$.
\end{enumerate}
Here $C$ is a universal constant that we don't care about.
\end{theorem}

We think of the Near Reduction Theorem and the Far Reduction
Theorem as geometric versions of our Descent Lemma II.
The Far Reduction Theorem is similar to the theoretical
argument we gave in order to reduce the Descent Lemma II
to a computer calculation, and the Near Reduction Theorem
is similar to this computer calculation. Indeed, the
proof of the Far Reduction Theorem is almost the
same as the proof of the theoretical part of the
Descent Lemma II.

\subsection{Discussion}
\label{tricky}

Recall that $\widehat A$ and $\widehat B$ are
partitioned into atoms such that the first
return map of $\widehat \Psi$ is entirely
defined on the interiors of these atoms.
Experimental evidence leads us to the
following conjecture.

\begin{conjecture}
\label{subtle}
No point on the boundary of an atom has a
well defined $\widehat \Psi$ orbit.
\end{conjecture}

Conjecture \ref{subtle} is subtle, because there
are boundary points on which at least
$2^{18}$ iterates of $\Psi$, both forwards
and backwards, are well-defined. In particular,
thousands of iterates of the first return map
$\widehat \Psi|\widehat A$ are well-defined
on such points.

Were the conjecture true, we could eliminate
the restriction to generic orbits.
The vertices of the atoms in the Renormalization
Theorem all lie in $\Z[\phi]^3$.  Given the
formula for the map $\widehat \Theta$,
 the preimage in $\Sigma$ of these
faces is a countable discrete set of 
non-horizontal lines defined over $\Z[\phi]$.
For this reason, $\Theta$ never maps a point
of a generic orbit into the boundary of one
of the atoms.  This is why we work with generic orbits.

\newpage

\section{Applications}
\label{apps}

\subsection{Proof of Theorem \ref{ULE}}
\label{qi}

We will first consider the case of Theorem \ref{ULE} in
which $y_2=R(y_1)$.  Here $R$ is the circle
renormalization map. After we finish this case, we will
deal with the general case.
Note that the winding number of an orbit coincides
with its $\Psi$-period. So, we will work with the
$\Psi$-period rather than the winding number.
We will often just say {\it period\/} in place
of $\Psi$-period.

Suppose that $O_1$ and $O_2$ respectively
are generic orbits having sufficiently high period.
We have already explained the renormalization
map $O_1 \leadsto O_2$.
This map is $3$-to-$1$ on the level
of orbits, and the inverse images of
a single orbit lie in different horizontal planes. 
So, in our case, we get a bijection between
a certain collection of orbits on
$\R^2_{y_1}$ and a certain collection
of orbits in $\R^2_{y_2}$. These collections
contain all generic infinite orbits and all
generic periodic orbits of sufficiently high
period.

\begin{lemma}
\label{coarse}
$O_1$ and $O_2$ are coarsely equivalent, and the
coarse equivalence constant is completely
uniform -- not even dependent on the parameters.
\end{lemma}

\startproof
We use the language of the Far Reduction Theorem.
We decompose $O_1$ into $\widehat A$-cores
$\{\alpha_k\}$.  At the same time we
decompose $O_2$ into $\widehat B$-cores
$\{\beta_k\}$.  We set things up so that
$\alpha_k \leadsto \beta_k$ for all $k \in \Z$.
Let $\lambda$ be either $1$ of $\phi^{-3}$.
According to the Far Reduction Theorem, the 
map $(x,y) \to (\lambda x,y)$ carries
$\alpha_k$ to $\beta_k$, up to a uniformly
bounded error, independent of $y_1$ and $y_2$.
Our result follows immediately from this.
\endproof

\noindent
{\bf Remark:\/}
We have established the existence of a coarse
equivalence constant that does not depend on
the parameters, and this is even stronger than what we are
claiming in Theorem \ref{ULE}.  The difference
is that the parameters in Theorem \ref{ULE} might
be related by a long chain of renormalizations,
and the coarse equivalence constants will
probably depend on the length of the chain.
In view of the freezing phenomenon discussed
in \S \ref{freeze1}, there couldn't possibly
be a uniform constant that worked for all
equivalent parameters.
\newline

Now we turn to the question of local equivalence.

\begin{lemma}
\label{local}
Corresponding orbits are locally equivalent.
\end{lemma}

\startproof
In the periodic case, each orbit is contained inside
a periodic tile and there is nothing to prove.
So, suppose that $O_1 \leadsto O_2$ and these
orbits are both infinite.
We think of the strip $\Sigma$ as a subset of $\widehat \Sigma$.
With this interpretation, the action of $\Psi$ on
$\Sigma$ coincides with the action of $\widehat \Psi$
on $\widehat \Sigma$.  So, in our proof we can think
of $O_1$ and $O_2$ as orbits of $\widehat \Psi$.
Every orbit is homogeneous, so it suffices to
produce points $p_j \in O_j$, a disk $\Delta_j$
containing $p_j$, and a similarity from
$\Delta_1$ to $\Delta_2$ which
maps $p_1$ and $p_2$ and (generically)
conjugates $\widehat \Psi|\Delta_1$ to
$\widehat \Psi|\Delta_2$.

We can just choose $p_1 \in \widehat A$ and
$p_2 \in \widehat B$ to be points are
related by $\widehat R$. We choose $\Delta_1$
small enough so that $\Delta_1 \subset \widehat A$,
and then we let $\Delta_2= \widehat R(\Delta_1)$.
By the Renormalization Theorem, $\widehat R$
conjugates  conjugates 
$\widehat \Psi|\Delta_1$ to 
$\widehat \Psi^{\pm 1}|\Delta_2$ on generic points.
This is exactly what we need.
\endproof

Now we turn to the general case of Theorem
\ref{ULE}. Let
${\cal O\/}(y)$ denote the set of generic
orbits in $\R^2_y$.
Say that a {\it full subset\/}
of ${\cal O\/}(y)$ is a set that
contains all infinite orbits and all
periodic orbits having sufficiently high
$\Psi$-period.

 Reflection in the $x$-axis
conjugates $\Psi$ to $\Psi^{-1}$, and
sets up an obvious and canonical bijection
between ${\cal O\/}(y)$ and
${\cal O\/}(-y)$.  (Here we think of
$-y$ as an element of $\R/2/Z$.)
So, for Theorem \ref{ULE},
it suffices to consider the case of
two parameters that equivalent under
the group $G_2^+$ studied in \S \ref{sdl}.

By the Descent Lemma II, two such parameters
$y_1$ and $y_2$ are such that
$$y=R^{n_1}(y_1)=R^{n_2}(y_2).$$
Let $y_1'=R(y_1)$ and $y_1''=R(y_1')$, etc.
The work above gives us a bijective
correspondence between full
subsets of ${\cal O\/}(y_1)$ and
${\cal O\/}(y_1')$.
Similarly, we get a bijection
between ${\cal O\/}(y_1')$ and
${\cal O\/}(y_1'')$.  And so on.
Composing all these bijections, we get
a bijection between
a full subset of ${\cal O\/}(y_1)$
and a full subset of ${\cal O\/}(y)$.
Call this {\it the first main bijection\/}.
From the lemmas proved above, corresponding
orbits are locally and coarsely equivalant,
and the coarse equivalence constant only
depends on the parameter.  

We get the same results for $y_2$ in place
of $y_1$.  Call this the {\it second main
bijection\/}.  Composing the first and
second main bijections, we get the
bijection between full subsets of
${\cal O\/}(y_1)$ and
${\cal O\/}(y_2)$ which has all the
properties advertised in Theorem \ref{ULE}.
This completes the proof.

\subsection{Bringing Orbits into View}
\label{renorm4proof}

Now we explore the consequences of the
Near and Far Reduction Theorems. We
use the same equivalence relation $(\to)$
as in the Near Reduction Theorem.

\begin{lemma}
\label{compress0}
Let $O_1$ be any generic infinite $\Psi$ orbit.
There is a sequence of orbits
$O_1 \leadsto ... \leadsto O_n$ such
that $O_n$ intersects the rectangle
$[-24,24] \times [-2,2]$.  If
$O_1$ intersects the region $[-N,N] \times [-2,2]$
then we can take $n \leq 4N$.
\end{lemma}

\startproof
We define the type of an orbit to be the
type of the $\widehat A$-cores comprising it.
Let $O_1$ be an infinite $\Psi$ orbit.
Let $O_2$ be such that $O_1 \leadsto O_2$.

We define $|S|_x$ exactly as in Equation
\ref{minx}, except that we use ``inf''
in place of ''min'' because we are
dealing with an infinite set.

Suppose first that $O_1$ has type $2$.
and $|O_1|_x>24$.
Let $\alpha$ be any
core and let
$\beta$ be such that $\alpha \leadsto \beta$.
Far Reduction Theorem tells us
that
$$|\beta|_x \leq \phi^{-3} |\alpha|_x+15.$$
When $|\alpha|_x>24$, a bit of
arithmetic tells us that
$|\beta|_x<|\alpha|_x-1$.
So, when $|O_1|_x>24$, we have
$|O_2|_x<|O_1|_x-1$.

Suppose now that $O_1$ has type $2$.
In this case, given the nature of $R$,
we have $O_1 \leadsto O_2 \leadsto O_3$, where
$O_2$ has type $2$.  The Compression
Theorem, combined with the same argument
as above, tells us that
$$|O_2|_x \leq |O_1|_x + 12; \hskip 30 pt
|O_3|_x=\phi^{-3}|O_2|_x+15.$$
Therefore
$$|O_3|_x \leq \phi^{-3} |O_1|_x + 12 \phi^{-3} + 15.$$
In this case, a bit of arithmetic tells us that
$|O_3|_x<|O_1|_x-1/2$ provided
that $|O_1|_x>24$.
In either case, at most $2$ renormalizations
the closest point of $O_1$ closer by $1/2$ units.
\endproof

Combining Lemma \ref{compress0} with the
Near Reduction Theorem, we obtain the
following result.

\begin{theorem}
\label{renorm4}
Let $O_1$ be any generic infinite orbit.
Then there is some $m$ such that
$O_1 \to ... \to O_m$, and $O_m$ intersects
the fundamental triangle $T$.
More precisely, if $O_1$ is an unbounded orbit that
intersects $[-2,2] \times [-N,N]$ then
we can take $n<4N+C_y$.
Here $C_y$ is a constant
that depends only on the value of
$y \in \R/2\Z$ such that $O_1 \subset \R^2_y$.
\end{theorem}

\subsection{Proof of Theorem \ref{penrose9}}

Theorem \ref{penrose9} deals
with the fundamental triangle $T$ and
the corresponding tiling $\cal T$.
Recall from  \S \ref{notation} that
${\cal T\/}^+={\cal T\/}$ and that
${\cal T\/}^-$ is the image of
${\cal T\/}$ under the reflection in the vertical
line $x=1$.  We really only care about the
tiling $\cal T$, but we find it convenient to
consider both ${\cal T\/}^+$ and
${\cal T\/}^-$ at the same time.

\begin{lemma}
\label{grow}
Let $p$ be a generic point contained in the
interior of a tile of ${\cal T\/}^{\pm}$. Then
$p$ has a periodic orbit.
There is an upper bound on the period of $p$ that depends
only on the depth of the tile.
Finally, the period of $p$ tends
to $\infty$ as the depth of the tile tends
to $\infty$.
\end{lemma}

\startproof
We use the language of the Fundamental Orbit Theorem.
We check directly that the tiles of
${\cal T\/}(2)$ are periodic tiles.
We handle the remaining tiles by induction.

Suppose that the first statement of the lemma holds
for all tiles having depth at most $n-1$.
Let $J_1$ be a tile of depth $n$ and
let $p \in J_1$ be a generic point.
There are indices $i$ and $j$, and some choice of
$\pm$, such that $J_1 \subset T_{ij}^{\pm}$ and
$J'=R_{ij}^{\pm}(J)$ is a tile of depth $n-1$.
Let $J_2=\psi^{-1}(J')$.  Then the lemma holds
for all generic points in $J_2$.
Let $s=\psi^{-1} \circ R_{ij}^{\pm}$ be the
similarity carrying $J_1$ to $J_2$.
By the Fundamental Orbit Theorem, the
orbit $p_1 \in J_1$ is periodic if and only
if the orbit $p_2=s(p_1) \in J_2$ is periodic.
By induction $p_2$ is a periodic point.
Hence, so is $p_1$.  This proves the first
statement of the lemma.

For the second statement, observe that the process
of renormalization shortens an orbit by at most
a factor of $812$.
So, the period of generic point in a tile of
depth $n$ is at most $812^n$. (This is a terrible estimate.)

For the third statement, we recall Statement 2
of the Renormalization Theorem:  The renormalization
operation shortens the orbit length.  So, if
$p$ lies in a tile of depth $n$, then $p$ has
period at least $n$.  (This is another terrible estimate.) 
\endproof

\begin{lemma}
\label{tile}
Suppose $p \in T$ has a well-defined orbit but
$p$ does not lie in the interior of a tile of
$\cal T$.  Then $p$ has an infinite orbit.
\end{lemma}

\startproof
If $p$ had a periodic orbit then some neighborhood
$U$ of $p$ would be such that all points in $U$
were periodic, with the same period. But $U$
necessarily contains infinitely many tiles
of $\cal T$. This contradicts Lemma \ref{grow}.
\endproof

\begin{lemma}
No point on the boundary of tile of $\cal T$ has a
well-defined orbit.
\end{lemma}

\startproof
Lemma \ref{tile} shows that $p$ cannot have a
periodic orbit.
On the other hand, suppose that $p \in \partial J$ has an
infinite orbit, where $J$ is a tile of $\cal T$.
  Then for any $n$, there is
an open neighborhood $U_n$ such that
$p \in U_n$ and the first $n$ iterates of
$\Psi$ are defined on $p$.  These iterates
all act by translation, so no point in $U_n$
has period less than $n$.  But $U_n$
contains some generic points of $J$, no matter
how large $n$.  This contradicts the second
statement of the previous lemma.
\endproof

\begin{corollary}
Every tile of $\cal T$ is a finite union of
orbit tiles.
\end{corollary}

\startproof
Let $J$ be a tile of $\cal T$.  The same argument
as in the first half of Lemma \ref{tile} shows that
there is a uniform bound on the period of $p$.
Let $N$ be this bound.
Since $\Psi$ is a polygon exchange map, there is
a {\it finite\/} set of lines in the plane such
that the first $N$ iterates of $\Psi$ are
defined in the complement of these lines.
These lines partition $J$ into finitely many
smaller convex polygons, and $\Psi$ is
well defined and periodic on the complement
of each of these lines. This shows that $J$ is
covered by a finite union of periodic tiles.
Since no point of $\partial J$ has a well-defined
orbit, none of these periodic tiles crosses
the boundary of $J$.
\endproof

Theorem \ref{penrose9} follows from the results above.
We prove one more related result in this section.

\begin{lemma}
\label{uniform}
Let $L$ be any horizontal line segment such that
$L \cap T^{\pm}$ contains contains no generic
points with infinite orbits.  Then there is a uniform upper
bound to the period of any orbit on $L \cap T^{\pm}$.
\end{lemma}

\startproof
When the height of $L$ lies in the interval $I_1$,
$L \cap T^{\pm}$ is contained in the closure of
the union of the two largest octagonal tiles.
The result is certainly true in this case.
In general, we either have $R^n(y)=[0]$ or
$R^n(y) \in I_1$ for some $n$.  In either
case, it follows from the symmetry of $\cal T$
that $L \cap T^{\pm}$ is
contained in the closure of the union of a finite
number of orbit tiles.  But then our bound
follows from Lemma \ref{grow}.
\endproof

\subsection{Unboundedness}

Theorems \ref{ULE} and \ref{renorm4} are the
workhorses in our overall proof.  They allow
us to transfer statements about orbits that
intersect $T$, the fundamental triangle,
to statements about orbits in general.
The results in this section illustrate the
technique.

We say that a {\it fundamental orbit\/} is an
infinite orbit that intersects the fundamental
triangle $T$.   

\begin{lemma}
\label{unbound}
Every fundamental orbit is unbounded in both directions.
\end{lemma}

\startproof
Every fundamental orbit is contained in the union of
two lines, and every line intersects $S$, the
fundamental fractal, in a nowhere dense set.
Now we apply our unboundedness criterion,
Lemma \ref{UNB}.
\endproof

Now we promote Lemma \ref{unbound} to a statement about
all orbits.

\begin{lemma}
\label{argue}
Every infinite orbit is unbounded in both directions.
\end{lemma}

Suppose first that $O$ is a generic infinite orbit.
By Theorem
\ref{renorm4}, we can write $O=O_1 \to ... \to O_n$, where $O_n$ is
an infinite orbit that intersects $T$ and is coarsely
equivalent to $O_1$.  But we have already seen that
$O_n$ is unbounded in both directions.  Hence $O_1$ is
also unbounded in both directions.

Now suppose that $O$ is an infinite (non-generic)
orbit that is bounded in, say,
the forward direction.  Then, by compactness, some point
of $O$ lies in the accumulation set of $O$.  But then
every point of $O$ lies in the accumulation set of $O$.
This implies that the set of accumulation points of $O$
is uncountable.  Indeed, this accumulation set
contains a Cantor set.  

The set of non-generic points on any horizontal line
is countable, and $O$ lies on $2$ horizontal lines.
(We take the $\Psi$-orbit, as usual.)  Hence, there is some
generic accumulation point $p$ of $O$. 
Note that $p$ necessarily has an
infinite orbit.  Hence $p$ has an orbit that
is unbounded in both directions.

In particular, for every $N$ there is some
$\epsilon$ such that any point within
$\epsilon$ of $p$ follows the orbit of
$p$ for $N$ steps in either direction.
Hence, $p$ cannot be the accumulation point
of an orbit that is bounded in one direction
or the other.  This contradiction shows
that $O$ is unbounded in both directions.
\endproof

\subsection{Ruling out Some Heights}

In this section we take care of a few annoying
technical details that will make our 
arguments in the next section go more smoothly.

\begin{lemma}
Let $y=\phi$.  Then $\R^2_y$ has no generic unbounded orbits.
\end{lemma}

\startproof
By the Descent Lemma I (or a direct calculation), the
point $\phi$ is $2$-periodic with respect to $R$.
We have $R(\phi)=2-\phi$. Let $L_1$ denote
the horizontal line $y=\phi$.
Let $L_2$ denote the horizontal line $y=2-\phi$.
Let $T$ be the fundamental triangle.
Note that $L_1 \cap T$ is just a single point,
and this point is non-generic.  At the same
time, $L_2 \cap T$ is contained entirely
in the periodic tile $K_0$.  Hence
neither line intersects $T$ in any point
that could have a generic unbounded orbit.

Suppose, on the other hand, that $\R^2_y$ has
a generic unbounded orbit $O_1$.  By 
Theorem \ref{renorm4} there is some $m$
such that $O_1 \to ... \to O_m$ and
$O_m$ intersects $T$.  The intersection
$O_m \cap T$ is a generic unbounded orbit
by Theorem \ref{ULE}.  The operation of
replacing an orbit by an associate simply
switches from $R^2_y$ to $R^2_{2-y}$.
Hence, $O_m \cap T$ is a point of
either $L_1 \cap T$ or $L_2 \cap T$.
This is a contradiction.
\endproof

\begin{corollary}
Suppose that $y \in C-C^{\#}$.
Then $\R^2_y$ has no generic unbounded orbits.
\end{corollary}

\startproof
We have $y \in \Z[\phi]$ by the construction of $C$.
By Lemma \ref{mod1}, we see
that $y = m + n\phi$, with $m$ even.
But then $y \sim 0$ or $y \sim \phi$.
Here $\sim$ denotes $G_2$-equivalence.
The map $\Psi$ is the identity on $\R^2_0$.
So, Theorem \ref{ULE} rules out the
possibility that $y \sim 0$.  On the
other hand, Theorem \ref{ULE} and the preceding
lemma rule out the possibility that $y \sim \phi$.

\begin{corollary}
\label{gap}
Any generic unbounded orbit that intersects the
fundamental triangle $T$ contains a point of
the form $(x,y)$ where $y \in C^{\#}$.
\end{corollary}

\startproof
By Theorem \ref{penrose9}, and the structure
of the fundamental fractal $S$, any such point
$(x,y)$ must have $y \in C$.  But our
previous result rules out the
possibility that $y \in C-C^{\#}$.

\subsection{Proof of Theorem \ref{self}}
\label{analy}

Now we turn to the question of self-accumulation.
Lemma \ref{qi1} tells us, in particular, that
there is a canonical bijection between the
$\Psi$ orbits in $\Sigma_+$ and the $\Psi$
orbits in $\Sigma_-$.  We call two such
orbits {\it partners\/}.  
Just for this section, we introduce the notation
$O_1 \Rightarrow O_2$ to mean that $O_1 \leadsto O_2'$,
where $O_2'$ is the partner of $O_2$.
The following result is just a reformulation
of the Fundamental Orbit Theorem.

\begin{corollary}
Let $p_1 \in T_{ij}^{\pm}-{\cal T\/}(2)$ be a generic point
and let 
$p_2 =R_{ij}^{\pm}(p_1)$.
Let $O_1$ and $O_2$ respectively be the orbits
of $p_1$ and $p_2$. Then $O_1 \Rightarrow O_2$.
\end{corollary}

\startproof
If $O_2'$ is an orbit in $\Sigma_-$, and
$\psi(O_2')$ intersects $\Sigma_+$, then
the partner of $O_2$ is the orbit of
any point of $\psi(O_2')$.
The corollary follows from this fact, and
from the Fundamental Orbit Theorem.
\endproof

Let $U_y$ denote the set of unbounded orbits
in $\R^2_y$.
We also have the following lemma,
which is just a reformulation of part
of Theorem \ref{ULE}.

\begin{lemma}[Bijection Principle]
Let $O_1$ and $O_1'$ be two generic orbits
in $U_y \cap \Sigma_+$ such that
$O_1 \Rightarrow O_2$ and
$O_1' \Rightarrow O_2$.  Then $O_1=O_1'$.
\end{lemma}

\startproof
We have $O_1' \leadsto O_2'$ and
$O_1 \leadsto O_2'$ where $O_2'$ is
the partner of $O_2$. 
On the level of infinite orbits, the
renormalization map is $3$-to-$1$.
At the same time, the $3$ preimages
of an orbit lie in different sets $U_y$.
Since $O_1$ and $O_1'$ lie in the
same set $U_y$, we have $O_1=O_1'$.
\endproof

\begin{lemma}
Let $p \in T$ be a generic fundamental orbit.
Then the orbit of $p$ is
self-accumulating.
\end{lemma}

\startproof
We already know that the orbit of $p$ is
unbounded in both directions.
By Lemma \ref{gap}, we have
$p=(x,y)$, where $y \in C^{\#}$.
Let $O(p)$ denote the orbit of $p$.
We refer to the notation and terminology in \S \ref{rset}.
Let $\Upsilon(p)$ denote the remormalization set for $p$.
Let $q \in \Upsilon(p)$.
Here $p$ and $q$ lie on the same horizontal line, and
\begin{equation}
p=p_0 \to p_1 \ldots \to p_n; \hskip 30 pt
q=q_0 \to q_1 \ldots \to q_n; \hskip 30 pt p_n=q_n.
\end{equation}
By the Corollary above,
\begin{equation}
O(p_0) \Rightarrow ... \Rightarrow O(p_n); \hskip 30 pt
O(q_{n-1})  \Rightarrow O(q_n)=O(p_n).
\end{equation}
Since $p_{n-1}$ and $q_{n-1}$ lie on the
same horizontal line, the Bijection Principle
says that $O(q_{n-1})=O(p_{n-1})$.
It now follows from induction on $n$ that
$O(p)=O(q)$.  

Since $q$ is an arbitrary point of
$\Upsilon(p)$ we now know that
\begin{equation}
\Upsilon(p) \subset O(p).
\end{equation}
By the Density Lemma,
$\Upsilon(p)$ is dense in $\Lambda=S \cap L$, where
$L$ is the horizontal line containing $p$.
By the Horizontal Lemma, $\Lambda$ is a
Cantor set. Hence $p$ is an accumulation point of
$O(p)$.  But, $O(p)$ is homogeneous.
Hence, every point of $O(p)$ is an accumulation
point of $O(p)$.   This proves that
$O(p)$ is self-accumulating.
\endproof

Now we prove Theorem \ref{self}.
Let $O_1$ be a generic infinite orbit.
By Theorem \ref{renorm4}, we have
$O_1 \to ... \to O_n$ where $O_n$
is an infinite orbit that intersects
the fundamental triangle.  But we
have already shown that $O_n$ is 
self-accumulating.  By Theorem
\ref{ULE}, the orbits
$O_1$ and $O_n$ are
locally similar.  Hence $O_1$
is also self-accumulating.

\subsection{Proof of Theorem \ref{penrose1}}

We've already proved that every orbit is either periodic
or unbounded in both directions.  It only remains to
show that the union $U$ of the unbounded orbits has
Hausdorff dimension $1$.
By Lemma \ref{DIM} the set of generic points in $S$
has Hausdorff dimension 1. All such points
have well-defined orbits and, by Theorem
\ref{penrose9}, all the orbits
are unbounded.  Hence $U$ contains a $1$ dimension set.
Hence $\dim(U) \geq 1$.

Let $G$ be the set of generic points in $\R^2$.
note that $\R^2-G$ has Hausdorff dimension $1$
because it is a countable set of lines.
Combining Theorems \ref{ULE} and \ref{renorm4}, we
see that every generic point $p \in U$ has a neighborhood
$\Delta$ such that $\Delta \cap U$ is similar to
a subset of $S$.  Hence $\dim(U \cap G \cap \Delta)\leq 1$.
But we can cover any compact subset of $U$ by
finitely many such neighborhoods.  Hence $\dim(U \cap G) \leq 1$.
Since $U-G$ is contained in a countable family of lines, we have
$\dim(U-G) \leq 1$.   Hence $\dim(U) \leq 1$.

\subsection{Proof of Theorem \ref{penrose2}}

First we establish part of Theorem \ref{non-generic}.
Let $U_y^*$ denote the set of generic unbounded
orbits in $\R^2_y$.

\begin{lemma}
\label{equiv}
$U_y$ is empty if
$U_y^*$ is empty.
\end{lemma}

\startproof
Let $R$ be the renormalization map.
Suppose $R^n(y)=0$ for some
$n$.  All orbits on the line $y=0$ have
period $1$, and renormalization decreases
periods by at most a factor of $812$.
Hence, there is a uniform bound to the
period of any generic point in
$U_y$. 

On the other hand, suppose that $y$ is not
in the inverse image of $0$. Let $p_1$ be
any generic point that is, say, within
$1$ unit of $y$.  Let $O_1$ be the orbit
containing $p_1$. Combining
Lemma \ref{compress0} and the second statement
of the Near Reduction Theorem, we see that
there is some uniform $m$ such that
$O_1 \to ... \to O_m$ and $O_m$ intersects
$T$.  But the horizontal line containing
$O_m \cap T$ has no generic points
with infinite orbits.  Hence, by
Lemma \ref{uniform}, there is a uniform
upper bound to the period of $O_m$.
Hence, there is a uniform upper bound
to the period of $O_1$.

Finally, we can take a sequence of generic
points approximating a supposed infinite orbit.
This is incompatible with the uniform upper
bound we have on the periods of this
approximating sequence.
\endproof

Now we prove Theorem \ref{penrose2}.
Let $Y \subset \R/2\Z$ denote those $y$
such that $\R^2_y$ contains unbounded orbits.
Let $y \in C^{\#}$.   The horizontal line of
height $y$ intersects $S^{\#}$ in an
uncountable set which must have points with
well-defined orbits.
By Theorem \ref{penrose9},
these points are not periodic.  By
Lemma \ref{unbound}, these points
have unbounded orbits.
Hence $C^{\#} \subset Y$.
Theorem \ref{ULE} now shows that $Y$ contains
all points $y$ such that $y \sim c$ and
$c \in C^{\#}$.  

Now for the converse. Suppose that $y_1 \in Y$.
By Lemma \ref{equiv}, we can assume
that $\R^2_{y_1}$ has a generic unbounded
orbit $O_1$.  By Theorem
\ref{renorm4}, we have
$O_1 \to ... \to O_m$ with
$O_m$ intersecting $T$.
Let $y_k$ be such that
$O_k$ is an orbit of 
$\R^2_{y_k}$.  Passing to
associates preserves the
$G_2$-equivalence class,
and so does $R$.  Hence $y_1 \sim y_m$.

By Theorem \ref{ULE}, the orbit
$O_m$ is generic and unbounded.
By Theorem \ref{penrose9}, we have
$y_m \in C$.
By Lemma \ref{gap}, we have
$y_m \not \in C-C^{\#}$.
Hence $y_m \in C^{\#}$
and $y_1 \sim y_m$, as
desired.  Finally, we compute easily
that $\dim(C)=\log(3)/\log(\phi^3)$.

\subsection{Proof of Theorem \ref{winding}}

As we remarked above, the concepts of
winding number and $\Psi$-period coincide.
We will work with $\Psi$-period, as usual.

\begin{lemma}
\label{linear}
Every horizontal line in $\Sigma$ contains
a dense set of periodic orbits.
\end{lemma}

\startproof
Let $U$ denote the union of unbounded orbits.
If $L \cap U$ is empty, there is nothing to prove.
Otherwise, by Theorem \ref{renorm4} and Theorem \ref{ULE},
the set $L \cap U$ is locally similar to
$L' \cap U$, where $L'$ is a line segment in
$T$, the fundamental triangle.  But $U \cap T$
is nowhere dense in each horizontal line.
Hence $U \cap L'$ is nowhere dense in $L'$. 
Hence $U \cap L$ is nowhere dense in $L$.
\endproof

\begin{lemma}
\label{partition1}
Let $y \in \R/2\Z$ be any value other than 
$$0; \hskip 30 pt 4-2\phi; \hskip 30 pt -2+2\phi; \hskip 30 pt 2$$
Then $\R^2_y$ contains a periodic orbit
that intersects $\widehat B$, the renormalization set.
\end{lemma}

\startproof
Let $y \in \R/2\Z$ be any value other than
the ones listed. Let $\Pi$ be the horizontal
plane of height $y$.  The plane $\Pi$
intersects $\widehat B$ in
an open set.  Let $L \subset \Sigma$ be
the horizontal line of height $y$.
Since $\Theta(L)$ is dense in $\Pi$, there
is an open subset $V \subset L$ such
that $\Theta(V) \subset \widehat B$.
By Lemma \ref{linear}, there is a dense subset
of $V$ consisting of periodic points.
Any periodic point in $V$ works for us.
\endproof

\begin{lemma}
\label{big}
Let $y \in \R/2\Z$. Suppose it never happens
that $R^n(y)=[0]$. Then $\R^2_y$ contains 
periodic points having arbitrarily high $\Psi$-period.
\end{lemma}

\startproof
Our proof refers to the terminology used in the
Renormalization Theorem.  Let $y_n=R^n(y)$.
For our analysis, we identify $\Sigma$
with a dense subset of $\widehat \Sigma$.
Define $\Sigma_n=\R^2_{y_n} \cap \Sigma.$
Let $\widehat \Sigma_n$ be the horizontal 
plane of height $y_n$.

 Referring to Lemma \ref{partition1},
the heights of the excluded horizontal planes all
lie in the inverse image of $0$.  For this reason,
$\widehat \Sigma_n$ contains a periodic tile $P$
whose orbit intersects $\widehat B$ in an
atom that is a translate of $P$.

Since $\Sigma_n$ is dense in $\widehat \Sigma_n$,
the periodic tile $P$ must intersect $\Sigma_n$ in
an open interval.  In particular
$\Sigma_n$ contains a generic periodic point
whose orbit $O_n$ intersects $\widehat B$.
By the Renormalization Theorem, we can find
a generic orbit $O_{n-1}'$ in $\Sigma_{n-1}$ such that
$O_{n-1}' \leadsto O_n$.  The period of
$O_{n-1}'$ is longer than the period of $O_n$.

$O_{n-1}'$ necessarily intersects
$\widehat A$. By the Renormalization Theorem,
$O_{n-1}'$ also intersects $\widehat B$.  So, we can
set $O_{n-1}=O'_{n-1}$ and repeat the above
argument to produce a generic periodic orbit
$O_{n-2}'$ in $\Sigma_{n-2}$ that intersects
$\widehat B$. And so on.
The final orbit $O_1$ has period at least $n-1$.
But $n$ is arbitrary.
\endproof

The full inverse image of $[0]$ is precisely $2\Z[\phi]$.
Hence, if $y \not \in 2\Z[\phi]$, then
$\R^2_y$ contains orbits of arbitrarily
high $\Psi$-period.

\begin{lemma}
Suppose that
$y=m+n\phi$ where both $m$ and $n$ are even.
Then there is a uniform bound on the $\Psi$-period
of any point of $\R^2_y$ with a well defined orbit.
\end{lemma}

\startproof
We consider generic points first.
Suppose $\R^2_y$ contains a generic orbit with an enormous
winding number -- either finite or infinite.
There is some $n$ such that
$R^n(y)=[0]$.  If the winding number of our orbit
is too large, then we can renormalize this orbit
more than $n$ times. But this contradicts the
fact that $\Psi$ is the identity on $\widehat \Sigma_0$,
and no such orbit on this plane has a renormalization.

So, every generic point of $\R^2_y$ is periodic and
we have a uniform bound on the periods.
But any point in $\R^2_y$ can be
approximated by a sequence of generic points.
The uniform bound on the approximating sequence
immediately gives the same uniform bound on the
limit, provided that the limit has a well-defined orbit.
\endproof

Let $y \in 2\Z[\phi]$ be some point
for which the following property $(*)$ holds:
\begin{equation}
\label{noimage}
R^k(y) \not = [1]; \hskip 30 pt k=0,...,n.
\end{equation}
Then the same argument as in Lemma \ref{big}
shows that $\R^2_y$ contains generic points
having $\Psi$-period at least $n-1$.  
Finally, for any fixed $n$, there are only
finitely many
points of $2\Z[\phi]$ which fail to satisfy
Equation \ref{noimage}.
This establishes Theorem \ref{winding}.

\subsection{Proof of Theorem \ref{stable}}

This is one of the more subtle results.  We will
try to break the proof down into small steps.

\begin{lemma}
\label{small period}
Let $y$ be one of the $4$ values
$$-1+\phi \hskip 18 pt -3+3\phi, \hskip 18 pt 3-\phi, \hskip 18pt 5-3\phi.$$
Let $\Lambda$ be the intersection of the
line $\R \times \{y\}$ with the fundamental
triangle $T$.  Then there is a uniform upper
bound to the period of any point that is
sufficiently close to $\Lambda$.
\end{lemma}

\startproof
Let $\Lambda_1,...,\Lambda_4$ be the $4$
line segments in question.
We check by direct inspection that there is some
$\epsilon_0>0$ with the following property.
The $\epsilon_0$ neighborhood of $\Lambda_k$
is contained in a finite union of tiles of
$\cal T$.  Indeed,
the two lines $\Lambda_2$ and $\Lambda_4$ are disjoint from $T$,
the interior of $\Lambda_1 \cap T$ is contained in
the interior of a single periodic tile, and
the interior of $\Lambda_3 \cap T$ is contained
in the interior of the closure of two periodic tiles.
\endproof

Let $y=m+n\phi$ where $m$ and $n$ are odd.
Let $L=\R \times \{y\}$. Let $N=N_y$ be as in
Lemma \ref{descent3}. 

\begin{lemma}
\label{descent5}
Let $k>N$ be a fixed integer.
There is some $\epsilon>0$ with the following
property.  Suppose that $O_1$ is an orbit that
contains a point within $\epsilon$ of $L$.
Suppose $O_1 \to ... \to O_k$ and that
$O_k$ intersects the
fundamental triangle.  Then there is a uniform
bound on the period of $O_1$.
\end{lemma}

\startproof
By Lemma \ref{descent3}, and the fact that
$R^k$ expands distances by at most
$\phi^{3k}$, the orbit
$O_k$ must intersect the fundamental triangle
in a point very close to one of the $4$ line
segments listed in Lemma \ref{small period},
provided that $\epsilon$ is small enough.
The point here is that, when we keep track of
which set $\R^2_y$ contains our successive orbits,
the renormalization operation implements the map $R$
and switching from an orbit to an associate
implements the map $S$.  

If $O_k$ is very close to one of the segments
in Lemma \ref{small period}, then there is a
uniform upper bound to the period of $O_k$.
But the renormalization operation, and the
operation of switching to an associate, only
can decrease the period by a uniformly bounded
factor.  Hence, there is a uniform upper bound
on the period of $O_1$.
\endproof

\begin{lemma} 
\label{descent4}
Let $p \in L$ be some point.  Let $M$ be any positive
integer.  Then there is some $\epsilon>0$ with the
following property. If $O_1$ is a generic periodic
orbit that comes within $\epsilon$ of $p$, then
there is some $k>M$ such that
$O_1 \to ... \to O_k$ and $O_k$ intersects the
fundamental triangle.  The same $k$ works for all
orbits satisfying the hypotheses.
\end{lemma}

\startproof
Combining Lemma \ref{compress0} and
the second statement of the Near Reduction Theorem,
we see that there is some $k$, which works for
all choices of $O_1$, such that
$O_1 \to ... O_k$ and $O_k$ intersects
the fundamental triangle. The only trouble
is that we might have $k<M$. However, we can
renormalize again and then, if necessary,
apply Theorem \ref{renorm4} again.  This
produces a larger value of $k$.  We keep
going like this until we arrive at some
$k>M$, and then we stop.
\endproof

\begin{corollary}
Let $p \in L$ be any point.  Then there is
some $\epsilon>0$ and some constant $Z$
 with the following property.
If $O_1$ is a generic periodic orbit that
comes within $\epsilon$ of $p$ then
$O_1$ has period at most $Z$.
\end{corollary}

\startproof
Choose $M=N$, the constant in Lemma \ref{descent5}.
If $\epsilon$ is small enough then we have
some fixed $k>M$ such that
$O_1 \to ... \to O_k$ and $O_k$ intersects
the fundamental triangle.  This fixed $k$
works for any orbit that comes within $\epsilon$ of $p$.
But now Lemma \ref{small period} applies to
$O_1$.
\endproof

Suppose that the conclusion of Theorem \ref{stable}
is false. Then we can find some $p \in L$ which is
the accumulation point of points having unbounded
orbits.  By Theorem \ref{penrose1}, we also 
have a sequence of generic periodic
points converging to a point of $L$ 
whose period tends to $\infty$.
This contradicts the corollary we have just proved.
This establishes This proves Theorem \ref{stable}.

We mention the following corollary of Theorem \ref{stable}.
\begin{corollary}
\label{corstable}
Let $y=m+n\phi$ with $m$ and $n$ odd.
Let $K$ be and compact subset of $\R^2$.
Then there is some $\epsilon>0$, depending
on $K,n,m$, such that
no point of $K$ within $\epsilon$ of
the line $L_y$ has an unbounded orbit.
\end{corollary}

\startproof
This is just an application of compactness.
\endproof

\subsection{Proof of Theorem \ref{outside}}

Let $Y_r$ be as in Theorem \ref{outside}.
As above, let $U_y$ denote the union
of unbounded orbits in $\R^2_y$.

\begin{lemma}
$Y_r$ is nowhere dense.
\end{lemma}

\startproof
Rather than work precisely with the set mentioned in
Theorem \ref{outside}, we define
$Y'_r$ to be the set $y \in Y$ such that
$U_y$ contains a generic orbit that
intersects the rectangle $[-r,r] \times [-2,2]$.
We have $Y_{r-1} \subset Y'_r \subset Y_{r+1}$, so
it suffices to prove Theorem \ref{outside} for
$Y'_r$ in place of $Y_r$.

The analysis in Lemma \ref{compress0} shows that
\begin{equation}
\label{expand}
R^{n}(Y'_r) \cup R^{n+1}(Y'_r) \subset Y'_{24}; \hskip 24 pt \forall n>4r.
\end{equation}

Suppose that $Y'_r$ is not nowhere
dense for some $r$.
Since $R^2$ is an expanding map,
and $Y'_r$ is not nowhere dense,
the left hand side of Equation \ref{expand} is
dense in $\R/2\Z$ for sufficiently large $n$.
Hence $Y'_{24}$ is dense in $\R/2\Z$.
But this contradicts Corollary \ref{corstable}.
\endproof

We have already computed that $\dim(C^{\#})=\log(3)/\log(\phi^3)$.
Observe that every neighborhood of the Penrose
kite vertex $(\phi^{-3},0)$ contains a similar copy
of $S$.  Hence $\dim(Y_r)=\dim(C^{\#})$.

\subsection{Proof of Theorem \ref{non-generic}}

When we proved Theorem \ref{penrose2}, we proved the first
statement of Theorem \ref{non-generic}.
Now we prove the second.

The set of non-generic points is
contained in a countable union of lines.
In particular, the set of non-generic
unbounded orbits is contained in a
countable union of sets of the form
\begin{equation}
L \cap (\R \times Y),
\end{equation}
where $Y$ is the set of $y$ such that
$\R^2_y$ has unbounded orbits.
By Theorem \ref{penrose2}, we have
$\dim(Y)=\log(3)/\log(\phi^3)$.
But then the union of non-generic
unbounded orbits is contained
in a countable union of sets
that have the same dimension as $Y$.
This completes the proof.

\subsection{Proof of Theorem \ref{penrose10}}

Almost every statement in Theorem \ref{penrose10}
is an immediate consequence of
Theorems \ref{penrose9}, \ref{ULE}, and
\ref{renorm4}. The one statement that is not 
immediate is that a neighborhood of the
dynamical tiling about a generic
point with unbounded orbit is {\it isometric\/}
to a small patch of $\cal T$. 

The results 
mentioned above only imply that the abovementioned
neighborhood is {\it similar\/} to a patch
of $\cal T$.  However, tracing through
our argument, we see that the similarity
factor is $\phi^{3k}$ for some integer $k$.
But $\cal T$ is a self-similar set with
expansion constant $\phi^3$.   So, we
can take $k=0$ in the similarity factor
between the neighborhood of interest to us
and a suitable patch of $\cal T$.

\subsection{Proof of Theorem \ref{ring}}

Let $y=m+n \phi$.
We already know that $1 \in C^{\#}$.  So,
by Theorem \ref{penrose2}, we know
that $\R^2_y$ has unbounded orbits
provided that $m$ is odd and $n$ is even.

Conversely, suppose that $\R^2_y$ has
unbounded orbits.
The set $C-C^{\#}$ contains points of
the form $m+n\phi$ with $m$ even and $n$
having either parity.  Combining
Theorem \ref{penrose2} and Corollary \ref{gap},
we see that $m$ cannot be even.
When $m$ is odd and $n$ is odd,
Theorem \ref{stable} shows that 
$\R^2_y$ has no unbounded orbits in this case.
The only case left is when $m$ is odd and $n$ is even,
as claimed.

\newpage

\section{Computational Methods}

\subsection{Golden Arithmetic}
\label{arithmetic}

A large percentage of our proofs involve 
exact computer calculations, done over
the number ring $\Z[\phi]$. We abbreviate
these kinds of computer calculations 
as {\it golden arithmetic\/}.  Here we
describe how the computer does golden
arithmetic, for the sake of making
our calculations completely reproducible.
There exist computer packages, such as
Pari, which perform these kinds of
calculations in extreme generality.  A
software package like Pari is vastly
more complicated than our own much
more limited collection of routines.

We represent the number $a_0+a_1\phi$ as the
integer pair $(a_0,a_1)$.  The computer can only
represent finitely many such numbers, but the 
large finite set of representable numbers is
suitable for our purposes.
\newline
\newline
{\bf The Ring Operations:\/}
We have the following obvious rules
$$
(a_0,a_1) \pm (b_0,b_1)=(a_0\pm b_0,a_1 \pm b_1);
$$ 
\begin{equation}
(a_0,a_1) \times (b_0,b_1)=(a_0b_0+a_1b_1,a_0b_1+a_1b_0+a_1b_1),
\end{equation}
which represent the ordinary ring operations in
$\Z[\phi]$.  As long as all integers stay less
than, say, $10^6$ in absolute value, the computer
adds and multiplies them correctly.
The Galois map $\tau(a+b \phi)= a-b/\phi$ is a ring automorphism
of $\Z[\phi]$. In terms of our representation, we have
\begin{equation}
\tau: \hskip 5 pt (a,b) \to (a+b,-b).
\end{equation}
\newline
{\bf Division:\/}
Since $\Z[\phi]$ is not a field, we cannot
generally perform division in $\Z[\phi]$.  However,
it occasionally happens that we know in advance
that $\gamma=\alpha/\beta$ lies in $\Z[\phi]$ and we want to
find $\gamma$ given $\alpha$ and $\beta$.
In this cases, we find $\gamma$ by computing
\begin{equation}
\gamma=\frac{\alpha \times \tau(\beta)}{\beta \times \tau(\beta)}.
\end{equation}
The number in the denominator is an integer, and we
find the quotient on the right hand side by dividing
the coefficients of the numerator by this integer.
Before returning the value, we verify that the
resulting element belongs to
$\Z[\phi]$ and satisfies the equation
$\alpha=\beta \times \gamma$.
\newline
\newline
{\bf Positivity:\/}
Now we explain how we check that
$a_0+a_1\phi$ is positive.  Our
method, applied to $a_0+a_1\phi$, 
returns a value of {\it true\/} only when
$a_0+a_1\phi>0$. We call
$a_0+a_1\phi$ {\it strongly positive\/} if
\begin{equation}
\label{fib}
a_0f_{100}+a_1f_{101}>0; \hskip 30 pt
a_0f_{101}+a_1f_{102}>0.
\end{equation}
Here $f_n$ is the $n$th Fibonacci number.
Since the successive quotients of Fibonacci
alternately over and under approximate 
$\phi$, the fact that both linear combinations
are positive guarantees that $a_0+a_1\phi$
is also positive.  In short, a strongly
positive element of $\Z[\phi]$ is positive.
Certainly there are positive numbers that
are not strongly positive.  However, we
do not encounter these numbers in our calculations.

In Equation \ref{fib},
we use the BigInteger class in Java, which
does exact integer arithmetic for integers
up to many thousands of digits long.
We might have used the BigInteger class
for all our calculations, but this would
make the calculations much slower and the
computer routines much more tedious to
program.  Thus, we use the BigInteger class
only when needed.

Equation \ref{fib} is the basis for all our
computations that involve inequalities between
elements of $\Z[\phi]$.  For instance,
to verify that $(a_0,a_1)$ represents a
number larger than $(b_0,b_1)$, we
apply the test to $(c_0,c_1)$, where
$c_j=a_j-b_j$. 
\newline
\newline
{\bf Golden Structures:\/}
The {\it GoldenReal\/}, a pair $(a_0,a_1)$ as
above, is our basic object.  We also define
more complicated objects based on the GoldenReal:
\begin{itemize}
\item A {\it GoldenComplex\/} is a pair
$x+iy$, where $x$ and $y$ are both GoldenReals.
\item A {\it GoldenVector\/} is a
tuple of GoldenReals.
\item A {\it GoldenPolytope\/} is a finite list
of GoldenVectors.
\end{itemize}
The usual operations for these objects are
done using the ring operations described above.
\newline
\newline
{\bf Remark:\/}
It is worth mentioning that our graphical user
interface mainly operates with floating point arithmetic,
for the purposes of speed.  We mainly use the
special arithmetic when we need to do rigorous
calculations for the purpose of making a proof.
The special structures are all isolated in
separate files, so as not to interfere with
the rest of the program.

\subsection{Overflow Error}
\label{overflow}

There is one computational issue that we
must face when we do exact integer arithmetic
calculations.  The computer
does not reliably perform the arithmetic operations
on very large integers.  For example, when
multiplying together ordinary integers, or
{\it ints\/}, our computer
tells us that 
$$100000 \times 100000 = 1410065408.$$
Each {\it int\/} is a length $32$ binary string.  
One of the
bits records the sign of the integer and the
remaining bits give the binary expansion.
The problem with the solution to 
$100000 \times 100000$ is that it requires
more than $32$ bits to express the answer.

To give us a bit more flexibility, we use
{\it longs\/}.  A {\it long\/} is a $64$-bit
representation of an integer.
The calculation above comes out right
when we use {\it longs\/} in place of
{\it ints\/}.

Here are some conservative bounds on what
the computer can do with {\it longs\/}
\begin{itemize}
\item The computer reliably computes
the value $a_1 \pm a_2$ provided that
we have $\max(|a_1|,|a_2|)<2^{60}$.
\item The computer reliably computes
the value $a_1 \times a_2$ provided that
we have $\max(|a_1|,|a_2|)<2^{60}$.
\end{itemize}
Whenever we perform one of the basic operations
above, we check that the conditions above hold.
The code is set up to interrupt the calculation
if one of the conditions above is not met.

The only other operations we perform are
Galois conjugation within $\Z[\phi]$ and
the guided form of division mentioned above.
For conjugation, we make sure that the 
coefficients of the GoldenReal are less than
$2^{60}$ in absolute value.

The division operation requires a special
explanation. When we perform the division
operation to {\it find\/} the quotient
$c=a/b$, we do not bound the sizes of the
coefficients.  Rather, when we then verify
that indeed $a=b \times c$, we check the
bounds of the coefficients of $b$ and $c$.
This means that the equality $a=b \times c$
has really been verified.

\subsection{Some Basic Operations and Tests}
\label{basictest}

Here we explain some basic operations we perform
on polygons and polyhedra.  When we have
golden polygons and golden polyhedra, the operations
are all done using the golden arithmetic
described above.
\newline
\newline
{\bf Positive Convexity:\/}
We say that a polygon $P$ is
{\it positively convex\/} if $P$ is convex, and
if the orientation of $\partial P$ given by the ordering
on the vertices is counterclockwise.  $P$ is
positively convex if and only if 
\begin{equation}
\label{posconvex}
{\rm Im\/}(\overline z_{21} z_{31}) \geq 0; \hskip 30 pt
z_{ij}=z_i-z_j.
\end{equation}
for  every triple $z_1,z_2,z_3$ of consecutive vertices of $P$.
Note that a positively convex polygon need not be strictly
convex.  
\newline
\newline
{\bf Containment Test for Polygons:\/}
Let $P$ be a positively convex polygon with vertices
$z_1,...,z_k$.  A point $w$ lies in the interior of
$P$ if and only if
\begin{equation}
{\rm Im\/}\Big((\overline{z_i-w})(z_{i+1}-w)\Big)>0;
\end{equation}
holds for every index $i$, with indices taken mod $k$.
\newline
\newline
{\bf Strict Convexity for Polyhedra:\/}
To check that a golden polyhedron is strictly convex, 
we exhibit a golden vector $V=V(P,v)$ such
that $v \cdot P>v' \cdot V$ for all other vertices
$v'$ of $P$.  In practice, we search for $V$ amongst
all golden vectors whose coefficients have real and
imaginary parts of the form $a+b\phi$ with $\max(|a|,|b|) \leq 2$.
In general, one would need to make a more extensive search
or else have some special information about $P$.
\newline
\newline
{\bf Golden Interior Points:\/}
The most natural interior point of
a polytope is its center of mass.
However, the center of mass of a 
golden polytope might not be a golden
vector.  Here we describe a less
canonical way of picking an interior
point of a golden polytope which
results in a golden vector.

Given two golden vectors $V_1$ and $V_2$, 
the vector
$a(V_1,V_2)=\phi^{-1} V_2+\phi^{-2} V_2$
is a golden vector that lies
on the line segment $\overline{V_1V_2}$.
Note generally, suppose we have 
a golden polytope with vertices
$V_1,...,V_n$.  We define
\begin{equation}
W_2=a(V_1,V_2); \hskip 30 pt
W_3=a(W_1,V_3); \hskip 30 pt \ldots
\hskip 30 pt
W_{n}=a(W_1,V_n).
\end{equation}
then $W_n$ is a golden vector in
the interior of the polytope.
\newline
\newline
{\bf Disjoint Interiors I:\/}
Here is one method we use to verify
that two golden polyhedra $P$ and $P'$ have
disjoint interiors.  
We exhibit a golden vector $W$
such that $V \cdot W<V' \cdot W$ for
all pairs $(V,V')$, where $V$ is a vertex
of $P$ and $V'$ is a vertex of $P'$.
In practice, we search for $W$ amongst all
vectors $E \times E'$, where $E$ is an
edge of $P$ and $E'$ is an edge of $P'$.
\newline
\newline
{\bf Disjoint Interiors II\/}
A second way we will know (without directly
testing) that two distinct polygons or polyhedra
$P_1$ and $P_2$ have disjoint interiors is
that there is
some map $f: P_k \to \R^2$ such that $f$ is
entirely defined and affine on the interior
of $P_k$, but $f$ does not extend to be
defined and affine in a neighborhood
of any point of $\partial P_k$.
If $P_1$ and $P_2$ did not have disjoint
interiors, then some point of $\partial P_1$
lies in the interior of $P_1$, and then
the map $f$ is then defined and locally
affine in a neighborhood of this boundary
point.  This is a contradiction.
\newline
\newline
{\bf Coplanarity Test:\/}
Let $V_1,...,V_k$ be finite list of golden vectors.
Consider the successive vectors
\begin{equation}
\label{normals}
N_i=(V_{i+1}-V_i) \times (V_{i+2}-V_i),
\end{equation}
with indices taken mod $k$.  The vectors
are coplanar if and only if $N_i$ and $N_{i+1}$
are parallel for all $i$.  
\newline
\newline
{\bf Raw Face Enumeration:\/}
Here is how we enumerate the faces of
a polyhedron that is given in terms of its vertices.
Let $V(P)$ be the vertex set of a polyhedron $P$.  For
each subset $S \subset V(P)$, having at least $3$ elements,
we first check if $S$ is a coplanar set of vertices,
as explained above.  If $S$ is not coplanar, we
eliminate $S$ from consideration.

Assuming that $S$ is coplanar, let
$V_1,...,V_k$ be the vertices of $S$.
Let $N_1$ be the vector from Equation \ref{normals}.
Let $V=V_1+...V_k$.
Also, let $kS$ denote the result
of scaling all vectors of $S$ by a factor of $K$.
Likewise define $kP$.
We try to check  that either $V \cdot N_1 \leq W \cdot N_1$
for all vertices $W$ of $kP$ or that
$V \cdot N_1 \geq W \cdot N_1$ for all vertices $W$ of $P$.
If one of these two things is true, then the center
of mass of the convex hull of $S$ is contained in
$\partial P$, and this means that $S$ is a subset of
a face of $P$.  The reason why we scale everything
by $k$ is that we want to work entirely in $\Z[\phi]$.
(To avoid scaling, we could have used the golden
interior point described above in place of $V$,
but we didn't.)
\newline
\newline
{\bf Edge Enumeration:\/}
Given the list of vertices and the raw list of
faces, we find the edges as follows.
A pair of distinct vertices $(v_1,v_2)$ is an edge of $P$
if any only if $P$ has $4$ distinct
vertices $v_1,v_2,w_1,v_2$ such that $(v_1,v_2,w_1)$
and $(v_1,v_2,w_2)$ are two distinct faces of $P$.
For each pair $(v_1,v_2)$, we try to find the
pair $(w_1,w_2)$ with the property just mentioned.
If we succeed, we add $(v_1,v_2)$ to the list of edges.
\newline
\newline
{\bf Polished Face List:\/}
There is one more step, which we perform for the sake
of tidiness, and also for the sake of making
some of our other computations go more quickly.
We would like the vertices of each face to be listed
so that they go in cyclic order around the face.
Our convention is that the vectors produced
by Equation \ref{normals} should be outward
normals.  
To order the vertices around a face, we first use the
raw face list to enumerate the edges of the face.
We then use the edges, in a fairly obvious way,
as a guide for placing the vertices.  
\newline
\newline
{\bf Slicing:\/} As an application of our edge and
face enumeration, we explain how we compute
the intersection of a golden polyhedron $P$ and a
horizontal plane $\Pi$ that contains no
vertices of $P$.  Every edge of $P$ will satisfy
Equation \ref{fibered} and $\Pi$ will have
height in $\Z[\phi]$, so $P \cap \Pi$ will be
a golden polygon.

Every edge of $P$ either
crosses $\Pi$ at an interior point or
else is disjoint from $\Pi$. 
We enumerate those edges $e_1,...,e_k$ of $P$
which cross $\Pi$. We order
the edges so that $e_i$ and $e_{i+1}$
lie in a common face.  Finally, we set
$v_i=\Pi \cap e_i$.  Then
$v_1,...,v_k$ are the vertices of the
polygon $P \cap \Pi$, and they
are cyclically ordered around the
boundary of $P \cap \Pi$.
\newline
\newline
{\bf Containment Tests for Polyhedra:\/}
Let $V$ be a point and let $P$ be a
polyhedron. Let $F_1,...,F_k$ be the
faces of $P$.  Let $W_1,...,W_k$ be
such that $W_j$ is a vertex of $F_j$
for each $j$.  Finally, let
$N_j$ be a vector that is perpendicular
to $F_j$.  We compute
$N_j$ by applying Equation \ref{normals} 
to the first $3$ vertices of $F_j$.
The point $V$ lies in the interior
of $P$ if and only if
$V \cdot N_j < V \cdot W_j$ for 
$j=1,...,k$.  
Similarly, $V$ lies in $\partial P$ if
and only if $V \cdot N_j \leq V \cdot W_j$
for all $j$, and there is equality for
at least one index $j$.
\newline
\newline
{\bf Remark:\/}
The tests we have described above look fairly
intensive, but for the most part we will be
applying them over and over to the same list
of $64$ polyhedra -- the ones
that define our polyhedron exchange map. For these,
we pre-compute everything and store the
results in a look-up table.

\subsection{Covering Tests}
\label{covertest}

Here we describe how we show that
a given polygon or polyhedron is
covered by a finite union of other
polygons or polyhedra.  We will
explain the polygon case first.
We then reduce the polyhedron
case to the polygon case by
a trick involving slicing.

Let $Q$ be a golden polygon and let
Let $P_1,...,P_n$ be a finite list of
golden polygons.  Here we explain the
computational test we
perform in order to prove that
$Q \subset \bigcup P_k$.
\newline
\newline
{\bf Eliminating Spanning Edges:\/}
Say that $P_k$ {\it has a spanning edge\/}
relative to $Q$
if some edge of $P_k$ intersects $Q$ in an interior
point, but has no vertices in the interior of $Q$.
In practice, we eliminate all spanning edges
by adding an extra ``vertex'' of $P_k$ to each 
edge.   
Since we want our polygons to
remain golden, we pick a golden interior
point for each edge. See above.

Our method is not guaranteed to eliminate
all spanning edges, but in practice it does.
In all cases but one, we also have
$P_k \subset Q$ for all $k$.  In this case, the only
spanning edges are those that have both
vertices on $\partial Q$. For such spanning
edges, the addition of a single vertex
does the job.  In the one remaining case, namely the
one involving the rectangle $Q=\Sigma^*$ in the
next chapter, there are no spanning edges at all.
\newline
\newline
{\bf The Matching Property:\/}
Let $v$ be an interior vertex of some $P_i$.
Let $P_1,...,P_k$ be the (re-indexed) union of all
the polygons $P_i$ such that $v \in \partial P_i$.
When $v$ is a vertex of $P_i$, we
let $e_i^+$ be the edge of $P_i$ such that $v$ is
the leading vertex, and we let $e_i^-$ be the
edge of $P_i$ such that $v$ is the trailing vertex.
When $p$ lies in the interior of an edge $e_i$ of
$P_i$, we let $e_i^+$ and $e_i^-$ be the two
oriented copies of $e_j$, labelled so that
the indices of the endpoints of $e_i^+$ are
cylically increasing.
We think of these $2k$ edges as vectors, all oriented
away from $v$.
For each $i \in \{1,...,k\}$ we check computationally
that there is exactly one $j \in \{1,...,k\}$ such that the two vectors
$e_i^-$ and $e_j^+$ are parallel. Call
this {\it the matching property\/}. 

\begin{lemma}
\label{coverX}
Suppose $P_1,...,P_k$ is a collection of 
positively convex polygons having pairwise disjoint
interiors.  Suppose also that each $P_k$
intersects $Q$ and has no spanning edges
with respect to $Q$. Finally, suppose that
the matching property holds for each vertex of $P_k$
that lies in the interior of $Q$.
Then $Q \subset \bigcup P_k$.
\end{lemma}

\startproof
Let $p$ be a vertex of $P_k$ that lies in the
interior of $Q$.
The matching property, the positive convexity,
and the matching property combine to show
that the link of $v$ is a circle. Hence,
$\bigcup P_k$ contains a neighborhood of $p$.
Let $Q'=\bigcup P_k$. If $Q$ is not a subset of $Q'$, then
we can find some point $q \in {\rm int\/}(Q)$ contained in
an edge $e$ of some $P_k$ such that
$Q'$ does not contain a neighborhood of $q$.
Since $e$ is not a spanning edge, we
can follow along $e$ in one direction until
we encounter a vertex of some $P_k$ that
lies in the interior of $Q$. Note that $Q'$ cannot
contain a neighborhood of any point on
$\overline{qp}$ because 
$e$ can only intesect other edges at
vertices. (This comes from the disjoint
interiors hypothesis.)
Now we can say that
$Q'$ does not contain a neighborhood
of $p$, which is a contradiction.
\endproof

Now we turn to the case of polyhedra.
Let $Q$ be a golden polyhedron and let
$P_1,...,P_n \subset Q$ be a finite collection of
strictly convex golden polyhedra having
pairwise disjoint interiors.  Suppose also that
every polyhedron in sight satisfies
Equation \ref{fibered}.   This means
that slicing the polyhedra by golden
horizontal planes results in golden polygons.
In this section, we explain the test we use to show that
$Q$ is partitioned by $\{P_k\}$.
\newline
\newline
{\bf The Slicing Method\/}
Let $P_h$ denote the intersection of $P$
with the horizontal plane $z=h$. As we
just mentioned, $P_h$ is a golden polygon for all
$h \in \Z[\phi]$.  
Let $S$ be the union of all the third coordinates
of all the vertices of our polyhedra, including $Q$.

We choose a subset $S^* \subset \Z[\phi]$ that
is {\it interlaced\/} with $S$. This means
that there is one point of $S^*$ between
each pair of consecutive points of $S$.
For each $h \in S^*$, we use the $2$ dimensional
covering test to show that
$Q_h$ is partitioned by $\{(P_k)_h\}$.
  
To see why this test suffices, consider a
vertex $v$ in the polygon partition at
height $h$.
We would like to see that the link of
$v$ remains a circle as we continuously
vary the height of the horizontal slice.
Suppose some polygons separate as we increase
the height of the slice.
Then these same polygons
overlap as we decrease the height from
$h$.  This cannot happen because of
the disjoint interior condition.

As we continuously vary the horizontal plane
away from $h$, every slice we encounter
satisfies the criterion for our covering test,
namely every interior vertex has a circle link.
So, for each open interval $I \in \R-S$, and
any $h \in I$, the polygons
$\{(P_k)_h\}$ partition $Q_k$. But this is
clearly sufficient to show that
$\{P_k\}$ partitions $Q$.

\subsection{Polyhedron Exchange Dynamics}
\label{pd}

Here we explain how we compute various
dynamical quantities related to the
polyhedron exchange dynamics.
\newline
\newline
{\bf Computing the Orbit\/}
Let $\widehat \Sigma=(\R/2\Z)$ as usual.
Let $P_1,...,P_{64}$ be the polyhedra
in the partition of $\widehat \Sigma$
corresponding to the dynamics.
Let $V_1,...,V_{64}$ be the corresponding
vectors associated to each of
these regions, so that $\widehat \Psi(p)=p+V_k$
for $p \in P_k$.

Here we describe one step of the polyhedron
exchange dynamics.
Given a point $p \in \widehat \Sigma$,
we use the interior point detection algorithm
to find $k$ such that $p \in P_k$.  We
then replace $p$ by $p+V_k$ and (if
necessary) subtract even integers from
the coordinates of $p+V_k$ until the result
$p'$ lies
in the chosen fundamental domain for $\widehat \Sigma$.
The map $p \to p'$ is one step of the
polyhedron exchange dynamics.  

In practice, we use floating point
arithmetic to produce the {\it canditate\/}
itinerary for a given point or region.  By this,
we mean the sequence of regions entered
by the orbit of that point or region.  Once we have the
candidate itinerary, we switch to
golden arithmetic to verify that the itinerary
really is the correct one, as explained below.

There is one fine point of our computer
calculations which seem worth mentioning.
The map $\Theta$ from the Compactification
Theorem involves division by $2$.  Since we
want to work entirely in $\Z[\phi]$ we
scale all vectors in $\R^3$ by a factor
of $2$ when we perform the golden
arithmetic calculations.  We mention this
for the benefit for a reader who wants
to survey our computer code.
\newline
\newline
{\bf Domain Verification:\/}
Let $P_1,...,P_{64}$
and $V_1,...,V_{64}$ be as above.
Suppose that $W$ is a vector and
$I=\{i_1,...,i_k\}$ is a sequence,
with each $i_j$ in the set 
$\{1,...,64\}$.  We say that the
pair $(W,I)$ is {\it feasible\/} if
\begin{itemize}
\item $W$ lies in the closure of $P_{i_1}$.
\item $W+V_{i_1}$ lies in the closure
of $P_{i_2}$.
\item $W+V_{i_1}+V_{i_2}$ lies in the closure
of $P_{i_3}$.
\item And so on.
\end{itemize}
Suppose that $Q$ is a polyhedron
with vertices $W_1,...,W_k$.  If
$(W_i,I)$ is feasible for all
$i=1,...,k$, then the first
$N$ iterates of $\widehat \Psi$
are defined on all interior points
of $Q$ and $I$ is the itinerary for
such points.  This works for the
following reason.  Since $Q$ lies
in the closure of $P_{i_1}$, each
interior point lies in the interior
of $P_{i_1}$.  Since $Q+V_{i_1}$
lies in the closure of $P_{i_2}$,
every interior point of $Q+V_{i_1}$
lies in the interior of $P_{i_2}$.
And so on.
\newline
\newline
{\bf Maximality Verification:\/}
Suppose we are given a polyhedron
$Q$ and some integer $N$.  Suppose
that we have already verified that
the first $N$ iterates are completely
defined on the interior of $Q$.
We would like to check that the
first $N$ iterates are nowhere defined
on $\partial Q$.  

We perform the following test
for each face $F$ of $Q$.  Let
$W$ be some point in the interior 
of $Q$, as described above.
We consider the orbit
$$W_1=W; \hskip 30 pt
W_k=W_{k-1}+V_{i_k}.$$
and check that there is some
$k<N$ such that $W_k \in \partial P_k$.

\subsection{Tile Creation}
\label{tc}

Suppose that $p \in \widehat \Sigma$ is some
point and $N$ is some integer. Here we explain
how we produce the maximal convex polyhedron $Q$
on which the first $N$ iterates of $\widehat \Psi$
are defined. From the viewpoint of making a
rigorous verification, it doesn't matter
how we produce $Q$.  Once we have produced
$Q$ {\it somehow\/}, we use the
domain verification algorithm to prove we are correct.
For what it is worth, we explain without proof how
matter how we produce $Q$.

We create $Q$ in $4$ steps.
\begin{enumerate}
\item We produce the first $N$ iterates
of $p$ using floating point calculations, as above.
\item We examine how each point of the
orbit sits inside polygon of $\widehat \Sigma$
that contains it.  In the end, we produce
a list of numbers, called the {\it tile data\/},
which carries the relevant information.
\item We construct a polyhedron $Q'$ from the initial
point $p$ and from the tile data.
\item We replace the floating point coordinates
of $Q'$ by the best golden real approximations.
The result is $Q$.
\end{enumerate}

We emphasize that there is no guarantee that our
construction will produce the correct answer.
It could fail at any stage.  However, in practice,
it works just fine, as later verified by golden
arithmetic calculations.

The second and third steps of the sketch above
need more explanation.  
\newline
\newline
{\bf Pseudo Distances\/}
There is a
list of $18$ planes $\Pi_1,...,\Pi_{18}$
such that each face of each $P_k$ is parallel to one
of these planes.  We have precomputed a certain
list $N_1,...,N_{18}$ of golden vectors such
that $\Pi_k$ is perpendicular to $N_k$ for all $k$.
Given a polygonal face $F$ of one of the $64$ polygons
$P_m$, and a point $q \in \R^3$, we define the
{\it pseudo distance\/} from $q$ to $F$ to be
the quantity
\begin{equation}
\label{pseudo}
\delta(q,F)=N_k \cdot (q-V).
\end{equation}
Here $k$ is such that $F$ is parallel to $\Pi_k$,
and $V$ is a vertex of $F$.  The
quantity $\delta(q,F)$ is independent of the
choice of $V$.  We choose the vectors $\{N_k\}$ so that
$p$ lies in the
interior of $P_m$ iff $\delta(q,F)>0$ for all
$F$ of $P_m$. This is another way of expressing
our interior point detection test.
\newline
\newline
{\bf Tile Data:\/}
Let $p_1,...,p_N$ be the forward iterates of $p$. 
There is a
list of $18$ planes $\Pi_1,...,\Pi_{18}$
such that each face of each $P_k$ is parallel to one
of these planes. 
Initially, let $D$ be the length $18$ list $\{100,...,100\}$.
(The number $100$ here is just some convenient large number.)
At the $k$th step of our process, we
compute the pseudo-distance from $p_k$ to each
of the faces in $V_{i_k}$.  If any of
the numbers we get is smaller than the
corresponding number in the list $D$,
we replace this number of $D$ by the
new, smaller one.  When we have done all
$N$ steps of the process, we have a final
list $D$.  We call $D$ the
{\it tile data\/}.
\newline
\newline
{\bf Tile Construction:\/}
Let $D$ be the tile data associated to the
pair $(p,N)$.  
 Call an index $i$ {\it active\/}
if $D_i<100$.  If $i$ is active, it means
that some face of some $P_{i_k}$
is parallel to $\Pi_i$.  
For each active index, we let
$\Pi'_i$ be the plane parallel to $\Pi$
such that the pseudo-distance from
$p$ to $\Pi'_i$ is $D_i$.
Let $H'_i$ be the half-plane which contains
$p$ and has $\Pi'_i$ as boundary.
Then $Q'$ is the intersection of all
the $H_i'$.  

In practice, we compute the intersection
$\bigcup H_i$ as follows. For each triple
of indices $(i,j,k)$ we find the
point $\Pi'_i \cap \Pi'_j \cap \Pi'_k$.
We call such a point {\it extraneous\/}
if some $\Pi'_m$ separates it from
the origin.  Once we eliminate all
extraneous intersections, the remaining
ones are the vertices of $Q'$.

\subsection{Tile Filling}
\label{tf}

Our proof of the Renormalization Theorem involves
partitioning the sets $\widehat A$ and $\widehat B$
into atoms.   From the point of view of giving
a proof, it doesn't matter how we produce these
partitions.  We simply have to verify that they
work. However, it seems worthwhile describing
what we actually did.
\newline
\newline
{\bf Generating the Individual Atoms:\/}
We generate the $\widehat B$ atoms in the
following way. Given a point $p \in \widehat B$,
we compute the points
$\widehat \Psi^k(p)$ for $k=1,2,3...$ until
the orbit returns to $\widehat B$. Then
we use the tile creation algorithm discussed
in the previous section.

We do not actually create the $\widehat A$ atoms
by this method.  Once we have created all the
$\widehat B$ atoms, we use the action of
$\widehat R$, the renormalization map, to 
define the $\widehat A$ atoms.  During
the proof of the Renormalization Theorem we
then use the domain verification method
discussed in \S \ref{pd} to
verify rigorously that the polyhedra
we have defined really are $\widehat A$-atoms.
\newline
\newline
{\bf Probing Vectors:\/}
We found by trial and error a list of $50$ vectors 
having the following property. Let $P$ be a special polyhedron
and let $v$ be a vertex of $P$. Then at least one of
our $50$ vectors, when based at $v$, points into $P$.
Given that our polyhedron exchange map is defined by
translations, the same property will hold with respect
to any polyhedron -- e.g. a periodic tile or an atom -- 
produced by the dynamics.  We call these vectors
the probing vectors.
\newline
\newline
{\bf Probing a List of Atoms:\/}
Suppose we have some list of 
$\widehat B$-atoms contained in 
some branch $Q$ (i.e. one of the $4$ big
convex polyhedra in each layer) of $\widehat B$.
We choose some vertex $v$ of one of the
$\widehat B$-atoms, or of $Q$.
Let $P_1,...,P_k$ be all the atoms
on our list that have $v$ as a vertex.
Say that a {\it new probe\/} is a
probing vector $W$ that
does not point into and of the $P_j$.

Given a new probe $W$,
we choose a point $v'$ so that
$v'-v$ is a small multiple of $W$,
and then we find the $\widehat B$-atom
that contains $v'$ and add it to our
list.  
\newline
\newline
{\bf Generating the Complete List:\/}
For each branch, we start with the
empty list of $\widehat B$ atoms.
(The initial list of vertices is
just the list of vertices of the
branch.)  We then iterate the
probing algorithm, cycling through
all vertices repeatedly, until we can find
no new probes.  This gives us our
list of atoms filling the branch.
We do this for all branches, and
this gives us our complete list
of $\widehat B$-atoms.

\subsection{Variations on a Theme}
\label{variation}

In \S \ref{pd} and \S \ref{tc} we explained
some constructions that we perform relative
to our polyhedron exchange map.  In the
next chapter, we will do similar
things with respect to the square
outer billiards map, and also with
respect to an auxilliary map which
we call the pinwheel map.
We will not describe the corresponding
operations in detail, because they are
so similar, but we will give one example.

Recall that $\Psi$ is the first return map
to the strip $\Sigma$, as described
in \S \ref{returnmap}.
Suppose we want to verify that
a given convex polygon $Q$ is a maximal domain
on which the first return map $\Psi$ is
entirely defined, and that the return
takes $n$ steps.  We first use floating
point arithmetic to compute a
(candidate) length $n$ itinerary $I$ for the center of mass of
$Q$.  We then check that $(V,I)$ is feasible,
in the sense of \S \ref{pd},
for all vertices $V$ of $Q$.
Next, we check that $\psi^n(V) \subset \Sigma$
for each vertex $V$ of $\Sigma$. Finally,
assuming that $Q$ has $k$ edges,
we generate a list of $k$ points
$v_1,...,v_k$, such that $v_j$ lies in the
interior of the $j$th edge of $P$ and
for each $j$ there is some
$m=m_j \in \{1,...,n-1\}$ such
that $\psi^m(v_j)$ lies in the boundary
of the relevant domain of definition
for $\psi$.

\newpage

\section{The Pinwheel Lemma}
\label{PINWHEEL}

\subsection{The Pinwheel Strips}
\label{pinstrip}

We are interested in the first return map
$\Psi: \Sigma_{\pm} \to \Sigma_{\pm}$ discussed
in \S \ref{returnmap}.  Our Pinwheel Lemma
below equates the map $\Psi$ with a new
map $\Pi$ that is simpler to manage.
As a step in the Compactification Theorem,
we prove a result which we call the Pinwheel
Lemma.  The Pinwheel Lemma equates
$\Psi$ with a map $\Pi$ which is easier
to analyze and ``compactify''.

We first explain how to associate $4$ pairs
$(S_i,V_i)$ for $i=0,1,2,3$ to the kite $K$.
Here $S_i$ is an infinite strip and $V_i$ is a
vector that points from one component of
$\partial S_i$ to the other.

\begin{center}
\resizebox{!}{4.8in}{\includegraphics{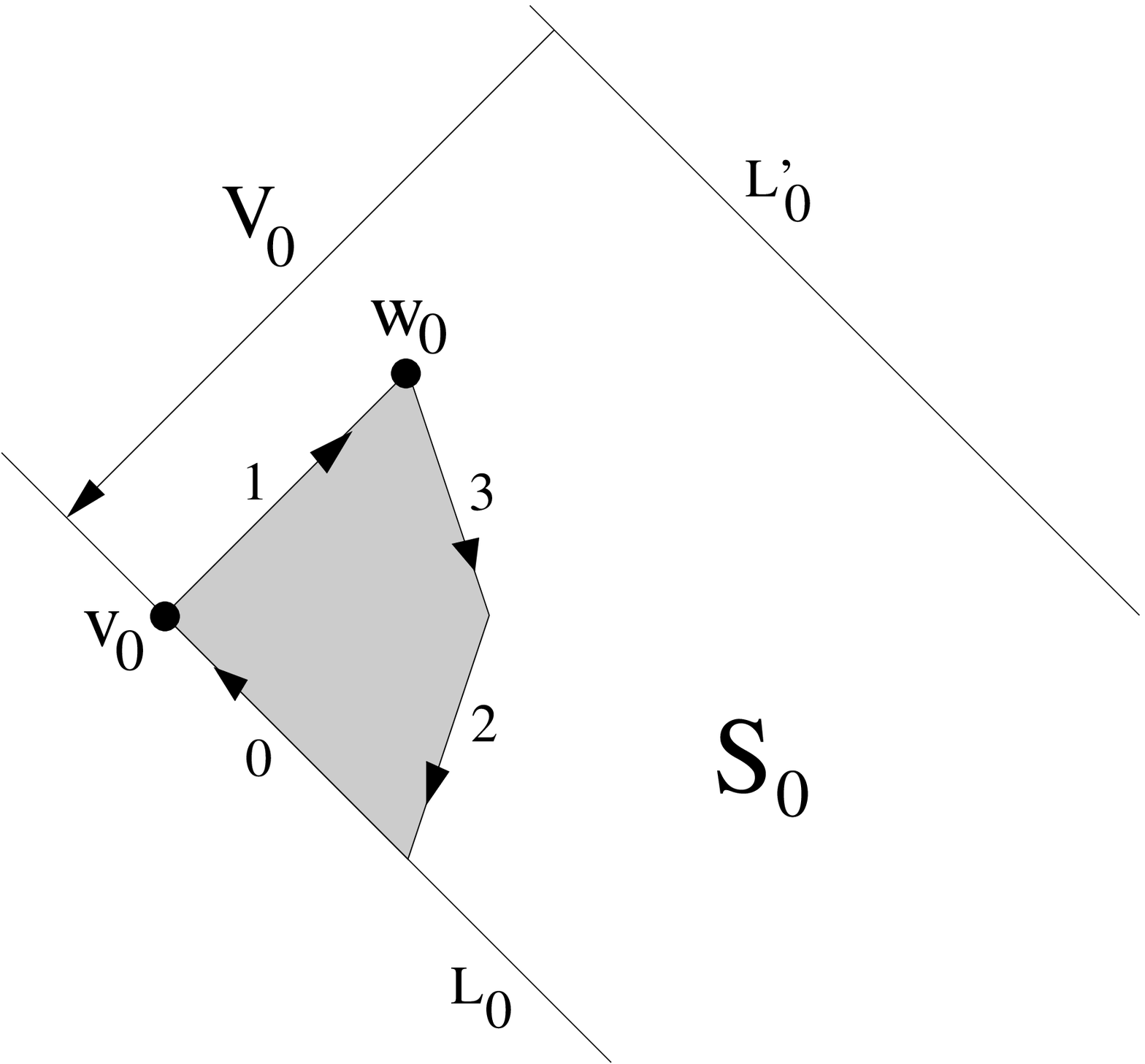}}
\newline
{\bf Figure 8.1:\/} The strip associated to $e_0$.
\end{center}  

We label the edges of $K$ as in Figure 8.1 
We call these edges $e_0,e_1,e_2,e_3$.  We orient
these edges so that they go clockwise around $K$.
  The labels
are such that the lines extending the edges
hit the ``circle at infinity'' (in the projective plane
that compactifies $\R^2$) in counterclockwise
cyclic order.

To the edge $e_j$ we associate the strip $S_k$
that has the following properties.  One boundary component of
$L_k$ of $S_k$ extends $e_k$.  The other boundary component is
$L_K'$ is such that the vertex $w_k$ of $K$ farthest from $L_k$
lies halfway between $L_k$ and $L_k'$.  So, $K$ extends
exactly up to the centerline of $S_k$.     We
define
\begin{equation}
V_k= \pm 2(v_k-w_k).
\end{equation}
Here $v_k$ is the head vertex of $e_k$.
We choose the signs so that the vectors are as shown
in Figure 8.2 below.
 Figure 8.1 shows
$S_0$ and $V_0$.

Let $A=\phi^{-3}$.
Here are the formulas for these strip pairs.
\begin{itemize}
\item $V_0=(2,2)$ and $S_0$ is bounded by the lines $x+y+1 \in \{0,4\}$.
\item $V_1=(0,4)$ and $S_1$ is bounded by the lines $x-y-1 \in \{0,4\}$.
\item $V_2=(2,-2)$ and $S_2$ is bounded by the lines $x-Ay-A\in \{0,-4/\phi\}$.
\item $V_3=(-2-2A,0)$ and $S_3$ is bounded by the lines $x+Ay-A \in \{0,-4/\phi\}$.
\end{itemize}

Figure 8.2 shows a rough picture of the strips and vectors relative
to $K(1/4)$.  The picture is very close to what one sees
for $K(\phi^{-3})$.
The black kite in the middle is $K(1/4)$.
Figure 8.2 also shows the strip $\Sigma=\R \times [-2,2]$, as well as the
regions $\Sigma_+$ and $\Sigma_-$.   These regions are denoted
$(+)$ and $(-)$ respectively.

To make our later definitions cleaner, we define
recursively define 
\begin{equation}
V_{4+i}=-V_i; \hskip 30 pt
S_{4+i}=S_i.
\end{equation}
This gives us pairs $(S_i,V_i)$ for all $i=0,1,2,3,4...$
The strips repeat with period $4$ and the vectors
repeat with period $8$.  However, note that
$(S,V)$ and $(S,-V)$ define the same strip maps.
So, the strip maps repeat with period $4$.  Nonetheless,
it is useful to distinguish between $V_1$ and $V_5=-V_1$, etc.
This way, the vectors $V_1,...,V_8$ can be drawn so that
they naturally circulate counterclockwise around $K$,
as in Figure 8.2.

\begin{center}
\resizebox{!}{4.8in}{\includegraphics{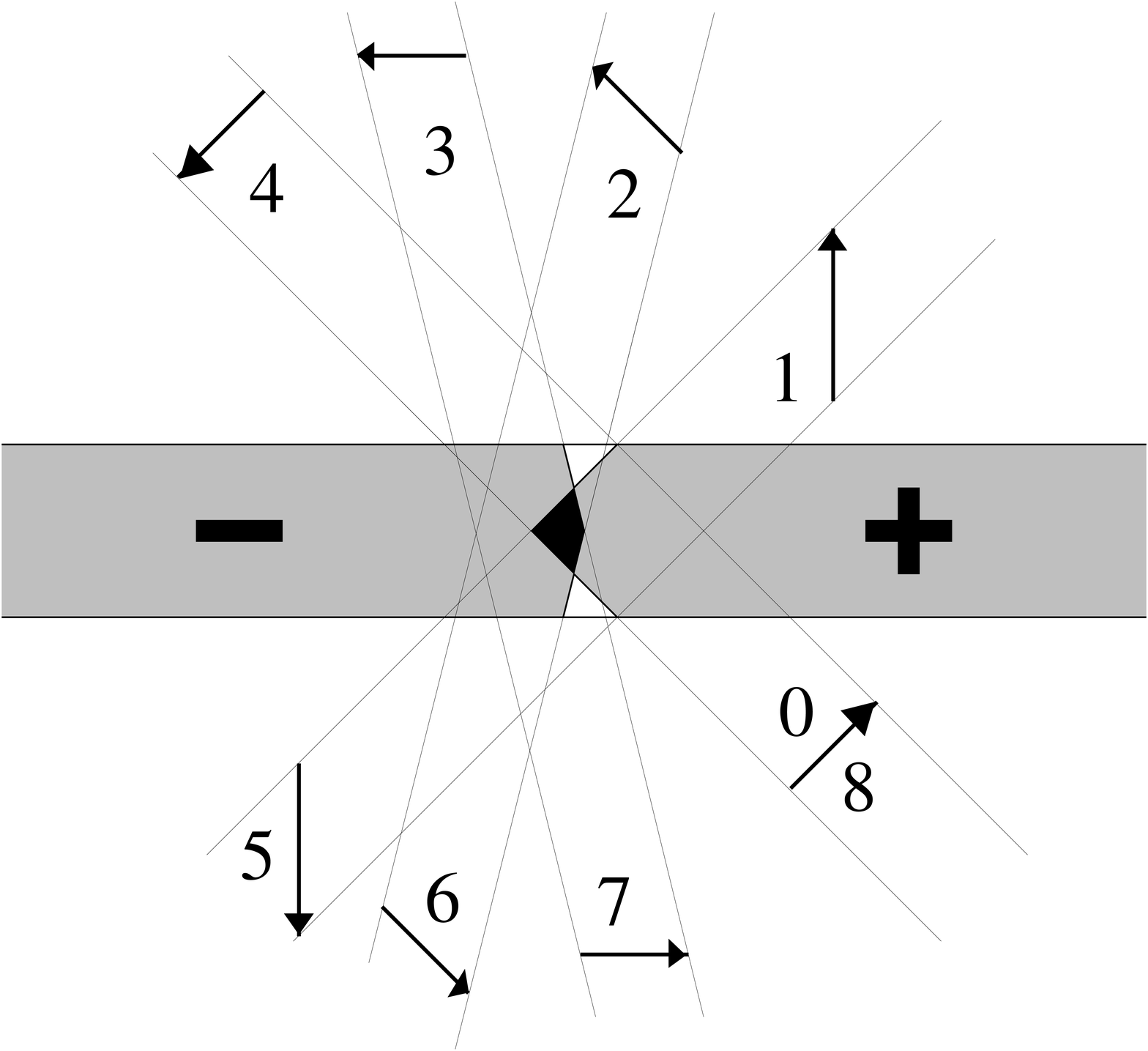}}
\newline
{\bf Figure 8.2:\/} The strips and vectors associated to $K$.
\end{center}  

\subsection{The Main Result}
\label{pinwheelmap}

Let $S$ be an infinite strip in the plane. 
We say that a {\it strip pair\/} is a pair $(S,V)$ where $S$ is a strip and
$V$ is a vector that points from one component of $\partial S$ to the other.
Given the strip pair, we have an associated {\it strip map\/}
$E: \R^2 \to S$, given by the formula
\begin{equation}
E(p)=p+nV \in S.
\end{equation}
Here $n=n(p)$ is chosen so that $E(p) \in S$.  The map
$E$ is well-defined on the complement of an infinite discrete
set of lines parallel to the boundary components of $V$.
Note that the vectors $V$ and $-V$ define the same map.
The pinwheel map is essentially the composition of strip maps.

Referring to the strip pairs $(S_i,V_i)$ defined in the previous
section, we recursively define $S_{4+i}=S_i$ and $V_{4+i}=V_i$.
This gives us pairs $(S_i,V_i)$ for all integers $i \geq 0$.
Of course, these pairs repeat with period $4$.
We let $E_i$ denote the strip map
associated to $(S_i,V_i)$.   This gives strip maps
$E_1$.  We let $\zeta$ be the strip map defined
relative to the pair $(\Sigma,(0,4))$.    We define the pinwheel map
\begin{equation}
\Pi = \zeta \circ E_8 \circ \ldots \circ E_1: \Sigma \to \Sigma.
\end{equation}
Recall that $\Psi$ is the first return map on
$\Sigma_+$ and likewise the first return map
on $\Sigma_i$.  

We proved versions of the Pinwheel Lemma in
[{\bf S1\/}], [{\bf S2\/}], and [{\bf S3\/}],
but unfortunately none of these results gives
us the exact statement we need.

\begin{theorem}[Pinwheel]
A point of $\Sigma_+ \cup \Sigma_-$ has a well-defined
$\Pi$-orbit if and only if it has a well-defined
$\Psi$-orbit, and $\Pi=\Psi$ on all points which have
well defined orbits.
\end{theorem}

Let $\Sigma_{\pm}^*$ denote those points in
$\Sigma_{\pm}$ of the form $(x,y)$ with
$|x| \leq 20$.  We see by direct inspection that
the Pinwheel Lemma holds for all points in 
$\Sigma_{\pm}-\Sigma_*$.  Both maps just circulate
the point around the kite, as shown in Figure 8.2.
See any of the papers
cited above for more details about this.

\begin{center}
\resizebox{!}{2.3in}{\includegraphics{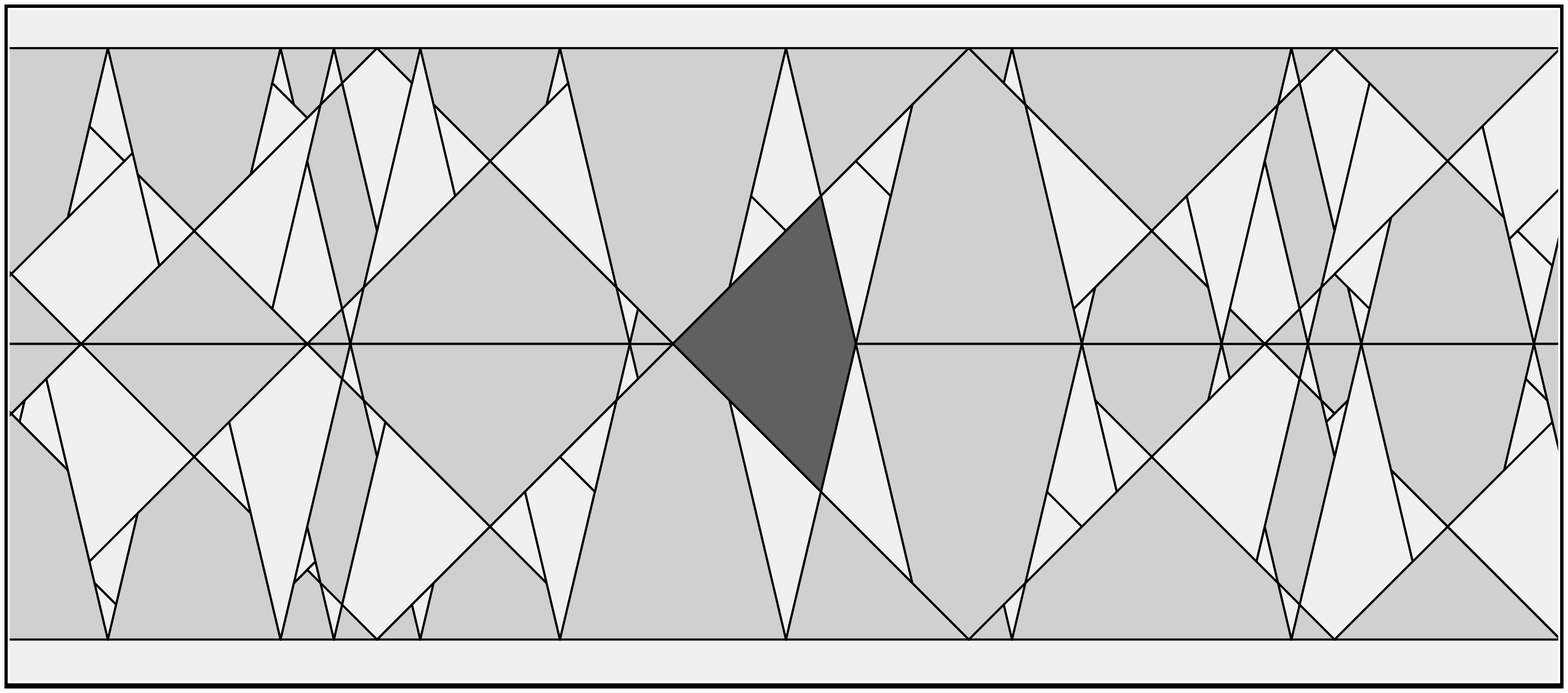}}
\newline
{\bf Figure 8.3:\/} Part of the covering
\end{center}

To deal with $\Sigma_{\pm}^*$, we make a direct
calculation.
We produce a list of $572$ golden convex polygons
$P_0,...,P_{571}$ whose union 
covers $\Sigma_{\pm}^*$.
Figure 8.3 shows some of these polygons.  The central dark
polygon is the Penrose Kite. The two triangles which
share the top and bottom vertex with $K$ belong to
$\Sigma-\Sigma_{\pm}$, and our dynamical statements
do not apply to them. We do not consider them to be
on our list.

We use the method described in \S \ref{variation}
to establish the following result for each 
polygon $P$ on our list.
\begin{itemize}
\item $\Psi$ is entirely defined on the interior of $P$.
\item $\Psi^2$ is nowhere defined on any non-horizontal edge of $P$.
\item $\Psi$ is the identity on the interior
of any horizontal edge of $P$.
\end{itemize}
As discussed in \S \ref{covertest}, this result implies
that our polyhedra have pairwise disjoint interiors.  We then
use the method in \S \ref{covertest} to show that
the rectangle $[-20,20] \times [-2,2]$ is covered
by the union of our polygons and the Penrose kite and
the two special triangles. This implies that
our union of polygons covers $\Sigma^*_+ \cup \Sigma^*_-$.

Next, we computationally establish exactly the same results for
$\Pi$.   Finally, we check that $\Pi=\Psi$ for one
interior point of each polygon.  Since $\Pi$ and
$\Psi$ are entirely defined on each polygon,
we see that $\Pi=\Psi$ on the interiors of
each of our polygons, and also on the interiors
of the horizontal edges.  None of the remaining
points have well-defined orbits for either map.

\subsection{Discussion}
\label{generate}

For the purposes of making a rigorous proof, it
doesn't really matter how we generate our polygons.
However, it seems worth saying
a word about how we do it.
First of all, we use the method described in
\S \ref{tc}, with $\Psi$ in place of
$\widehat \Psi$, to generate the individual tiles.

To generate the collection of polygons,
we let $N$ be the set of points of the form
\begin{equation}
\bigg(\frac{i}{100}+\epsilon,\frac{j}{100}+\epsilon\bigg);
\hskip 30 pt i=-2500,...,2500; \hskip 30 pt j=1,...,199.
\end{equation}
The small ``fudge factor'' $\epsilon=10^{-5}$ is present to
guarantee that we choose points on which the map
$\Pi$ is defined.
For each $p \in N$ we
compute $P(p)$, and then we weed out
redundancies and remove the ones that are
disjoint from $\Sigma^*_{\pm}$. 

\newpage
\section{The Compactification Theorem}

\subsection{Discussion}

Recall that $\Sigma=\R \times (-2,2]$ and
$\Pi: \Sigma \to \Sigma$ is the pinwheel map.
In light of the Pinwheel Lemma,
it suffices to prove the Compactification Theorem for
$\Pi$ in place of $\Psi$. 
We first compare the Compactification Theorem with the
corresponding results in [{\bf S1\/}] and [{\bf S2\/}].

Let $\Sigma_1 \subset \Sigma$ denote the union
of two horizontal lines $\R \times \{\pm 1\}$.
Let $\widehat \Sigma_1$ denote the 
$2$-torus slice
of the $3$-torus $\widehat \Sigma$, obtained
by interecting $\widehat \Sigma$ with the plane
of height $1$.  See Figure 4.6.
The map $\Pi$ carries $\Sigma_1$ to itself.

Our result [{\bf S1\/}, Arithmetic Graph Lemma] is equivalent to the
restricted version of the Compactification Theorem, when
$\Sigma_1$ replaces $\Sigma$, and
$\widehat \Sigma_1$ replaces $\widehat \Sigma$.
However, the result in [{\bf S1\/}] has a
somewhat different emphasis.  There we were concerned with how 
$\widehat \Sigma_1$ controls the structure of the
corresponding arithmetic graphs.  
In [{\bf S2\/}, Master Picture Theorem] we proved a
result that is 
equivalent to [{\bf S1\/}, Arithmetic Graph Lemma] for all kite
parameters.  The result in [{\bf S2\/}] is stated in 
the same general terms as the Compactification Theorem.

The proof in [{\bf S2\/}] is different than the proof
in [{\bf S1\/}]. It is more conceptual and also much
easier to generalize. In spite of this, we
will prove the Compactification Theorem following the
ideas in [{\bf S1\/}].  We do this partly because
we think of this paper as a sequel to [{\bf S1\/}], and
also because we would like to take the opportunity
to revisit that proof and give a cleaner exposition.
The proof we give is specially adapted to the
Penrose kite, and would not easily generalize.

\subsection{The Covering Property}
\label{cover}

The first order of business is to show that the
$64$ golden 
polyhedra defining the map $\widehat \Psi$ really
do partition the torus $\widehat \Sigma$
into strictly convex polyhedra. (This would be
the first step in any approach we took to
proving the main result.)  We use the
tests described in \S \ref{basictest} to show
that the polyhedra involved are strictly
convex and have pairwise disjoint interiors.
We also verify that they satisfy Equation \ref{fibered}.
Finally, we use the slicing method described
in \S \ref{covertest} to show that our
polyhedra really do partition $\widehat \Sigma$.

\subsection{The Proof Modulo Computations}

Let $\widehat F$ be the fundamental domain
for $\widehat \Sigma$ that we defined
in \S \ref{partitionpix}.
 Recall that
\begin{equation}
\Pi=\zeta \circ E_8 \circ \ldots E_1
\end{equation}
where $E_1,...,E_8$ and $\zeta$ are the strip maps defined in
\S \ref{pinwheelmap}.  In the next section, we
define the following objects.
\begin{itemize}
\item An affine map $\widetilde \Theta: \R^2 \to \R^4$. See Eq. \ref{tildetheta}
\item A locally affine map $\widetilde \zeta: \R^4 \to \widehat \Sigma$.
See Eq. \ref{tildezeta}.
\item An affine embedding $I_t: \widehat F \to \R^4$, for any $t \in \Z[\phi]$.
  See Eq. \ref{iota1}.
\item A set $\widetilde L_k \subset \R^4$ of parallel hyperplanes, for $k=1,...,8$.
See Eqs. \ref{hyp1}--\ref{hyp2}
\item Piecewise affine maps $\widetilde E_1,...,\widetilde E_8: \R^4 \to \R^4$.
See Eqs. \ref{aff0}--\ref{aff1}.
\end{itemize}

Associated to $E_j$ is an infinite family
$L_j$ of parallel lines, such that
$E_j$ is defined exactly in the complement of
the lines. At the same time,
$\widetilde E_k$ is defined precisely in the 
complement of $\widetilde L_k$.  We write
$V \sim V$ is $V-W \in (4\Z)^4$.

\begin{lemma}
\label{tech1}
The following is true.
\begin{enumerate}
\item $\widetilde \zeta \circ \widetilde \Theta = \Theta \circ \zeta$.
\item $\widetilde \zeta \circ T = \widetilde \zeta$.
\item $\widetilde \zeta \circ I_t={\rm Identity\/}$ for all $t$.
\item $\widetilde E_1 \circ I_t$ is independent of $t$.
\item $T(\widetilde L_k)=\widetilde L_k$.
\item $\widetilde \Theta^{-1}(\widetilde L_k)=L_k$.
\item $\widetilde E_k \circ T \sim \widetilde E_k$.
\item $\widetilde \Theta \circ E_k \sim \widetilde E_k \circ \widetilde \Theta$.
\end{enumerate}
Here $T$ is and element of $(4\Z)^4$, considered
as a map of $\R^4$, and $k=1,...,8$.
\end{lemma}

For the next result, we interpret $\Theta(p)$ as a
point of $\widehat F$.  When the point lies in
the interior of $\widehat F$, it is the unique
representative of $\Theta(p)$. Otherwise, we make
some choice amongst the possibilities.

\begin{lemma}
\label{tech2}
There is a $4$-element set $\Omega \subset \Z[\phi]$ with the following
property. For any $p \in \Sigma$ there is some $t \in \Omega$ such that
that $I_t \circ \Theta(p) \sim \widetilde \Theta(p)$. 
\end{lemma}

Define the following piecewise (locally) affine map.

\begin{equation}
\label{major}
\widetilde \Pi = \widetilde \zeta \circ \widetilde E_8 \circ \ldots
\widetilde E_1: \R^4 \to \widehat \Sigma.
\end{equation}

\begin{lemma}
\label{tech3}
Let $P$ be a special polygon.
$\widetilde \Pi \circ I_t$ is everywhere
defined and locally affine on the interior
of $P$, but nowhere defined on the boundary of $P$.
\end{lemma}

We prove the above $3$ results in the
sections following this one.

\begin{lemma}
\label{defined}
Suppose that $p$ lies in the interior of
$\Sigma$ and $\Theta(p)$ lies
in the interior of a special polyhedron. Then
$\Pi$ is well-defined on $p$.
\end{lemma}

\startproof
Let $t \in \Z[\phi]$ be the value
guaranteed by 
Lemma \ref{tech2}, and let $I=I_t$.
We define
\begin{equation}
\widetilde q_0=I \circ \Theta(p); \hskip 30 pt
\widetilde p_0=\widetilde \Theta(p).
\end{equation}
By Lemma \ref{tech2} we have
$\widetilde p_0 \sim \widetilde q_0$.
We inductively define
$\widetilde q_k=\widetilde E_k(\widetilde q_{k-1})$
for $k=1,...,8$.
Applying Lemma \ref{tech3} to the open special
polygon that contains $\Theta(p)$, we see that the
points $\widetilde q_1,...,\widetilde q_8$ are
well defined, meaning that
$\widetilde q_{k-1} \not \in \widetilde L_k$.

Since $\widetilde L_1$ is $(4\Z^4)$-invariant
and $\widetilde p_0 \sim \widetilde q_0$, we
have $\widetilde p_0 \not \in \widetilde L_1$.
Hence, we may define
$\widetilde p_1 = \widetilde E_1(\widetilde p_0)$.
By Statement 7 of Lemma \ref{tech1}, we have
$\widetilde p_1 \sim \widetilde q_1$.
Shifting the indices and repeating the same
argument $7$ more times, we find that
we can inductively define points
$\widetilde p_k=\widetilde E_k(\widetilde p_{k-1})$
for $k=1,...,8$.  That is, $\widetilde \Pi$ is
defined on $p$.  Put another way,
$\widetilde p_{k-1} \not \in \widetilde L_k$
for $k=1,...,8$.

Let $p_0=p$.  We have
$\widetilde p_0 \not \in \widetilde L_1$.  But then
$\widetilde \Theta(p_0) \not \in \widetilde L_1$
because  $\widetilde L_1$ is $(4\Z)^4$-invariant and
$\widetilde \Theta(p_0) \sim \widetilde p_0$.
By Lemma \ref{tech1}, we have
$p_0 \not \in L_1$.
Hence $E_1$ is defined on $p_0$.  This allows us to define
Let $p_1=E_1(p_0)$.
Shifting the indices and repeating the same
argument $7$ times, we can inductively
define the points $p_k=E_k(p_{k-1})$ for
$k=1,....,8$.  
In short, the composition
$E_8 \circ ... \circ E_1$ is defined on $p$.

Finally, the set of lines $\R \times (2\Z)$ is
invariant under all our strip maps.  Since
$p$ does not lie on these lines, neither
does $p_8$.  Hence $\zeta$ is defined on $p_8$.
Hence, $\Pi$ is defined on $p$ and
$\Pi(p)=\zeta(p_8)$.
\endproof

\begin{lemma}
\label{diagram}
$\widetilde \Pi \circ \widetilde \Theta=\Pi \circ \Theta$
whenever both sides are defined.
\end{lemma}

\startproof
We write $V_1 \sim V_2$ when $V_1-V_2 \in (4\Z)^4$.
Let $p \in \Sigma$.  Repeatedly using Statements 7 and 8 of Lemma
\ref{tech1}, we have
\begin{equation}
\label{ind}
\widetilde E_8 \circ \ldots \circ \widetilde E_1 \circ \widetilde \Theta(p) \sim
\widetilde \Theta \circ E_8 \circ \ldots \circ E_1(p).
\end{equation}
By Statement 2 of Lemma \ref{tech1}, we get equality above when we apply
$\widetilde \zeta$ to both sides of Equation \ref{ind}. That is,
\begin{equation}
\widetilde \Pi \circ \widetilde \Theta(p)=
\widetilde \zeta \circ \widetilde \Theta \circ E_8 \circ \ldots E_1(p)=^*
\Theta \circ \zeta \circ E_8 \circ \ldots \circ E_1(p)=
\Theta \circ \Pi(p).
\end{equation}
The starred equality is Statement 1 of Lemma \ref{tech1}.
\endproof

\begin{lemma}
\label{conjugacy}
Suppose $p_1, p_2$ lie in the interior of $\Sigma$
and $\Theta(p_1), \Theta(p_2)$ 
lie in the interior of the same
special polyhedron.  Then 
$\Theta \circ \Pi(p_j)-\Theta(p_j)$
is the same for $j=1$ and $j=2$.
\end{lemma}

\startproof
Let $P$ be the special polygon whose interior contains
$\Theta(p_1)$ and $\Theta(p_2)$.  
\begin{equation}
\label{indep1}
\widetilde \Pi \circ I_{t_1} \circ \Theta(p_1)-\Theta(p_1)=
\widetilde \Pi \circ I_{t_1} \circ \Theta(p_2)-\Theta(p_2)=
\widetilde \Pi \circ I_{t_2} \circ \Theta(p_2)-\Theta(p_2).
\end{equation}
Here $t_j=t(p_j)$, in the sense of Lemma \ref{tech2}, and
the subtraction takes place in the abelian group
$(\R/2\Z)^3$. The first equality comes from the fact that
$\widetilde \Pi \circ I_{t_1}$ is locally affine and
entirely defined on $P$. 
The second equality has the following explanation.
By Statement 4 of Lemma \ref{tech1}, the map
$\widetilde E_1 \circ I_t$ is independent of $t$.
Hence $\widetilde \Pi \circ I_t$ is independent of $t$.

Noting that $\zeta(p_j)=p_j$, we use the identities in
Lemmas \ref{tech1} and \ref{diagram} to conclude that
$$
\Theta \circ \Pi(p_j)-\Theta(p_j)=
\Theta \circ \Pi(p_j)-\Theta\circ \zeta(p_j)=
\widetilde \Pi \circ \widetilde \Theta(p_j)-\widetilde \zeta \circ \widetilde \Theta(p_j)=
$$
\begin{equation}
\label{indep2}
\widetilde \Pi \circ I_{t_j} \circ \Theta(p_j)-\widetilde \zeta \circ I_{t_j} \circ \Theta(p_j)=
\widetilde \Pi \circ I_{t_j} \circ \Theta(p_j)-\Theta(p_j),
\end{equation}

Combining Equations \ref{indep1} and \ref{indep2}, we get the
equation we seek.
\endproof

Recall that $\widehat \Psi$ is a piecewise translation on
$\widehat \Sigma$.  For each of the $64$ polyhedra $P$
that partition $\widehat \Sigma$, we produce a point
$p \in \Sigma$ such that 
\begin{equation}
\Theta(p) \in {\rm interior\/}(P); \hskip 30 pt \Theta \circ \Pi(p)=\widehat \Psi(p).
\end{equation}
Indeed, this is how we defined $\widehat \Psi$ in the first place.
If $q \in \Sigma$ is any other point such that
$\Theta(q)$ lies in the interior of $P$ then, by Lemma
\ref{conjugacy} we have
$$\Theta \circ \Pi(q)-\Theta(q)=\Theta \circ \Pi(p)-\Theta(p)=
\widehat \Psi(p)-p = \widehat \Psi(q)-q.
$$
This shows that
$\Theta \circ \Pi = \widehat \Psi \circ \Theta$
for all points $p \in \Sigma$ such that $\Theta(p)$ lies in the interior of
a special polyhedron.

To finish the deduction of the Compactification Theorem, we just have
to understand what happens in exceptional cases.

\begin{lemma} Suppose that $p$ lies in the interior of $\Sigma$ and
$\Theta(p)$ does not lie in the interior of a special polygon.
Then $\Pi$ is not defined on $p$.
\end{lemma}

\startproof
We will revisit the proof of Lemma \ref{defined} and
use the notation set up in that proof.  We will suppose
that $\Pi$ is defined on $p$ and then derive a contradiction.
Since $\Pi$ is defined on $p$, the compositions
$E_1$, $E_2 \circ E_1$, etc. are defined on $p$.  Hence,
the points $p_1,...,p_8$ are well defined.  But then
the same argument as in the proof of Lemma \ref{defined}
(but done in reverse) shows that the points
$\widetilde p_1,...,\widetilde p_8$ are
all defined.  Hence, the maps $\widetilde E_1,...,\widetilde E_8$
are all defined on $\widetilde E_0$.  Finally, the map
$\widetilde \zeta$ is everywhere defined on $\R^4$.
Hence $\widetilde \Pi$ is defined on $\widetilde p_0$, a point
in the boundary of $\widetilde P_0$.  This contradicts the
second statement of Lemma \ref{tech3}.
\endproof

It remains only to deal with the points on $\partial \Sigma$.
The only issue is the action of the map $\zeta$, and
this only comes up at the end of our proof
of Lemma \ref{defined}. The problem is that $\zeta$
is undefined on the set
$\R \times 2\Z$. These lines are invariant under
the action of $E_1,...,E_8$.  Hence, if $p \in \partial \Sigma$,
then $p_8 \in \R \times 2\Z$ and
$\zeta$ is not defined on $p_8$.  However, recall that
$\Sigma=(-2,2] \times \R$.  The bottom
boundary is left off.  Given this convention, we define
\begin{equation}
\zeta(p_8)=(x_8,2) \in \partial \Sigma.
\end{equation}
Here $x_8$ is the first coordinate of $p_8$.  Once we
make this definition, our proofs above go through
without a hitch.

\subsection{Proof of Lemma \ref{tech1}}
\label{equations}
\label{techproof1}

Recall that
\begin{equation}
\Theta(x_1,x_2)=\bigg[\bigg(1,\frac{1}{2},0\bigg)+\bigg(\frac{x_1}{\phi},\frac{x_1-x_2}{2},x_2\bigg)\bigg]_2
\hskip 30 pt
\zeta(x_1,x_2)=(x_1,[x_2]_4).
\end{equation}
Here $[V]_2$ means that we take the image of $V$ in $(\R/2\Z)^3$, and
$[x_2]_4$ means representative
of $x_2$ in $\R/4\Z$ that lies in $(-2,2]$.

We define
\begin{equation}
\label{tildetheta}
\widetilde \Theta(x_1,x_2)=\Bigg(x_1+x_2,x_1-x_2,\frac{x_1+x_2}{\phi},\frac{x_1-x_2}{\phi}\Bigg).
\end{equation}
This is exactly the same map we used in [{\bf S1\/}].
Next, we define
\begin{equation}
\label{tildezeta}
\widetilde \zeta(x_1,x_2,x_3,x_4)=\bigg(1,\frac{1}{2},0\bigg)+
\frac{1}{2}\big(x_3+x_4,x_2,x_1-x_2\big)
\in (\R/2\Z)^3.
\end{equation}

For any $t \in \Z[\phi]$ we define
We define
$$
I_t(x_1,x_2,x_3)=
(2y_2+y_3,2y_2+y_3,y_1,y_1)+(y_3,-y_3,\frac{y_3}{\phi}-t,-\frac{y_3}{\phi}+t),
$$
\begin{equation}
\label{iota1}
(y_1,y_2,y_3)=(x_1,x_2,x_3)-\bigg(1,\frac{1}{2},0\bigg).
\end{equation}

Now we define the families $\widetilde L_k$ of parallel
hyperplanes.  We introduce the functions $g_k: \R^4 \to \R$ as follows.
\begin{equation}
\label{hyp1}
g_1=x_2+1; \hskip 15 pt
g_2=x_2+x_3+2+\phi; \hskip 15 pt
g_3=x_1+x_4+2+\phi; \hskip 15 pt
g_4=x_1+1.
\end{equation}
We define
\begin{equation}
\label{hyp2}
\widetilde L_{k+4}=\widetilde L_k=g_k^{-1}(4\Z).
\end{equation}

\noindent
{\bf Statements 1-4:\/}
\begin{enumerate}
\item A direct calculation shows
$\widetilde \zeta \circ \widetilde \Theta=\Theta$
and $\Theta=\Theta \circ \zeta$.
Hence $\widetilde \zeta \circ \widetilde \Theta=
\Theta \circ \zeta$.
\item It is immediate from the formula that $\widetilde \zeta \circ T=\widetilde \zeta$
for any $T \in (4\Z)^4$.
\item A direct calculation shows that
$\widetilde \zeta \circ \widetilde I_t$ is the identity.
\item It is immediate from the equations that
$\widetilde L_k$ is $(4\Z^4)$-invariant.
\end{enumerate}

For $x \in \R-4\Z$ let $\tau(x)$ denote the greatest element in
$4\Z$ that is less than $x$.  Define
\begin{equation}
\label{aff0}
\gamma_k(X) = \phi \times \bigg(\tau \circ g_k(X)\bigg) \in 4\phi \Z.
\end{equation}
The $(\times)$ in Equation \ref{aff0} means multiplication.
Here $X=(x_1,x_2,x_3,x_4)$ is a point of $\R^4-\widetilde L_k$.
Finally, we define

\begin{equation}
\label{aff1}
\matrix{
\widetilde E_1&=&
\left(\matrix{x_1 \cr x_2 \cr -x_2/\phi+x_3+x_4 \cr x_2/\phi}\right)&+ &
\left(\matrix{0 \cr 0 \cr \gamma_1 \cr -\gamma_1}\right) \cr \cr \cr
\widetilde E_2&=&\left(\matrix{x_1 \cr x_2 \cr x_3 \cr -x_1+x_2/\phi+x_3 \phi}\right)&+
&\left(\matrix{0\cr 0 \cr 0\cr -\gamma_2}\right) \cr \cr \cr
\widetilde E_3&=&
\left(\matrix{
x_1 \phi - x_2-x_3 +x_4\phi \cr
x_1/\phi + 0 x_2 - x_3 + x_4 \phi \cr
-x_1/\phi+x_2+2x_3 - x_4\phi \cr
-x_1/\phi + x_2 + x_3 - x_4/\phi}\right)&+&
\left(\matrix{-\gamma_3 \cr -\gamma_3 \cr \gamma_3 \cr \gamma_3}\right) \cr \cr \cr
\widetilde E_4 &=&
\left(\matrix{x_1\cr x_2 \cr x_1/\phi \cr x_4}\right) &+&
\left(\matrix{0 \cr 0 \cr -\gamma_4 \cr 0}\right)}
\end{equation}

Here we have simplified our notation slightly by writing $\widetilde E_1=\widetilde E_1(X)$
and $\gamma_1=\gamma_1(X)$, etc.   Each $\widetilde E_k$ is the sum of a linear
transformation and a piecewise constant vector-valued function.  The linear part
is defined on all of $\R^4$, but the vector-valued function changes valued
when one passes through a hyperplane of $\widehat L_k$.  We define
$\widetilde E_{4+k}=\widetilde E_k$.
\newline
\newline
{\bf Statement 5:\/}
The maps $I_0$ and $I_t$ agree
up to post-composition with the translation 
form $V \to V+(0,0,-t,t)$.   Inspecting the formula for
$\widetilde E_1$, we see that the only expression involving
the third and fourth coordinates is $x_3+x_4$.
Hence $\widetilde E_1 \circ \widetilde I_t$ is
independent of $t$.  Hence $\widetilde \Pi \circ I_t$
is independent of $t$.
\newline
\newline
{\bf Statement 6:\/}
We can deduce the formulas for the lines $L_k \subset \R^2$ by
looking at the formulas for the strip pairs given in
\S \ref{pinstrip}.
We have $L_{4+k}=L_k$, so we just need to consider
$L_k$ for $k=1,2,3,4$.  Let $A=\phi^{-3}$.
\begin{itemize}
\item $L_1$ consists of those lines satisfying $x_1-x_2-1 \in 4\Z$.
\item $L_2$ consists of those lines satisfying $x_1-Ax_2-A \in 4\phi^{-1}\Z$.
\item $L_3$ consists of those lines satisfying $x_1+Ax_2-A \in 4\phi^{-1}\Z$.
\item $L_4$ consists of those lines satisfying $x_1+x_2+1 \in 4\Z$.
\end{itemize}
The fact that $\widetilde \Theta^{-1}(\widetilde L_k)=L_k$ is obvious for $k=1,4$.
For $k=2$ we compute 
$$g_2 \circ \widetilde \Theta(x_1,x_2)=
\phi x_1-\phi^{-2} x_2 +4 -\phi^{-2}=\phi(x_1-Ax_2-A).$$
The result is obvious from here. The proof for
$k=3$ is similar.
\newline
\newline
{\bf Statement 7:\/}
Let $e_1,e_2,e_3,e_4$ be the standard basis vectors of $\R^4$.
let
\begin{equation}
u(j,k,X)=\widetilde E_j(X+4e_k)-\widetilde E_j(X).
\end{equation}
Statement $7$ is equivalent to the statement that
$u(j,k,X) \in (4\Z)^4$ for all $j,k$ and $X$.

We see directly that
$u(1,k,X) \in (4\Z)^4$ for $k=1,3,4$, and we compute
$$u(1,2,X)=(0,0,-4/\phi,4/\phi)+(0,0,4\phi,-4\phi)=(0,0,-4,4).$$
This proves the result for $\widetilde E_1$. 
The proof for $\widetilde E_4$ is similar.

We see directly that
$u(2,k,X) \in (4\Z)^4$ for $j=1,4$, and we compute
$$u(2,2,X)=(0,0,0,4/\phi) + (0,0,0,-4\phi) = (0,0,0,-4),$$
$$u(2,3,X)=(0,0,0,4\phi)+(0,0,0,-4\phi)=(0,0,0,0).$$
This proves the result for $\widetilde E_2$.  

For $\widetilde E_3$, we compute
$$u(3,1,X)=(4\phi,4/\phi,-4/\phi,-4\phi)+(-4\phi,-4\phi,4\phi,4\phi)=
(0,-4,4,4)$$
$$u(3,2,X)=(-4,0,4,4); \hskip 30 pt
u(3,3,X)=(-4,-4,8,4).$$
$$u(3,4,X)=(4\phi,4\phi,-4\phi,-4/\phi)+(-4\phi,-4\phi,4\phi,4\phi)=(0,0,4,4).$$
\newline
\newline
{\bf Statement 8:\/}
We want to prove that
$$\widetilde \Theta \circ E_k \sim \widetilde E_k \circ \widetilde \Theta,$$
where $\sim$ denotes equivalence mod $(4\Z)^4$.
Recall that the strip map $E_k$ is defined in terms of a
pair ($S_k,V_k)$, where $S_k$ is a strip and
$V_k$ is a vector. 
The two lines of $\partial S_k$ are consecutive lines of
$L_k$.  Hence, by Lemma \ref{tech1}, the set
$\widetilde S_k=\widetilde \Theta(S_k)$ lies between two
hyperplanes in $\widetilde L_k$.  Hence
$\widetilde E_k$ is an affine map on $\widetilde S_k$.
Indeed, $\widetilde E_k$ is linear on $\widetilde S_k$
because 

$$\gamma(0,0,0,0)=0; \hskip 30 pt
(0,0,0,0) =\widetilde \Theta(0,0) \in \widetilde \Theta(S_k)=\widetilde S_k$$
A direct calculation shows that
the linear part of $\widetilde E_k$ is the identity on the $2$-plane
$\widetilde \Theta(\R^2)$.  Hence
$\widetilde E_k$ is the identity on $\widetilde \Sigma$. 
In summary
\begin{equation}
\widetilde \Theta(p) = \widetilde E_k \circ \widetilde \Theta(p)
\hskip 30 pt \forall p \in {\rm interior\/}(\Sigma).
\end{equation}

Let $\widetilde V_k=\widetilde \Theta(V_k)$.
Given the piecewise affine nature of our maps,
the value of
$\widetilde E_k(X+\widetilde V_k)-\widetilde E_k(X)$
is independent of the choice of $X \in \R^4-\widetilde L_k$.
Compare the proof of Statement 7.  A direct calculation
shows that $\widetilde E_k(\widetilde V_k) \in (4\Z)^4$.
Therefore
$$
\widetilde E_k(X+\widetilde V_k) \sim \widetilde E_k(X);
\hskip 30 pt
\forall X \in \R^4-\widetilde L_k.
$$
Iterating this formula, we get
\begin{equation}
\widetilde E_k(X+\widetilde nV_k) \sim \widetilde E_k(X);
\hskip 30 pt
\forall X \in \R^4-\widetilde L_k.
\end{equation}
Here $n$ is an integer.

Every point $p' \in \R^2-L_k$ has the form $p+nV_k$ for
some $p$ in the interior of $S_k$.  Note that $E_k(p')=p$.
We have
$$\widetilde \Theta \circ  E_k(p')=
\widetilde \Theta(p)=\widetilde E_k \circ \widetilde \Theta(p) \sim
$$
$$
\widetilde E_k(\Theta(p)+n\widetilde V_k)=
\widetilde E_k \circ \widetilde \Theta(p+nV_k)=
\widetilde E_k \circ \widetilde \Theta(p').$$
This proves Statement 8.

\subsection{Proof of Lemma \ref{tech2}}
\label{techproof2}

Let $p=(x_1,x_2) \in \Sigma$. When we interpret $\Theta(p)$ as
a point of the fundamental domain $\widehat F$, we
are choosing integers $A_1,A_2,A_3$ such that
\begin{equation}
\label{rep}
\Theta(p)=\bigg(2A_1+1,2A_2+\frac{1}{2},2A_3\bigg)+
\bigg(\frac{x_1}{\phi},\frac{x_1-x_2}{2},x_2\bigg) \in \widehat F.
\end{equation}
Note that 
\begin{equation}
A_3 \in \{0,1\}
\end{equation}
 because $x_2 \in (-2,2]$.
 Using
Equation \ref{rep} and the definition of the map
$I_t$, we compute
\begin{equation}
I_t \circ \Theta(p)=
\widetilde \Theta(p)+\bigg(4A_2+4A_3,4A_2,2A_1+\frac{2A_3}{\phi}+t,2A_1-\frac{2A_3}{\phi}-t\bigg).
\end{equation}
The appropriate choice of $t$ from amongst the values
\begin{equation}
0; \hskip 30 pt 2; \hskip 30 pt
-2/\phi; \hskip 30 pt 2-2/\phi
\end{equation}
leads to $$I_t \circ \Theta(p)\sim\widetilde \Theta(p)$$
This proves Lemma \ref{tech2}.

\subsection{Proof of Lemma \ref{tech3}}

The arguments we give here are very similar to what
we discussed in \S \ref{pd}, but the setting is
different enough that we will go through the
details.

Given the $64$ special polyhedra $P_1,...,P_{64}$, we
define the $256$ polyhedra
\begin{equation}
P_{ij}=I_{t_i}(P_j) \subset \R^4; \hskip 30 pt
i \in \{1,2,3,4\} \hskip 30 pt
j \in \{1,...,64\}.
\end{equation}

\begin{lemma}
$\widetilde \Pi$ is defined on the interior of $P_{ij}$
for all $i=1,2,3,4$ and all $j=1,...,64$.
\end{lemma}

\startproof
Let $q$ be the center of mass of $P_{ij}$.
Using floating point arithmetic, we associate to
$q$ a certain $8$-tuple of integers, as follows.
Let $q_0=q$ and inductively define $q_k=\widetilde E_k(q_{k-1})$.
We define 
\begin{equation}
n_k=\phi^{-1} \gamma_k(p_{k-1}).
\end{equation}
This integer locates
the component of $\R^4-\widetilde L_k$ that contains
$q_{k-1}$.  The $8$-tuple of interest to us is $(n_1,...,n_8)$.
For the purposes of our proof, it doesn't matter how
we generated the sequence $(n_1,...,n_8)$.
This sequence just serves as a guide for our
rigorous calculations.

Guided by the sequence we have produced,
we define 
the map $\widetilde E_k^*$ to be the same map as
$\widetilde E_k$, except that $4\phi n_k$ replaces 
$\gamma_k$ in the formula.  We check, using golden
arithmetic, that the composition
$\widetilde E_8^* \circ \ldots \widetilde E_1^*$
is defined on all the vertices of $P_{ij}$.
For each such vertex $v_0$ we inductively define
$v_k=\widetilde E_k^*(v_{k-1})$.  We check that
\begin{equation}
g_k(v_k) \in [n_k,n_k+4].
\end{equation}
Here $g_k$ is as in Equation \ref{hyp1}.

The same argument we gave when we discussed
domain verification in \S \ref{pd} now says
that $\widetilde E_8 \circ \widetilde E_1$ is
defined on the interior of $P_{ij}$.  Indeed,
we can now say rigorously that the
sequence $(n_1,...,n_8)$ ``works'' for all
points in the interior of $P_{ij}$.
\endproof

\begin{lemma}
$\widetilde \Pi$ is not defined on any boundary point of $P_{ij}$
for any indices $i =1,2,3,4$ and any $j=1,...,64$.
\end{lemma}

\startproof
We use the notation from the previous proof.
For each face $F$ of $P_{ij}$ we exhibit an
index $k$ such that one of the following two equations
simultaneously holds relative to any vertex $v_0$ of $F$:
\begin{itemize}
\item $g_k(v_{k-1}) = n_k$.
\item $g_k(v_{k-1}) = n_{k+4}$.
\end{itemize}
This shows that 
$$F_{k-1}=\widetilde E_{k-1} \circ \ldots \widetilde E_1(F)$$
lies entirely in one hyperplane of $\widetilde L_k$.
Hence $\widetilde E_k$ is nowhere defined on $F_{k-1}$.
Hence $\widetilde \Pi$ is nowhere defined on $F$.
Since $F$ is an arbitrary face, $\widetilde \Pi$ is
nowhere defined on the boundary of $P_{ij}$.
\endproof

\newpage

\section{Proof of the Renormalization Theorem}

\subsection{Some Terminology}
\label{renormoverview}

We use the notation from the
Renormalization Theorem.
Here we establish some terminology which will
be useful to us in this chapter and in
later ones.  
\newline
\newline
{\bf Forward Atoms:\/}
Recall that $\widehat A$ 
is partitioned into atoms.
We call these atoms {\it forward\/}
$\widehat A$-atoms.  We define the
forward $\widehat B$-atoms in the 
same way.
\newline
\newline
{\bf Neutral and Backward Atoms\/}
Let $P$ be a forward $\widehat A$ atom.
There is some $n>1$, depending on $P$, such that
$\widehat \Psi^k(\widehat P)$ is disjoint from
$\widehat A$ for $k=1,...,n-1$ but
$\widehat \Psi^n(P) \subset \widehat A$.
We call the polyhedra
$\widehat \Psi^k(P)$ {\it neutral\/} 
$\widehat A$-atoms when $k=1,...,n-1$.
We call $\widehat \Psi^n(P)$ a
{\it backward\/} $\widehat A$-atom.
$\widehat A$ is also partitioned into
backward $\widehat A$-atoms, and
the first return map of
$\widehat \Psi^{-1}$ is entirely
defined and a translation on each one.
\newline
\newline
{\bf Chains:\/} 
Let $P$ be a forward $\widehat A$-atom.  We
call the sequence
$P, ... ,\widehat \Psi^{n-1}(P)$ a
forward $\widehat A$-chain.  So, each
forward $\widehat A$-chain is a finite union
of $\widehat A$-atoms, the first of
which is forward and the rest of
which are neutral.  We make the same definitions for
$\widehat B$.
\newline
\newline
{\bf Backward Chains:\/}
We define {\it backward chains\/} just as
we defined the forward chains, except that we use
backward atoms in place of forward
atoms, and we use $\widehat \Psi^{-1}$
in place of $\widehat \Psi$.  We think
of the backward atom as being the last
atom in the backwards chain.
For each forward chain there is a backward chain
that shares all its atoms except the first one.
Likewise, for every backward chain there is a
forward chain that shares all of its atoms except
for its last one.
\newline
\newline
{\bf Types:\/}
We distinguish two types of layers of $\widehat A$.
The {\it type-1\/} layers are the ones on which
the restriction of $\widehat R$ to that layer
is a translation.  The other layers we call
type-2.  Note that this definition is
consistent with the notion of type defined
in association with the Reduction Theorems.
We define the type of a $\widehat A$-atom
and a $\widehat A$-chain according to the
type of the layer that contains it.

\subsection{Statements 1, 2, and 5}

We compute the list of $\widehat A$ and
$\widehat B$ atoms using the method
described in \S \ref{tf}.
Figure 10.1 shows a slice of the set of
$\widehat B$-atoms at a particular height.

\begin{center}
\resizebox{!}{3.7in}{\includegraphics{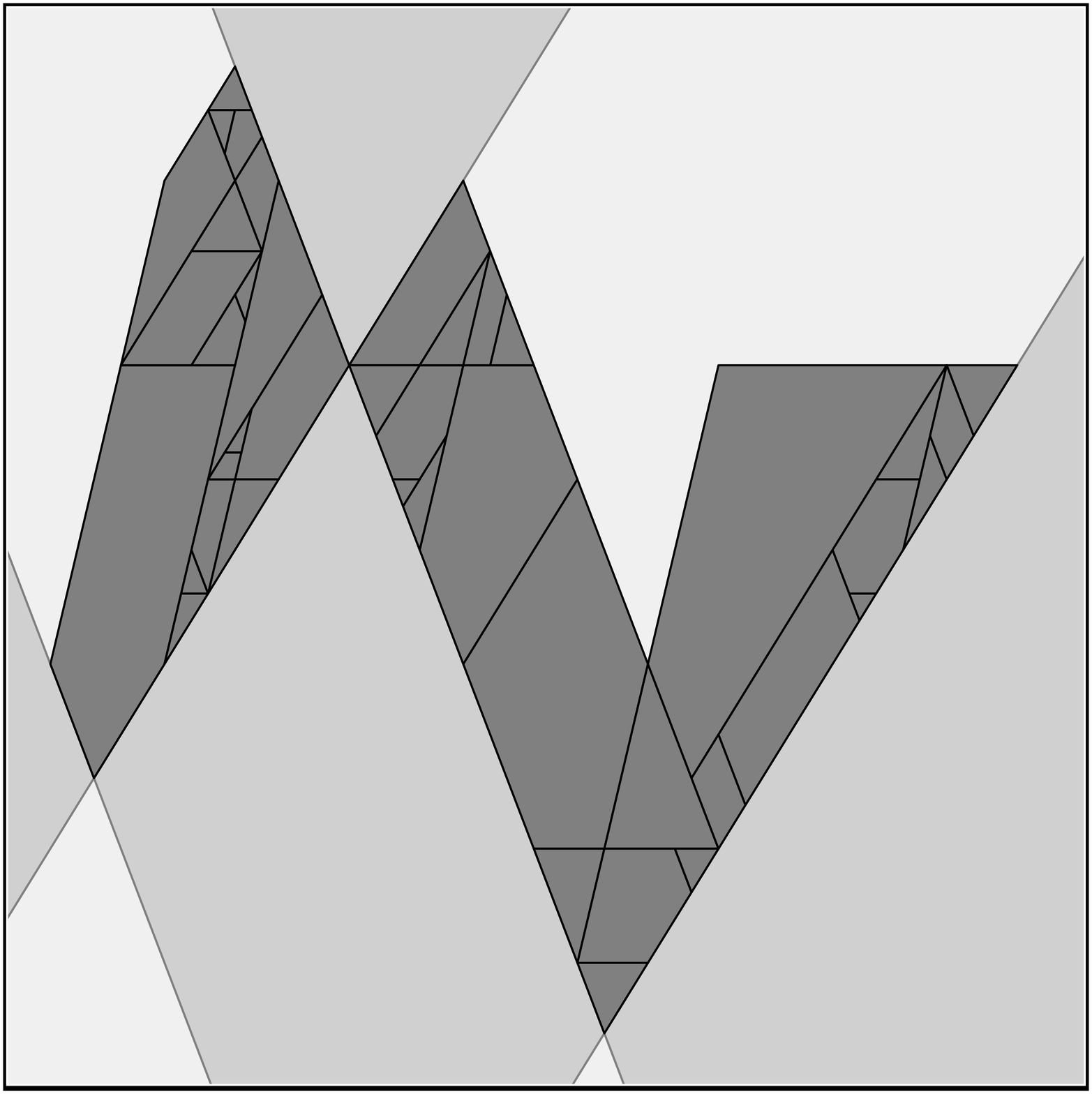}}
\newline
{\bf Figure 10.1:\/} The $\widehat B$-atoms sliced at
$-7+4\phi \approx .5278$
\end{center}  

We use the method described in \S \ref{pd} to
check that each of the polyhedra on our candidate
list of $\widehat A$-atoms is indeed a $\widehat A$-atom.
We do the same thing for the $\widehat B$-atoms.
In the case of the $\widehat B$-atoms, we also
verify that each face $F$ of a $\widehat B$-atom has
one of $3$ properties.
\begin{itemize}
\item The first return map $\widehat \Psi|\widehat B$
is not defined on $F$. That is, some small power of
$\widehat \Psi$ maps $F$ into one of the faces of
one of the $64$ special polyhedra that define
$\widehat \Psi$.
\item $F$ lies in $\partial \widehat B$.
\item The first return of $F$ to $\widehat B$ lies in $\partial \widehat B$.
\end{itemize}
Since these properties hold for every edge, all the
$\widehat B$ atoms are maximal.
\newline
\newline
{\bf Remark:\/} At some point we also checked that all the
$\widehat A$-atoms are maximal, but we did not save the
computer code for this.  This check is not necessary for
our arguments.
\newline

Given that the $\widehat B$-atoms are all maximal,
they have pairwise disjoint interiors.  We also
check that they are strictly convex.  We now apply
the slicing method of \S \ref{covertest} to
each branch of each layer of $\widehat B$. We
conclude that the $\widehat B$-atoms partition
$\widehat B$. That is, we have obtained the
complete list.

As we mentioned in \S \ref{tf}, we obtain the
candidate list of $\widehat A$ atoms by pulling
back the $\widehat B$ atoms by the action of the
renormalization map $\widehat R$.  From this
way of constructing things, we automatically
know that the $\widehat A$-atoms are strictly
convex and have pairwise disjoint interiors.
(We don't need the maximality for this.)
We now apply the slicing method to
$\widehat A$ just as we did for
$\widehat B$. We conclude that
$\widehat A$ is partitioned into the
$\widehat A$ atoms on our list.

We have set things up so that, on each
type-1 layer, $\widehat R$ establishes
a bijection between the forward $\widehat A$-atoms
and the forward $\widehat B$-atoms.
We check that this bijection is
compatible with action of the first return maps
$\widehat \Psi|\widehat A$ and
$\widehat \Psi|\widehat B$.  Since we have already
verified that everything in sight is an
atom, we just need
to check the action of the relevant
maps on one interior point of each
atom.  For the check, we choose interior golden
points as in \S \ref{basictest}.

We have set things up so that, on each
type-2 layer, $\widehat R$ establishes
a bijection between the forward $\widehat A$-atoms
and the backward $\widehat B$-atoms.
In the same way as in the type-1 case, we 
check that the action of $\widehat R$ is
compatible with action of the first return maps
$\widehat \Psi|\widehat A$ and
$\widehat \Psi^{-1}|\widehat B$.

These calculations prove Statement 1 of the Renormalization
Theorem.  

Statement 2 of the Renormalization Theorem is really just
a description of the action of $\widehat R$.  We have
constructed $\widehat R$ so that it acts on the
horizontal slices as $R$ acts on $\R/2\Z$.

Statement 5 is purely combinatorial. 
We simply check, in each case,
that each $\widehat B$-chain has a shorter length than
the corresponding $\widehat B$-chain.
\newline
\newline
{\bf Remark:\/} For the record, we mention that there are
$678$ $\widehat B$ atoms and $3 \times 678=2034$ $\widehat A$-atoms.
We also compute that the longest $\widehat A$-chain has
length $703$ and the longest $\widehat B$-chain has length
$109$.  In other words, for generic points in $\widehat A$,
the first return map is defined in at most $704$ iterates
and for generic points in $\widehat B$ the first return
map is defined in at most $110$ iterates.

\subsection{Statement 3}

When we lift $\Theta(\Sigma)$ to $\R^3$ we see that it consists
precisely in the planes spanned by the vectors
$$(2\phi^{-1},1,0); \hskip 30 pt
(0,-2,1)$$
which intersect the line
$\{0\} \times \R \times \{0\}$ in points of the form
$(0,t,0)$ where $t \in \Z[\phi]$.   This follows from
a routine calculation. Another calculation
shows that each choice of $\widehat R$ mentioned
above preserves this set.  

\subsection{Statements 4 for $\widehat B$}

Recall that $\widehat B$ is partitioned into $6$ layers,
as discussed in \S \ref{renormstruct}.  It turns out that
some periodic tiles intersect some layers of $\widehat B$
but not others.  To make our proof go cleanly, we
introduce the notation of a $\widehat B$-periodic
tile.  A $\widehat B$ periodic tile is an intersection
of the form
\begin{equation}
P \cap \Big((\R/2\Z) \times (\R/2/Z) \times J_{\lambda}\Big)
\end{equation}
where $J_{\lambda}$ is the $6$ intervals
$J_1,...,J_6$ in the 
partition defining the layers of $\widehat B$.

Using essentially the same method as in \S \ref{generate},
we generate a candidate list of $\widehat B$ periodic
tiles whose orbits avoid $\widehat B$.  In practice,
we simply enumerate all the periodic tiles whose
orbits avoid at least one layer of $\widehat B$, and
then we keep track of which layers each orbit intersects.
We check that the period of the longest tile on the list is $33$.

Let $\lambda \in \{1,2,3,4,5,6\}$ be an integer that describes
a layer of $\widehat B$ currently of interest to us.
Let $i \in \{1,...,64\}$ be an index that describes
one of the $64$ special polyhedra that partition
$\widehat \Sigma$.

Let ${\cal B\/}(i,\lambda)$ be the union of all the
$\widehat B$-tiles and forward/neutral $\widehat B$-atoms
that are contained in the region
\begin{equation}
\label{puck}
P(i,\lambda)=P_i \cap \Big(\R/2\Z) \times (\R/2/Z) \times J_{\lambda})\Big)
\end{equation}
By construction, these tiles have pairwise disjoint interiors.
Say that $(i,\lambda)$ is {\it good\/} if
the tiles in  ${\cal B\/}(i,\lambda)$ partition
the region in Equation \ref{puck}.

\begin{lemma}
\label{caseB}
Statement 4 holds for $\widehat B$ provided that
every pair $(i,\lambda)$ is good.
\end{lemma}

\startproof
Let $p_1,...,p_{109}$ be a (non-periodic)
 generic orbit portion of length
$109$.  Consider the point $p=p_{109}$.  There is some
pair $(i,\lambda)$ such that $p \in P(i,\lambda)$.
Since $(i,\lambda)$ is good and $p$ is generic,
$p$ lies in the interior of either a
$\widehat B$ periodic tile or in the interior of a
forward/neutral $\widehat B$-atom.  In the
former case, $p$ has period at most $33$.
This is a contradiction.  If $p$ lies in
the interior of a forward $\widehat B$-atom,
then $p \in \widehat B$, by definition, and
we are done.

So, suppose $p$ lies in the interior of
neutral $\widehat B$-atom.  Since the
longest $\widehat B$-chain has length
$109$, there is some $m \leq 109$ such
that $\widehat \Psi^{-m}(p) \in \widehat B$.
But then $p_{109-m} \in \widehat B$.
This covers all the cases.
\endproof

Our computer code is equipped to verify that
all indices are good, but we can get away
with doing less computational work.
The indices $(i,\lambda)$ are automatically good
for $i \in \{1,...,22\}$.  These correspond to
the special polyhedra on which $\widehat \Psi$
acts trivially. So, we only have to consider
indices $i>22$.

We can eliminate more indices by considering the
dynamics of $\widehat \Psi$.  For instance,
$\widehat \Psi(P_{26})=P_{25}$.  Hence, if
$(25,\lambda)$ is good, then so is $(26,\lambda)$.
In Lemma \ref{transition}, we prove that the list of
indices 
\begin{equation}
I=\{23,25,32,40,41,44,46,53,61,62\}
\end{equation}
is such that every well-defined and nontrivial
orbit intersects $P_i$ for some $i \in I$.
For this reason, we only have to check
that $(i,\lambda)$ is good for
$i \in I$ and $\lambda \in \{1,2,3,4,5,6\}$.

There is one more savings we can make.  
Everything in sight is invariant under the
involution
\begin{equation}
(x,y,z) \to (2,2,2)-(x,y,z).
\end{equation}
This involution permutes the special polyedra,
commutes with $\widehat \Psi$, preserves
$\widehat A$ and $\widehat B$, and permutes
the list of atoms. For this reason, it
suffices to prove that the indices
$(i,\lambda)$ are good for
$i \in I$ and $\lambda \in \{1,2,3\}$.
Using the slicing method of \S \ref{covertest}
we check that this is the case.  This establishes
Statement 4 for $\widehat B$.  Incidentally,
our calculaton show that our candidate list
of $\widehat B$-periodic tiles is correct and
complete.
\newline
\newline
{\bf Remark:\/} The calculation we make is
fairly massive, but the reader can use our
program and see the individual slices plotted
as they are tested.  The reader can also select
any index $(i,\lambda)$ and plot the corresponding
partition.

\subsection{Statement 4 for $\widehat A$}
\label{PROOFA} 

In principle, we could prove Statement 4 for
$\widehat A$ just as we proved it for
$\widehat B$. This time there are $18$
layers.  The problem with this approach
is that it leads to a massive calculation.
For instance, the partition associated
to the index pair $(23,7)$ has $2562$
polyhedra.   

Were we to take this approach (and succeed)
we would know than any generic orbit portion of
length $703$ intersects $\widehat A$.
The bound of $703$ is the sharp constant.
There is an orbit portion of length $702$
that avoids $\widehat A$. 
We will take a different approach and get
the slightly worse constant $812$.
The benefit we get from the other approach 
is that the calculation we make is vastly
shorter.  The idea is to use the result
we have already proved for $\widehat B$
and only calculate ``the difference''
between $\widehat B$ and $\widehat A$,
so to speak.

Recall that each layer
of $\widehat B$ is divided into $4$
branches.  Moreover, each layer of
$\widehat B$ is partitioned into finitely
many smaller pieces, corresponding to
the layers of $\widehat A$.  We call
these pieces of $\widehat B$ the
{\it sub-layers\/} of $\widehat B$.
Precisely,
\begin{itemize}
\item Layer $1$ of $\widehat B$ has $5$ sublayers.
\item Layer $2$ of $\widehat B$ has $1$ sublayers.
\item Layer $3$ of $\widehat B$ has $3$ sublayers.
\item Layer $4$ of $\widehat B$ has $3$ sublayers.
\item Layer $5$ of $\widehat B$ has $1$ sublayers.
\item Layer $6$ of $\widehat B$ has $5$ sublayers.
\end{itemize}
The palindromic nature of the list comes from
the involutive symmetry mentioned in the
previous section.

Each layer of $\widehat B$ is divided into
$4$ branches. Each branch is a golden
polyhedron.  By slicing each branch with 
relevant horizontal planes, we can say that
each sublayer of $\widehat B$ is also
divided into $4$ branches.  In summary,
$\widehat B$ has $18$ sublayers and
each sublayer has $4$ branches.
We index the sublayers by $\lambda \in \{1,...,18\}$
and the branches by $\beta \in \{1,2,3,4\}$.

Finally we discuss $\widehat A$.
We define the $\widehat A$-periodic tiles just as
we defined the $\widehat B$-periodic tiles.
We find a candidate list of all $\widehat A$-periodic
tiles as we did in the $\widehat B$ case.
The longest period in this case is $213$.
We say that a pair
$(\lambda,\beta)$ is {\it good\/} if
the polyhedron corresponding to $(\lambda,\beta)$
is partitioned into a union of $\widehat A$-periodic
tiles and forward/neutral $\widehat A$-atoms.

\begin{lemma}
Statement 4 holds for $\widehat A$ provided that
every index $(\lambda,\beta)$ is good.
\end{lemma}

\startproof
the same argument as in Lemma \ref{caseB} proves that
any orbit portion of length $703$ that
starts in $\widehat B$ intersects $\widehat A$.
Combining this result with Lemma \ref{caseB}, 
we see that every generic orbit
portion of length $812=703+109$ intersects
$\widehat A$. 
\endproof

To finish the proof of Statement 4 for $\widehat A$ we
just have to prove that every index
$(\lambda,\beta)$ is good.  Using the involutive
symmetry discussed in the previous section,
it suffices to prove this result for
indices $\lambda \leq 9$.  
Using the slicing method of \S \ref{covertest}
we verify that the indices $(\lambda,\beta)$
are good for $\lambda \in \{1,...,9\}$
and $\beta \in \{1,2,3\}$. This completes
the proof.
\newline
\newline
{\bf Remark:\/} This time, our calculations do
not verify that our list of
periodic $\widehat A$ tiles is complete.
Rather, all we know is that each
polyhedron on the list is indeed a
$\widehat A$-periodic tile, and that
our list of $\widehat A$-periodic tiles
whose orbits intersect $\widehat B$
is complete.  We can still say something
useful, even with an incomplete list.
Statement $4$ for $\widehat A$ implies that
$\widehat \Sigma$ is partitioned into
the forward/neutral $\widehat A$-atoms
and a finite union of $\widehat A$ periodic
tiles, some of which perhaps are not on our list.
If we have a generic orbit portiont that is
not contained in a periodic orbit, then it
must intersect one of the forward or
neutral $\widehat A$-atoms. But the longest
$\widehat A$-chain has length $703$.
this establishes that
a generic orbit portion of length $703$
intersects $\widehat A$ provided it
does not lie on a periodic orbit.
The constant $703$ is sharp because
of the $\widehat A$-chain of length $703$.

\subsection{Statement 6}

Any generic orbit that avoids $\widehat B$ is
contained in the orbit of one of the
$\widehat B$-periodic tiles on our list.
We check that all these orbits avoid $\widehat A$.
Hence, any orbit that avoids $\widehat B$
also avoids $\widehat A$.

\newpage

\section{Critical Renormalization Calculations}

\subsection{Discussion}
\label{renormdyn}

The goal of this chapter is to prove the
Fundamental Orbit Theorem and a second result
which we call the Fixed Point Theorem.
Both results involve fixed points
of the renormalization operator.  They are
the cleanest applications of the Renormalization
Theorem.  In the next chapter, we will give
some messier applications, including the proof
of the Near Reduction Theorem.

We first discuss the general nature of our calculations.
We use the notation from the previous chapter.
Let $\cal A$ be the set of $\widehat A$-chains.
Likewise define $\cal B$.   
The Renormalization Theorem gives us a canonical
$3$-to-$1$ map from $\cal A$ to $\cal B$.
This map does not depend on any conventions
concerning the boundaries of the 
layers of $\widehat A$.   Let
$\chi: {\cal A\/} \to {\cal B\/}$ denote the
map just described.

Since $\Sigma$ can be considered as a dynamically
invariant subset of $\widehat \Sigma$, we can transfer
the structure to $\Sigma$.  For instance, an
$A$-{\it atom\/} of $\Sigma$ is a connected intersection
of $\Sigma$ with the $\widehat A$-atoms.   Every
polygon in sight is a golden polygon.  We deal
with these polygons and polyhedra using
the methods from \S 6.

Given an $A$-atom $X$ and a $B$-atom $Y$, we write
$X \to Y$ if $\chi([X])=[Y]$.  Here $[X]$ is
the forward $A$-chain containing $X$ and $[Y]$ is the
$B$-chain containing $Y$.  We take the $B$-chains
to be forwards when they are discussed
in connection with type-1 $A$-chains, and
backwards when discussed in connection with
type-2 $A$-chains.
\newline
\newline
{\bf Important Observation:\/}
By construction, a generic point never intersects
the boundary of an $A$-atom or the boundary
of a $B$-atom.

\subsection{Proof of the Fundamental Orbit Theorem}

Referring to the notation in the Fundamental Orbit Theorem,
we find a partition of each $T_{ij}^{\pm}$ by a
union of $A$-atoms and periodic tiles such that
\begin{equation}
\label{instance}
\tau \to \psi^{-1}(R_{ij}^{\pm}(\tau)),
\end{equation}
for every $A$-atom $\tau$. The periodic
tiles all lie in ${\cal T\/}(2)$.

For future reference, we refer to the
above mentioned partitions as the
$A$-{\it partitions\/}. The left side of
Figure 11.1 shows the partition for $T_{11}^+$.

\begin{center}
\resizebox{!}{5in}{\includegraphics{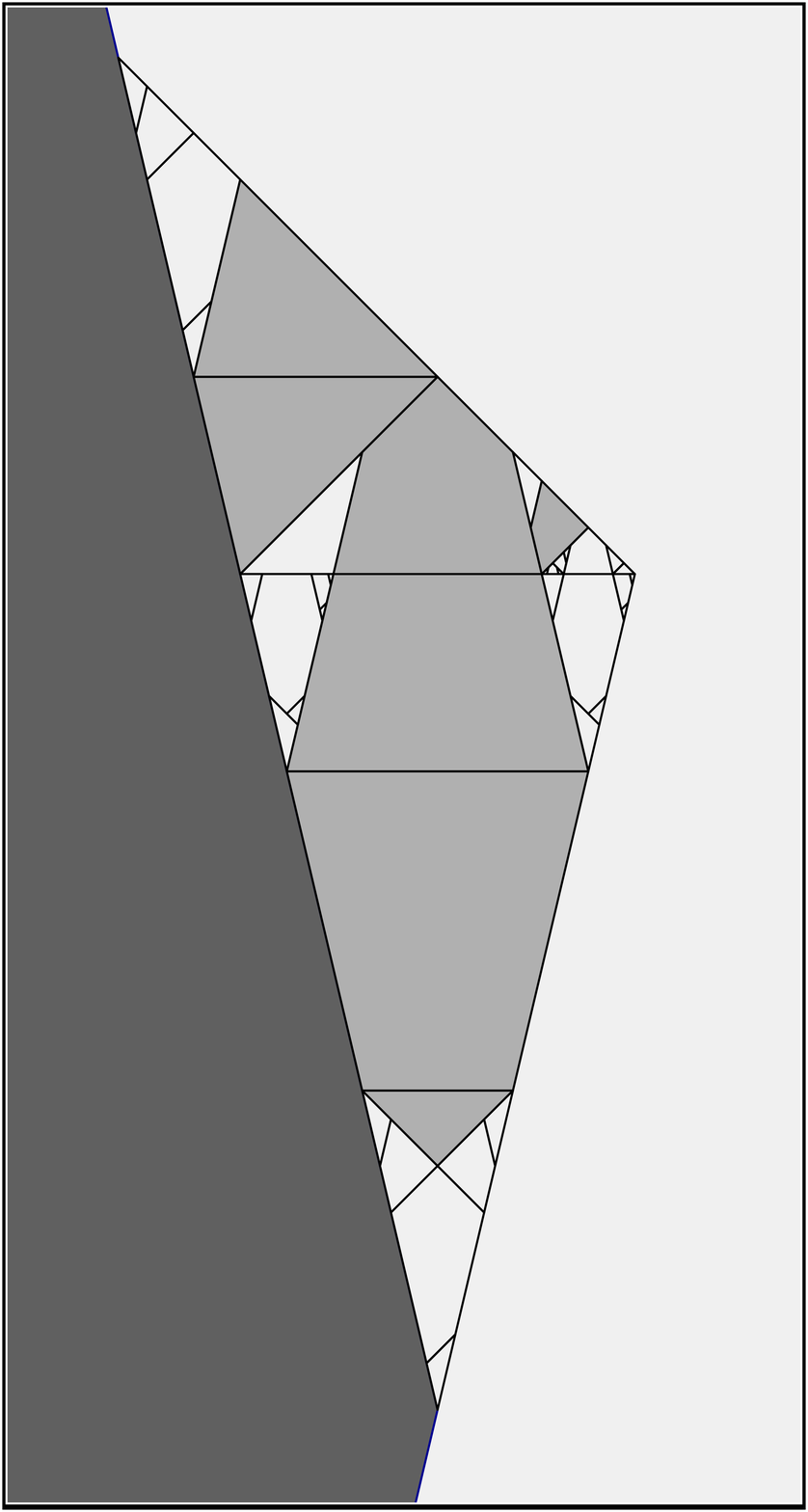}}
\newline
{\bf Figure 11.1:\/} The partition of $T_{11}^+$ and
an edge-covering.
\end{center}

The $A$-partition of $T_{11}^+$ is 
a union of $41=34+7$ $A$-atoms and periodic tiles.  The
periodic tiles are all shaded. The dark polygon on the side
is part of the Penrose kite.  The horizontal
lies through the periodic tiles correspond to
the tops and bottoms of the
layers of $\widehat A$.  
Thus, for instance, even though
the large octagon is a single periodic tile,
we treat it as $4$ separate tiles in our proof.
Each of the $3$ triangles
$T_{12}^+$, $T_{31}^+$ and $T_{41}^+$ are
also partitioned into $41=43+7$ polygons.  These partitions
look exactly like the one in Figure 11.1, up to
reflection in a vertical line.
The triangle $T_{41}^+$ is partitioned into $4=3+1$ polygons
and $T_{41}^-$ is partitioned into $3=2+1$ polygons.
It turns out that we do not need to consider
the $4$ remaining triangles, as we explain below.
The reader can see all the partitions, and
in color, using our applet.

\begin{lemma}
\label{finite}
Suppose that Equation \ref{instance} holds for all the
$A$-atoms in the $A$-partitions.  Then the
Fundamental Orbit Theorem is true.
\end{lemma}

\startproof
Let $O_1$ be a generic orbit that intersects
$T_{ij}^{\pm}$ but avoids
set ${\cal T\/}(2)$.  Suppose first
that $O_1$ is not a distinguished orbit.
Since all the
periodic tiles in the $A$-partitions lie
in ${\cal T\/}(2)$, some point $p_1 \in O_1$
intersects an $A$-tile $\tau$ of the $A$-partition.
Since $p_1$ is generic, it cannot intersect
the boundary of an $A$-atom.
Hence $p_1$ intersects the interior of some
$A$-atom $\tau$ of the $A$-partition.
But then, by Equation \ref{instance},
we have $O_1 \leadsto O_2$, where
$O_2$ is the $\Psi$-orbit containing
$p_2=\psi^{-1}(T_{ij}^{\pm}(p_1)$.

When $O_1$ is a distinguished orbit, we attach
it to one of the slabs or the other.  The
argument above is the same, except that
what we get is that $p_1$ must lie
in the interior of one of the horizontal
edges of the $A$-partition.
But then Equation \ref{instance} gives
us the same result as for the ordinary
orbits.  Were we to attach the special
orbit to the other slab, we could
still use Equation \ref{instance}, but
it would involve a different $A$-atom.
\endproof

To prove the Fundamental Orbit Theorem,
it only remains to check Equation \ref{instance}
for each of the $A$-atoms in the $A$-partitions.
Each individual check is a finite calculation,
similar to what we did in the last chapter.
We now discuss this calculation.

First of all, we reduce the calculation from
$10$ triangles to $6$ triangles, as follows.
We have the relations
\begin{equation}
\Psi(T_{11}^-)=T_{12}^+; \hskip 15 pt
\Psi(T_{12}^-)=T_{11}^+; \hskip 15 pt
\Psi(T_{31}^-)=T_{32}^+; \hskip 15 pt
\Psi(T_{32}^-)=T_{31}^+.
\end{equation}
Using these relations, we see right away that the truth
of the Fundamental Orbit Theorem for
each triangle on the right hand side of the relation
implies the truth of the Fundamental Orbit Theorem
for the triangle on the left hand side.
Thus, we only need to check the $5$ triangles
$T_{ij}^+$ and the triangle $T_{41}^-$.

First of all, we check that  each relevant triangle
$T_{ij}^{\pm}$ is indeed
partitioned by the polygons we have listed.
The check is similar to what we did in \S 6.
As for verifying Equation \ref{instance}, we
have to make the same verification
$141=4 \times 34+3+2$ times.

Let $P$ and $Q$ respectively be the polygons on the left and
right hand sides of Equation \ref{instance}.
To check that $P \to Q$, in each of the $141$ cases,
we do the following calculation.
\begin{itemize}
\item We check that $\widehat P=\Theta(P)$ is contained in a $\widehat A$-atom.
This amounts to producing the two integers $m$ and $n$,
with $m \leq 0<n$, such that
$\widehat \Psi^k$ is entirely defined on the interior of
$\widehat P$ for all $k \in [m,n]$, and
$\widehat \Psi^k(\widehat P) \subset \widehat A$
for $k=m,n$ but not for $k \in (m,n)$.
We verify these things using the methods discussed in \S \ref{pd}.
Define $\widehat P_0=\widehat \Psi^m(\widehat P)$.
\item 
We check that $\widehat Q=\Theta(Q)$ is contained in a $\widehat B$-atom.
This is the same kind of check as in the previous case.
Let $\widehat Q_0=\widehat \Psi^m(\widehat Q)$ be the
iterate that is contained in $\widehat B$.  (Depending
on the layer, $m$ is either non-positive or non-negative.)
\item 
We check that $\widehat R(\widehat P_0)=\widehat Q_0$.
\end{itemize}

We simply perform the $141$ calculations and observe
that they all work.  
This completes the proof of the Fundamental Orbit Theorem.
\newline
\newline
\noindent
{\bf Remark:\/}
Note that we are not quite verifying that
$\widehat P$ and $\widehat Q$ are atoms, but
only that they are contained in atoms.  
Similarly, we are only verifying that $P$ and $Q$
are contained in atoms, and that their
boundaries are all contained in lines
defined over $\Z[\phi]$.  This suffices for our
purposes.  We mention, however, that we produce
our $A$-partitions by slicing the $\widehat A$-atoms
and periodic tiles by the image of $\Sigma$.
Thus, were to make the extra check, it would work out.
\newline

The reader can use our applet to survey the calculations.
Here we give an example illustrating the rough size and
complexity of the calculation.
There is a small isosceles $A$-tile in the
partition of $T_{11}^+$ that contains the rightmost vertex
of $T_{11}^+$.  This triangle is perhaps too small to see
in Figure 11.1.  The vertices of this triangle are
$$(997-616 \phi,-236+146 \phi) \hskip 20 pt
(1984 - 1226 \phi,-3+2 \phi) \hskip 20 pt
(10-6\phi,-3+2\phi).$$
In this case
$$\widehat P_0=\widehat \Psi^{-508}(\widehat X); \hskip 30 pt
\widehat Q_0=\widehat \Psi^{87}(\widehat Q).$$
Things don't get much worse than this, because the maximum
chain length is $703$.

\subsection{The Fixed Point Theorem}
\label{fpt}

It turns out that there is a second region that plays a role
that is similar to the fundamental triangles.
Let $Q_0$ be the quadrilateral with vertices
\begin{equation}
(6-2\phi,\phi^{-3}); \hskip 20 pt
(3,0); \hskip 20 pt
(2\phi,\phi^{-3}); \hskip 20 pt
(3,2\phi^{-1})
\end{equation}
$Q_0$ is a small kite whose line of symmetry is
the vertical line $x=3$.  The bottom point of
$Q$ is $(3,0)$.  See Figure 11.2 below.
Let $D$ denote the homothety that
fixes $(3,0)$ and expands distances
by a factor of $\phi^{3}$.  Define
\begin{equation}
Q_n=D^{-n}(Q).
\end{equation}
The family $\{Q_n\}$ is a shrinking family
of kites that limits to $(3,0)$.
Let $\rho$ denote the reflection in
the vertical line $x=3$.

\begin{theorem}[Fixed Point]
Let $O_1$ be any sufficiently long and generic orbit that
intersects $Q_1$.  Then $O_1 \leadsto O_2$,
and $O_2$ is such that 
$$\psi'(\overline O_2) \cap Q_0=\rho \circ D(O_1 \cap Q_1).$$
Here $O_2$ is the reflection of $O_2$ in the $x$-axis
and $\psi'$ is the outer billiards map.  In particular
$O_1 \to O_2$, in the sense of Theorem \ref{renorm4}.
\end{theorem}

\startproof
The proof here is essentially the same as the proof we gave
for the Fundamental Orbit Theorem.  
Let $\Delta: \R^2 \to \R^2$ be the dilation by
$\phi^3$ such that $\Delta(3,0)=(-3,-2)$.
We have the equation
\begin{equation}
\label{identity1}
\psi' \circ \overline{\Delta(p)}=\rho \circ D(p)
\end{equation}
for all $p \in Q_1$.

We exhibit a tiling of $Q_1$ by a union of $48$ polygons,
$42$ of which are $A$-atoms and $6$ of which are
periodic tiles.
For each $A$-tile $\tau$ in the
partition, we verify,
by the same means as above, 
the equation
\begin{equation}
\label{identity2}
\tau \to \Delta(\tau),
\end{equation}
Our result follows from Equation \ref{identity1},
Equation \ref{identity2}, and the same argument as
in Lemma \ref{finite}.
\endproof

\begin{center}
\resizebox{!}{5.6in}{\includegraphics{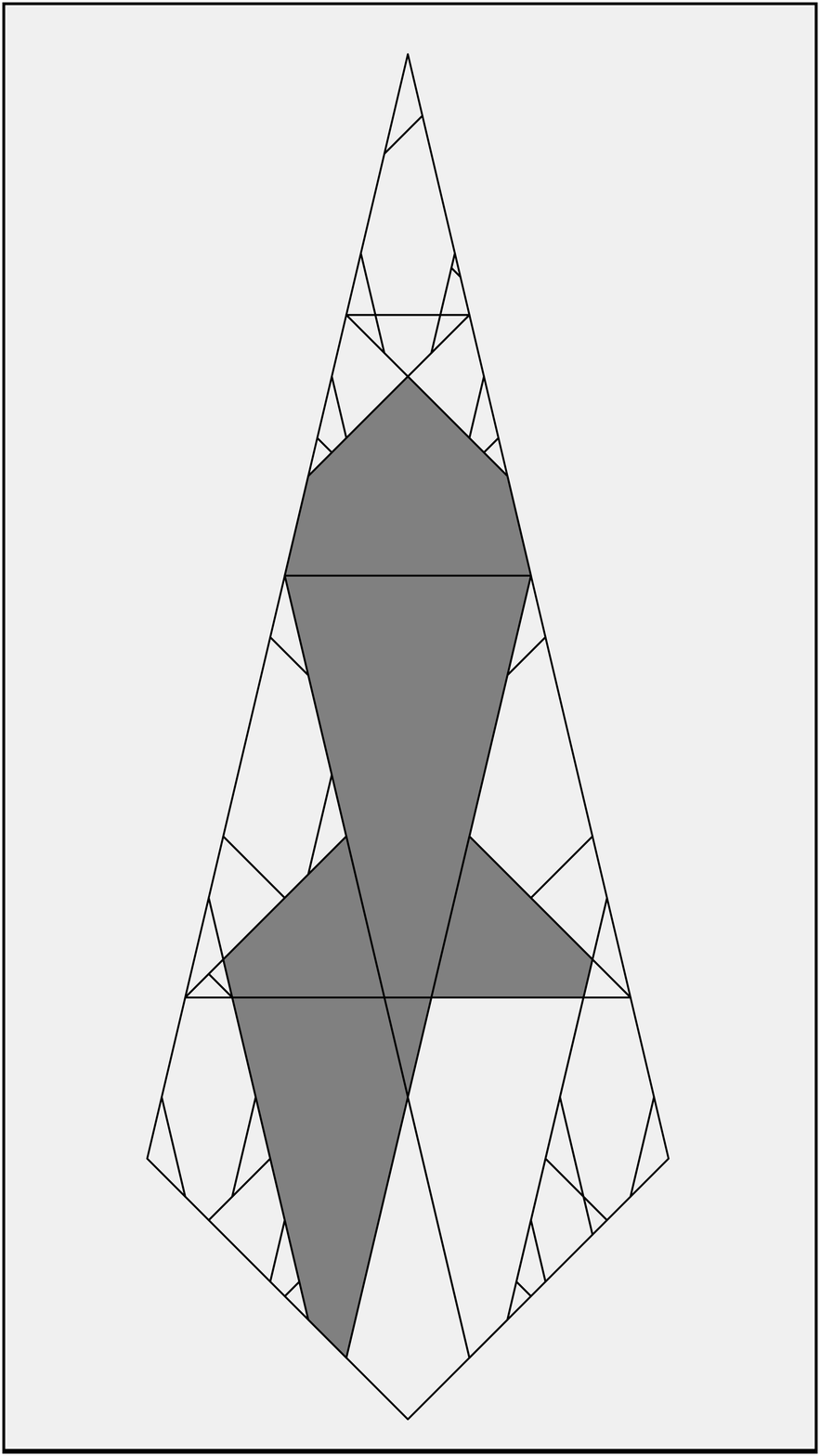}}
\newline
{\bf Figure 11.2:\/} $Q_1$ and its partition into
atoms and periodic tiles.
\end{center}

It is a consequence of the Fixed Point Theorem that the
point $(3,0)$ is a point of dilation symmetry for the
dynamical tiling associated to outer billiards on the penrose
kite.  Figure 1.3 shows what is going on.  The white
regions in Figure 1.3 are small copies of portions the
tiling $\cal T$ and the shaded regions are increasingly
small periodic tiles.  The reader can probably
spot $Q_1$ in Figure 1.3:  It is flanked by the two largest
dark tiles.

\newpage

\section{Proof of the Near Reduction Theorem}

\subsection{Discussion}

The main structural relation between
the $\widehat A$-atoms and the
$\widehat B$-atoms is the $3$-to-$1$
map from the set of $\widehat A$-chains
to the set of $\widehat B$-chains. Now
we mention a second piece of structure.
Though we did not formally check this in general,
in the
cases of interest it turns out that each
$\widehat B$-atom is partitioned into a finite union
periodic tiles and 
$\widehat A$-atoms.  Thus, we have a kind of
auxilliary dynamical system defined on the
level of $\widehat A$-atoms.  Starting with
a $\widehat A$-atom $\tau$, we can consider
all those $\widehat B$-atoms -- there are
finitely many, and all in the same chain -- such
that $\tau \to \tau'$.  Next, we subdivide
each such $\tau'$ back into $\widehat A$-atoms.
And so on.  All our proofs in this chapter
amount to tracing through this auxilliary
dynamical system.

Our computer program is designed so that
the user can trace through all the
calculations manually, plotting the
relevant tiles.

The way to view this chapter is that we have
one main calculation, contained in \S \ref{bigg},
that takes care of almost all the orbits that
come fairly close to the origin.  The remaining
calculations, such as those in \S \ref{techX},
deal with the several small regions that
the big calculation does not cover.

\subsection{Points Very Near the Kite}
\label{techY}

In our next result, the notation $O_1 \to O_2$ is as in
Theorem \ref{renorm4}.
Let
\begin{equation}
S_{r}=[0,r] \times [0,2].
\end{equation}

\begin{lemma}
\label{TECH1}
Let $r=2+\phi^{-3}$. 
Let $O_1$ be any generic infinite orbit
that intersects $S_r$. 
either $O_1$ intersects $T^+$ or $O_1 \to O_2$ and
$O_2$ intersects $T^+$.
\end{lemma}

\startproof
Our proof refers to
Figure 12.1.  Let $r=2+\phi^{-3}$.   Let
$\rho$ denote complex conjugation -- i.e. reflection in the $x$-axis.
Figure 12.2 shows a covering of
$S_r$ by $13$ tiles.  The left black tile is the
Penrose kite, and the right vertex of the rightmost dark
tile lies on the vertical line $y=r$.  The two
white tiles comprise $T$, the fundamental triangle.
Let $X_1$ be the lower of the two black tiles and let
$L$ be the smallest light tile -- namely, the little kite.
Note that $X_1$ and $L$ share an edge and
$X_1 \cup L$ is a triangle.

Short calculations reveal the following.
\begin{itemize}
\item $\Psi$ is the identity on all the light tiles. 
\item The top black tile $X_2$ is such that $(\psi')^{-1} \circ \Psi^2(X_2) \subset T$.
\item The leftmost dark tile $D_1$ is such that $\Psi(D_1) \subset T$.
\item The bottom dark tile $D_2$ is such that $\rho \circ \Psi(D_2)=X_1$.
\item The rightmost dark tile $D_3$ is such that $\rho \circ \Psi(D_3)=X_1 \cup L$.
\end{itemize}

\begin{center}
\resizebox{!}{2.7in}{\includegraphics{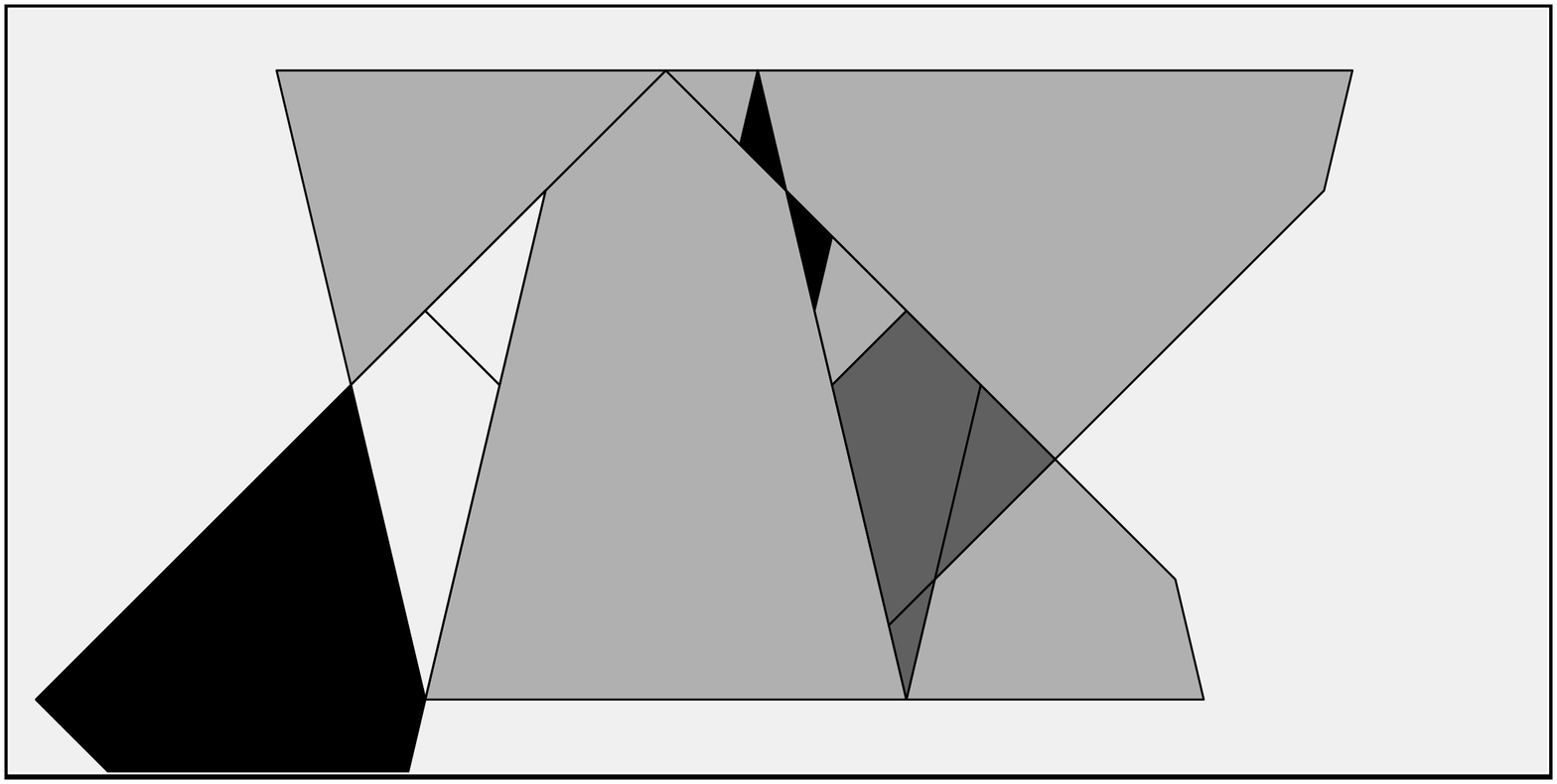}}
\newline
{\bf Figure 12.1:\/}  A covering of $\Sigma_r$.
\end{center}

The region $X_1$ is the triangle $T_{41}^-$
analyzed in the proof of the Fundamental Orbit Theorem.
By the Fundamental Orbit Theorem, the
orbits intersecting $X_1$ renormalize to orbits
that intersect $T$.
Lemma \ref{TECH1} now follows directly from case-by-deductions.
We analyze of the cases in detail.  The remaining cases are similar, and
in fact easier. (We chose the most elaborate case.)
Suppose that $O_1$ intersects $D_3$, the rightmost
dark tile.  Then an associate $O_1'$ intersects $L \cup X_1$.
Being an infinite orbit, $O_1'$ does not intersect the interior of $L$.
Being a well-defined orbit, $O_1'$ does not intersect
$\partial L \cup \partial X_1$.  Hence $O_1'$ intersects
the interior of $X_1$. But then $O_1' \leadsto O_2$, where
$O_2$ intersects $T$.  This shows that $O_1 \to O_2$ where
$O_2$ intersects $T$.
\endproof

\subsection{Two Trouble Spots}
\label{techX}.

Let $\sigma$ be the triangle with vertices
\begin{equation}
(-15+11\phi,2-\phi) \hskip 30 pt
(-36+24 \phi,-19+12 \phi) \hskip 30 pt
(53-31 \phi,2-\phi).
\end{equation}
Let $\sigma^*$ be the triangle with vertices
\begin{equation}
(21-11\phi,2-\phi) \hskip 30 pt
(42-24 \phi,-19+12 \phi) \hskip 30 pt
(-47+31 \phi,2-\phi).
\end{equation}
$\sigma$ and $\tau$ are two small $A$-atoms contained
in $Q_0-Q_1$.  Reflection in the vertical line
$x=3$ interchanges them.  Let $\psi'$ denote the
outer billiards map.

\begin{lemma}
Let $O_1$ be any generic infinite orbit that
intersects $\sigma$.  Then we have
$O_1 \leadsto O_2 \leadsto O_3$
and $\psi'(O_3)$ intersects $T^-$.
\end{lemma}

\startproof
Let $\sigma_1=\sigma$.
Using the notation above, and the same kind of calculations,
we check that $\sigma_1 \leadsto \sigma_2$, where $\sigma_2$
is the triangle (a translate of $\sigma_1$) with vertices
$$(-11+9\phi,-2+\phi) \hskip 30 pt
(-32+22 \phi,-23+14 \phi) \hskip 30 pt 
(57-33 \phi,-2+\phi).$$
We check that $\sigma_1$ is a union of one periodic
tile and $2$ $A$-atoms, $\sigma_2$ and $\sigma_3$.
Next, we check that $\sigma_2 \to \sigma_4$
and $\sigma_3 \to \sigma_5$, where
$\sigma_6=\sigma_4 \cup \sigma_5$ is the triangle with vertices
$$(5-4\phi,-6+4\phi) \hskip 30 pt
(18-12 \phi,7-4\phi) \hskip 30 pt
(-3+\phi,2-\phi).$$
A direct calculation shows that $\psi'(\sigma_6) \subset T^-$.
An argument just like the one
given in Lemma \ref{finite} finishes the proof.
\endproof

\begin{lemma}
Let $O_1$ be any generic infinite orbit that intersects
$\sigma^*$. Then we have $O_1 \leadsto O_2 \leadsto O_3 \leadsto O_4$
and $\psi'(O_4)$ intersects $T^+$.
\end{lemma}

\startproof
The argument here is very similar to what we did in the
previous lemma.  We omit the details.
\endproof

\subsection{The Big Calculation}
\label{bigg}

Now we prove a result similar to the ones in the
previous section, except that the computational
part of the proof is much more intensive.

Given any $p=(x,y)$, define $|p|_x=|x|$.
In our next result, the significance of the
number $2\phi^{-3}$ is that it is precisely the
horizontal width of $Q_0$.

\begin{lemma}
\label{near1}
Let $r=2+\phi^{-1}$.
Let $p_1 \in S_{24}-S_r-Q_1-\sigma-\sigma^*$ be a point
with a generic infinite orbit $O_1$.  Then
$O_1 \leadsto O_2 \leadsto O_3$,
and there is some $k \in \{1,2,3\}$ and some
$q_k \in O_k$ such that $|q_k|_x \leq |p_k|_x-2\phi^{-3}.$
\end{lemma}

\startproof
Say that a {\it symmetric return tile\/} is a 
maximal polygon in the main strip $\Sigma$ on
which both $\Psi$ and $\Psi^{-1}$ are defined.
For each point $p$ in a symmetric return tile $\tau$,
the expressions
$$K_{\pm}(\tau)=|\Psi^{\pm 1}(p)|_x-|p|_x$$
are independent of the choice of $p$.  
If $f(\tau)=0$ it means that $\Psi$ is
the identity on $\tau$.  Otherwise, by
Equation \ref{pregraph}, we have
$\min|K_{\pm}(\tau)| \leq 2\phi^{-3}$.
Accordingly, we only have to consider points
in symmetric return tiles $\tau$ such that
$K_+(\tau)$ and $K_-(\tau)$ are both positive.
Call such tiles {\it positive symmetric
return tiles\/}, or PSRT's for short.

\begin{center}
\resizebox{!}{2.9in}{\includegraphics{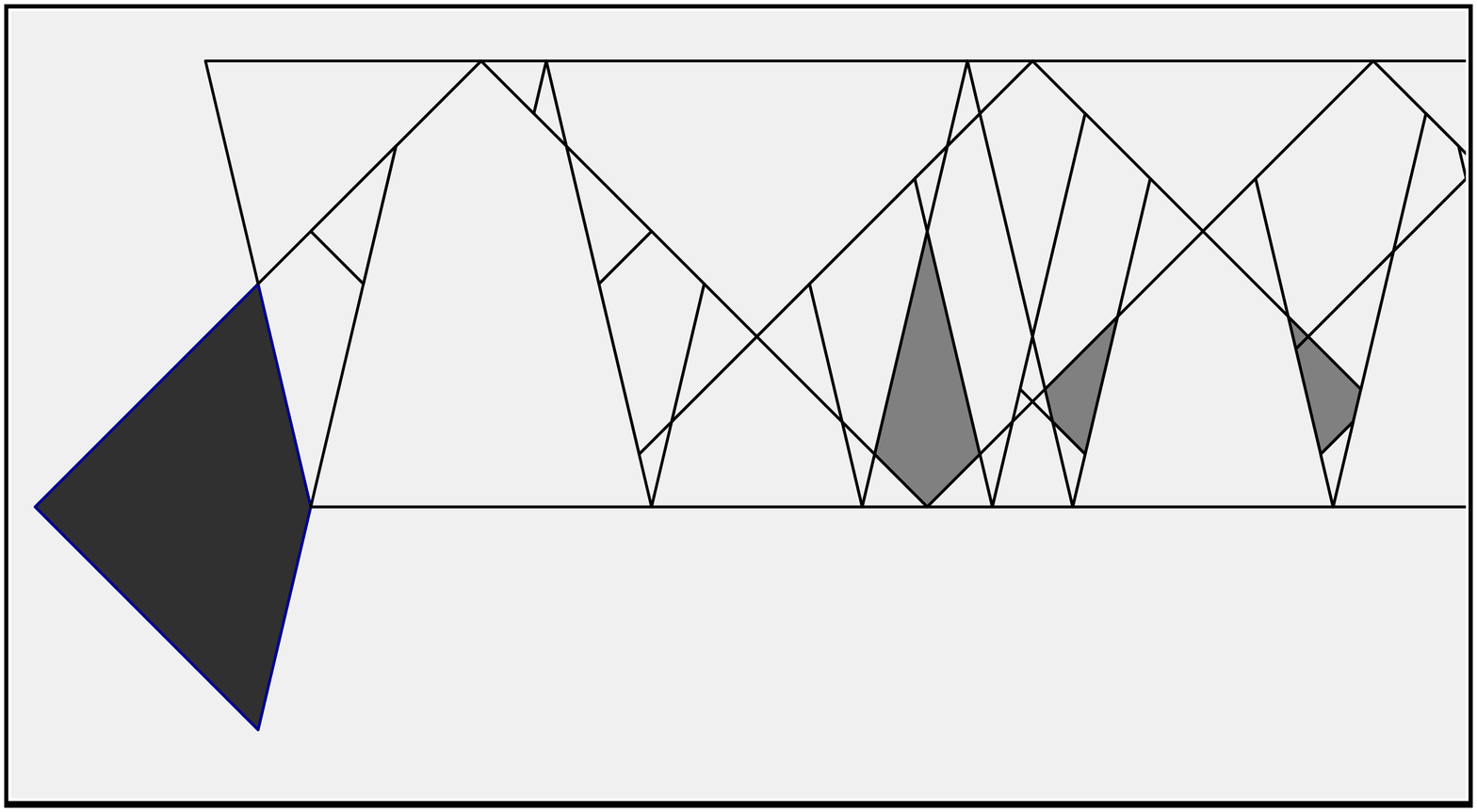}}
\newline
{\bf Figure 12.2:\/} Four of the PSRT's.
\end{center}

It turns out that there are $23$ PSRT's
that intersect $S_{24}$.
The PSRT closest to the origin is $Q_0$.
Figure 12.2 shows the first $4$ PSRT's,
shaded, as well as a
portion of the partition of the top half
of $\Sigma$ into symmetric return tiles.
(The righmost shaded piece is a union of two tiles.)
The reader can see the remaining tiles
using our applet.

We partition the union of $23$ PSRT's
by $848$ polygons, each of which is either
an $A$-atom or a periodic tile.  Since we
are only interested in infinite orbits,
we ignore the periodic tiles.  Let
$\tau$ be an $A$-atom in our partition,
and let $[\tau]$ be the chain containing $\tau$.
Say that $\tau$ is {\it slack\/} if
there is some other $A$-atom $\tau' \subset [\tau]$
such that
\begin{equation}
\label{taut}
\max_{p' \in \tau'} |p'|_x \leq
\max_{p \in \tau} |p|_x-2\phi^{-3}.
\end{equation}
Otherwise, we call $\tau$ {\it taut\/}.
A taut atom is one that is nearly as close
as possible to the origin within its chain.
Note that $\tau'$ need not be in our
partition.  

Ignoring the atoms in $Q_1$, which we do not
need to consider for this lemma, we find that
there are $123$ taut atoms.
With $12$ exceptions, we take each taut atom
$\tau$ and exhibit a $B$-atom $\tau'$ such
that $\tau \to \tau'$ and the pair
$(\tau,\tau')$ satisfies Equation \ref{taut}.
Call the $12$ exceptional atoms {\it super taut\/}.
The $12$ super taut atoms cannot be
brought substantially closer to the origin
by a single renormalization.  

To deal with the $12$ super taut atoms, we
trace through the second renormalization,
as follows.  Let $\tau$ be a super taut atom.
\begin{itemize}
\item We find a $B$-atom $\tau'$ such that
$\tau \to \tau'$.  In practice, we take
$\tau'$ to be the atom in its chain which
is closest to the origin.
\item The atom $\tau'$ is partitioned into
a finite union $\tau_1',...,\tau_k'$ of
$A$-atoms and periodic tiles.  The value $k$ depends on
$\tau$.  The largest value is $27$.
\item For each periodic tile $\tau'_k$ we
exhibit a $B$-atom $\tau''_i$
such that $\tau'_i \to \tau''_i$ and
$(\tau''_i,\tau)$ satisfies Equation \ref{taut}.
Again, we take $\tau''_k$ so that it minimizes
the distance to the origin within its chain.
\end{itemize}
The analysis above handles the points
in $10$ of the $12$ super-taut atoms.

The two exceptional super taut atoms are
$\sigma$ and $\sigma^*$, which we have
deliberately excluded in our hypothesis.
\endproof

\subsection{Proof of the Near Reduction Theorem}

The Near Reduction Theorem concerns 
generic infinite orbits that
intersect the rectangle  $\Sigma_{24}=[-24,24] \times [-2,2]$,
but our results above concern the smaller
region $S_{24}=[0,24] \times [0,2]$.  However,
replacing an orbit $O_1$ by on of its $3$
associate orbits, we can assume that
the orbit intersects $S_{24}$.  Let
$O_1$ be such an orbit.  

Since $O_1$ is a generic infinite orbit, we have
$O_1 \leadsto O_2' \leadsto O_3'...$.
These orbits intersect the strip
$\Sigma=\R \times [-2,2]$ but they
might not intersect the smaller
region $S=\R^+ \times [0,2]$.   However,
replacing our orbits with associates,
if necessary, we get a new sequence
$O_1 \to O_2 \to O_3...$ in which each
orbit intersects $S$.
\newline
\newline
\noindent
{\bf Remark:\/}
We allow ourselves the liberty of passing
to associates mainly for convenience.  We
might have proved the Near Reduction Theorem
for the original sequences of renormalizations,
but it is more computational work.
\newline

Say that an orbit $O$ is {\it passed\/} $Q_1$ if
there is some $p \in O$ such that 
$$|p|_x \leq \min_{q \in Q_1} |q|_x.$$

\begin{lemma}
\label{prenrt}
The Near Reduction Theorem holds for a
generic orbit $O_1$
that is passed $Q_1$.
\end{lemma}

\startproof
Iterating Lemma \ref{near1}, we find that
there is some $m$ such that one of three
things is true.
\begin{enumerate}
\item $O_m$ intersects $S_r$ for $r=2+\phi^{-3}$.
\item $O_m$ intersects $\sigma \cup \sigma^*$, the union of two special $A$-atoms.
\end{enumerate}
In either case, the lemmas in the preceding section
combine to show that $O_{m+k}$ intersects $T^+$ for
some $k \in \{0,1,2,3\}$.
\endproof

Now let's consider the general case when $O_1$ is not
necessarily past $Q_1$.  We can still iterate
Lemma \ref{near1}, but now we arrive at the possibility
that $O_m$ intersects $Q_1$ for some $m$. 
The map $\Psi$ is the identity on the line $y=0$, so $O_m$
does not intersect this line.  Hence, there is some
$\epsilon>0$ such that $O_m$ avoids the $\epsilon$ neighborhood
of the point $(3,0)$. But then there is some $n$
such that $O_m$ intersects $Q_n$ but not $Q_{n+1}$.
Applying the Fixed Point Theorem, we see that $O_{m+1}$
intersects $Q_{n-1}$ but not $Q_{n}$, and
$O_{m+2}$ intersects $Q_{n-2}$ but not $Q_{n-1}$.  
Continuing in this way, we see that
$O_1 \to ... \to O_m \to ... \to O_h$, where
$O_h$ intersects $Q_0$ but not $Q_1$.

If $O_h$ intersects $S_r \cup \sigma \cup \sigma^*$
then we have the same situation as in Lemma \ref{prenrt}.
Otherwise, we can apply Lemma \ref{near1} one last
time to produce an orbit $O_{h+k}$ which is passed
$Q_1$.  This works thanks to the following three facts.
\begin{itemize}
\item $O_h$ intersects $Q_0$ in an interior point (because the orbit is well-defined).
\item $|O_{h+k}|_x \leq |O_h|_x-2\phi^{-3}$.
\item $Q_0$ has width $2\phi^{-3}$.
\end{itemize}
We finish the proof by applying Lemma \ref{prenrt} to $O_{h+k}$.
\newline

For each orbit $O_1$, the Near Reduction Theorem
produces some integer $m=m(O_1)$ such that
$m$ such that $O_1 \to ... \to O_m$, and $O_m$ intersects
$T^+$.  There is no uniform $m$ which works for all orbits.
The difficulty is that orbits that come sufficiently near
$(3,0)$ do not intersect $T^+$, and iterated renormalization
only gradually moves such orbits away from $(3,0)$,
as a consequence of the Fixed Point Theorem.  On the
other hand, our analysis above yields the second statement
of the Near Reduction Theorem:

\begin{lemma}
\label{eff0}
Let $y \in (0,2)$ be any point such that $\{R^n(y)\}$ does
not contain $0$.  Then there is some $N$ with the
following property.  If $O_1$ contains a point within
$1/N$ of the segment $[-24,24] \times \{y\}$, then
$O_1 \to ... \to O_m$ and $O_m$ intersects $T^+$ and
$m<N$.
\end{lemma}

\startproof
Let's consider iteratively applying Lemma \ref{near1} to $O_1$.
After we apply Lemma \ref{near1} at most $24/\phi^{-3}$ times,
the resulting orbit is either passed $Q_1$ or else intersects
$Q_1$.  In the latter case, there is some $\epsilon$ such
that every point of the resulting orbit is at least 
$\epsilon$ from $(3,0)$.  The number $\epsilon$ works
uniformly as long as $O_1$ contains a point sufficiently
close to our line segment.  But then, an additional
$C \log(\epsilon)$ renormalizations produce an orbit
that intesects $T^+$. Here $C$ is some uniform constant.
\endproof

\newpage

\section{Proof of the Far Reduction Theorem}

\subsection{Deviation Estimates}

Let $\Theta: \Sigma \to \widehat \Sigma$ be as in the
Compactification Theorem.
Recall that an $A$-{\it core\/} is a finite portion
$\alpha=\{p_0,...,p_n\}$ of a $\Psi$-orbit with the following
properties.
\begin{enumerate}
\item $\Theta(p_0) \in \widehat A$, the first of the renormalization sets.
\item $\Theta(p_1),...,\Theta(p_n)$ do not lie in $\widehat A$.
\item $\Theta(p_{n+1}) \in \widehat A$.
\end{enumerate}
Up to translation, there are only finitely many $A$-cores, one
per $\widehat A$-chain.  To the
$\widehat A$-core $\alpha$ we associate the
{\it deviation interval\/}
\begin{equation}
\langle \alpha \rangle=\left[\min_k \pi_1(p_k-p_0),\max_k \pi_1(p_k-p_0)\right].
\end{equation}
Here $\pi_1$ is projection onto the first coordinate.
That is, we translate the $A$-core so that its initial point
is at the origin, and we measure how far the sequence deviates
to the left and right of the origin.

The sequence of points $\{p_k-p_0\}$ associated to $\alpha$
is combinatorial in nature:  In view of
Equation \ref{pregraph}, there is a finite sequence 
$\{(a_k,b_k)\}$ such that 
\begin{equation}
\pi_1(p_k-p_0)=2Aa_n+2b_k; \hskip 30 pt k=1,...,n.
\end{equation}
Indeed, we can form a connected lattice path by connecting
the consecutive integers
$(0,0)=(a_0,b_0),...,(a_n,b_n)$.
We call this path an $A$-{\it strand\/}. An
$A$-strand is nothing more than an arc of the arithmetic
graph, translated to that the initial point is the origin.
The $A$-strands are not essential for our analysis, but
we introduce them because they help us visualize what
is going on.  

For the purposes of drawing
pictures, we will continue the strands to the point
$(a_{n+1},b_{n+1})$, so that the endpoints of
the strand correspond to points in $\widehat A$.
In this way, the arithmetic graph component
corresponding to an infinite orbit is partitioned
into uniformly short arcs, with the endpoints
being special points that map into $\widehat A$.
The distribution of these special
{\it marker points\/} controls the coarse geometry
of the arithmetic graph, and thereby the coarse
geometry of the orbit itself.

We define $B$-cores and $B$-strands in the same way, except
that in the type-2 case, the last point of the $B$-core
maps into $\widehat B$. The Renormalization Theorem
gives us a $3$-to-$1$- map from the set of $A$-cores
to the set of $B$-cores. As in the Far Reduction Theorem,
the notation $\alpha \leadsto \beta$ indicates that
the $B$-core $\beta$ is the image of the $A$-core
$\alpha$ under this map.

Let $N_k(S)$ denote the $k$-tubular neighborhood of
a set $S$.  Let $\phi^{-3}I$ denote the interval
obtained by scaling $I$ down by a factor of
$\phi^{-3}$ about the origin.

\begin{lemma}[Trivial]
There is a universal constant $U$ with the following property.
Suppose that $\alpha \to \beta$.
\begin{enumerate}
\item If $\alpha$ has type 1 then $\langle \alpha \rangle \subset
N_U(\langle \beta \rangle)$.
\item If $\alpha$ has type 2 then $\phi^{-3} \langle \alpha \rangle \subset
N_U(\langle \beta \rangle)$.
\item If $\alpha$ has type 1 then $\langle \beta \rangle \subset
N_U(\langle \alpha \rangle)$.
\item If $\alpha$ has type 2 then $\langle \beta \rangle \subset
N_U(\phi^{-3}\langle \alpha \rangle)$.
\end{enumerate}
\end{lemma}

\startproof
Up to translation there are only finitely many cores.
Hence, all deviation intervals are contained in a
single compact subset of $\R$.
\endproof

\noindent
{\bf Remark:\/}
There is nothing special about $\phi^{-3}$ in the
Trivial Lemma.  Were we to pick any other constant, we
would get a similar result.
However, we place $\phi^{-3}$ in this result for the
sake of comparison with the more effective result that
we establish below.   For the more effective result,
the choice of $\phi^{-3}$ is the natural choice.
\newline

Items 3 and 4 serve our needs perfectly. They correspond
to the constant $C$ in the Far Reduction Theorem that we
don't care about.  
However, we do need a more effective version of
Items 1 and 2.

\begin{lemma}[Deviation]
Suppose that $\alpha \to \beta$.
\begin{enumerate}
\item If $\alpha$ has type 1 then $\langle \alpha \rangle \subset
N_r(\langle \beta \rangle)$.  Here $r=-20+24 \phi<12-2\phi$.
\item If $\alpha$ has type 2 then $\phi^{-3} \langle \alpha \rangle \subset
N_5(\langle \beta \rangle)$.
\end{enumerate}
\end{lemma}

\startproof
We simply enumerate all the integer sequences and check that
the bounds hold.  We use golden arithmetic to make sure that
everything is computed correctly.
\endproof

\noindent
{\bf Remark:\/} We could get similar improvements
on Items 3 and 4 of the Trivial Lemma, but
we don't care about those cases.
\newline

\subsection{An Estimate for the Markers}

The Trivial Lemma and the Deviation Lemma reduce the
Far Reduction Theorem to a statement that just
concerns points that map into $\widehat A$ and
$\widehat B$, as we now explain.  Say that an
$A$-{\it marker\/} is a point $p \in \Sigma$
such that $\Theta(p) \in \widehat A$.  We define
the type of the point just as we defined the
type of an $A$-core:  It depends on which
layer of $\widehat A$
contains the image.   

We write $a \leadsto b$ if $a$ is an $A$-marker and
$b$ is a $B$-marker, and 
\begin{equation}
\label{entail}
\widehat R \circ \Theta(a)=b.
\end{equation}
By definition $a \leadsto b$ if and only if
$\alpha \leadsto \beta$, where $\alpha$ is
the $A$-chain containing $a$ and $\beta$ is
the $B$-chain containing $b$.  Below we
will prove the following result.

\begin{lemma}[Marker]
Suppose that $a \leadsto b$.
\begin{enumerate}
\item If $a$ has type 1 then $|b|_x \leq |a|_x+2\phi$.
\item If $a$ has type 2 then $|b|_x<\phi^{-3}|a|_x+10$.
\end{enumerate}
\end{lemma}

The lower bounds in the Far Reduction Theorem are
immediate consequences of the triangle inequality,
the Marker Lemma, and Items 3 and 4 of the
Deviation Lemma.  Let's establish the upper 
bound $|\beta|_x<|\alpha|_x+10$ in the type-$1$ case.
Let $a$ and $b$ respectively be the markers
in $\alpha$ and $\beta$.
Choose the point $a' \in \alpha$ that realizes
$|\alpha|_x$.  Set $r=12-2\phi$. We have
\begin{equation}
|a'|_x-|a|_x \in \langle \alpha \rangle \subset N_r(\langle \beta \rangle).
\end{equation}
Therefore, there is some $b' \in \beta$ such that
\begin{equation}
|b'|_x-|b|_x \leq |a'|_x-|a|_x + (12-2\phi).
\end{equation}
By the Marker Lemma,
\begin{equation}
|b'|_x-|a|_x < |b'|_x-|b|_x+2\phi.
\end{equation}
Therefore 
\begin{equation}
|b'|_x-|a|_x < |a'|_x-|a|_x + 12.
\end{equation}
Adding $|a|_x$ to both sides, we get
\begin{equation}
|\beta|_x \leq |b'|_x < |a'|_x+12=|\alpha|_x+12,
\end{equation}
as desired.  The proof of the upper bound
in the type-2 case works the same way.

The remainder of the chapter is devoted to
proving the Marker Lemma.  The type-1 case
turns out to be fairly trivial, and the
type-2 case has an arithmetical feel to it.

\subsection{The Type-$1$ Case}

Let us recall the information given in \S \ref{structural}.
First of all, we have
\begin{equation}
\Theta(x,y)=\bigg(1,\frac{1}{2},0\bigg)+\bigg(\frac{x}{\phi},\frac{x-y}{2},y\bigg)
\end{equation}
The image of $\Theta$ is a certain fundamental domain of the form $F \times [0,2]$,
where $F$ is described in \S \ref{partitionpix}.  For our purposes, we can
deal with the simpler region
\begin{equation}
\Omega=[1/2,3/2] \times [0,2] \times [0,2] \subset F \times [0,2]
\end{equation}
because (by direct inspection) this region contains both
$\widehat B$ and $\widehat A$.

In the type 1 case, namely in the slabs corresponding to the
partition intervals $I_2$ and $I_4$, the formula for
$\widehat R$ is 
\begin{equation}
\widehat R(x,y,z)=(x,y,z \pm 2\phi^{-1}).
\end{equation}
In this case, we can restrict $\widehat R$
to a smaller domain $S \subset \Omega$
such that $\widehat R(S) \subset \Omega$ as well.

Given the explicit equation for
$\widehat R$, we get the following result.
If $a=(x,y)$ and $a \to b$, then
\begin{equation}
\label{ab1}
b=(x \mp 2\phi,y \pm 2 \phi^{-1} \pm 2).
\end{equation}
What makes this equation work is that
$2\phi-2\phi^{-1} \equiv 0$ mod $2\Z$.
Once we see that this formula works,
it must be the uniquely correct formula,
because $\Theta$ is an injective map.

Equation \ref{ab1} immediatly implies the type-1 case
of the Marker Lemma.

\subsection{The Type-2 Case}

The type-2 case involves the
$7$ maps $\widehat R_{ijk}$, for
various choices of $i,j,k$.
These maps are listed in \S \ref{renormstruct}.
Define
\begin{equation}
S_{ijk}=\widehat R_{ijk}^{-1}(\Omega).
\end{equation}
Here $S_{ijk}$ is a small rectangular solid
that either contains or abuts the fixed
point $(i,j,k)$ of $\widehat R_{ijk}$.
In the Renormalization Theorem, the map
$\widehat R_{ijk}$ is only applied to 
polyhedra that are entirely contained
in $S_{ijk}$.

Since $\widehat R_{ijk}$ expands distances
by a factor of $\phi^3$, we see that the
sidelength of $S_{ijk}$ in the
$x$-direction is $\phi^{-3}$ and the
sidelengths in the other directions
are $2\phi^{-3}$.

For the sake of notational convenience,
we will give our proof in the case
of $\widehat R_{111}$.  The other cases
have essentially the same proof.

Let $p=(p_1,p_2)$ be some point.
Let $\Theta_k(p)$ denote the $k$th coordinate
of $\Theta(p)$.   We first observe that
$\Theta_1(p)=1$ if and only if
$p_1=2N \phi$ for some integer $N$.
So, we can always write 
\begin{equation}
p_1=2N_p\phi+\epsilon_p,
\end{equation}
where $\epsilon_p$ is an error term of size
at most $\phi$. 

When $p=a$, we have a better estimate.
Since $\Theta_1(a)$ lies within
$\phi^{-3}/2$ of $1$, the error term
$\epsilon_a$ is at most $\phi^{-2}/2$.
We also observe that $a_2$ either lies
within $\phi^{-3}$ of $1$ or within
$\phi^{-3}$ of $-1$.   We will treat the
former case.  The latter case has the same
treatment.  Summarizing the discussion, we have
\begin{equation}
a=(2N_a\phi+\epsilon_a,1+\delta_a); \hskip 30 pt
|\epsilon_a| \leq \phi^{-2}/2; \hskip 30 pt
|\delta_a| \leq \phi^{-3}.
\end{equation}
From the equation
\begin{equation}
\label{repel}
\Theta_k(b)-1=\Big(\Theta_k(a)-1\Big) \times \phi^3,
\end{equation}
we get
\begin{equation}
b=(2N_b \phi + \phi^3 \epsilon_a,1+\phi^3 \delta_a).
\end{equation}

In the formula for $\Theta$, we take the listed
coordinates and suitably translate them by
even integers until these coordinates all lie in $[0,2]$.
Hence, there is an even integer $M_a'$ such that
$$\Theta_2(a)=N_a\phi+M'_a+\epsilon_a/2-\delta_a/2 \in [0,2].$$
We write $M'_a=1+M_a$, so that
\begin{equation}
\label{coeff2}
\Theta_2(a)=1+N_a\phi+M_a+\epsilon_a/2-\delta_a/2
\end{equation}
Similarly,
\begin{equation}
\label{coeff22}
\Theta_2(b)=1+N_b\phi+M_b+\phi^3(\epsilon_b/2-\delta_b/2).
\end{equation}
Given our bounds on $\epsilon_a$ and $\delta_a$,
we have
\begin{equation}
\label{bound}
M_a \in [-N_a \phi - 2, -N_a \phi + 2].
\end{equation}

Combining the $k=2$ case of Equation \ref{repel} 
with Equations \ref{coeff2}
and \ref{coeff22}, we get
\begin{equation}
N_b\phi+M_b=\phi^3(N_a \phi + M_a).
\end{equation}
Expanding out the right hand side of this last equation, we get
$$
N_b \phi + M_b=
\phi^4 N_a + \phi^3 M_a =$$
$$
(2+3\phi) N_a + (1+2\phi)M_b=$$
\begin{equation}
(3N_a+2M_b)\phi+(2N_a+N_b).
\end{equation}
Equating coefficients, we get
\begin{equation}
N_b=3N_a+2M_a.
\end{equation}
Combining this last equation with Equation \ref{bound}, we find that
$$
N_b \in [(3-2\phi)N_a-2,(3-2\phi)N_a+2] = [\phi^{-3}N_a-2,\phi^{-3}N_a+2].
$$
In short,
\begin{equation}
\label{bound2}
\phi^{-3}|N_a| \leq |N_b|+2.
\end{equation}
Multiplying through by $2\phi$, we
get
\begin{equation}
\phi^{-3} \times \bigg(2\phi|N_a|\bigg) \leq 2 \phi|N_b|+4\phi.
\end{equation}
Finally, we mention that
\begin{equation}
\label{bound3}
\bigg{|}|a|_x-2\phi|N_a|\bigg{|}=\epsilon_a \leq \phi^{-2}/2; \hskip 30 pt
\bigg{|}|b|_x-2\phi|N_b|\bigg{|}=\epsilon_b \leq \phi/2; \hskip 30 pt
\end{equation}
By the triangle inequality, we get
\begin{equation}
b_x \leq \phi^{-3} a_x + 4\phi + \phi/2 + \phi^{-3} \times (\phi^{-2}/2)
<10.
\end{equation}
This completes the proof of the type-$2$ case of the Marker Lemma.

\newpage

\section{Coordinates}
\label{appendix}

\subsection{The Polyhedron Exchange Map}

\subsubsection{Notation and Conventions}

The polyhedron exchange map is defined in terms of
$64$ polyhedra $P_0,...,P_{63}$.  The polyhedra
$P_0,...,P_{21}$ are what we call {\it inactive\/},
in the sense that $\widehat \Psi$ acts as the identity on
these polyhedra.   The remaining polyhedra are what
we call {\it active\/}.  The $42$ active polyhedra are such that
\begin{equation}
\label{involute}
P_{j+21}=I(P_j); \hskip 30 pt j=22,...,42.
\end{equation}
Here $I$ is the involution
\begin{equation}
I(x,y,z)=(2,2,2)-(x,y,z).
\end{equation}
The map $I$ commutes with $\widehat \Psi$.  Thus $\Psi$
has the same action both the lefthand and the righthand polyhedra
listed in Equation \ref{involute}.  Accordingly, to save space,
we will just list $P_1,...,P_{42}$.

Our notation is such that
\begin{equation}
\left[\matrix{a_0 & a_1 \cr b_0& b_1 \cr c_0 & c_1}\right]=
(a_0+a_1 \phi,b_0+b_1\phi,c_0+c_1\phi).
\end{equation}
We list each polyhedron by its vertices.  The
vertices are not given in any special order.

\scriptsize

\subsubsection{The Inactive Polyhedra}

\begin{eqnarray}
\nonumber
P0=\left[\matrix{3&-1\cr 0&1\cr 2&-1}\right]\left[\matrix{-3&3\cr 2&0\cr 2&-1}\right]\left[\matrix{0&1\cr 2&0\cr 2&0}\right]\left[\matrix{1&0\cr 1&0\cr 2&0}\right]\left[\matrix{3&-1\cr 0&0\cr 2&0}\right]\left[\matrix{0&1\cr 0&0\cr 2&0}\right]\\ \nonumber \left[\matrix{-3&3\cr 2&-1\cr 2&0}\right]\left[\matrix{-1&2\cr 2&0\cr 2&0}\right]\left[\matrix{0&1\cr 2&0\cr 0&0}\right]
\end{eqnarray}

\begin{eqnarray}
\nonumber
P1=\left[\matrix{-1&1\cr 2&0\cr 2&-1}\right]\left[\matrix{5&-3\cr 0&1\cr 2&-1}\right]\left[\matrix{2&-1\cr 0&0\cr 2&0}\right]\left[\matrix{1&0\cr 1&0\cr 2&0}\right]\left[\matrix{-1&1\cr 2&0\cr 2&0}\right]\left[\matrix{2&-1\cr 2&0\cr 2&0}\right]\\ \nonumber \left[\matrix{5&-3\cr 0&1\cr 2&0}\right]\left[\matrix{3&-2\cr 0&0\cr 2&0}\right]\left[\matrix{2&-1\cr 2&0\cr 0&0}\right]
\end{eqnarray}

\begin{eqnarray}
\nonumber
P2=\left[\matrix{3&-1\cr 0&0\cr 0&1}\right]\left[\matrix{-3&3\cr 2&-1\cr 0&1}\right]\left[\matrix{0&1\cr 0&0\cr 2&0}\right]\left[\matrix{0&1\cr 2&0\cr 0&0}\right]\left[\matrix{1&0\cr 1&0\cr 0&0}\right]\left[\matrix{3&-1\cr 0&0\cr 0&0}\right]\\ \nonumber \left[\matrix{0&1\cr 0&0\cr 0&0}\right]\left[\matrix{-3&3\cr 2&-1\cr 0&0}\right]\left[\matrix{-1&2\cr 2&0\cr 0&0}\right]
\end{eqnarray}

\begin{eqnarray}
\nonumber
P3=\left[\matrix{-1&1\cr 2&-1\cr 0&1}\right]\left[\matrix{5&-3\cr 0&0\cr 0&1}\right]\left[\matrix{2&-1\cr 0&0\cr 2&0}\right]\left[\matrix{2&-1\cr 0&0\cr 0&0}\right]\left[\matrix{1&0\cr 1&0\cr 0&0}\right]\left[\matrix{-1&1\cr 2&0\cr 0&0}\right]\\ \nonumber \left[\matrix{2&-1\cr 2&0\cr 0&0}\right]\left[\matrix{5&-3\cr 0&1\cr 0&0}\right]\left[\matrix{3&-2\cr 0&0\cr 0&0}\right]
\end{eqnarray}

\begin{eqnarray}
\nonumber
P4=\left[\matrix{1&0\cr 0&0\cr 1&0}\right]\left[\matrix{6&-3\cr -1&1\cr 0&0}\right]\left[\matrix{1&0\cr 0&0\cr 0&0}\right]\left[\matrix{2&-1\cr 0&0\cr 0&0}\right]\left[\matrix{1&0\cr 1&0\cr 0&0}\right]
\end{eqnarray}

\begin{eqnarray}
\nonumber
P5=\left[\matrix{1&0\cr 2&0\cr 1&0}\right]\left[\matrix{-4&3\cr 3&-1\cr 2&0}\right]\left[\matrix{1&0\cr 2&0\cr 2&0}\right]\left[\matrix{0&1\cr 2&0\cr 2&0}\right]\left[\matrix{1&0\cr 1&0\cr 2&0}\right]
\end{eqnarray}

\begin{eqnarray}
\nonumber
P6=\left[\matrix{1&0\cr 1&0\cr 1&0}\right]\left[\matrix{-4&3\cr 3&-1\cr 0&0}\right]\left[\matrix{1&0\cr 2&0\cr 0&0}\right]\left[\matrix{0&1\cr 2&0\cr 0&0}\right]\left[\matrix{1&0\cr 1&0\cr 0&0}\right]
\end{eqnarray}

\begin{eqnarray}
\nonumber
P7=\left[\matrix{1&0\cr 1&0\cr 1&0}\right]\left[\matrix{6&-3\cr -1&1\cr 2&0}\right]\left[\matrix{1&0\cr 0&0\cr 2&0}\right]\left[\matrix{2&-1\cr 0&0\cr 2&0}\right]\left[\matrix{1&0\cr 1&0\cr 2&0}\right]
\end{eqnarray}

\begin{eqnarray}
\nonumber
P8=\left[\matrix{6&-3\cr 2&0\cr -1&1}\right]\left[\matrix{-4&3\cr 3&-1\cr -1&1}\right]\left[\matrix{1&0\cr 1&0\cr 1&0}\right]\left[\matrix{1&0\cr 2&0\cr 1&0}\right]\left[\matrix{6&-3\cr -2&2\cr 3&-1}\right]\left[\matrix{-4&3\cr 3&-1\cr 3&-1}\right]\\ \nonumber \left[\matrix{1&0\cr 1&0\cr 2&0}\right]\left[\matrix{1&0\cr 2&0\cr 0&0}\right]
\end{eqnarray}

\begin{eqnarray}
\nonumber
P9=\left[\matrix{-4&3\cr 4&-2\cr -1&1}\right]\left[\matrix{6&-3\cr -1&1\cr -1&1}\right]\left[\matrix{1&0\cr 0&0\cr 1&0}\right]\left[\matrix{1&0\cr 1&0\cr 1&0}\right]\left[\matrix{6&-3\cr -1&1\cr 3&-1}\right]\left[\matrix{-4&3\cr 0&0\cr 3&-1}\right]\\ \nonumber \left[\matrix{1&0\cr 0&0\cr 2&0}\right]\left[\matrix{1&0\cr 1&0\cr 0&0}\right]
\end{eqnarray}

\begin{eqnarray}
\nonumber
P10=\left[\matrix{6&-3\cr 0&0\cr -1&1}\right]\left[\matrix{1&0\cr 0&0\cr 0&0}\right]\left[\matrix{6&-3\cr 0&0\cr 0&0}\right]\left[\matrix{11&-6\cr -3&2\cr 0&0}\right]\left[\matrix{6&-3\cr -1&1\cr 0&0}\right]
\end{eqnarray}

\begin{eqnarray}
\nonumber
P11=\left[\matrix{-4&3\cr 3&-1\cr -1&1}\right]\left[\matrix{1&0\cr 2&0\cr 0&0}\right]\left[\matrix{-4&3\cr 2&0\cr 0&0}\right]\left[\matrix{-9&6\cr 5&-2\cr 0&0}\right]\left[\matrix{-4&3\cr 3&-1\cr 0&0}\right]
\end{eqnarray}

\begin{eqnarray}
\nonumber
P12=\left[\matrix{11&-6\cr 0&0\cr -3&2}\right]\left[\matrix{6&-3\cr 0&0\cr 0&0}\right]\left[\matrix{11&-6\cr -3&2\cr 0&0}\right]\left[\matrix{3&-1\cr 0&0\cr 0&0}\right]
\end{eqnarray}

\begin{eqnarray}
\nonumber
P13=\left[\matrix{-9&6\cr 5&-2\cr -3&2}\right]\left[\matrix{-4&3\cr 2&0\cr 0&0}\right]\left[\matrix{-9&6\cr 5&-2\cr 0&0}\right]\left[\matrix{-1&1\cr 2&0\cr 0&0}\right]
\end{eqnarray}

\begin{eqnarray}
\nonumber
P14=\left[\matrix{5&-2\cr 0&0\cr -3&2}\right]\left[\matrix{0&1\cr 0&0\cr 0&0}\right]\left[\matrix{-3&3\cr 2&-1\cr 0&0}\right]\left[\matrix{5&-2\cr 0&0\cr 0&0}\right]
\end{eqnarray}

\begin{eqnarray}
\nonumber
P15=\left[\matrix{-3&2\cr 5&-2\cr -3&2}\right]\left[\matrix{2&-1\cr 2&0\cr 0&0}\right]\left[\matrix{5&-3\cr 0&1\cr 0&0}\right]\left[\matrix{-3&2\cr 2&0\cr 0&0}\right]
\end{eqnarray}

\begin{eqnarray}
\nonumber
P16=\left[\matrix{1&0\cr 0&0\cr 2&0}\right]\left[\matrix{6&-3\cr 0&0\cr 2&0}\right]\left[\matrix{11&-6\cr -3&2\cr 2&0}\right]\left[\matrix{6&-3\cr -1&1\cr 2&0}\right]\left[\matrix{6&-3\cr -1&1\cr 3&-1}\right]
\end{eqnarray}

\begin{eqnarray}
\nonumber
P17=\left[\matrix{1&0\cr 2&0\cr 2&0}\right]\left[\matrix{-4&3\cr 2&0\cr 2&0}\right]\left[\matrix{-9&6\cr 5&-2\cr 2&0}\right]\left[\matrix{-4&3\cr 3&-1\cr 2&0}\right]\left[\matrix{-4&3\cr 2&0\cr 3&-1}\right]
\end{eqnarray}

\begin{eqnarray}
\nonumber
P18=\left[\matrix{6&-3\cr 0&0\cr 2&0}\right]\left[\matrix{11&-6\cr -3&2\cr 2&0}\right]\left[\matrix{3&-1\cr 0&0\cr 2&0}\right]\left[\matrix{11&-6\cr -3&2\cr 5&-2}\right]
\end{eqnarray}

\begin{eqnarray}
\nonumber
P19=\left[\matrix{-4&3\cr 2&0\cr 2&0}\right]\left[\matrix{-9&6\cr 5&-2\cr 2&0}\right]\left[\matrix{-1&1\cr 2&0\cr 2&0}\right]\left[\matrix{-9&6\cr 2&0\cr 5&-2}\right]
\end{eqnarray}

\begin{eqnarray}
\nonumber
P20=\left[\matrix{0&1\cr 0&0\cr 2&0}\right]\left[\matrix{-3&3\cr 2&-1\cr 2&0}\right]\left[\matrix{5&-2\cr 0&0\cr 2&0}\right]\left[\matrix{5&-2\cr -3&2\cr 5&-2}\right]
\end{eqnarray}

\begin{eqnarray}
\nonumber
P21=\left[\matrix{2&-1\cr 2&0\cr 2&0}\right]\left[\matrix{5&-3\cr 0&1\cr 2&0}\right]\left[\matrix{-3&2\cr 2&0\cr 2&0}\right]\left[\matrix{-3&2\cr 2&0\cr 5&-2}\right]
\end{eqnarray}

\subsubsection{Half the Active Polyhedra}

\begin{eqnarray}
\nonumber
P22=\left[\matrix{-4&3\cr 3&-1\cr -1&1}\right]\left[\matrix{1&0\cr 1&0\cr 1&0}\right]\left[\matrix{1&0\cr 2&0\cr 0&0}\right]\left[\matrix{-4&3\cr 3&-1\cr 0&0}\right]
\end{eqnarray}

\begin{eqnarray}
\nonumber
P23=\left[\matrix{5&-3\cr 0&1\cr 2&-1}\right]\left[\matrix{3&-2\cr 0&0\cr 1&0}\right]\left[\matrix{3&-2\cr 0&0\cr 2&0}\right]\left[\matrix{2&-1\cr 0&0\cr 2&0}\right]\left[\matrix{2&-1\cr 2&0\cr 0&0}\right]\left[\matrix{5&-3\cr 0&1\cr 0&0}\right]\\ \nonumber 
\end{eqnarray}

\begin{eqnarray}
\nonumber
P24=\left[\matrix{3&-1\cr 0&1\cr 2&-1}\right]\left[\matrix{-5&4\cr 0&0\cr 5&-2}\right]\left[\matrix{3&-1\cr 0&0\cr 2&0}\right]\left[\matrix{0&1\cr 0&0\cr 2&0}\right]\left[\matrix{-5&4\cr 5&-2\cr 0&0}\right]\left[\matrix{0&1\cr 2&0\cr 0&0}\right]\\ \nonumber 
\end{eqnarray}

\begin{eqnarray}
\nonumber
P25=\left[\matrix{-1&1\cr 0&0\cr 2&-1}\right]\left[\matrix{1&0\cr 1&0\cr 1&0}\right]\left[\matrix{1&0\cr 0&0\cr 2&0}\right]\left[\matrix{2&-1\cr 0&0\cr 2&0}\right]\left[\matrix{2&-1\cr 0&0\cr 0&0}\right]\left[\matrix{-1&1\cr 2&-1\cr 0&0}\right]\\ \nonumber 
\end{eqnarray}

\begin{eqnarray}
\nonumber
P26=\left[\matrix{3&-1\cr 0&0\cr 2&-1}\right]\left[\matrix{11&-6\cr -3&2\cr 5&-2}\right]\left[\matrix{3&-1\cr 0&0\cr 2&0}\right]\left[\matrix{6&-3\cr 0&0\cr 2&0}\right]\left[\matrix{11&-6\cr -3&2\cr 0&0}\right]\left[\matrix{6&-3\cr 0&0\cr 0&0}\right]\\ \nonumber 
\end{eqnarray}

\begin{eqnarray}
\nonumber
P27=\left[\matrix{8&-5\cr 0&0\cr -1&1}\right]\left[\matrix{3&-2\cr 0&0\cr 1&0}\right]\left[\matrix{3&-2\cr 0&0\cr 0&0}\right]\left[\matrix{8&-5\cr -1&1\cr 0&0}\right]
\end{eqnarray}

\begin{eqnarray}
\nonumber
P28=\left[\matrix{5&-2\cr -3&2\cr 5&-2}\right]\left[\matrix{5&-2\cr 0&0\cr 2&0}\right]\left[\matrix{0&1\cr 0&0\cr 2&0}\right]\left[\matrix{0&1\cr 0&0\cr 3&-1}\right]
\end{eqnarray}

\begin{eqnarray}
\nonumber
P29=\left[\matrix{1&0\cr 0&0\cr 5&-2}\right]\left[\matrix{1&0\cr 0&0\cr 2&0}\right]\left[\matrix{6&-3\cr 0&0\cr 2&0}\right]\left[\matrix{6&-3\cr -1&1\cr 3&-1}\right]
\end{eqnarray}

\begin{eqnarray}
\nonumber
P30=\left[\matrix{-9&6\cr 7&-4\cr -3&2}\right]\left[\matrix{-1&1\cr 0&0\cr 2&-1}\right]\left[\matrix{-4&3\cr 0&0\cr 4&-2}\right]\left[\matrix{-9&6\cr 0&0\cr 4&-2}\right]\left[\matrix{-1&1\cr 2&-1\cr 0&0}\right]\left[\matrix{-4&3\cr 4&-2\cr 0&0}\right]\\ \nonumber 
\end{eqnarray}

\begin{eqnarray}
\nonumber
P31=\left[\matrix{5&-2\cr 0&0\cr -3&2}\right]\left[\matrix{-3&3\cr 2&-1\cr 2&-1}\right]\left[\matrix{0&1\cr 0&0\cr 4&-2}\right]\left[\matrix{5&-2\cr 0&0\cr 4&-2}\right]\left[\matrix{-3&3\cr 2&-1\cr 0&0}\right]\left[\matrix{0&1\cr 0&0\cr 0&0}\right]\\ \nonumber 
\end{eqnarray}

\begin{eqnarray}
\nonumber
P32=\left[\matrix{-7&5\cr 6&-3\cr 0&1}\right]\left[\matrix{3&-1\cr 0&0\cr 0&1}\right]\left[\matrix{1&0\cr 1&0\cr 5&-2}\right]\left[\matrix{-5&4\cr 0&0\cr 5&-2}\right]\left[\matrix{-7&5\cr 4&-2\cr 2&0}\right]\left[\matrix{3&-1\cr 0&0\cr 2&0}\right]\\ \nonumber \left[\matrix{1&0\cr 1&0\cr 0&0}\right]\left[\matrix{-5&4\cr 5&-2\cr 0&0}\right]
\end{eqnarray}

\begin{eqnarray}
\nonumber
P31=\left[\matrix{5&-2\cr 0&0\cr -3&2}\right]\left[\matrix{-3&3\cr 2&-1\cr 2&-1}\right]\left[\matrix{0&1\cr 0&0\cr 4&-2}\right]\left[\matrix{5&-2\cr 0&0\cr 4&-2}\right]\left[\matrix{-3&3\cr 2&-1\cr 0&0}\right]\left[\matrix{0&1\cr 0&0\cr 0&0}\right]\\ \nonumber 
\end{eqnarray}

\begin{eqnarray}
\nonumber
P32=\left[\matrix{-7&5\cr 6&-3\cr 0&1}\right]\left[\matrix{3&-1\cr 0&0\cr 0&1}\right]\left[\matrix{1&0\cr 1&0\cr 5&-2}\right]\left[\matrix{-5&4\cr 0&0\cr 5&-2}\right]\left[\matrix{-7&5\cr 4&-2\cr 2&0}\right]\left[\matrix{3&-1\cr 0&0\cr 2&0}\right]\\ \nonumber \left[\matrix{1&0\cr 1&0\cr 0&0}\right]\left[\matrix{-5&4\cr 5&-2\cr 0&0}\right]
\end{eqnarray}

\begin{eqnarray}
\nonumber
P33=\left[\matrix{1&0\cr 1&0\cr 5&-2}\right]\left[\matrix{1&0\cr 1&0\cr 2&0}\right]\left[\matrix{-7&5\cr 4&-2\cr 2&0}\right]\left[\matrix{-7&5\cr 6&-3\cr 0&1}\right]
\end{eqnarray}

\begin{eqnarray}
\nonumber
P34=\left[\matrix{-9&6\cr 7&-4\cr -3&2}\right]\left[\matrix{-4&3\cr 4&-2\cr -1&1}\right]\left[\matrix{-4&3\cr 0&0\cr 4&-2}\right]\left[\matrix{-9&6\cr 0&0\cr 4&-2}\right]\left[\matrix{1&0\cr 0&0\cr 1&0}\right]\left[\matrix{-4&3\cr 0&0\cr 3&-1}\right]\\ \nonumber \left[\matrix{1&0\cr 1&0\cr 0&0}\right]\left[\matrix{-4&3\cr 4&-2\cr 0&0}\right]
\end{eqnarray}

\begin{eqnarray}
\nonumber
P35=\left[\matrix{5&-2\cr 0&0\cr 4&-2}\right]\left[\matrix{0&1\cr 0&0\cr 4&-2}\right]\left[\matrix{0&1\cr 0&0\cr 3&-1}\right]\left[\matrix{-3&3\cr 2&-1\cr 2&0}\right]\left[\matrix{5&-2\cr 0&0\cr 2&0}\right]\left[\matrix{-3&3\cr 2&-1\cr 2&-1}\right]\\ \nonumber 
\end{eqnarray}

\begin{eqnarray}
\nonumber
P36=\left[\matrix{-9&6\cr 2&0\cr 5&-2}\right]\left[\matrix{-9&6\cr 5&-2\cr 2&0}\right]\left[\matrix{-1&1\cr 2&0\cr 2&0}\right]\left[\matrix{-1&1\cr 2&0\cr 0&1}\right]
\end{eqnarray}

\begin{eqnarray}
\nonumber
P37=\left[\matrix{1&0\cr 2&0\cr -3&2}\right]\left[\matrix{-4&3\cr 3&-1\cr -1&1}\right]\left[\matrix{-4&3\cr 2&0\cr 4&-2}\right]\left[\matrix{1&0\cr 2&0\cr 4&-2}\right]\left[\matrix{-9&6\cr 5&-2\cr 1&0}\right]\left[\matrix{-4&3\cr 3&-1\cr 3&-1}\right]\\ \nonumber \left[\matrix{-9&6\cr 5&-2\cr 0&0}\right]\left[\matrix{-4&3\cr 2&0\cr 0&0}\right]
\end{eqnarray}

\begin{eqnarray}
\nonumber
P38=\left[\matrix{-4&3\cr 2&0\cr 3&-1}\right]\left[\matrix{-9&6\cr 5&-2\cr 2&0}\right]\left[\matrix{-4&3\cr 3&-1\cr 2&0}\right]\left[\matrix{-9&6\cr 5&-2\cr 1&0}\right]
\end{eqnarray}

\begin{eqnarray}
\nonumber
P39=\left[\matrix{-4&3\cr 3&-1\cr 3&-1}\right]\left[\matrix{-4&3\cr 3&-1\cr 2&0}\right]\left[\matrix{1&0\cr 1&0\cr 2&0}\right]\left[\matrix{1&0\cr 2&0\cr 1&0}\right]
\end{eqnarray}

\begin{eqnarray}
\nonumber
P40=\left[\matrix{-11&8\cr 2&0\cr 4&-2}\right]\left[\matrix{-6&5\cr 2&0\cr 4&-2}\right]\left[\matrix{-6&5\cr 3&-1\cr 3&-1}\right]\left[\matrix{-3&3\cr 2&-1\cr 2&0}\right]\left[\matrix{-11&8\cr 4&-2\cr 2&0}\right]\left[\matrix{-3&3\cr 2&0\cr 2&-1}\right]\\ \nonumber 
\end{eqnarray}

\begin{eqnarray}
\nonumber
P41=\left[\matrix{-1&2\cr 2&0\cr 1&0}\right]\left[\matrix{-6&5\cr 2&0\cr 3&-1}\right]\left[\matrix{-6&5\cr 3&-1\cr 3&-1}\right]\left[\matrix{-11&8\cr 4&-2\cr 2&0}\right]\left[\matrix{-6&5\cr 3&-1\cr 2&0}\right]\left[\matrix{-11&8\cr 2&0\cr 4&-2}\right]\\ \nonumber \left[\matrix{-6&5\cr 2&0\cr 4&-2}\right]
\end{eqnarray}

\begin{eqnarray}
\nonumber
P42=\left[\matrix{-9&6\cr 5&-2\cr 1&0}\right]\left[\matrix{-4&3\cr 3&-1\cr 3&-1}\right]\left[\matrix{-4&3\cr 2&0\cr 3&-1}\right]\left[\matrix{1&0\cr 2&0\cr 2&0}\right]\left[\matrix{-4&3\cr 3&-1\cr 2&0}\right]\left[\matrix{1&0\cr 2&0\cr 4&-2}\right]\\ \nonumber \left[\matrix{-4&3\cr 2&0\cr 4&-2}\right]
\end{eqnarray}

\normalsize

Using our computer program, the reader can survey all these polyhedra.
The program is set up so that the vertices of the polyhedron are
displayed whenever the polyhedron is selected.

\subsection{The Map}

Now we describe the action of the map $\widehat \Psi$ on the polyhedron $P_j$.
To each integer $j$ we assign a pair of integers $N_j$ according to the
following lookup table.
\begin{equation}
\left[\matrix{
0-21 & (0,0) \cr
22-29 & (1,0) \cr
30-32 & (1,1) \cr
33-40 & (0,1) \cr
41-42 & (-1,1) \cr
43-50 & (-1,0) \cr
51-53 & (-1,-1) \cr
54-61 & (0,-1) \cr
62-63 & (1,-1)} \right]
\end{equation}
For instance, if $j=41$ then $N_j=(-1,1)$.
The numbers $N_j=(m,n)$ are such that

\begin{equation}
\label{pregraph2}
\Psi(p)-p=(2m + 2n \phi^{-3},*)
\end{equation}
for any point $p$ such that $\Theta(p) \in P_j$.
As in Equation \ref{pregraph}, the third
coordinate, which lies in $\{-2,0,2\}$ depends on
the parity of $m+n$.  When $m+n=0$ the third coordinate
in Equation \ref{pregraph2} is $0$. Otherwise, it
is $\pm 2$, depending on which value yields a point
in the strip $\Sigma$.

The action of $\widehat \Psi$ on $P_j$ is expressed
in terms of the pair $(m,n)$ as follows.  
\begin{equation}
\widehat \Psi(\widehat p)-\widehat p =
\left[\matrix{
-2m + 10 n & 2m-6n \cr
-2n & 2n \cr
0&0} \right]
\end{equation}
Thus, for instance, for $N_{41}$, the pair $(-1,1)$ yields the vector
$$(12-8\phi,-2+2\phi,0).$$

Notice that $N_{21+j}=-N_j$ for $j=22,...,41$.  This is
consistent with the action of the involution $I$ defined
in the previous section.

Now we discuss the allowable transitions for our map.
We write $a \to b$ is there exists a point
$p$ in the interior of $P_a$ such that $\widehat \Psi(p) \in P_b$.
For $a<22$ we only have $a \to a$.  We also have the general
symmetry
\begin{equation}
a \to b \hskip 40 pt \Longleftrightarrow \hskip 40 pt (a+21) \to (b+21)
\end{equation}
which holds for all $a \geq 22$.
For this reason, we just list the transitions for $a=22,...,42$.

\begin{equation}
\matrix{
22& \to & 61 & 23 & 56   \cr
23& \to & 63 & 58 & 29 & 32 & 26 & 24   \cr
24& \to & 25 & 53 & 37   \cr
25& \to & 31 & 35 & 28 & 27 & 23 & 62 & 61   \cr
26& \to & 25   \cr
27& \to & 63   \cr
28& \to & 25   \cr
29& \to & 23   \cr
30& \to & 23   \cr
31& \to & 32   \cr
32& \to & 25 & 30 & 47 & 36 & 34 & 37 & 42 & 38 & 39   \cr
33& \to & 32   \cr
34& \to & 46 & 26 & 32   \cr
35& \to & 44 & 40 & 41 & 48   \cr
36& \to & 39   \cr
37& \to & 32   \cr
38& \to & 32 & 33   \cr
39& \to & 32   \cr
40& \to & 48 & 23   \cr
41& \to & 43 & 32 & 46   \cr
42& \to & 44 & 40 & 41} 
\end{equation}

\begin{lemma}
\label{transition}
Every nontrivial generic orbit intersects $P_j$ for
one of the indices $j=23,25,32,40,41,44,46,53,61,62$.
\end{lemma}

\startproof
Our list of indices is invariant under the involution
$i \to i+21$, which corresponds to the involution of
the system, which we have mentioned several tiles.
For that reason, it suffices to solve the following
modified problem. We take each number $b>42$ on the
above list and replace it by $b-21$.  We when show
for every $a \in [22,42]$ that $a=a_0 \to ... \to a_k$,
such that  $a_k \in L=\{23,25,32,40,41\}$.  This is equivalent
to the original lemma.

Here is the modified list.

\begin{equation}
\matrix{
22& \to & 40 & 23 & 35   \cr
23& \to & 42 & 37 & 29 & 32 & 26 & 24   \cr
24& \to & 25 & 32 & 37   \cr
25& \to & 31 & 35 & 28 & 27 & 23 & 41 & 40   \cr
26& \to & 25   \cr
27& \to & 42   \cr
28& \to & 25   \cr
29& \to & 23   \cr
30& \to & 23   \cr
31& \to & 32   \cr
32& \to & 25 & 30 & 26 & 36 & 34 & 37 & 42 & 38 & 39   \cr
33& \to & 32   \cr
34& \to & 25 & 26 & 32   \cr
35& \to & 23 & 40 & 41 & 27   \cr
36& \to & 39   \cr
37& \to & 32   \cr
38& \to & 32 & 33   \cr
39& \to & 32   \cr
40& \to & 27 & 23   \cr
41& \to & 22 & 32 & 25   \cr
42& \to & 23 & 40 & 41} 
\end{equation}

Note that $a \to b \in L$ when 
$a$ is one of $26,28,29,30,31,33,37,39,42$.
But $27 \to 42$ and $36 \to 39$. So, we
can eliminate these as well.
Erasing the eliminated cases, and also erasing
the cases corresponding to elements of $L$, we have

\begin{equation}
\matrix{
22& \to & 40 & 23 & 35   \cr
24& \to & 25 & 32 & 37   \cr
34& \to & 25 & 26 & 32   \cr
35& \to & 23 & 40 & 41 & 27   \cr
38& \to & 32 & 33}
\end{equation}

Note that $34 \to 26$ or $36 \to b \in L$.
Since we have eliminated $26$, we
can eliminate $34$.
Since $35 \to 27$ or $35 \to b \in L$, and
we have eliminated $27$, we can eliminate $35$.
Since $22 \to 35$ or $22 \to b \in L$, we
eliminate $22$. 
Since $24 \to 37$ or $24 \to b \in L$, and
we have eliminated $37$, we eliminate $24$.
Since $38 \to 33$ or $38 \to b \in L$,
and we have eliminated $33$, we eliminate $38$.
\endproof

\subsection{The set $\widehat B$}
\label{B}

The set $\widehat B$ is divided into $6$ layers.
These layers correspond to $6$ intervals
$J_1',...,J_6'$, which are just re-indexed versions of the
levels corresponding to $\widehat B$.  That is
\begin{equation}
J_1'=J_1; \hskip 15 pt
J_2'=J_2; \hskip 15 pt
J_3'=J_{31}; \hskip 15 pt
J_4'=J_{32}; \hskip 15 pt
J_5'=J_4; \hskip 15 pt
J_6'=J_5.
\end{equation}
The intervals $J_1',...,J_6'$ are defined by
the partition
\begin{equation}
\label{BBpart}
0< \phi^{-2}<2  \phi^{-2}<1<2-\phi^{-2}<2-2  \phi^{-2}<2.
\end{equation}

Each layer is decomposed into $4$ convex
polyhedra, which we call branches.  The first $3$ layers
are contained in the bottom half of $\widehat \Sigma$,
namely the region
\begin{equation}
\label{bottomhalf}
(\R/2\Z) \times (\R/2\Z) \times [0,1].
\end{equation}
The remaining $3$ layers are the images of
the first $3$ layers under the involution
$I$ discussed above.  For this reason, we just
list the $12$ polyhedra in the first $3$ layers.
Our notation is such that $P_{ij}$ is the
$j$th branch of the $i$th layer.

\scriptsize

\subsubsection{Layer $1$}

\begin{eqnarray}
\nonumber
B11=\left[\matrix{1&0\cr 1&0\cr 0&0}\right]\left[\matrix{-7&5\cr 6&-3\cr 0&0}\right]\left[\matrix{6&-3\cr -2&2\cr 5&-3}\right]\left[\matrix{-7&5\cr 4&-2\cr 2&-1}\right]\left[\matrix{1&0\cr 1&0\cr 2&-1}\right]\left[\matrix{-7&5\cr 6&-3\cr 2&-1}\right]\\ \nonumber \left[\matrix{6&-3\cr 1&0\cr 2&-1}\right]
\end{eqnarray}

\begin{eqnarray}
\nonumber
B12=\left[\matrix{9&-5\cr -4&3\cr 2&-1}\right]\left[\matrix{1&0\cr 1&0\cr 2&-1}\right]\left[\matrix{9&-5\cr -2&2\cr 2&-1}\right]\left[\matrix{-4&3\cr 1&0\cr 2&-1}\right]\left[\matrix{1&0\cr 1&0\cr 0&0}\right]\left[\matrix{-4&3\cr 3&-1\cr 0&0}\right]\\ \nonumber 
\end{eqnarray}

\begin{eqnarray}
\nonumber
B13=\left[\matrix{9&-5\cr 0&1\cr 0&0}\right]\left[\matrix{22&-13\cr -3&3\cr 0&0}\right]\left[\matrix{-4&3\cr 3&-1\cr 5&-3}\right]\left[\matrix{9&-5\cr -10&7\cr 10&-6}\right]\left[\matrix{-25&16\cr 11&-6\cr 10&-6}\right]\left[\matrix{-25&16\cr 11&-6\cr 2&-1}\right]\\ \nonumber \left[\matrix{-4&3\cr 3&-1\cr 2&-1}\right]\left[\matrix{-12&8\cr 3&-1\cr 2&-1}\right]\left[\matrix{22&-13\cr -5&4\cr 2&-1}\right]\left[\matrix{9&-5\cr 0&1\cr 2&-1}\right]
\end{eqnarray}

\begin{eqnarray}
\nonumber
B14=\left[\matrix{1&0\cr 2&0\cr -3&2}\right]\left[\matrix{1&0\cr 2&0\cr 0&0}\right]\left[\matrix{-4&3\cr 2&0\cr 0&0}\right]\left[\matrix{9&-5\cr 0&1\cr 2&-1}\right]\left[\matrix{-4&3\cr 0&1\cr 2&-1}\right]\left[\matrix{-12&8\cr 5&-2\cr 2&-1}\right]\\ \nonumber 
\end{eqnarray}

\subsubsection{Layer $2$}

\begin{eqnarray}
\nonumber
B21=\left[\matrix{-7&5\cr 4&-2\cr 2&-1}\right]\left[\matrix{1&0\cr 1&0\cr 2&-1}\right]\left[\matrix{-7&5\cr 6&-3\cr 2&-1}\right]\left[\matrix{-15&10\cr 9&-5\cr 2&-1}\right]\left[\matrix{-15&10\cr 4&-2\cr 7&-4}\right]\left[\matrix{6&-3\cr -1&1\cr -1&1}\right]\\ \nonumber \left[\matrix{6&-3\cr -1&1\cr 4&-2}\right]\left[\matrix{1&0\cr 1&0\cr 4&-2}\right]
\end{eqnarray}

\begin{eqnarray}
\nonumber
B22=\left[\matrix{1&0\cr 1&0\cr 4&-2}\right]\left[\matrix{9&-5\cr -4&3\cr 4&-2}\right]\left[\matrix{-4&3\cr 4&-2\cr -1&1}\right]\left[\matrix{1&0\cr 1&0\cr 2&-1}\right]\left[\matrix{9&-5\cr -4&3\cr 2&-1}\right]\left[\matrix{-4&3\cr 1&0\cr 2&-1}\right]\\ \nonumber \left[\matrix{9&-5\cr -2&2\cr 2&-1}\right]
\end{eqnarray}

\begin{eqnarray}
\nonumber
B23=\left[\matrix{6&-3\cr -1&1\cr 4&-2}\right]\left[\matrix{-15&10\cr 4&-2\cr 4&-2}\right]\left[\matrix{6&-3\cr 4&-2\cr -1&1}\right]\left[\matrix{-15&10\cr 4&-2\cr 7&-4}\right]
\end{eqnarray}

\begin{eqnarray}
\nonumber
B24=\left[\matrix{-10&7\cr 4&-2\cr 4&-2}\right]\left[\matrix{-15&10\cr 4&-2\cr 4&-2}\right]\left[\matrix{6&-3\cr 4&-2\cr -1&1}\right]\left[\matrix{-15&10\cr 4&-2\cr 7&-4}\right]\left[\matrix{-15&10\cr 9&-5\cr 2&-1}\right]\left[\matrix{6&-3\cr 1&0\cr 2&-1}\right]\\ \nonumber \left[\matrix{-10&7\cr 6&-3\cr 2&-1}\right]\left[\matrix{-15&10\cr 6&-3\cr 2&-1}\right]
\end{eqnarray}

\subsubsection{Layer $3$}

\begin{eqnarray}
\nonumber
B31=\left[\matrix{6&-3\cr -2&2\cr 4&-2}\right]\left[\matrix{1&0\cr 1&0\cr 4&-2}\right]\left[\matrix{1&0\cr 1&0\cr 1&0}\right]\left[\matrix{14&-8\cr -4&3\cr 1&0}\right]\left[\matrix{6&-3\cr 1&0\cr 1&0}\right]
\end{eqnarray}

\begin{eqnarray}
\nonumber
B32=\left[\matrix{1&0\cr 1&0\cr 1&0}\right]\left[\matrix{22&-13\cr -12&8\cr 9&-5}\right]\left[\matrix{1&0\cr 1&0\cr 4&-2}\right]\left[\matrix{22&-13\cr -7&5\cr 4&-2}\right]
\end{eqnarray}

\begin{eqnarray}
\nonumber
B33=\left[\matrix{6&-3\cr -1&1\cr 9&-5}\right]\left[\matrix{-7&5\cr 4&-2\cr 4&-2}\right]\left[\matrix{-15&10\cr 4&-2\cr 4&-2}\right]\left[\matrix{6&-3\cr -1&1\cr 1&0}\right]\left[\matrix{-15&10\cr 4&-2\cr 1&0}\right]\left[\matrix{-7&5\cr 4&-2\cr 1&0}\right]\\ \nonumber 
\end{eqnarray}

\begin{eqnarray}
\nonumber
B34=\left[\matrix{-15&10\cr 4&-2\cr 4&-2}\right]\left[\matrix{11&-6\cr -2&2\cr 4&-2}\right]\left[\matrix{-15&10\cr 4&-2\cr 1&0}\right]\left[\matrix{-2&2\cr 1&0\cr 1&0}\right]\left[\matrix{11&-6\cr 1&0\cr 1&0}\right]\left[\matrix{-2&2\cr 4&-2\cr 1&0}\right]\\ \nonumber 
\end{eqnarray}

\normalsize

\subsection{The set $\widehat A$ and the Renormalization Map}
\label{AR}

The set $\widehat A$ is divided into $18$ layers,
$9$ of which are contained in the bottom half
of $\widehat \Sigma$.  The involution $I$ maps
the first $9$ layers to the second $9$.
For this reason, we will just explain the
first $9$ layers.  Again, each layer has $4$
branches, and so we have $36$ polyhedra in all.

Each layer of $\widehat A$ corresponds to an interval
of the form
\begin{equation}
I_{ij}=I_i \cap R^{-1}(J'_j)
\end{equation}
Here $I_1,...,I_5$ is the partition relative to
which the circle renormalization map $R$ is
defined.  These intervals are defined by the
partition 
\begin{equation}
\label{AApart}
0< \phi^{-2}<2  \phi^{-2}<2-\phi^{-2}<2-2  \phi^{-2}<2.
\end{equation}

The notation 
$A_{ijk}$ denotes the $k$th branch of the $(i,j)$th layer of $\widehat A$.
Rather than list the coordinates of the $A_{ijk}$, we will define these
polyhedra in terms of the renormalization map $\widehat R$.  This
approach simultaneously defines $\widehat A$ and $\widehat R$.
The map $\widehat R: \widehat A \to \widehat B$ is a
piecewise golden affine transformation of $\R^3$.  We
encode golden affine transformations by $8$-tuples of
integers.  The $8$-tuple $(a_1,...,a_8)$ corresponds
to the map $T(V)=rV+W$, where
\begin{equation}
r=a_1+a_2 \phi; \hskip 40 pt
W=(a_3+a_4 \phi, a_5+a_6 \phi,a_7+a_8 \phi).
\end{equation}
The maps we actually list are $\widehat R^{-1}$.

There are $6$ special cases.
\begin{itemize}
\item $A_{314}=\widehat R^{-1}(B_{14})$, where $\widehat R^{-1}=(-3,2,4,-2,8,-4,4,-2)$.
\item $A_{364}=\widehat R^{-1}(B_{64})$, where $\widehat R^{-1}=(-3,2,4,-2,0,0,4,-2)$.
\item $A_{11k}=\widehat R^{-1}(B_{13})$, where $\widehat R^{-1}=(-3,2,4,-2,8,-4,0,0)$ for $k=3,4$
\item $A_{56k}=\widehat R^{-1}(B_{63})$, where $\widehat R^{-1}=(-3,2,4,-2,0,0,8,-4)$ for $k=3,4$
\end{itemize}
Otherwise, we have
\begin{enumerate}
\item $A_{1jk}=\widehat R^{-1}(B_{jk})$, where  $\widehat R^{-1}=(-3,2,4,-2,4,-2,0,0)$
\item $A_{2jk}=\widehat R^{-1}(B_{jk})$, where  $\widehat R^{-1}=(1,0,0,0,0,0,2,-2)$
\item $A_{3jk}=\widehat R^{-1}(B_{jk})$, where, $\widehat R^{-1}=(-3,2,4,-2,4,-2,4,-2)$.
\item $A_{4jk}=\widehat R^{-1}(B_{jk})$, where, $\widehat R^{-1}=(1,0,0,0,0,0,-2,2)$
\item $A_{5jk}=\widehat R^{-1}(B_{jk})$, where, $\widehat R^{-1}=(-3,2,4,-2,4,-2,8,-4)$
\end{enumerate}

\newpage

\section{References}

\noindent
[{\bf DeB\/}] N. E. J. De Bruijn, {\it Algebraic theory of Penrose's nonperiodic tilings\/},
Nederl. Akad. Wentensch. Proc. {\bf 84\/}:39--66 (1981).
\newline
\newline
[{\bf D\/}] R. Douady, {\it These de 3-eme cycle\/}, Universit\'{e} de Paris 7, 1982.
\newline
\newline
[{\bf DF\/}] D. Dolyopyat and B. Fayad, {\it Unbounded orbits for semicircular
outer billiards\/}, Annales Henri Poincar\'{e}, to appear.
\newline
\newline
[{\bf DT1\/}] F. Dogru and S. Tabachnikov,  {\it Dual billiards\/}, 
Math. Intelligencer {\bf 26\/}(4):18--25 (2005).
\newline
\newline
[{\bf F\/}] K. J. Falconer, {\it Fractal Geometry: Mathematical Foundations
 and Applications\/},
John Wiley and Sons, New York (1990).
\newline
\newline
[{\bf G\/}] D. Genin, {\it Regular and Chaotic Dynamics of
Outer Billiards\/}, Pennsylvania State University Ph.D. thesis, State College (2005).
 \newline 
\newline
[{\bf GS\/}] E. Gutkin and N. Simanyi, {\it Dual polygonal
billiard and necklace dynamics\/}, Comm. Math. Phys.
{\bf 143\/}:431--450 (1991).
\newline
\newline
[{\bf Ko\/}] Kolodziej, {\it The antibilliard outside a polygon\/},
Bull. Pol. Acad Sci. Math.
{\bf 37\/}:163--168 (1994).
\newline
\newline
[{\bf M1\/}] J. Moser, {\it Is the solar system stable?\/},
Math. Intelligencer {\bf 1\/}:65--71 (1978).
\newline
\newline
[{\bf M2\/}] J. Moser, {\it Stable and random motions in dynamical systems, with
special emphasis on celestial mechanics\/},
Ann. of Math. Stud. 77, Princeton University Press, Princeton, NJ (1973).
\newline
\newline
[{\bf N\/}] B. H. Neumann, {\it Sharing ham and eggs\/},
Summary of a Manchester Mathematics Colloquium, 25 Jan 1959,
published in Iota, the Manchester University Mathematics Students' Journal.
\newline
\newline
[{\bf S1\/}] R. E. Schwartz, {\it Unbounded Orbits for Outer Billiards\/},
J. Mod. Dyn. {\bf 3\/}:371--424 (2007). 
\newline
\newline
[{\bf S2\/}] R. E. Schwartz, {\it Outer Billiards on Kites\/},
Annals of Mathematics Studies {\bf 171\/}, Princeton Press, 2009.
\newline
\newline
[{\bf S3\/}] R. E. Schwartz, {\it Outer Billiards and the Pinwheel Map\/},
preprint 2009.
\newline
\newline
[{\bf T1\/}] S. Tabachnikov, {\it Geometry and billiards\/},
Student Mathematical Library 30,
Amer. Math. Soc. (2005).
\newline
\newline
[{\bf T3\/}] S. Tabachnikov, {\it Billiards\/}, Soci\'{e}t\'{e} Math\'{e}matique de France, 
``Panoramas et Syntheses'' 1, 1995
\newline
\newline
[{\bf VS\/}] F. Vivaldi and A. Shaidenko, {\it Global stability of a class of discontinuous
dual billiards\/}, Comm. Math. Phys. {\bf 110\/}:625--640 (1987).

\end{document}